\documentclass[fleqn,leqno]{article}
\usepackage{amssymb,latexsym}
\usepackage{amsfonts}
\usepackage[noload]{qtree}
 \usepackage[dvips]{color}

\newcommand{\bi}{\bigskip}
\newcommand{\sm}{\smallskip}

\newcommand{\wh}{\widehat}
\newcommand{\ee}{\end{equation}}
\newcommand{\eea}{\end{eqnarray}}
\newcommand{\bean}{\begin{eqnarray*}}
\newcommand{\eean}{\end{eqnarray*}}
\newif\ifpctex
\newcommand{\noi}{\noindent}

\newcommand{\ve}{\varepsilon}

\newcommand{\wt}{\widetilde}

\newtheorem{theorem}{Theorem}
\newtheorem{remark}{Remark}
\newtheorem{example}{Example}
\newcommand{\D}{\displaystyle}

\newcommand{\qad}{{}\hfill \square}

\newtheorem{proposition}{Proposition}[section]
\newtheorem{definition}[proposition]{Definition}
\newtheorem{lemma}[proposition]{Lemma}
\newtheorem{corollary}[proposition]{Corollary}

\newtheorem{conjecture}[proposition]{Conjecture}
\newenvironment{proof}{\par\noindent{\bf Proof\ }}

 \renewcommand{\theequation}{\thesection.\arabic{equation}}

 \newcommand{\rand}[1]{}

\newcommand{\La}{\Longrightarrow}

\newcommand{\Rand}[1]{\marginpar{#1}}
     \renewcommand{\Rand}[1]{}

\newcommand{\be}[1]{\Rand{\vspace{0,6cm}\tt #1}\begin{equation}\label{#1}}
\newcommand{\bea}[1]{\Rand{\vspace{0,7cm}\tt
#1\vspace{-0,7cm}}\begin{eqnarray}\label{#1}} \marginparwidth2.5cm

\newcommand{\beL}[1]{\Rand{\vspace{0,6cm}\tt #1}\begin{lemma}\label{#1}}
\newcommand{\beD}[1]{\Rand{\vspace{0,6cm}\tt #1}\begin{definition}\label{#1}}
\newcommand{\beT}[1]{\Rand{\vspace{0,6cm}\tt #1}\begin{theorem}\label{#1}}
\newcommand{\beP}[1]{\Rand{\vspace{0,6cm}\tt #1}\begin{proposition}\label{#1}}
\newcommand{\beC}[1]{\Rand{\vspace{0,6cm}\tt #1}\begin{corollary}\label{#1}}
\newcommand{\beCj}[1]{\Rand{\vspace{0,6cm}\tt #1}\begin{conjecture}\label{#1}}

\newcommand{\ep}{\end{proposition}}

\newcommand{\suml}{\sum\limits}
\newcommand{\intl}{\int\limits}
\newcommand{\liml}{\lim\limits}
\newcommand{\supl}{\sup\limits}

\newcommand{\limsupl}{\limsup\limits}
\newcommand{\prodl}{\prod\limits}

\newcommand{\B}{\mathbb{B}}
\newcommand{\F}{\mathbb{F}}

\newcommand{\K}{\mathbb{K}}
\newcommand{\R}{\mathbb{R}}
\newcommand{\N}{\mathbb{N}}

\newcommand{\I}{\mathbb{I}}
\newcommand{\Z}{\mathbb{Z}}

\newcommand{\Q}{\mathbb{Q}}
\newcommand{\text}{\mbox}

\newcommand{\CA}{{\mathcal A}}
\newcommand{\CB}{{\mathcal B}}

\newcommand{\CE}{{\mathcal E}}
\newcommand{\CF}{{\mathcal F}}

\newcommand{\CI}{{\mathcal I}}
\newcommand{\CL}{{\mathcal L}}

\newcommand{\CM}{{\mathcal M}}
\newcommand{\CN}{{\mathcal N}}

\newcommand{\CP}{{\mathcal P}}
\newcommand{\CU}{{\mathcal U}}
 
\newcommand{\CW}{{\mathcal W}}

\newcommand{\CS}{{\mathcal S}}

\newcommand{\tto}{{_{\D \Longrightarrow \atop t \to \infty}}}

\newcommand{\ttooo}{{_{\D \longrightarrow \atop t \to -\infty}}}

\newcommand{\ntoo}{{_{\D \longrightarrow \atop n \to \infty}}}
\newcommand{\Nto}{{_{\D \Longrightarrow \atop N \to \infty}}}

\newcommand{\Ntoo}{{_{\D \longrightarrow \atop N \to \infty}}}

\newcommand{\la}{\longrightarrow}

\begin{document}


\title{Invasion by rare mutants in a spatial two-type Fisher-Wright system with selection}
\author{{\bf Donald A. Dawson$^{1,2}$ \mbox{ } Andreas Greven$^{2,3}$}}

\date{\small  \today}

\maketitle

\begin{abstract}

We consider a meanfield system of interacting Fisher-Wright diffusions with selection and rare
mutation on the geographic space $\{1,2,\cdots,N\}$. The type 1 has fitness 0, type 2 has
fitness 1 and (rare) mutation occurs from type  1 to 2 at rate $m \cdot N^{-1}$, selection
is at rate $s>0$. The system starts in the state concentrated on type 1, the state of low fitness.
We investigate this system for $N \to \infty$ on the original and large time scales.

We show that for some $\alpha \in (0,s)$ at times $\alpha^{-1} \log N+t, t \in \R, N \to \infty$
the emergence of type 2 (positive global type-2 intensity) at a global level occurs, while
at times $\alpha^{-1} \log N+t_N$, with $t_N \to \infty$ we get fixation on type 2 and on the
other hand with $t_N \to -\infty$ as $N \to \infty$ asymptotically only type 1 is present.

We describe the transition from emergence to fixation in the time scale
$\alpha^{-1} \log N+t, t \in \R$ in the limit $N \to \infty$ by a McKean-Vlasov random entrance law. This entrance law
behaves for $t \to -\infty$ like $^\ast\CW e^{-\alpha |t|}$ for a positive random variable $^\ast \CW$.
The formation of small droplets of type-2 dominated sites in times
$o(\log N)$, or $\gamma \cdot \log N, \gamma \in (0,\alpha^{-1})$ is described in the limit
$N \to \infty$ by a measure-valued process following a stochastic equation driven by
Poissonian type noise which we identify explicitly.
The total mass of this limiting $(N \to \infty)$ droplet process grows like $\CW^\ast e^{\alpha t}$ as $t \to \infty$.
We prove that exit behaviour from the small time scale equals
the entrance behaviour in the large time scale, namely
 $\CL[^\ast\CW] =\CL[\CW^\ast]$.

\end{abstract}

\noi {\bf Keywords:} Interacting Fisher-Wright diffusions, mutation,
selection, rare mutation, punctuated equilibrium,
random McKean-Vlasov equation, random entrance law, droplet formation,
atomic measure-valued processes.\\
\vspace*{2cm}

\footnoterule
\noi
\hspace*{0.3cm}{\footnotesize $^1$ {School of Mathematics and Statistics,
Carleton University, Ottawa K1S 5B6, Canada, e-mail:
ddawson@math.carleton.ca}}\\
\hspace*{0.3cm}{\footnotesize $^2$
{Research supported by NSERC, DFG-Schwerpunkt 1033 and later DFG-NWO
Forschergruppe 498, and D. Dawson's Max
Planck Award for International
Cooperation.}}\\
\hspace*{0.3cm}{\footnotesize $^3$
{Department Mathematik, Universit\"at Erlangen-N\"urnberg,
Bismarckstra{\ss}e 1 1/2,
D-91054 Erlangen, Germany, e-mail: greven@mi.uni-erlangen.de}}
\newpage

\tableofcontents

\newpage

\newcounter{secnum}
\setcounter{secnum}{\value{section}}
\setcounter{section}{0}
\setcounter{equation}{0}
\renewcommand{\theequation}{\mbox{\arabic{secnum}.\arabic{equation}}}

\section {The model and the results}
\label{s.model}

\setcounter{secnum}{\value{section}} \setcounter{equation}{0}

\subsection{Background and motivation}
\label{ss.background}

Mutation and selection both play essential roles in the
evolution of a population. Mutation increases genetic diversity,
while selection reduces this diversity pushing it towards states
concentrated on the fittest types. The balance between these two
forces is influenced heavily by the effect of migration. Randomness
enters through resampling (pure genetic drift).

To give a precise setting for this study we model the
population dynamics according to a stochastic process arising as
the diffusion limit of many individuals incorporating the basic
mechanisms for a spatial population traditionally employed in
genetics and evolutionary biology. In this model {\it random
fluctuations} are included  in contrast to other infinite
population limits which are ruled by deterministic differential
equations (for example, see \cite{Bu01}). More precisely the model
we use arises from the particle model driven by migration of particles
between colonies, and in each site by resampling of types,
mutation and selection. Increase the number of particles as
$\ve^{-1}$ and give them mass $\ve$. As $\ve\to 0$, a diffusion
limit of interacting multitype diffusions results if the
resampling rate is proportional to the number of pairs and there
is weak selection, i.e. selection occurs at a rate decreasing with
the inverse of the number of particles per site. Otherwise with
strong selection, i.e. selection at a fixed rate, we get the
deterministic limit, often referred to as infinite population
limit.

These stochastic models
which describe populations of possibly infinitely many types  are
called interacting Fleming-Viot processes and have been
intensively studied in recent years (for example, see \cite{EG},
\cite{D}, \cite{EK1} - \cite{EK3}, \cite{DGV}, \cite{DG99}).  The basic mechanisms are
\emph{migration} between colonies in the diffusion limit this takes the
form of deterministic mass flow between sites and then within each colony
\emph{resampling (i.e. pure genetic drift)} which gives the diffusion term,
\emph{selection} of
haploid type based on a fitness function on the space of  types
and \emph{mutation} between types occuring at rates
given by a transition kernel on the
space of types, the latter two in the diffusion limit are drift terms
giving a deterministic massflow between types.

\bi

The purpose of this paper is to investigate how new types brought in by rare mutation
invade and conquer the whole population through the composed effect of resampling,
mutation, selection and migration.

The first part of the picture we  focus on is a  scenario (\emph{punctuated
equilibrium}) which has received considerable attention in the biological
literature (see for example \cite{EG1},\cite{EG2} and the references below).
In this scenario a population
remains for long periods of time in a nearly stationary state
(\emph{quasi-equilibrium or stasis}) until the emergence of one or several rare mutants
of higher fitness that then take over the population relatively
quickly. The explanation of this phenomenon has two aspects
tunnelling between local fitness maxima and the effect of the
spatial structure.

For example one explanation of this  tunnelling was given by Newman, Cohen and Kipnis
\cite{CKN} in the context of Wright's theory of evolution in an
adaptive landscape (see \cite{Wr1},\cite{Wr2},\cite{W}). In this
explanation the long times between transitions from one
quasi-equilibrium to another corresponds to ``tunnelling'' between
one adaptive peak and another in a situation in which local
stabilizing selection maintains a population near a local peak. In
this setting, the Ventsel-Freidlin theory of small random
perturbations of dynamical systems  asserts that the actual
transition, when it occurs, takes place rapidly.  Thus this
scenario would be consistent with the observations of
paleontologists of periods  of rapid evolution separated by very
long periods of ``stasis''. This is a stochastic  effect  not
present in the deterministic version, that is,  the infinite
population limit, but must be formulated in terms of a stochastic
diffusion limit, which we consider here.

A second part of the picture is the role of geographic space and its strong
effects on the qualitative behaviour. Namely typically a large population is subdivided into
small subpopulations (colonies) occupying different geographic
regions. Then tunnelling between adaptive peaks can occur in these
subpopulations and then  spatial migration of  individuals from a
colony corresponding to a higher adaptive peak can spread and result in the
take-over of  the entire population which then fixates on the  fitter type.
\bi

Our goal is to develop a rigorous framework to discuss a model incorporating these features
of a punctuated equilibrium (by rare mutation) but with a different driving mechanism
for the tunnelling starting from rare mutants
than in the literature quoted above. In this context our goal is also more generally to give a framework which
provides insight into the  respective roles of
mutation, selection and spatial migration in the emergence and
spatial spread of rare mutants corresponding to a higher adaptive
peak (quasi-equilibrium).

In this paper we develop the scenario in the {\em two-type} case, one type of low and one
type of high fitness.
The full picture meaning (1) an infinite hierarchy of fitness levels, (2) with
$M \geq 2$ (instead of $M=1$ here) types on the successive levels and (3)
all this embedded on a countable geographic space instead of $\{1,2,\cdots,N\}$,
which approximates $\Z^2$ is developed in \cite{DGsel}.

Our goal in this paper is to exhibit how a rare mutant which is fitter than the
types in the current population invades the population, by starting to take
over a sparse set of colonies by selective advantage and then spreads as
an increasing droplet in space till the complete population is finally taken
over. The basic technique to do so is to consider the limit $N \to \infty$.

Indeed using this technique of taking $N \to \infty$ we can produce
limiting objects, which have a simple and fairly explicit
description.  However due to the interaction between the
migration, selection, mutation and resampling mechanisms, we have
to develop some new ideas to adapt the method of multi-scale
analysis to resolve some delicate conceptional and mathematical problems that arise
in this case. (In previous work \cite{DG99} the qualitative analysis
of the longtime behaviour did not include mutation.)

A further point is that we give here a rigorously treated example for two-scale phenomena
and develop as mathematical tool the notion
of random nonlinear evolutions, more specifically solutions with \emph{random McKean-Vlasov
entrance laws} from time $-\infty$. The technique of random entrance laws allows
to connect the separating time scales. This phenomenon is of relevance in  many other
applications.

In order to describe the early occurrence
of rare mutants we work with a description of the sparse set of colonies
conquered by the rare mutants via {\em atomic measure-valued processes} which
in the limit $N \to\infty$ can be described by stochastic processes driven by
(inhomogeneous) Poisson random measures and in order to calculate the intensity
measure we are using excursion theory.

We have developed in \cite{DGsel} a technique to extend this picture to infinite
geographic space, multiple types on each level of fitness and infinite hierarchy
of fitness levels. There we get a good set-up to study the scenario of punctuated
equilibria.
However this analysis still suffers some shortcomings.
The restriction of our method is that we assume
that the different time and space scales {\em separate} in the limit of
large times and large space scales which is of course only an
approximation but allows us to make precise the notion of a
quasi-equilibrium as equilibrium of a limiting non-linear McKean-Vlasov
dynamic.  Moreover the techniques can be extended to the multitype case and then allow to define
the concept of a {\em quasi-equilibrium} in a mathematically rigorous
fashion as equilibria of certain limiting dynamics in the various
time scales.

\subsection{The model and key functionals }
\label{ss.functionals}

We begin here below with the description of the problem for the precisely
stated model, define key functionals and formulate the problem.

\subsubsection{The model }
\label{sss.model}

We define a system of $N$-exchangeably interacting  finite populations ($d>0$) with two
 types $\{1,2\}$ with fitness 0 and 1 respectively,  selection rate  $s>0$ and
 migration rate $c>0$.

This means that we study the stochastic process $(X^N(t))_{t \geq 0}$, which is uniquely
defined as follows:
\be{angr1}
X^N(t) =((x^N_1 (i,t),x^N_2 (i,t)); i=1,\dots,N),
\ee
\be{angr1b}
x^N_1 (i,0) =1, \quad i=1,\cdots,N,
\ee
satisfying the wellknown SSDE where $i \in \{1,2,\cdots,N\}$:
\bea{angr2} d x^N_1 (i,t) = c(\bar x^N_1 (t) - x^N_1 (i,t)) dt
 &-& s\,x^N_1(i,t) x^N_2 (i,t) dt
 - \frac{m}{N} x^N_1 (i,t)dt \\
 &+& \sqrt{d \cdot x^N_1(i,t)x^N_2(i,t)} dw_1(i,t), \nonumber
 \eea

\bea{angr3}
d x^N_2 (i,t) = c(\bar x^N_2 (t) - x^N_2 (i,t)) dt
 &+& s\,x^N_1(i,t) x^N_2 (i,t) dt
 + \frac{m}{N} x^N_1 (i,t)dt \\
 &+& \sqrt{d \cdot x^N_1 (i,t)x^N_2 (i,t)} dw_2 (i,t),\nonumber
 \eea
where $\{(w_1(i,t))_{t \geq 0},
i=1,\dots,N\}$ are i.i.d. Brownian motions and $w_2, w_1$ are coupled via
 $w_2 (i,t) = - w_1(i,t)$,
\be{angr3c}
\wh {x}_\ell (t) =  \suml^{N}_{i=1} x^N_\ell(i,t)
   \mbox{ with } \ell=1,2
\ee and
\be{angr3b} \bar {x}_\ell (t) = \frac{1}{N}\wh {x}_\ell (t)
\mbox{ with } \ell=1,2. \ee
Without loss of generality we can  assume for the resampling rate $d$
that $d=1$ and will do this except where we
wish to indicate the roles of the parameters $c,s,d,m$.
The simple collection of $N$ one-dimensional stochastic
differential equations in (\ref{angr1})- (\ref{angr3}) allows us to focus first
on the key features of
emergence, which are already complicated enough as we shall see.

{\bf Step 2}  {\em (Key functionals and main objectives)}

The description of the system
proceeds  via four objects, namely either locally by a
{\em sample of tagged sites} or globally by the {\em empirical measure}
over the whole spatial collection and furthermore we use as well an {\em atomic measure}
to globally describe sparse (in space) spots of type two mass respectively the
{\em Palm measure} to describe this locally by formalizing the concept of a
typical type-2 site. The asymptotic analysis of these
four objects leads to propagation of chaos results, respectively, nonlinear
McKean-Vlasov limiting dynamics and atomic measure-valued processes driven
by Poisson processes and formulas for Palm measures based on this.

Since it often suffices to
consider the tagged site 1, i.e. $x^N_\ell (1,t)$,  in the sequel
to designate a  tagged site we set\be{angr4} x^N_\ell (t) =
x^N_\ell (1,t) \quad ; \ell=1,2.
\ee

If we are interested in the {\em local} picture we consider for some $L \in \N$ the tagged sites
\be{angr5++}
(x^N_\ell(1,t), \cdots, x^N_\ell(L,t)), \quad \ell=1,2.
\ee
We abbreviate the marginal law of the type-2 mass as
\be{add133}
\mu^N_t := \CL [x^N_2(1,t)] \in \CP([0,1]).
\ee

The empirical measure process gives a {\em global} picture and is
defined by: \be{angr5} \Xi_N (t) := \frac{1}{N} \suml^N_{i=1}
\delta_{(x^N_1 (i,t), x^N_2(i,t)})\in
\mathcal{P}(\mathcal{P}(\{1,2\}))\ee and the empirical measure
process of either type is defined by
\be{angr5+}
\Xi_N(t,\ell):=\frac{1}{N} \sum_{i=1}^N \delta_{x_\ell^N(i,t)}\in
\mathcal{P}([0,1]),\; \ell=1,2. \ee Note that since $x^N_1 (i,t) +
x^N_2 (i,t) =1$ it suffices in the case of two types to know one
component of the pair in (\ref{angr5+}), the other is then determined by this
condition. Then we can for two types effectively replace $\CP(\CP(\{1,2\}))$ by $\CP([0,1])$
as states of the empirical measure.

Next we have to describe the sparse sites where type 2 appears.
First consider $N$ sites starting only with type  $1$ mass at time
$0$ and evolving in the  time interval $[0,T_0]$. As a result of the rare
mutation from type $1$ to type $2$ at rate $\frac{m}{N}$ mutant
mass will appear and then form non-negligible colonies
at a sparse set of the total of $N$ sites, but the typical site
will have mass of order $o(1)$ as $N \to \infty$.
We can now again take a {\em global} perspective, where we describe
the sparse set of colonised (by type 2) sites or a {\em local}
perspective, where we consider a typical colonised site via
the concept of the Palm measure.
We introduce both viewpoints successively.

Turn to the global perspective. In order to
keep track of the sparse set of sites at which nontrivial mass appears we
will give a random label to each site and define the following
{\em atomic-measure-valued process}.

 We assign independent of the process a point $a(j)$ randomly in $[0,1]$ to each site
$j\in\{1,\dots,N\}$, that is, we define the collection
\be{dd29a}
\{a(j), \quad j=1,\cdots, N \} \mbox{  are i.i.d. uniform on } [0,1].
\ee
We then associate with our process and a realization of the
random labels a measure-valued process on $\CP( [0,1])$, which we denote by
\be{dd29}
\gimel^{N,m}_t=\sum_{j=1}^N x^N_2(j,t) \delta_{a(j)}.
\ee

This description by $\gimel^{N,m}_t$ is complemented by a local perspective,
arising by zooming in on the random location where we actually find mass by studying a
{\em typical} type-2 dominated site.  In order to make
precise the notion of a site as seen by a typical type-2 mass,
i.e. a site seen from a randomly (among the total population in all $N$ sites)
chosen individual of type 2, we
use the concept of the {\em Palm-distribution}, which we indicate
by a hat. Define
\be{Y12d}
\wh \mu^N_t = \wh \CL [x^N_2 (1,t)],
\mbox{ by } \wh \mu^N_t(A) = \left(\intl_A x \mu^N_t
(dx)\right)/E[x^N_2 (1,t)], \ee
where $\mu^N_t =  \CL [x^N_2 (1,t)]$.
The law $\wh \mu^N_t$ describes now for an i.i.d. initial state the law of a {\em site, typical} for
the type-2 population at time $t$ in the limit $N \to \infty$. Namely we can think of this as picking
a type-2 individual at random and then determine the site where it sits. Note that
\be{tsi1}
\int x\wh\mu_t^N(dx)=\frac{E[(x^N_2(1,t))^2]}{E[x^N_2(1,t)]}.
\ee
\bi

{\bf The scenario}

The basic objective of this paper is to describe the emergence
of the more fit type, i.e. type $2$, using the above functionals of the process.
This emergence occurs in {\em three regimes}, two are of short duration
separated by a third of long duration and arise as follows in an asymptotic
description as $N \to \infty$.

{\bf (1)}  In the first regime
during times of order $O(1)$ rare mutants arise at a small set of sites
(droplet) where they locally make
up substantial part of the population and then develop mutant
``droplets'' by migration. We will see later in Subsection
\ref{ss.bounds} in all detail that the mutant population during this stage
develops as a type of measure-valued branching process. We refer to this
stage as {\em droplet formation}.

{\bf (2)} The next
stage, namely, the stage at which type $2$ has very small but
 $O(1)$ frequency at a typical (i.e. randomly chosen site) site, or equivalently at
the macroscopic level has $O(1)$ density, occurs at time
$O(\log N)$, we refer to this as {\em emergence}.

{\bf (3)} The last regime at which the macroscopic density of type 2 reaches
a frequency arbitrarily close to 1 is taking a further piece of time but is
again only of order $O(1)$ and is referred to
as {\em fixation}. For the regime between emergence and fixation we give a
description by a limiting evolution in the form of an entrance law
 evolving from
time $-\infty$ and type-1 states to type-2 states at time $+\infty$ and hence
altogether with time index $t \in \R$.

In the sequel we will analyse these three regimes
of droplet formation, emergence and fixation
by studying  the behaviour as $N \to \infty$ of
\be{agrev80}
\Xi_N(t), \quad \{x^N_\ell (i,t);  \ell = 1,2; i=1,\cdots,L\}
\mbox{ and } \gimel^{N,m}_t, \quad \wh \mu^N_t,
\ee
in finite times as well as in larger $N$-dependent time scales thus capturing
each of the three regimes described above.

\subsection{Statement of results}
\label{ss.6results}

We have three tasks, first to determine the
\emph{macroscopic emergence time scale}, that is the time scale at
which the level two type appears at a typical site or equivalently has positive intensity, second, to
describe the pre-emergence picture, in particular the {\em droplet
formation} respectively the {\em limiting droplet dynamics}  and thirdly to describe the \emph{limiting dynamic of
fixation} after the emergence with which the take-over by type 2 then occurs.

We begin the analysis with introducing
the two different limit dynamics in the subsequent two subsubsections and then
subsequently the emergence, fixation and convergence results are stated in
a further subsubsection.

\subsubsection{The limiting random McKean-Vlasov dynamics and random entrance laws}
\label{sss.MVdyn}

We start with the initial condition $\Xi_N(0,\{1\})=\delta_1$ (recall (\ref{angr5+})), that is
only type $1$ appears. The objective is to establish the emergence
and fixation of type $2$ in times of the form  $t_N
=C\cdot \log N+t$, to identify the constant $C$,  and to
identify the limit of $\Xi_N(t_N, \cdot))$ as a function of $t$ in
the process of emergence and determine the dynamics in $t$ of the fixation (takeover)
which leads to the concentration of the complete population in type-2 as $N\to \infty$.
In the next subsubsection we shall then discuss the
behaviour at early times, i.e. times $t_N$ with
$\limsup ((\log t_N)/\log N)<C$, in the stage of droplet formation.

The key ingredients in the limiting processes in times $t_N$ for the local
and global description are:
\begin{itemize}
\item the limiting {\em  (nonlinear) McKean-Vlasov dynamics},
\item the {\em entrance law} from $-\infty$ of the McKean-Vlasov dynamic and
\item {\em random} solutions to the McKean-Vlasov equation.
\end{itemize}
In order to define these ingredients we proceed in three steps,
we recall first in Step 1 the ``classical'' McKean-Vlasov limit
(associated with a nonlinear Markov process)
before in Step 2 and Step 3 we introduce the two new objects.

{\bf Step~1} Consider the above system (\ref{angr1})-(\ref{angr3})
of $N$ interacting sites with type space $\mathbb{K}=\{1,2\}$ starting at
time $t=0$ from a product measure (that is, i.i.d. initial values
at the $N$ sites). The basic McKean-Vlasov limit (cf. \cite{DG99},
Theorem 9) says that if we start initially in an i.i.d.
distribution, then \be{angr6} \{\Xi_N(t)\}_{0\leq t\leq T} \Nto
\{\mathcal{L}_t\}_{0\leq t\leq T}, \ee where the
$\CP(\CP(\mathbb{K}))$-valued path $\{\mathcal{L}_t\}_{0\leq t
\leq T}$ is the law of a {\em nonlinear} Markov process, namely the unique weak
solution of the {\em McKean-Vlasov equation}: \be{angr7b} \frac {d
\mathcal{L}_t}{dt}= ({L}_t^{\mathcal{L}_t})^\ast \mathcal{L}_t,
\ee where for $\pi \in \mathcal{P}(\mathcal{P}(\mathbb{K}))$,  and $F(\mu)=f(\mu(2))$,
\be{}L^{\pi}F = c[\int y\pi(2,dy)-x]\frac{\partial f}{\partial x}+sx(1-x)\frac{\partial f}{\partial x}+\frac{d}{2}x(1-x)\frac{\partial^2 f}{\partial x^2}
\ee
and the $\ast$ indicates the adjoint
of an operator mapping from a dense subspace of
$C_b(E, \R)$ into $C_b(E,\R)$
w.r.t. the pairing of $\CP(E)$ and $C_b(E,\R)$
given by the integral of the function with respect to the measure. Similar equations have been
studied extensively in the literature (e.g. \cite{Gar}).
The process $(\CL_t)_{t \geq 0}$ corresponds to a nonlinear Markov process since
$\CL_t$ appears also in the expression for the generator $L$.

As pointed out above in (\ref{angr5+}), in the special case
$\mathbb{K}=\{1,2\}$, we can simplify by  considering the
frequency of type 2 only and by reformulating (\ref{angr7b})
living on $\CP (\CP(\{1,2,\}))$ in
terms of $\mathcal{L}_t(2) \in \mathcal{P}[0,1]$.  This we carry
out now.

Namely we note that given the mean-curve
\be{angr7b0a}
m(t) = \intl_{[0,1]} y\,\mathcal{L}_t(2)(dy), \ee the process
$(\CL_t(2))_{t \geq 0}$ is the {\em law} of the strong solution of
(i.e. the unique weak solution) the SDE:
\be{angr7b0b} dy(t) = c(m(t)-y(t)) dt + sy(t)(1-y(t))dt +
\sqrt{dy(t)(1-y(t)} dw(t).
\ee
Then informally $(\CL_t)_{t \geq 0}$ corresponds to the
solution of the nonlinear diffusion equation. Namely for $t >0$,
$\CL_t(2)(\cdot)$ is absolutely continuous and for
\be{angr7b1}
 \mathcal{L}_t(2)(dx)=u(t,x)dx \;\in \mathcal{P}([0,1])
\ee
the evolution equation of the density $u(t,\cdot)$ is given by:
\be{angr7b0}
 \frac{\partial}{\partial t} u(t,x) =
  - c\frac{\partial}{\partial x}\{[\intl_{[0,1]} y u(t,y) dy - x]{u(t,x)}\}
-s\frac{\partial}{\partial x}(x(1-x)u(t,x))
  + \frac{d}{2}\frac{\partial^2}{\partial x^2}(x(1-x)u(t,x)).
\ee

{\bf Step~2} In order to describe the emergence and invasion
process via $\Xi_N$, we introduce in this step and the next two
extensions of the nonlinear McKean-Vlasov dynamics which describes
only the limiting evolution over finite time stretches, $t\in
[t_0,t_0+T]$ given an initial condition $\mathcal{L}_{t_0}$ which
we then follow as $t_0\to -\infty$.  Hence we have to define the
dynamics for $t\in (-\infty,\infty)$ in terms of an {\em entrance
law} at $-\infty$.

Since we consider the limits of systems observed in the interval
$const\cdot \log N +[-\frac{T}{2},\frac{T}{2}]$ with $T$ any
positive number, that is setting $t_0(N) = const \cdot \log N
-\frac{T}{2}$,  we need to identify entrance laws for the process
from $-\infty$ (by considering $T \to \infty$) out of the state
concentrated on type 1 with certain properties.

\beD{D.entlaw} {(Entrance law from $t = -\infty$)}

We say in the two-type case that a probability measure-valued function
$\CL : \R \to \CP([0,1])$, is an entrance law at $-\infty$
starting from type 1 if
$(\CL_t)_{t \in (-\infty,\infty)}$ is
such that $\CL_t$ solves the McKean-Vlasov equation (\ref{angr7b}) and
$\CL_t \to \delta_1$ as $t \to -\infty$.

In the case of more than two
types we work with maps $\CL: \R \to \CP(\CP(\K))$ where we require
as $t \to -\infty$ that $\CL_t$ converges to measures $\delta_\mu$
with $\mu \in \CP(\K)$ such that $\mu$ is a measure concentrated on
the lower level types.
$\qquad \square$
\end{definition}

The existence of such an object is obtained in part (c) of the
proposition below. \sm

{\bf Step~3} The usual formulation of the McKean-Vlasov limit
requires that we start in an i.i.d. initial configuration. This is
not sufficient for us since even though we consider systems
observed in finite time stretches  we do so only {\em after large
times}, namely,  after the time $const \cdot \log N$ with $N \to \infty$.
Hence in our context the McKean-Vlasov limit is valid in a fixed
time scale but if the system is viewed in time scales that depend
on $N$ this can (and will) break down since we only know that the initial state
is exchangeable and this then leads to a random solution.

We will indeed establish that the emergence of rare mutants gives
rise to \emph{``random''} solutions of the McKean-Vlasov dynamics.
In particular we will show that the limiting empirical measures at
times of the form $C \log N +t$ are {\em random} probability measures on
$[0,1]$ and therefore given by exchangeable sequences of $[0,1]$-
valued truely exchangeable random variables which are \emph{not} i.i.d., that is,
the exchangeable $\sigma$-algebra is not trivial.
This means that the empirical mean turns out to be a {\em random variable}
and since this is the term driving via the migration the local evolution
of a site in the McKean-Vlasov limit, the non-linearity of the evolution
equation comes seriously into play. However once we condition on the
exchangeable $\sigma$-algebra, we then get for the further evolution
again a deterministic limiting equation for the empirical measures,
namely the McKean-Vlasov equation. The reason for this is the fact
that conditioned on the exchangeable $\sigma$-algebra we obtain on
asymptotically (as $N \to \infty$) i.i.d. configuration to which
the classical convergence theorem applies. Using the Feller property
of the system, which is a direct consequence of the duality, we
get our claim.

This leads to
the task of identifying an entrance law in terms of a {\em random
initial condition at time $-\infty$}.
A consequence of this scenario is that when we use the duality from
time $T_N$ to $T_N+t$, we apply it to a random initial state in the
limit and we therefore have to use because of the non-linearity of the
evolution the appropriate formulas.

The above discussion shows that we need to
introduce  the notion of a truly random McKean-Vlasov entrance law
from $-\infty$.

\beD{D.ranMKV}{(Random solution of McKean-Vlasov)}

We say that the  probability measure-valued process
$\{\CL(t)\}_{t\in\mathbb{R}}$ is a {\em random solution of the
McKean-Vlasov equation} (\ref{angr7b}) if
\begin{itemize}
\item $\{\CL_t:t\in \R\}$ is a.s. a solution to (\ref{angr7b}),
that is, for every $t_0$ the distribution of $\{\CL_t:t\geq t_0\}$
conditioned on $\mathcal{F}_{t_0}=\sigma\{\CL_s:s\leq t_0\}$ is
given by $\delta_{{\{\mu_t}\}_{t\geq t_0}}$ where $\mu_t$ is a
solution of the McKean-Vlasov equation with $\mu_{t_0}=\CL_{t_0}$,
\item the time $t$ marginal distributions of $\{\CL_t : t \in
\R\}$ are truly random.$\qquad \square$
\end{itemize}
\end{definition}
\bigskip

The key result of this subsubsection on the objects introduced in the
previous three steps is now the following existence and uniqueness results
on the solution of the McKean-Vlasov equation (\ref{angr7b}) and of its
(random) entrance laws from $t = -\infty$:

\beP{P.MVentrlaw} {(McKean-Vlasov entrance law from $-\infty$)} \\
(a) Given the initial state $\mu_0 \in \CP([0,1])$
there exists a unique solution
\be{angr8b} \mathcal{L}_t(2)(dx) = \mu_t(dx),\;\; t\geq t_0 \ee to
(\ref{angr7b}) with initial condition $\mathcal{L}_{t_0}(2)=\mu_0$. \\
(b) If $s>0$ and $ \int_{[0,1]} x\mu_{t_0}(dx) >0$, then this
solution satisfies: \be{angr8c}
   \lim_{t\to\infty}
\mathcal{L}_t(2)(dx) =\delta_1(dx).\ee
\noindent
(c) There exists a solution
$(\mathcal{L}_t^{\ast \ast}(2))_{t \in \R}$ to  equation
(\ref{angr7b})
satisfying the conditions:
\bea{Ex1b}
\lim_{t\to-\infty} \mathcal{L}^{\ast\ast}_t(2) &= &\delta_0,\\
\lim_{t\to\infty} \mathcal{L}^{\ast\ast}_t(2) &= &\delta_1 \nonumber \\
\int_{[0,1]}x\mathcal{L}^{\ast\ast}_0 (2,dx) &=
&\frac{1}{2}.\nonumber \eea This solution is called an entrance
law from $-\infty$ with mean
$\frac{1}{2}$ at $ t=0$.\\
(d) We can obtain a solution in (c) such that:
\be{Ex2}
\exists\; \alpha
\in (0,s)\mbox{ and  }A_0\in(0,\infty)\mbox{ such that }
\lim_{t\to-\infty}e^{\alpha |t|}
 \int_{[0,1]} x \mathcal{L}^{\ast\ast}_t(2,dx) =A_0.
\ee (e) The solution of (\ref{angr7b}) also satisfying
(\ref{Ex2}) for prescribed $A_0$ is unique and if $A_0 \in (0, \infty)$
then $\alpha$ is necessarily uniquely determined.

For any deterministic solution
\be{Ex2a}
\mathcal{L}_t,\; t \in \R
\ee
 to (\ref{angr7b}) with
\be{Ex2b} 0 \leq \limsup_{t\to -\infty}e^{\alpha |t|}\int_{[0,1]}
x \mathcal{L}_t(2,dx) <\infty, \ee the limit
$A=\lim_{t\to-\infty}e^{\alpha |t|} \int_{[0,1]} x
\mathcal{L}_t(2,dx)$ exists.

If $A>0$, then $\{\mathcal{L}_t, t \in \R\}$ is
given by a time shift of the then unique
$\{\mathcal{L}^{\ast \ast}_t, t \in \R\}$ singled out in (\ref{Ex2}), i.e.
\be{Ex2c}
\CL_t = \CL^{\ast\ast}_{t+\tau}, \quad \tau =
\alpha^{-1} \log \frac{A}{A_0}.
\ee
For future reference we define
$(\CL^\ast_t)_{t \in \R}$ to be the unique solution satisfying
\be{L.ast}\lim_{t\to-\infty}e^{\alpha |t|}
 \int_{[0,1]} x \mathcal{L}^\ast_t(2,dx)=1
\mbox{ for some } \alpha \in (0,s).
\ee

(f) Any random solution $(\wt {\mathcal{L}}_t)_{t \in \R}$ to (\ref{angr7b})
such  that
\be{Ex2d}
\limsup_{t\to -\infty}e^{\alpha |t|}E[\int_{[0,1]} x \mathcal{L}_t(2,dx)] <\infty,\quad
\liminf_{t\to -\infty}e^{\alpha |t|}[\int_{[0,1]} x
\mathcal{L}_t(2,dx)] >0 \mbox{ a.s.},
\ee is a random time shift of
$(\mathcal{L}^{\ast\ast}_t)_{t \in \R}$ (and of $(\CL^\ast_t)_{t \in \R}).
 \qad$
\end{proposition}

\begin{example}  Let $\mathcal{L}^\ast_t$ be a solution
satisfying (\ref{angr7b0}), (\ref{Ex2}) and for a given value of
$A$  let $\tau$ be a true real-valued random variable.  Then
$\{\mathcal{L}^\ast_{t-\tau}\}_{t\in\mathbb{R}}$ is a truly random
solution. This can also be viewed as saying that we have a
solution with an exponential growth factor $A$ which is truly
random.
\end{example}

\begin{remark} We shall derive later on a bound from above on the
cardinality of our dual process which will imply that the first relation in (\ref{Ex2d})
must be satisfied for a limiting  dynamic (as $N \to \infty$) of our
empirical measures $\Xi_N$ and furthermore we
get a lower bound implying that $A$ is a.s. positive for this
limiting dynamic, (see (\ref{grev90.a}) and Subsubsection \ref{sss.collreg}).
Then we are able to use the identification of the solution given in
(\ref{Ex2c}) or of the one satisfying (\ref{L.ast}) to identify the limiting dynamic of the process of
empirical measures.
\end{remark}

\subsubsection{Limiting dynamics of sparse sites with substantial type-2 mass: Droplets}
\label{sss.dynsparse}

We next  describe the limiting dynamic of the sparse set of sites
which have been colonized by type 2 prior to the onset of
emergence. There are two time regimes, first an initial finite
time horizon of order $O(1)$, then large times $t_N \to \infty$ as
$N \to \infty$ but which remain $<< \alpha^{-1} \log N$ so that
we are in the preemergence regime with still a global density of
type 2 which is asymptotically zero.

It is very unlikely that a randomly chosen site   has at time O(1) mass of type
2 at least $\ve
>0$  as $N \to \infty$.  In fact an explicit
calculation (see (\ref{agrev1b})) shows that the transition density
decays like $const N^{-1}$.  On the other hand the number of sites
increases with $N$  so that here we have  a compensation provided we do
not look at only one single tagged site  but the complete population at all sites and the
result is that considering all sites there is a  {\em finite random number} of sites with
substantial type 2 mass.  Newly added  sites in this set of sites arise
from a process
 entering from the state 0. In fact we obtain here a Poisson
distribution in the limit. This whole scenario  will be made
precise using entrance laws from state 0 which we derive using
diffusion theory and excursion law theory.

Therefore before formulating the main results of this section we
first consider entrance laws of a single site with no mutation
from type $1$ to type $2$ which provides us with a key ingredient,
the {\em excursion measure $\Q$}. We state the wellknown fact:

\beL{L.Pa2} {(Single site: entrance and excursion laws)}

(a)
 Let $c>0, d>0,\; s>0$.  Then $0$ is an exit boundary for the
the Fisher-Wright diffusion
\be{onedim}
dx(t)=-c x(t)dt
+sx(t)(1-x(t))dt +\sqrt{d \cdot x(t)(1-x(t))}dw(t),
\ee
which then has a $\sigma$-finite entrance law from state $0$
at time 0, the $\sigma$-finite excursion law
\be{dd300}
{ \mathbb{Q}} = \mathbb{Q}^{c,d,s}\ee on
\be{dd30} W_0:=\{w\in C([0,\infty),\mathbb{R}^+),\;w(0)=0,\;
w(t)>0\mbox{ for  } 0<t<\zeta \mbox{ for some }
\zeta\in(0,\infty)\}.\ee

(b) Moreover, denoting by $P^\ve$ the law of the process started
with $w(0)=\ve$ and $\ve> 0, \Q$ is given by:
\be{dd31a}
 \mathbb{Q}(\cdot) =\lim_{\ve\to 0}\frac{P^\ve(\cdot)}{S(\ve)},
 \ee
where $S(\cdot)$ is the scale function of the diffusion
(\ref{onedim}), defined by the relation,
\be{dd31ab}
P_\ve(T_\eta<\infty)=\frac{S(\ve)}{S(\eta)},\qquad
0<\ve<\eta<\infty,
\ee
where $T_\eta$ is the first hitting time of $\eta$.

For the Fisher-Wright diffusion $S$ is given by (cf. \cite{RW},
V28)) the initial value problem:
\be{dd31b} S(0)=0,
\frac{dS}{dx}=\frac{e^{-2sx}}{(1-x)^{2c}}, \ee so that \be{dd31c}
\lim_{\ve\to 0}\frac{S(\ve)}{\ve}=1. \ee

(c) The measure $\Q$ is $\sigma$-finite, namely
for any $\eta >0, \zeta$ as in (\ref{dd30}),
\be{dd32x} \mathbb{Q} (\{w:\zeta(w)>\eta\})<
\infty,\ee
\be{dd32b2}
\mathbb{Q}(\sup_t(w(t))>\eta)=\mathbb{Q}(T_\eta
<\infty)= \frac {1}{S(\eta)} \la \infty \mbox{ as } \eta \to 0,
\ee and \be{dd32b} \int_0^1 x
\mathbb{Q}(\sup_t(w(t))\in dx)=\infty.\qquad \square \ee
\end{lemma}

\begin{proof}
It is wellknown that $0$ is an exit boundary  for (\ref{onedim}).
In this case (a) - (c) then  follows immediately from  the results of Pitman
and Yor \cite{PY}, Section 3.

\end{proof}
\bi

We next  identify the limit in distribution  of $\gimel^{N,m}_t$
(defined in (\ref{dd29}) ) as $N\to\infty$. Let
\be{ddx1}
\CM_a ([0,1]) = \mbox{ the set of finite atomic measures on } [0,1],
\ee
equipped with the weak
atomic topology (see \cite{EK4} and in this present work in
(\ref{ag41c})-(\ref{EKcriterion}) below for details on this topology)
which forms a Polish space.

We will need Poisson random measures in four variables, the actual
time called $s \in [0,t]$, the location parameter $a \in [0,1]$,
the mutation and immigration potential called $u \in [0,\infty)$ and
the path of an excursion called $w \in W_0$. This allows us to obtain
the droplet dynamic $(\gimel^m_t)_{t \in [0,\infty)}$ and its
properties:

\beP{CSB}{(A continuous atomic-measure-valued Markov process)}

Let $\gimel^\ast_0(t) =\sum_i y_i(t)\delta_{a_i}\in \mathcal{M}_a([0,1])$
where the $y_i(t)$ are independent solutions of the SDE (\ref{onedim}) (describing
the initial droplet and its evolution) and let
${ N(ds,da,du,dw)}$ be a Poisson random measure on (recall
(\ref{dd30}) for $W_0$) \be{agrev81} [0,\infty)\times [0,1]\times
[0,\infty)\times W_0, \ee with intensity measure
\be{dd32c}
 ds\, da\, du\,\mathbb{Q}(dw),
 \ee
where $\mathbb{Q}$ is the single site excursion law defined in (\ref{dd31a}) in Lemma
\ref{L.Pa2}.

Then the following three properties hold.

(a) The stochastic integral equation for $(\gimel^m_t)_{t \geq 0}$ a process with values
in $\CM_a([0,1])$ is given as
\be{ZL2m}
\gimel^m_t=\gimel^\ast_0 (t) +\int_0^t\int_{[0,1]}\int_0^{q(s,a)}\int_{W_0}w(t-s)\delta_a
{ N(ds,da,du,dw)}, \quad t \geq 0, \ee
where  $\delta_a(\cdot) \in \CM_1 ([0,1])$ and where $q(s,a)$ denotes the
non-negative predictable function \be{dd32d} q(s,a):= (m+ c
\gimel^m_{s-}([0,1])), \ee has a unique continuous
$\mathcal{M}_a([0,1])$-valued solution, which we call \be{sol1a}
(\gimel^{m}_t)_{t \geq 0}. \ee

(b)  $(\gimel^{m}_t)_{t \geq 0}$ is a $\CM_a([0,1])$-valued strong
Markov process.

 (c) The process $({\gimel^m_t})_{t \geq 0}$ has the following properties:
\begin{itemize}
\item the mass of each atom observed from the time of its creation
follows an excursion from zero generated from the excursion law
$\mathbb{Q}$ (see (\ref{dd31a})), \item new excursions are produced
at time $t$ at rate \be{agrev3} m+ c\gimel^m_t([0,1]), \ee \item
each new excursion produces an atom located at a point  $a \in
[0,1]$ chosen according to the uniform distribution on $[0,1]$,
\item at each $t$ for $\ve >0$ there are at most finitely many
atoms of size $\geq \ve$, \item $t\to \gimel^m_t([0,1])$ is a.s.
continuous. $\qquad \square$
\end{itemize}
\end{proposition}

\begin{remark}
The process $(\gimel^m_t)_{t \geq 0}$ can be viewed as a continuous state
analogue of the Crump-Mode-Jagers branching process with
immigration (see Subsubsection \ref{sss.dualcfr}, Step 2 for a review and more information
on this type of processes). We shall see in (\ref{CSB3}) that the total mass grows
exponentially as in a supercritical branching process. $\qquad \square$
\end{remark}

\begin{remark}
A necessary and sufficient condition for extinction
of the analogue of $\gimel^{\ast}_t$ for a general class of
one-dimensional diffusions was obtained by M. Hutzenthaler
(\cite{Hu}). In the Fisher-Wright case with $c>0,\;s>0 $ the fact
that the probability of non-extinction is  non-zero follows from
the Proposition \ref{CSB-longtime}. $\qquad \square$
\end{remark}

\begin{proof}
(a) and (b)  The existence and uniqueness of the solution to
(\ref{ZL2m}) and the strong Markov property follows as in \cite{D-Li}
and \cite {FL}.

(c) follows directly from the construction via the Poisson
measure $N$ given in (\ref{ZL2m}) .

\end{proof}

\begin{remark}
We can enrich the process $(\gimel^m_t)_{t \geq 0}$ to include the
genealogical information, namely which mass results from which
mutation. For that purpose we have in addition to the location record the birth
times of atoms due to migration (successful colonization). This
means we have to split the rate at which excursions are created
into mutation at the site and  immigration from other sites. At
time 0 we start with \be{xg1} \wh \gimel^\ast_0(t) \ee describing
the initial atoms we consider to be present at time 0.

The genealogical enrichment denoted
\be{agrev3b}
(\wh\gimel^m_t)_{t \geq 0}
\ee
is a measure-valued process on a richer set $E$ and is obtained as follows.

Let
\be{agrev3c}
E= (\left(\mathbb{R}_+\times[0,1])\cup\O\right)^{\N}
\ee
and for $y \in E$ define $\tau(y)=min\{n:y(n)=\O\}-1$.

 Let  ${ N(ds,da,du,dw)}$ be a Poisson random measure on
 $[0,\infty)\times [0,1]\times [0,\infty)\times W_0$
with intensity measure ${ ds\, da\, du\,\mathbb{Q}(dw)}$ where
$\mathbb{Q}$ is the single site excursion law. Then the following stochastic
integral equation  has a unique continuous solution, $\wh\gimel^{m}_t$:
\be{ZL2b}
\wh\gimel^m_t= \wh \gimel^\ast_0(t)+
\int_0^t\int_{E}\int_{[0,1]}\int_0^{q(s,y,a)}\int_{W_0}w(t-s)
\delta_{y \diamond (s,a)} { N(ds,da,du,dw)},\ee
where $(s,a)$ is short for $((s,a),\O,\O,\cdots)$,
\be{agrev2b}
q(s,y,a)=m \mbox{ if }
\tau(y)=0,\quad q(s,y,a):= c \wh\gimel^m_{s-}(y),\;\tau(y)\ne 0, \ee
and
\be{agrev2c}\begin{array}{l}
\diamond:E\times E \to E ,\;y_1\diamond
y_2=(y^\prime_1,y^\prime_2, \underline{\O})\quad \mbox{ (concatenation) with },\\
 y^\prime= (y^\prime(1), \cdots, y^\prime(\tau(y^\prime))) \quad ,
 \underline{\O}=(\O,\O, \cdots).
\end{array}\ee

Here $\wh\gimel^m_t$ is a measure-valued process on $E$ and $ y
\diamond (s,a)$ denotes the offspring of $y$ with birth time $s$
and location $a\in [0,1]$. Moreover,
$\wh\gimel^m_t(\{y:\tau(y)=\infty\})=0$ a.s..

\end{remark}

{\bf Some topological facts concerning atomic measures.}
We finally recall the definition and some facts on the {\em weak atomic topology with
metric} $\rho_a$ on the space of finite atomic measures $\CM_a
([0,1])$ due to Ethier and Kurtz. Recall that we have the topology
on $\CM([0,1])$ induced by the weak topology which is induced by the Prohorov
metric $\rho$. Next choose a function $\Psi$,
where  $\Psi:[0,\infty)\to[0,1]$ is continuous, nonincreasing
and $\Psi(0)=1,\;\Psi(1)=0$.
Then for $\nu, \mu \in \CM_a ([0,1])$ one defines
\bea{ag41a2}
&&
\rho_a (\nu, \mu):= \rho(\nu,\mu)\\
&&
  +\sup_{0<\ve\leq 1}\left|\int_{[0,1]}\int_{[0,1]}\Psi(\frac{|x-y|}{\ve}\wedge 1)\mu(dx)\mu(dy)
  -\int_{[0,1]}\int_{[0,1]}\Psi(\frac{|x-y|}{\ve}\wedge 1)\nu(dx)\nu(dy)\right|.\nonumber
\eea
We refer to $\rho_a$ as the {\em Ethier-Kurtz metric}.

The space $(\CM_a ([0,1]), \rho_a)$ is a
Polish space and the topology, in other words convergence,
does not depend on the choice of $\Psi$ (the geometry of the space of course does).

The following lemma collects what we need on the relation between
the weak and the weak atomic topologies, which are different.

\beL{L.EK} {(Weak atomic topology and weak topology)}

(a) A sequence of
random finite atomic measures on $[0,1]$, $(\mu_n)_{n \in \N}$, converges to
$\mu$ in the weak atomic topology if and only if
\be{ag41c}
\mu_n \mbox{ converges weakly to } \mu
\ee
and
\be{ag41d}
\mu^\ast_n([0,1]) \ntoo \mu^\ast([0,1]), \mbox{ where} \quad
\mu^\ast ([0,1]) =\suml_{x \in [0,1]} \mu(\{x\})^2\delta_x.
\ee
(b) If the following three properties hold for
$\mu_n, \mu \in \CM_a ([0,1])$ (here $\La$ denotes weak convergence as
$n \to \infty$):
\be{ag41e}
\mu_n\Rightarrow \mu, \mbox{ the ordered atom sizes converge
and the set of atom locations converges},
\ee
then:
\be{ag41f}
\rho_a(\mu_n,\mu)\to 0.
\ee
(c) A continuous $\mathcal{M}_f$-valued process with a.s. continuous
(in the weak topology) sample paths of the form $\suml_i a_i (t) \delta_{x_i(t)}$
such that $\sum a^2_i(t)$ is a.s. continuous, has also sample
paths a.s. in $C ([0,\infty), (\CM_f ([0,1], \rho_a))$, where $\rho_a$
is the Ethier-Kurtz metric.\\
(d) Consider a sequence $\{Z_N, N \in \N\}$ of atomic measure-valued
processes with cadlag paths in the weak atomic topology,
so that $Z_N$ has the form
\be{ag41f2}
Z_N (t)= \suml_i a_{N,i}(t) \delta_{x_{N,i}(t)}
\ee
for suitable functions (of $t$) $a_{N,i}, x_{N,i}$. Assume furthermore that
\be{xg2}
\{\CL (Z_N), N \in
\N\}$ in $\CP (D ([0,\infty), (\CM_f ([0,1], \rho))) \mbox{ is
relatively compact}. \ee

Then the compact containment condition will hold also in
$(\mathcal{M}_a([0,1]),\rho_a)$, if and only if for each $T>0$ and
$\delta>0$ there exists  $\ve>0$ such that \be{EKcriterion}
\inf_NP\left[\sup_{t\leq
T}\left(\int_{[0,1]}\int_{[0,1]}\Psi(\frac{|x-y|}{\ve} \wedge
1)Z_N(t,dx)Z_N(t,dy)-\sum_i a_{N,i}^2(t)\right)\leq
\delta\right]\geq 1-\delta.\ee Here $\Psi$ is as in
(\ref{ag41a2}).
\end{lemma}
\begin{proof}
This follows from work by Ethier and Kurtz \cite{EK4}, namely (a)
- Lemma 2.2, (b) -Lemma 2.5, (c) - Lemma 2.11, (d) - Remark 2.13.
\end{proof}

\subsubsection{Emergence time, droplet formation, fixation dynamic
for $N$-interacting sites:\\ Statement of results }
\label{sss.Ninteract}

We now have  the ingredients and the background to continue the
analysis and state all the results of the finite population model with two types
$\mathbb{K}=\{1,2\}$ with fitness ${0}$ and $1$, respectively,
$s>0$ and $d>0$   and with $N$ exchangeable sites and with $c>0$.
We start with only type $1$ present initially.
The exposition has three parts, the {\em emergence}, the {\em preemergence } (droplet formation) and the {\em fixation}.

{\bf Part 1: Emergence}

The goal is to describe in mathematical precise form the initial
formation of germs for the expression of the fitter type 2 by
mutation which subsequently expand which then leads finally to the global emergence of the
fitter type, and over the period of order $O(1)$ the increase of the mass of this type
continues until it takes over almost the entire population (fixation). We state the results on
this scenario in four main propositions.

First on emergence are Proposition \ref{P6.1b} which
shows that (global) emergence occurs at times of order $const\cdot
\log N$ and which identifies the constant and then Proposition
\ref{P.measure} and Proposition \ref{P.chaos}  which identify the
limiting dynamics in $t$ after global emergence of type 2, meaning
we study the system observed  in times $const \cdot \log N+t$ with $t \in \R$. The
Proposition \ref{P.DF} explains the emergence behaviour by
describing the very early formation of droplets of ``type-2
colonies''.

In addition to
the four main statements mentioned above the exact
properties in the early stage of droplet formation are given in Proposition \ref{P.lt} and
Proposition \ref{CSB-longtime}. Namely droplet formation  occurs in a random manner in the
very beginning followed by a deterministic expansion of the
droplet size leading to emergence on a global level.

A remarkable fact we state after the part on fixation in Proposition \ref{P.DF2}, which is establishing
a close connection between exit behaviour in short scale and entrance
behaviour in the same scale but placed much later such that both time intervals
are separated by a long time stretch.

We begin with the emergence times.
\bi

\beP{P6.1b}{(Macroscopic emergence and fixation times)}

\noi{\bf (a)} {(Emergence-time)}

There exists a constant $\alpha$ with:
\be{ang3}
0<\alpha<s, \ee
such that if
$\; T_N = \frac{1}{\alpha}\log N$, then for $t \in \R$ and asymptotically as
$N \to \infty$ type 2 is present at times $T_N+t$, i.e. there exists an $\ve >0$
such that
\be{Y12b-}
  \liminf_{N \to \infty}  P[ x^N_2 (T_N+t)>\ve]>0,  \ee
and type 2 is not present earlier, namely for   $1>\ve >0$:
\be{Y12b}
 \lim_{t\to -\infty} \limsupl_{N \to \infty}
 [P( x^N_1 (T_N+t)<1-\ve]=0.  \ee

\noi{\bf (b)} (Fixation time)

After emergence the fixation occurs in times $O(1)$ as $N \to \infty$,
i.e. for any $\ve>0$  \be{Y12c}
\lim_{t\to\infty}\limsup_{N\to\infty}P[x^N_1(T_N+t)>\ve]=0.\ee

{\bf (c)} The constant $\alpha$ can be characterized as the Malthusian parameter for a
Crump-Mode-Jagers branching process denoted  $(\wt K_t)_{t \geq
0}$ which is explicitly defined below in (\ref{ang1ab1}).
The constant $\alpha$ can alternatively be introduced as exponential growth rate
for the limiting droplet growth dynamic $\gimel^m_t$, see (\ref{CSB2}) or in terms of the
excursion measure of a diffusion as specified below in (\ref{dd34}).
$\qquad \square$
\end{proposition}

\beC{C.emfix}{(Emergence and fixation times of spatial density)}

The relations (\ref{Y12b-}), (\ref{Y12b}) and (\ref{Y12c})
hold for $\bar x^N_2$ respectively $\bar x^N_1$ as well.
$\qquad \square$
\end{corollary}

{\bf Part 2: Droplet formation}

This raises the question how the global emergence of type 2
actually came about and what is the role of $\alpha$ in the
forward dynamics (instead of the view back from emergence)
and why do we have a random element in the
emergence. We will now demonstrate that rare mutation and
subsequent selection with the help of migration produce  in times
$O(1)$ a cloud of sites where the type-2 is already manifest with
substantial mass at a given time and this cloud then starts
growing as time increases. This growth of total type-2 mass is exponential with a
random factor arising as usual in the very beginning (i.e. at
times growing arbitrarily slow with $N$) and hence in particular
the growth process is given by a randomly shifted exponential. We
call the growing cloud of type-2 sites a droplet (which in fact in
the euclidian geographic space is literally accurate).
The droplet will be described using the atomic measure
$\gimel^{N,m}_t$ introduced in (\ref{dd29}) and by the
Palm measure explained in (\ref{Y12d}).

For the purpose of making our three regime scenario precise we
first need to investigate how for a finite time horizon the site
of a typical type-2 mass looks like in the limit $N \to \infty$
and if we can show that such a site exhibits nontrivial type-2
mass. Then we have to see how fast the number of such sites grows
and reaches size $\ve N$ at time $\alpha^{-1} \log N + t(\ve)$.

One might expect that a growth of the droplet at exponential rate
$\alpha$ would be the appropriate scenario.  This scenario will
also show that some randomness created initially, i.e. times
$O(1)$ as $N \to \infty$, remains in the system up to fixation of
type 2. This randomness arises since up to some finite time $T$
there will be among the $N$-sites as $N \to \infty, O(1)$ such
sites where the mass of type 2 exceeds some level $\delta
>0$. In the limit this will result in a compound Poisson number of
germs for expansion. Only descendants of this early mass will make
up at much later times the bulk of the type-2 mass, since we have
exponential growth.

Therefore we study
\begin{itemize}
\item the configuration of type-2 mass on the sparse set of sites
and its limiting law as $N \to \infty$ for a finite time horizon,
i.e. $\gimel^{N,m}_t$ and $\wh \mu^N_t$ for $t \in [0,T]$,
\item
the configuration in a typical type-2 site as time and $N$ tends
slowly to infinity, i.e. $\wh \mu^N_{t_N}$, with $t_N \uparrow \infty$
but $t_N= o(\log N)$, \item the growth of the droplet as time goes
to infinity but slower than the time needed for global emergence,
i.e. $\gimel^{N,m}_{t_N}$ and $\gimel^{N,m}_{t_N} ([0,1])$, with
$t_N$ as above.
\end{itemize}

To verify our scenario we want to show that $\wh \mu^N_{t_N}$ converges
for $t_N\uparrow \infty$ to a limit, say $\wh \mu^\infty_\infty$
which is nontrivial, i.e. $\wh \mu^\infty_\infty ((0,1))=1$. Then
we want to prove that the number of sites looking like $\wh
\mu^\infty_\infty$-realisations grows like $\CW^\ast e^{\alpha^\ast
t_N}$ for $N\to \infty$, with $t_N \to \infty, t_N= o(\log N)$ for
some non-degenerate random variable $\CW^\ast$ and a number $\alpha^\ast \in (0,\infty)$.
And finally we have to identify $\alpha^\ast$ as $\alpha$.

To verify that, we have to
establish either that for some suitable finite and not identically zero
random variable $\CW^\ast$ and number $\alpha^\ast$ we have
\be{Y12d1} \wh
x^N_2(t_N):=\suml^N_{i=1} x^N_2 (i,t_N) \sim \CW^\ast e^{\alpha^\ast t_N}
\ee
as $N\to\infty$, or alternatively that we need a spatial window of size
$Ne^{-\alpha^\ast t_N}$ to find such a type-2 colonised site with nontrivial probability
in that window.

In order to explain the origin of the randomness in
$\CW^\ast$ we need the behaviour of the law of
the random variable
\be{angr40} \suml^N_{i=1}x^N_2(i,t)
=\gimel^{N,m}_t([0,1])\;,\quad  \quad t \in [0,T] \ee and the
localization of this total mass on different sites for $N \to
\infty$ in a suitable description, which is given by the full
atomic measure
$(\gimel^{N,m}_{t})_{t \geq 0}$.

In particular if we have established the above scenario,
we can conclude that it
requires  time $(\alpha^\ast)^{-1}  \log N +t$ for some appropriate $t \in
\R$ to reach intensity $\ve$ of type 2 in the whole collection of
sites and then we can conclude that $\alpha^\ast =\alpha$.

We now want to show  three things to describe the droplet growth:
(1) There is a limiting law for the size of the type-2 population
in a typical colonized site. (2) Identify the exponential growth rate of the number of
colonized type-2 sites. (3) Identify the limiting dynamic as $N \to \infty$ of
$(\gimel^{N,m}_t)_{t \geq 0}$ exactly as the dynamic $(\gimel^m_t)_{t \geq 0}$.

Addressing first point (1), (2) (for (3) see below), we
shall prove the following on droplet formation.
We use here the notational convention
\be{Y12d2} a_N = \wh o(bf(N)) \mbox{ if }
\limsup\limits_{N\to\infty} (a_N/f(N))< b. \ee

\beP{P.DF}{(Microscopic emergence and evolution: droplet formation)}

a) The Palm distribution stabilizes, i.e.
\be{pm1} \wh \mu^N_{t_N} \Nto \wh \mu^\infty_\infty
\mbox{ for } t_N \uparrow \infty \mbox{ with } t_N=\wh o(\alpha^{-1} \log N)
\ee
and the limit has the property
\be{pm2}
\wh \mu^\infty_\infty ((0,1))=1.
\ee
The law $\wh \mu^\infty_\infty$ will be identified in terms of the excursion
measure in (\ref{elam2b}).

b) The total type-2 mass $\gimel^{N,m}_{t_N} ([0,1])$ grows at exponential rate
$\alpha^\ast$, i.e. (\ref{Y12d1}) holds and
\be{pm3}
\alpha^\ast=\alpha, \mbox{ with $\alpha$ from Proposition \ref{P6.1b}}.
\ee
Furthermore we have a random growth factor: \be{rgf1} \CL \left[ \wh x^N_2(t_N)
e^{-\alpha t_N}\right] \Nto \CL [ \CW^\ast], \quad \mbox{ for }
t_N \uparrow \infty \mbox{ with } t_N= \wh o(\alpha^{-1} \log N).
\ee The random variable $\CW^\ast$ is non-degenerate. $\qquad
\square$
\end{proposition}

\begin{remark}
This result is obtained by showing that for sufficiently large
$t_1$ and $t_2$:
\be{agrev84}
 \lim_{N\to\infty} P\left( |e^{-\alpha t_1}\wh
x^N_2(t_1)-e^{\alpha t_2}\wh x^N_2(t_2)|>\ve\right)<\ve.\ee
\end{remark}

\begin{remark}
The mutant population making up the total population at time
$\alpha^{-1} \log N+t$ of emergence dates back to a rare mutant
ancestor which appeared at a time in $[0,T(\ve)]$ with probability
at least $1-\ve$ and $T(\ve) < \infty$ for every $\ve>0$.
\end{remark}

\begin{remark}
Growth in the spatial Fisher-Wright and the spatial branching model
are different in nature.
In the branching model we have growth of the total type-2 mass at
exponential rate $s$, which is due to essentially a random number
of families growing at that rate and hence with few sites
with a large population. In the Fisher-Wright case we get an
exponentially (rate $\alpha <s$) growing number of sites with a
type-2 population of at least $\ve$. Macroscopic emergence occurs
once this growing droplet has volume $const \cdot N$.
\end{remark}

In order to understand the  structure of the quantity $\CW^\ast$
arising from (\ref{rgf1}), we have to investigate the behaviour
also in the earlier  time horizon
$t \in [0,T_0]$ of $\gimel^{N,m}_t$ as $N \to \infty$.
And addressing our point (3) from above (\ref{Y12d2}) we obtain the following convergence result.

\beP{P.lt}{(Limiting droplet dynamic)}

As $N\to\infty$ \be{agrev2}
\CL [(\gimel^{N,m}_t)_{t \geq 0}]
\Nto \CL[({\gimel^m_t})_{t \geq 0}]
\ee in the sense of convergence of
continuous $\mathcal{M}_a([0,1])$-valued processes where
$\mathcal{M}_a([0,1])$ is equipped with the weak atomic topology.
\hfill$\qquad \square$
\end{proposition}

In the growth behaviour of the limit dynamics as $N \to \infty$
we obtained above, we recover the quantity $\alpha^\ast$ and $\CW^\ast$ in
the longtime behaviour of the limit dynamic $(\gimel^m_t)_{t \geq 0}$
in the following proposition.

\beP{CSB-longtime}{(Long-time growth behaviour of $\gimel^m_t$)}

Assume that either $m>0$ or $\gimel^0_0([0,1])>0$.  Then the following
growth behaviour of $\gimel^m$ holds.\\
(a) There exists $\alpha^\ast$ such that the following limit exists
\be{CSB1}
\lim_{t\to\infty} e^{-\alpha^\ast t}E[\gimel^m_t([0,1])] \in (0,\infty).
\ee
We have
\be{CSB1b}
\alpha^\ast =\alpha,
\ee
where $\alpha$ is defined as in Proposition \ref{P6.1b} (see also \ref{ang1ab1}).

(b) Recall (\ref{rgf1}) for the law of $\CW^\ast$. Then:
\be{CSB3} e^{-\alpha t}
\gimel^m_t([0,1]) \tto {\mathcal{W}}^\ast \mbox{ in probability}.
\ee

(c) The growth factor in the exponential is truly random:
\be{CSB2}
0< Var[\lim_{t\to\infty} e^{-\alpha t}\gimel^m_t([0,1])]<\infty. \qquad \square
\ee

\end{proposition}

\begin{remark}
We shall see in the proof that the random variable
$\CW^\ast$
reflects the growth of $\gimel^m_t([0,1])$ in the beginning, as
is the case in a supercritical branching process and hence
$\CE^\ast = \alpha^{-1} \log \CW^\ast$ can be viewed as the random
time shift of that exponential $e^{\alpha t}$ which matches
the total mass of $\gimel^m_t$ for large $t$.

Given $\gimel^m_t$, let  $\nu_t(a,b],\; 0<a<b\leq 1$ denote the number of atoms in the interval $(a,b]$. Then the analogue of the size distribution for the discrete CMJ process is given by
\be{add48} \frac{\nu_t(\cdot)}{\gimel^m_t([0,1])}\;\in \mathcal{P}((0,1])).  \ee
It is reasonable to expect that there is a limiting stable size distribution in the $t\to\infty$ limit similar to that in  the discrete case but we do not follow-up on this here.

\end{remark}

{\bf Part 3: Fixation}

We now understand the preemergence situation
and the time of emergence.
In order to continue the program to describe the dynamics of
macroscopic fixation, we  consider the limiting
distributions of the empirical measure-valued processes in a
second time scale $\alpha^{-1} \log N+t, \quad t \in \R$. Define
\be{angr9z}
\Xi^{\log,\alpha}_N(t):=\frac{1}{N}
\sum_{i=1}^N \delta_{(x_1^N(i,\frac{\log
N}{\alpha}+t),x_2^N(i,\frac{\log N}{\alpha}+t))},\quad t \in
\R,\quad(\Xi^{\log,\alpha}_N(t)\in \CP(\mathcal{P}(\mathbb{K}))) \ee
and then the two empirical marginals are given as
\be{angr9}
\Xi^{\log,\alpha}_N(t,\ell):=\frac{1}{N} \sum_{i=1}^N
\delta_{x_\ell^N(i,\frac{\log N}{\alpha}+t)}, \; \ell=1,2 \mbox{ and } t
\in \R, \quad(\Xi^{\log, \alpha}_N(t,\ell)\in \mathcal{P}({[0,1]})).
\ee Note that for each $t$ and given $\ell$ the latter is a random
measure on $[0,1]$. Furthermore we have the representation of the
empirical mean of type 2 as follows:
\be{angr10}
\bar x^N_2(\frac{\log
N}{\alpha}+t) =\int_{[0,1]} x\;\Xi^{\log,\alpha}_N(t,2)(dx).\ee

The third main result of this section is on the fixation process in type 2,
saying that $\Xi^{\log,\alpha}_N$
converges as $N\to\infty$ and the limit can be explicitly
identified as a random McKean-Vlasov entrance law starting from time $-\infty$ in type 1.

\beP{P.measure} {(Asymptotic macroscopic fixation process)}

(a) For each $-\infty<t<\infty$ the empirical measures converge
weakly to a random measure:
\be{angr11}
\mathbf{\CL} [\{\Xi^{\log,\alpha}_N(t,\ell)\}] \Nto \mathbf{\CL}[\{ \CL_t(\ell)\}] =P^\ell_t \in
\CP(\CP({[0,1]})), \mbox{ for } \ell =1,2 . \ee
In addition we have path convergence:
\be{angr13}
 w-\liml_{N \to \infty}\mathbb{\CL}[(\Xi^{\log,\alpha}_N(t))_{t \in \R}]=P
 \in \CP [C((-\infty,\infty),\CP(\CP(\mathbb{K})))].
\ee
A realization of $P$ is denoted $(\CL_t)_{t \in \R}$
respectively its marginal processes $(\CL_t(1))_{t \in \R}, (\CL_t(2))_{t \in \R}$.
\sm

(b) The laws  $P^\ell_t, \ell=1,2$ are  non-degenerate probability measures on
the space $\CP({[0,1]})$. In particular
\be{LD5}
E\left[\left(\int_{[0,1]}x\mathcal{L}_t(2)(dx)\right)^2\right]>\left(E[\int_{[0,1]}x\mathcal{L}_t(2)(
dx)]\right)^2. \ee

(c) The process $(\CL_t)_{t \in \R}$ describes the emergence and
fixation dynamics, that is, for $t\in\mathbb{R}$, and $\ve
>0$, \be{LD2} \lim_{t\to-\infty}
\mbox{Prob}[\mathcal{L}_t(2)((\ve,1])>\ve]=0,\ee
\be{LD3}
\lim_{t\to\infty}
\mbox{Prob}[\mathcal{L}_t(2)([1-\ve,1])<1-\ve]=0,\ee
with
\be{LD3b}
\CL_t (2) ((0,1)) > 0\quad , \quad \forall t \in \R, \mbox{ a.s.}.
\ee
(d) For every $t\in\mathbb{R}$, always both type 1 and  type 2 are present:
\be{LD1}
Prob[\mathcal{L}_t(2)(\{0\})=0)]=1. \ee

(e) The limiting dynamic in (\ref{angr13}) is identified as follows:

The probability measure $P$ in (\ref{angr13}) is such
that the canonical process is a  random solution  (recall
Definition (\ref{D.ranMKV})) and entrance law from time $-\infty$ to the McKean-Vlasov equation
(\ref{angr7b}). \sm

(f) The limiting dynamic in (\ref{angr13}) satisfies with $\alpha$ as in (\ref{ang3}):
\be{angr14a2}
\CL [e^{\alpha |t|}\int_{[0,1]} x \CL_t(2)
(dx)]\Rightarrow \CL[^\ast \CW] \mbox{ as } t \to -\infty,
\ee
 and we explicitly identify the
random element generating $P$ in (\ref{angr13}), namely $P$ arises
from random shift of a deterministic path:
\be{angr14c}
P = \CL [\tau_{^\ast\CE} \CL^\ast]
 \quad , \quad ^\ast\CE =(\log ^\ast\CW)/\alpha, \quad
 \tau_r \mbox{ is the time-shift of path by } r,
\ee
where $\CL^\ast$ is the unique and deterministic entrance law of the McKean-Vlasov
equation (\ref{angr7b}) satisfying (\ref{Ex2b}) and with projection $\wt \CL_t(2)$ on the type 2
coordinate satisfying:
\be{angr14d} e^{\alpha |t|} \int_{[0,1]} x \wt
\CL_t(2)(dx) \longrightarrow 1,  \mbox{ as } t \to -\infty.
 \ee

The random variable $^\ast\CW$ satisfies \be{angr14a}
0< {^\ast{\CW}} < \infty \mbox{ a.s.},\quad E[^\ast\CW]<\infty,\quad
0<\mbox{Var}(^\ast\CW)<\infty. \ee

(g) We have for $s_N \to \infty$ with $s_N =o(\log N)$ the approximation
property for the growth behaviour of the limit dynamic by the finite $N$ model,
namely:
\be{angr14a3}
\CL [e^{\alpha s_N} \bar x^N_2 (\frac{\log N}{\alpha} - s_N)]
\Nto \CL[^\ast\CW]. \qquad \square \ee
\end{proposition}

\beC{C.ld5}{(Scaled total mass process convergence)}

The above implies for the total mass process that:
\be{xg3}
\CL [\{\bar x^N_2(\frac{\log N}{\alpha}+ t)\}_{t\in\mathbb{R}}] \Nto
\CL[\{\int_{[0,1]}x \mathcal {L}_t(2)(dx)\}_{t\in\mathbb{R}}],
\ee
and we have convergence in distribution on path space. $\qquad \square$
\end{corollary}

\begin{remark} The random variable $\CW^\ast$  of (\ref{rgf1}) can be viewed as the
{\em final value} of the microscopic development of the mutant
droplets and $^\ast \CW$ of (\ref{angr14a2}) can be viewed as an
{\em entrance value} of the macroscopic emergence of the mutant
type.
\end{remark}

This result of course immediately raises the question how the macroscopic
{\em entrance} from 0 is related to the {\em final} value in the smaller original time scale.
We have already stated that the two exponential growth rates from these two
directions are the same.
But more is true also the two growth constants have the same law:

\beP{P.DF2} (Relation between microscopic and macroscopic regimes)

The entrance value and the final value are equal in law:
\be{agrev83}
\CL[\CW^\ast]=\CL[^\ast\CW].\qquad \square
\ee
\end{proposition}

\begin{remark}
We note that this result cannot be derived from some soft argument since
we relate two time scales $t$ and $\alpha^{-1} \log N+t$ which separate
as $N \to \infty$ and indeed the limits $N \to \infty$
and $t \to \infty$ or $t \to - \infty$ cannot be interchanged in general.
\end{remark}

\begin{remark}\label{R63}
One could go further studying the relation between final and entrance
configuration and ask how the limiting Palm distribution
$\wh \mu^\infty_\infty$, which describes the typical sparse type-2
sites is related to the limit $\CL_{-\infty}$ as $t \to -\infty$
of the Palm distribution $\wh \CL_t$ of $\CL_t$, i.e.
\be{agrev83b} \wh \CL_t
(2) (dx) = \frac{x \CL_t(2)(dx)}{\int x \CL_t(2)(dx) }, \ee which
describes the typical site ''immediately'' after emergence.

Note that this new transformed measure is now deterministic  in the  $t \to -\infty$
limit, since in the ratio the same random
time shift is used in both numerator and denominator.

Here we have to observe that due to the exponential growth of
type-2 sites the state depends on when they were colonised and
that time point is typically a {\em finite} random time back.
Define for  the process given in (\ref{onedim})
\be{xg3c}
\mu_s^\ast(dx) = \mbox{ probability that  an excursion
of age }s \mbox{  has size in } dx.
\ee
Note that for
any time the total mass $\gimel_t([0,1])$ is finite but there are
countably many non-zero atoms in $\gimel_t$.  At time $t$ let
\be{add49}
\Gamma_t([a,b)) = \mbox{ number of atoms of age in
}[a,b),\quad 0<a<b\leq 1\ee and
then define the corresponding size-biased distribution
on the product space $\R^+ \times [0,1]$ of age and size given by
\be{b80}
\wh \Gamma_t (ds,dx) = \frac{\Gamma_t(ds) x
 \mu^\ast_s (dx)} {\int_0^t\Gamma_t (ds)\int_0^1 x  \mu^\ast_s
(dx)ds}. \ee
 By analogy with
the Crump-Mode-Jagers theory, we would expect that this age and
size distribution $\wh \Gamma_t(ds,dx)$  stabilizes as
$t\to -\infty$ with stable age and size distribution $\wh \Gamma
(ds,dx)$.

Then we would expect to have the relation:

\be{agrev83c} \wh \CL_{-\infty}(2)(dx) = \wh
\Gamma((0,\infty),dx). \ee
\end{remark}

 Corresponding to the limit dynamics of Proposition
\ref{P.measure}  we also have a  {\em stochastic propagation of
chaos result}, which characterizes the asymptotic behaviour of a
tagged sample of sites as follows:

\beP{P.chaos}{(Fixation dynamic for tagged sites)}

There is a unique limiting entrance law for the process of tagged components:
for fixed $L\in\mathbb{N}$, the sequence in $N$ of laws
\be{angr15}
 \CL [(\{x^N_2(i,\frac{\log N}{\alpha}+t))\}_{i=1,\dots,L})_{t\in\mathbb{R}}]
 \ee
 converges as $N \to \infty$ weakly to a law
\be{dg29}
\CL[(X_2(t))_{t \geq 0}] =
\CL [((x_2(1,t), \cdots, x_2 (L,t))_{t \in \R}],
\mbox{ with } X_2(t) \la \underline{0} \mbox{ as } t \to -\infty,
\ee
which is if restricted to $t \geq t_0$ for every $t_0 \in \R$ the weak solution of the SSDE:

 \bea{Y777}
 d x_2 (i,t)&& =
 c\left(\int_{[0,1]} x\;\mathcal{L}_t(2)(dx)  -  x_2 (i,t)\right)dt
 +  s  x_2(i,t) (1-x_2(i,t)) dt\\&&
 + \sqrt{d\cdot  x_2(i,t)(1-x_2(i,t))} \, dw_i (t),\quad i=1,\dots,L, \nonumber
\eea
where $w_1,\dots,w_L$ are independent standard Brownian motions which are also independent of the
process $\{\mathcal{L}_t\}_{t \in \R}$ and the latter is given by (\ref{angr14c}).
The behaviour at $-\infty$ in (\ref{dg29}) and the equation (\ref{Y777}) uniquely determine the law
for a given mean curve.
$\qquad \square $\end{proposition}

\begin{remark}
Note that $x_2 (i,t)$ will hit 0 at random times until $\bar
x_2(t)$ exceeds a certain level depending on the other parameters
in particular $c/d$, only after that time where $\bar x_2(t)$
remains large enough the component remains a.s. occupied with
positive type-2 mass for all future time. Before that time the
path will have zeros at an uncountable set of time points with
Lebesgue measure zero, the set of these zeros will be stochastically
bounded from above. This is a consequence of standard diffusion theory
(see \cite{RW}).
\end{remark}

\section{Duality}
\label{s.dual}

\setcounter{secnum}{\value{section}} \setcounter{equation}{0}

The key tool to obtain our results is {\em duality}
between our process and a particle system.

A duality relation for a Markov process $Z$ with state space $E$ is given as follows.
There is a process $Z^\prime$ with state space $E^\prime$ such that for a
function $H:E \times E^\prime \to \R$, which is bounded and continuous we have
\be{dual0}
E_{Z_0} [H(Z_t, Z^\prime_0)] = E_{Z^\prime_0} [H(Z_0,  Z^\prime_t], t \geq 0.
\ee

The duality we use is a special
case of the duality theory developed in \cite{DGsel} and \cite{DGDuality}.
We begin constructing the dual particle system and then state the duality.

\subsection{The dual process}
\label{ss.dualpr}

We construct now a dual process. This dual process is based on a particle system
$(\eta_t)_{t \geq 0}$, whose dependence on $N$ we suppress in the notation.

We need the following ingredients to define the process $\eta$:

\begin{itemize}
\item $N_t$ is a non-decreasing $\N$-valued process. $N_t$ is the number of {\em individuals} present in the dual process  with $N_0=n$, the number of initially tagged individuals,
and with $N_t-N_0\geq 0$ given by the number of individuals born during the interval $(0,t]$.

\item $\zeta_t =\{1,\cdots,N_t\}$ is an {\em ordered particle system} where the individuals are given an assigned order and the remaining particles are ordered by time of birth.

\item \begin{tabular}{ll} \label{dualcomp}
\hspace{-0.2cm}$\eta_t=$ & $(\zeta_t, \pi_t, \xi_t)$ is a trivariate process consisting of the above $\zeta$ and\\
& $- \pi_t$:  {\em partition} $(\pi^1_t,\cdots,\pi^{|\pi_t|}_t)$ of $
\zeta_t$, i.e. an ordered family of subsets, where the \\&\qquad
{\em index} of a partition element is the smallest element (in the ordering of individuals)\\&\qquad of the partition element.\\
& $- \xi_t: \pi_t \longrightarrow \Omega^{|\pi_t|}_N$,  giving {\em locations} of the
partition elements. Here $\xi_0$ is the vector \\&\qquad of the prescribed space points
where the initially tagged particles sit.
\end{tabular}
\item $\CF_t$ is for given  $\eta_t=( \zeta_t, \pi_t,\xi_t)$ is a function in
$L_\infty (\I^{|\pi_t|})$.
\end{itemize}

Therefore the state of $\eta$ is an element in the set
\be{dg30.1}
\CS = \bigcup\limits^\infty_{m=1} \{1,2, \cdots, m\}
\times \mbox{Part}^<(\{1,\cdots,m\}) \times \{\xi: \mbox{Part}^< (\{1,\cdots,m\})
\la \Omega_N\},
\ee
where Part$^< (A)$ denotes the set of ordered partitions of a set $A \subseteq \N$
and the set $\{1,\cdots, m\}$ is equipped with the natural order.

It is often convenient to order the partition elements according to their
indices
\be{dg30}
\mbox{{\em index} } (\pi^i_t)= min (k \in \N : k \in \pi^i_t)
\ee
and then assign them ordered
\be{l1}
\mbox{ {\em labels} } 1,2,3, \cdots, |\pi_t|.
\ee

\begin{remark} \quad
Note that
\be{3.2b2}
N_t = |\pi^1_t| +\cdots + |\pi_t^{|\pi_t|}|.
\ee
\end{remark}

\begin{remark} \label{R.eta2} \quad
The process $\eta$ can actually be defined starting with countably many individuals located in $\Omega_N$. This is due to the quadratic death rate mechanism at each site which implies that the number of individuals at any site will have jumped down into $\N$ by time $t$ for any $t >0$.
\end{remark}

The dynamics of $\{\eta_t\}$ is that of a pure Markov jump process with the following transition mechanisms which correspond to resampling,
migration and selection in the original model
(in this order):
\begin{itemize}
\item \emph{Coalescence} of partition elements:  any pair of  partition
elements which are at the same site  coalesce at a constant
rate; the resulting enlarged partition element is given the
smaller of the indices of the coalescing elements as index and the
remaining elements are reordered to preserve their original order.
\item \emph{Migration} of partition elements: for $j\in
\{1,2,\dots,|\pi_t|\}$, the partition element $j$ can migrate
after an exponential waiting time to another site in $\Omega_N$ which means we have a jump of
$\xi_t(j)$.
\item \emph{Birth} of new individuals by each partition element
after an exponential waiting time: if there are currently $N_t$ individuals then a newly
born individual is given index $N_{t-}+1$ and it forms a {\em partition
element consisting of one individual}  and this partition element is given a label
$|\pi_{t-}|+1$. Its location mark is the same as that of its parent
individual  at the time of birth,
\item {\em Independence}: all the transitions and waiting times described in the previous
points occur independently of each other.
\end{itemize}

Thus the  process $( \eta_t)_{t \ge 0}$ has the form
\be{3.4}
\eta_t = (\zeta_t, \; \pi_t, (\xi_t(1), \xi_t(2), \cdots, \xi_t(|\pi_t|))),
\ee
where $\xi_t(j) \in\Omega_N,\;j=1,\dots,|\pi_t|,$ and $\pi_t$ is an
ordered partition (ordered tuple of subsets) of the current basic set of the form
$\{1,2,\dots,N_t\}$ where $N_t$ also grows  as a random process,
more precisely as a pure birth process.

We order partition elements by their smallest elements and then
{\em label} them by $1,2,3,\dots$ in increasing order. Every partition
element (i.e. the one with label $i$) has at every time
a location in $\Omega_N$, namely the i-th partition
element has location $\xi_t(i)$. (The interpretation is that we
have $N_t$ individuals which are grouped in $|\pi_t|$-subsets and
each of these subsets has a position in $\Omega_N$).

Furthermore denote by
\be{3.2c}
\pi_t (1), \cdots, \pi_t (|\pi_t|)\ee
the index of the first, second etc. partition element.
In other words the map gives the index of a partition element
(the smallest individual number it contains) of a label
which specifies its current rank in the order.

For our concrete situation we need to specify in addition the
parameters appearing in the above description, this means that we define
the dual as follows.

\beD{D.Could}{( The dual particle system: $\eta_t$)}

(a) The initial state $\eta_0$ is of the form $(\zeta_0, \pi_0,\xi_0)$ with:
\be{dg31}
\zeta_0=\{1,2,\cdots,n\}, \quad \pi_0 =\{\{1\}, \cdots,\{n\}\}
\ee
\be{dg32}
\xi_0 \in (\Omega_N)^n.
\ee
(b)
The evolution of $( \eta_t)$ is defined as follows:
u\begin{itemize}
\item[(i)] each pair of partition elements which
 occupy the same site in $\Omega_N$ \emph{coalesces} during the
 joint occupancy of a location into one partition element after a \emph{rate $d$}
 exponential waiting time,
 \item[(ii)] every partition element performs, independent of the
 other partition elements, a \emph{continuous time random walk} on
 $\Omega_N$ with transition \emph{rates $c \bar a(\cdot, \cdot)$}, with $\bar a(\xi,\xi^\prime)=a(\xi^\prime,\xi)$,
 \item[(iii)] after a rate $s$ exponential waiting time each partition element
 gives \emph{birth} to a new
 particle   which forms a new (single particle)
 partition element at its location, and this new partition element
 is given as label $|\pi_{t-}|+1$ and the new particle the index $N_{t-}+1$.
 \end{itemize}
 All the above exponential waiting times are \emph{independent} of each other.
 $\qquad \square$
 \end{definition}
 Note that this process is well-defined since the total number of
 partition elements is stochastically dominated by an ordinary linear birth
 process which is well-defined for all times. \bi

\subsection{The duality relation}
\label{ss.dualrela}

This gives in our case  (see \cite{DGDuality})  a duality relation between $X^N$ given in (\ref{angr1})and $\eta$ given by (\ref{3.4}).

The duality relation reads in our case for $i_1, \cdots, i_n \in \{1,\cdots,N\}$
all different:
\be{dual1}
E[\prodl^n_{k=1} x^N_{i_k} (t)] = E[exp (-\frac{m}{N} \intl^t_0 \Pi_u^{N,n,1} du)],
\ee
where  $\Pi^{N(n,1)}_u$ denotes the number of dual particles at time $u$ starting with the
initial state (\ref{dg31}), (\ref{dg32}).

We also need the duality relation for the {\em McKean Vlasov process}
$\{(x^\infty_\ell (0,t), x^\infty_\ell(\ast,t)), \quad \ell =1,2\}$
which is given by $x_1=1-x_2$ and

\bea{kvp1}
d x^\infty_2 (0,t) = c (\theta - x^\infty_2 (0,t)) dt& + s(x^\infty_2 (0,t)
(1 - x^\infty_2(0,t))dt\\
& + \sqrt{d \cdot x^\infty_2 (0,t) (1-x^\infty_2 (0,t)} dw (t) \nonumber
\eea

\be{kvp1a}
x^\infty_2 (\ast,t) \equiv \theta.
\ee

The dual process for the McKean-Vlasov processes arises if we replace
\be{kvp2}
\{1, \cdots,N\} \mbox{ by } \{0,\ast\}
\ee
such that all rates of the dual particle on $\ast$ are put equal to zero and
the jump at rate $c$ from 0 is always to $\{\ast\}$. Otherwise the system
remains.

The duality relation now reads:
\be{dvp3}
E[(x^\infty_2(t,0))^k] = E[\theta^{|\Pi^{\infty,k}_t|}].
\ee

\section{Proof of results}
\label{s.proofresult}

\setcounter{secnum}{\value{section}} \setcounter{equation}{0}

We begin in Section \ref{s.dual} to state the duality theory and continue in
Subsubsections \ref{ss.singles}-\ref{ss.chaos} with the proofs. We start in Subsection \ref{}
 with a warm-up exhibiting some of the
methods in simpler cases. (This is not needed for the rigorous argument). Then we continue
with the formal proofs, namely proving
\begin{itemize}
\item
in Subsection \ref{ss.proofemfix} assertions related to the large time scale,
in particular Proposition \ref{P.MVentrlaw}, Proposition \ref{P6.1b} and
Proposition \ref{P.measure},
\item
in Subsection \ref{ss.bounds} assertions concerning the small time scale
and droplet formations, in particular Proposition \ref{CSB},
Proposition \ref{P.DF} and Proposition \ref{P.lt},
\item
in Subsections \ref{ss.relationsW} and \ref{ss.highermom} assertions concerning the
connection between the two time scales, in particular Proposition \ref{P.DF2},
\item
in Subsections \ref{ss.chaos}-\ref{ss.non-critical} additional
assertions where different methods are used,
in particular Proposition \ref{P.chaos} in Subsection \ref{ss.chaos}
and extensions in Subsection \ref{ss.non-critical}.
\end{itemize}

\subsection{A warm-up: The case of a single site model }
\label{ss.singles}

Before we start proving all the statements made in Subsection
\ref{ss.6results} we look at some simpler cases and already give
in the arguments a flavour of the methods used later on.
This subsection may be skipped from a purely logical point of view.
These simpler systems are single site models.

In the {\em absence of migration}, $c=0$,  the individual sites
evolve independently. In this case it suffices to study a system
\be{dg4} X(t) = (x_1(t), x_2 (t)) \ee at a single isolated site,
for example use
$x_1(t)=x_1(1,t), x_2(t) = x_2 (1,t)$.
Recall (\ref{angr1}-\ref{angr3b}).

We begin by showing in this case that emergence in $\log N$ time
scale does occur in the infinite population, that is,
deterministic case, i.e. $d=0$ but {\em does not} occur in the
finite population case, that is, random case, where $d>0$. This
demonstrates  the effect of the bounded component (mass
restriction) and the important role played by space and migration.

The proof of this result gives us the opportunity to
introduce some basic ideas and methods in a simple setting. We
consider here successively the deterministic and stochastic case
and collect the proofs with the dual calculations in the third
Subsubsection \ref{sss.proofsdual}.

\subsubsection {Emergence of rare mutants in the single site deterministic
population model}\label{raremutinf}

We first consider the  deterministic, i.e. $d=0$ (often called
infinite population) model without migration, i.e.  $c=0$. In the case of small
mutation rate, $\frac{m}{N}$, we will show with two different
methods that asymptotically as $N\to \infty$ emergence occurs at
times $\frac{\log N}{s}+t, t \in \R$ with fixation as $t \to
\infty$: \beL{L.determ} {} Assume that the mutation rate is
$\frac{m}{N}$, and  $d=0,\; s>0,\; c=0$. Then we have:
\be{grev75b}
\liml_{N\to \infty} x_2\left(\frac{\log N}{s}+t\right)
=1- \frac{s}{s + m e^{st}} \left\{
\begin{array}{ll}
    \to 1, & \mbox{as} \quad t \to\infty \\
    \to 0 , & \mbox{as}\quad t\to -\infty. \\
\end{array}
\right.\qquad \square
\ee
\end{lemma}

\begin{remark} \quad This means that in the deterministic case we have emergence
and fixation in the ``critical'' time scale $(s^{-1} log~N)+t$ as
in the branching approximation.
\end{remark}

\subsubsection {Emergence of rare mutants in the single site stochastic
population model} \label{raremutfin}

The above analysis can easily be extended to the stochastic case
(often called finite population), that is, $d>0$, but now we get a
qualitatively different conclusion. We again take mutation rate
$\frac{m}{N}$. Namely we observe that
\beL{L.emernonspat0} {}

Let $d>0$ and $c=0$, then
\be{grev76a}
\CL[x^N_2(a_N)] \Nto \delta_0 \mbox{ whenever } a_N/N \to 0 \mbox{ as }
N \to \infty.
\ee
Therefore emergence and takeover does not occur in time scale
$\log N$ if there is no migration, that is,   $c=0$. $\qquad
\square$
\end{lemma}

\begin{remark}
This means that the ``finite population sampling'', i.e. $d>0$,
seriously slows down the emergence of rare advantageous mutants.
This agrees with the corresponding result in the branching case
where the fluctuation term is responsible for the emergence on
$O(N)$ scale if $c=0$. This effect is removed as soon as we have
migration, which is therefore the key mechanism for the emergence
if $d>0$.
\end{remark}

\subsubsection{Proofs using the dual}
\label{sss.proofsdual}

{\bf Proof of Lemma \ref{L.determ}}. The Lemma \ref{L.determ} is
now proved in two ways with and without duality. We consider first
an auxiliary object $x_2(t)$ for which if we consider the case
where the parameter $m$ is replaced by $\frac{m}{N}$ we get our
$x^N_2$.\sm

\emph{(i) Analytic approach} \quad The quantity
$x_2(t)$ satisfies the ODE
\be{grev60}
\frac{dx_2}{dt}=m(1-x_2)+sx_2(1-x_2);\quad x(0)=0\ee
with solution
\be{grev61}
 x_2(t)=\frac {m(e^{(m+s)t}-1)}   {s+me^{(m+s)t}}.\ee

We can now study this explicit  solution to the ODE  for our
choices of the parameters and the appropriate time scales where
$t$ is replaced by $\frac{1}{s}\log N+t$. (See
point (iii) below).

\emph{(ii) Duality approach.} \quad
As a prelude to the analysis of the stochastic case we show
how the formula (\ref{grev61}) can be obtained using the \emph{dual
representation}.

We assume that $Y_0= \delta_{1}$, $s>0$ and we choose
$f_0=1_{\{1\}}$. The mutation part of the dual consists of the
simple transition corresponding to a mutation from type 1 to type
2. This means  for the dual that we have
jumps $1_{\{2\}} \to 1$, $1_{\{1\}}\to 0$ both at rate $m$. \bi

Let $\{p_k(t)\}_{k\in\mathbb{Z}_+}$ be the distribution of a pure
birth process $N_t$ starting with one individual (and each individual
independently gives birth at rate $s$ at time $t$) and let
\be{grev63}
 \tau = \mbox{ time of first mutation jump}.\ee

Let $\tilde p_k(t)$ denote the probability that there were $(k-1)$ births in $[0,t]$
and that no mutation jump has occurred.

Then the random state of $^{FK}\CF^{+}_t$ at time $t$ is given by:
\be{grev62}
 ^{FK}\CF^{+}_t= (1_{\{1\}})^{\otimes k}1_{\{\tau >t\}},
\ee with probability $\tilde p_{k}(t)$ and 0 with probability
$1-\sum _{k=1}^{\infty}\tilde p_k(t)$.
Then
\be{grev64} x_1(t) =
E[<X(t), 1_{\{1\}}> ]=E[\int\, ^{FK}\CF^{+}_0 d X_t]
  = E[\int \; ^{FK}\CF^{+}_td X_0^{\otimes |\pi_t|}] = P[\tau >t].
  \ee
We have: \be{grev65}
 P(\tau >t)=\sum_{k=1}^\infty \tilde p_k(t).\ee

Noting that the rate of a jump to zero is $mn$ when there are $n$
factors it follows that the collection $\{\tilde p_n(t), n
\in \N\}$ satisfies the system of ODE's:
\bea{grev66}
&& \frac{d\tilde p_{n}(t)}{dt}   =(n-1)s\ \tilde p_{n-1}(t)-n(s+m)\tilde p_{n}(t)\\
&& \tilde p_{1}(0)    =1. \nonumber \eea To determine the $\{\tilde
p_k(t), \quad k \in \N\}$ we introduce  the Laplace transform:
 \be{grev67}
R(\theta,t)    =\sum_{n=1}^{\infty}e^{n\theta}\tilde p_{n}(t). \ee
Then using (\ref{grev66})  it can be verified that
$R(\cdot,\cdot)$ satisfies \bea{grev67b}
R(\theta,0)  &&  =e^{\theta}\\
\frac{dR(\theta,t)}{dt}
  &&=\sum_{n=1}[se^{n\theta}(n-1) \tilde p_{n-1}(t)-(s+m)\sum ne^{n\theta}\tilde p_{n}(t)] \nonumber\\
  && =se^{\theta}\frac{\partial}{\partial\theta}R(\theta,t)
       -(s+m)\frac{\partial}{\partial\theta}R(\theta,t) \nonumber \\
  && =(se^{\theta}-(s+m))\frac{\partial}{\partial\theta}R(\theta,t).\nonumber
\eea

The solution to this PDE is obtained by the method of characteristics as
follows. In symbolic notation: \be{grev68}
\frac{dt}{1}=\frac{d\theta}{-(se^{\theta}-(s+m))}=\frac{dR}{0}, \ee
so that we obtain the characteristic curve \be{grev69}
e^{(s+m)t}(s-(s+m)e^{-\theta}) =\text{constant}. \ee Hence the
general solution is \be{grev70}
 R(\theta,t)   =\Psi(e^{(s+m)t}(s-(s+m)e^{-\theta})),\ee
 where $\Psi$ is an arbitrary function. Using the initial condition
 \be{grev71}
e^{\theta}    =\Psi(s-(s+m)e^{-\theta})\ee we get that \be{grev72}
\Psi(u)    =\frac{s+m}{s-u}. \ee It follows that
 \be{grev73}
R(\theta,t)   =\frac{s+m}{s-e^{(s+m)t}(s-(s+m)e^{-\theta})}. \ee
Therefore
\be{grev74}
 1-x_2(t) =\sum_{n=1}^{\infty} \tilde p_n(t) =R(0,t)= \frac{s+m}{s+me^{(s+m)t}}.\ee

\emph{(iii) Large $N$ asymptotic.} \quad We now consider the
emergence of a rare mutant with $m$ replaced by $\frac{m}{N}$ and
consider the time scale $(C\log N)+t, \quad t \in \R$. We get \be{grev75}
1-x_2(C\log N +t)
=\frac{s+\frac{m}{N}}{s+\frac{m}{N}e^{(s+\frac{m}{N})(C \log
N+t)}} = \frac{s+\frac{m}{N}}{s+\frac{m}{N}e^{(s+\frac{m}{N})t}
N^{C(s+\frac{m}{N})}}, \ee so that $1-x_2(C\log N+t)\to
1\;\mathrm{or}\;0$ for all $t \in \R$ depending on
$C<\frac{1}{s}$ or $C>\frac{1}{s}$. In other words the
rare mutant emerges and takes over in time scale $s^{-1} \log N$
as $N\to\infty$. We shall derive more precise results later in
Subsection \ref{ss.bounds} for example in (\ref{grev75c}).
Finally, for $C=\frac{1}{s}$, \be{grev75aa}
 1-x_2(C\log N
+t)=\frac{s+\frac{m}{N}}{s+me^{(s+\frac{m}{N})t}N^{\frac{m}{sN}}}\to\frac{s}{s+me^{st}},\quad
\text{as  }N\to\infty.\ee \bi

\par\noindent{\bf Proof of Lemma \ref{L.emernonspat0}\ }

We consider the dual process $(\eta^N_t, FK \CF^+_t)_{t \geq 0}$ and define
\be{grev76}
 \tau_N =\inf\{t:^{FK}\CF^{+}_t \equiv 0\}.\ee
For the deterministic case, $d=0$, the random time $\tau_N$ was of
order $s^{-1} \log N$, but this now changes to a much larger order
of magnitude. This is made precise in the following Lemma that
together with the dual identity completes the proof of Lemma
\ref{L.emernonspat0}.

\beL{L.emernonspat}{} Let $d>0$ and $c= 0$.

If  $\lim\frac{a_N}{N}\to 0$, then

\be{grev78}
  P(\tau_N < a_N ) \to 0
 \mbox{ as } N \to \infty \mbox{ for all } t>0.\qquad \square
 \ee
\end{lemma}

\begin{proof}{\bf of Lemma \ref{L.emernonspat}}

Start the dual $\CF^{+}$ in $1_{\{1\}}$. In this case up to the
time of the first mutation the number of factors $1_{\{1\}}$ in
(\ref{grev62}) is equal to or greater than $1$. Therefore for each
$N$, $P(\tau_N <\infty)=1$.

This implies that for every {\em fixed} $N$
\be{grev77}
  x^N_2(t)\to 1,\;\text{a.s. as }\;t\to \infty. \ee

On the other hand the number of factors  is stochastically bounded  due to the coalescence (which
 occurs at quadratic rate as soon as $d>0$).
To verify this we can calculate the expected number of factors and
show that it remains bounded by obtaining an upper bound given by
the solution of an ODE as follows.

Consider the birth and death process $\Pi^{(n)}_t = \Pi^{N,(n,1)}_t$ representing
the number of factors of which $\CF^{+}_t$ consists of at time $t$
if we are starting at $\Pi^{(n)}_0=n$ at one site. The dynamics of this object
follows immediately from the evolution of $(\CF^{+}_t)_{t\geq
0}$. The forward Kolmogorov equations of the process $\Pi^{(n)}_t$
are therefore given by \be{KE1}
p^\prime_{ij}(t)=s(j-1)p_{i,j-1}(t)+\frac{d}{2}(j+1)jp_{i,j+1}(t)-(sj+\frac{d}{2}j(j-1))p_{ij}(t),
\ee

\be{grev76b} p_{n,n}(0)=1.\ee Letting $E[\Pi^{(n)}_t]^k
=m_k(t)=\sum_{j=0}^\infty j^kp_{n,j}(t)$, we obtain with $\hat
s=s+\frac{d}{2}$: \be{grev76b2} \frac{dm_1(t)}{dt}= \hat
sm_1(t)-\frac{d}{2}m_2(t)\leq \hat
sm_1(t)-\frac{d}{2}m_1^2(t),\quad m_1(0)=n.\ee Therefore

 \be{grev78b} E[  \Pi^{(n)}_t ] \leq
\frac{2\hat sne^{\wh st}}{nd(e^{\hat st}-1)+2\wh s},\quad \ee so
that \be{grev78c} E[\int_{0}^{t} \Pi^{(n)}_u du] \leq \frac{2}{d}
\log(\frac{2\wh s}{d}+{n}(e^{\hat st}-1)).\ee

This implies that if we let $1> \varepsilon >0$ and choose
$\alpha$ so that $1- e^{-\alpha}=\varepsilon$,  then using
(\ref{grev78c}) we get:
\be{grev78d}
P[\tau_N <a_N]\leq
\varepsilon +P[\frac{m}{N}\int_0^{a_N} \Pi^{(n)}_sds
> \alpha]\leq \varepsilon + \frac{2m}{d \cdot \alpha}\frac{\text{const}+[\log
n +\hat sa_N]}{N} \ee Hence $\lim_{N\to\infty} P(\tau_N <a_N) \leq
\varepsilon$ for arbitrary $\varepsilon >0$ and the proof is
complete.
\end{proof}

\subsection{Proofs of Propositions \ref{P.MVentrlaw},
\ref{P6.1b}, and \ref{P.measure}.} \label{ss.proofemfix}

In this section we consider the large time scale
$\alpha^{-1} \log N +t, \quad t \in \R$ and prove Proposition \ref{P.MVentrlaw} on
properties of the limit dynamics in this {\em emergence} and {\em fixation time
scale} and the convergence results stated in Subsubsection
\ref{sss.Ninteract}, in particular, on the time needed for
emergence and the form of the process of fixation, i.e. Propositions
\ref{P6.1b}, \ref{P.measure}. (Recall the remaining assertions requiring analysis
in different time scales, in particular Propositions \ref{P.DF} and
\ref{CSB-longtime} on droplet formation are proved later in Subsection
\ref{ss.bounds} and Proposition \ref{P.chaos} on tagged sample convergence
which uses different methods will be proved in Subsection \ref{ss.chaos}.)

\subsubsection{Outline of the strategy of proofs for the asymptotic analysis}
\label{sss.outline}

The proofs for the large time scale proceeds in altogether thirteen subsubsections.
 Therefore we describe here in \ref{sss.outline} first the general strategy
 and list at the end what happens in the various subsubsections.

The proofs of the results on the $N \to \infty$ asymptotics
for the process in large time scales involve calculating the mean of $x^N_1(t)$ or $x^N_2(t)$
respectively of higher moments and hence
depend in an essential way on an asymptotic analysis of the  dual
process which we introduced in Section \ref{s.dual}.

In the case of two types the required dual calculation involves a
system of particles (corresponding to factors) located at sites in
$\{1,\dots,N\}$ with linear birth rate and quadratic death rate at
each site and migration between sites. Assuming that the initial
population consists of only  one individual and the initial
function $f$ is the indicator of type 1 at site 1, then we get a growing
string of factors of such indicators, until the rare mutation
operator (at rate $\frac{m}{N}$) acts on it  at which time it and
therefore the product becomes zero.  Hence the main point is to prove
that a mutation event somewhere in
the string will occur in the $\log N$ time scale.
For this to be the case we have to prove that the total
number of factors reaches $O(N)$ in that time scale.

Consider now the dual particle system $\eta$ starting with  only finitely
many particles. As long as the
number of occupied sites remains $o(N)$ the migration of individuals
leads in the limit $N \to \infty$ always to an unoccupied site.
More precisely, in this regime the proportion
of sites at which a collision occurs, that is a migration to a
previously occupied site, is negligible. For that reason
the collection of {\em
occupied sites} can be viewed (asymptotically as $N\to\infty$) as
a {\em Crump-Mode-Jagers branching process} (a notion we shall recall later on).
This process grows exponentially
and the identification of
the emergence time corresponds to the Malthusian parameter $\alpha$ of this
CMJ-branching process.

In order to investigate the dynamics in the next phase between emergence and fixation,
i.e. in times $C \log N+t$ with $ t\in \R$ the new time parameter, we must
on the side of the dual particle system
now consider the role of {\em collisions} of dual particles once the number of occupied sites
is of order $O(N)$ and the proportion of occupied sites at which
collisions occur is asymptotically non-zero. This can be handled
by identifying a dynamical law of large numbers described by a
{\em nonlinear} evolution equation.
Nonlinear due to the immigration term at a site, resulting from immigration
from other sites and which in the limit depends on their law. On the side
of the McKean-Vlasov limit, this
makes the evolution of the frequency of type 2 at single site depend on a global signal, namely
the current spatial intensity of type 2.

The delicate issue of the {\em transition}
between these two regimes with and without collision must be addressed  to fully describe the
emergence and to show that this leads to  the McKean-Vlasov
dynamics with random initial condition. For that purpose the
system is considered in the time scale
$\alpha^{-1} \log N+t$, with $t \in \R$ the macroscopic time parameter. The key element in the
analysis is the determination of the leading term in an asymptotic
expansion of the type-2 mass in time scale $\alpha^{-1} \log N+t$ with $t$
varying, once we have taken $N \to \infty$, as $t\to -\infty$.
On the side of the dual process this will correspond to an expansion of the
normalised number of occupied sites as $N \to \infty$ and then $t\to\infty$ in
the above time scale.
However also in the discussion on droplet formation later on in subsection \ref{ss.bounds} it is necessary to
keep track of the collisions since we need higher-order terms in the expansion
of the total type-2 mass and these terms involve collisions of dual particles.
This will be dealt with in
Subsubsection \ref{sss.dropletgrow} and the sequel.

To carry out this program  we must study the behaviour of the dual
particle system in two time regimes, namely,
\be{partsys1} s(N)+t
\mbox{ where } s(N)\to \infty \mbox{  but } s(N)=o(\log N) \ee and
for suitably chosen $\alpha$,

\be{partsys2} \frac{ \log N}{\alpha}+t,\quad -\infty<t<\infty, \ee
and their asymptotics as $N\to\infty$. The first we call the {\em
collision-free regime} or {\em pre-emergence regime}, the second
the {\em collision regime}.

 A key step in the dual process
analysis is a coupling of the pre-emergence and collision regimes
in order to understand the {\em transition} between the two regimes
and to be able to consider the time scale
\be{partsys3}
\frac{1}{\alpha} \log N +t_N \mbox{ with } t_N \to -\infty,
\quad t_N = o(\log N).
\ee
\bi

{\bf Outline of Subsection \ref{ss.proofemfix}.}

This subsection focusses
on the proofs of Proposition \ref{P.MVentrlaw} (Subsubsection \ref{sss.outline})
on the limiting McKean-Vlasov dynamics, Proposition \ref{P6.1b} (Subsubsections
\ref{sss.outline} -\ref{sss.collreg}) and Proposition \ref{P.measure}
on emergence respectively convergence to the limiting dynamics of fixation
(Subsubsections \ref{sss.growth}-\ref{sss.complproof}).

The Subsubsection \ref{sss.structure} is devoted to the proof of
Proposition \ref{P.MVentrlaw}, which specifies the limiting laws
in short and large time scales.

In Subsubsections
\ref{sss.dualcfr}-\ref{sss.collreg} we begin the proof of
Proposition \ref{P6.1b} on the order of magnitude of the emergence
time. In a first step we exhibit the branching structure for the
number of sites in \ref{sss.structure} and then analyse the
collision free, respectively the collision regime,
\ref{sss.dualcfr} to \ref{sss.collreg}.

The proof of Proposition \ref{P.measure} on the limit dynamics is
more involved and is broken in six preparatory steps which are
entirely concerned with the dual process at times
$\frac{1}{\alpha} \log N+t$ and a final argument returning to our
process, each of these is in one of the  separate Subsubsections
\ref{sss.growth}-\ref{sss.complproof}. Namely, we
study in Subsubsections \ref{sss.growth}, \ref{sss.asycolreg},  \ref{sss.prep2}
and \ref{sss.limrandmarg} the structure of the
dual process more carefully, show in \ref{sss.proofldy}
convergence of the process to the limiting dynamics and examine in
Subsubsection \ref{sss.proofmeasure} the random growth constant in
the emergence. Then we are able in Subsubsection
\ref{sss.complproof} to complete the proof of Proposition
\ref{P.measure} on the limiting dynamics of fixation.

\subsubsection{Proof of Proposition \ref{P.MVentrlaw} (Properties McKean-Vlasov dynamic)}
\label{sss.proofkv}

In order to prove the results on the McKean-Vlasov dynamic recall
the duality in Section \ref{s.dual}.
We now go through the different parts of the proposition to prove,
but here the hard part is assertions d and e.

{\bf (a)} This is the wellknown McKean-Vlasov limit for the interacting
Fleming-Viot model with migration and selection (compare, for example, \cite{DG99}, Theorem 9).

{\bf (b)} We want to show that $E[x_1(t)] \to 0$ as $ t \to \infty$ if
$E[x_1(0)]<1$. This follows by a simple application of the
mean-field dual $(\wt \eta_t, \CF^+_t)_{t \geq 0}$
started with $\CF^+_0 = 1_{\{1\}}$ and $\eta_0$
having one particle located at site 0. As time proceeds in the
dual process the selection operator keeps being applied. But the
$\ell$-fold application of the selection operator, which now in
the two-type case is simply obtained as $f \to f\otimes$
$1_{\{1\}}$, produces as surviving terms for $\CF^{+}_t$ a string
of the form (note coalescence only changes here the number of factors
belonging to distinct variables)
\be{angr8d} \Pi^{\ell(t)}_{i=1}
1_{\{1\}}(u_{\xi_i(t)}) \ee and since the total number of factors
diverges as soon as the migration rate $c$ is positive,
i.e. $\ell(t)\to\infty$,  the result follows from the
assumption $E[x_1(0)]<1$.

{\bf (c)} For existence it suffices to obtain a particular
solution $(\CL^\ast_t(2))_{t \in \R}$ where
$\CL^\ast_0(2)\in\mathcal{P}([0,1])$. Take the solution
$(\wt\CL^n_t(2))_{t \geq 0}$ with
$\wt{\mathcal{L}}^{n}_0(2)=\mu_n$ with
\be{ag50}
a_n := \int_0^1 x \mu_n(dx), a_n\to 0 \mbox{ as } n\to \infty.
\ee
In particular we can take here
 $\mu_n=\delta_{a_n}$. Choose now $(a_n)_{n \in \N}$ with  $0<a_n <\frac{1}{2}$
 and $a_n \to 0$ as $n \to \infty$.
 Then first of all for every fixed $n$
\be{11201}\lim_{t\to\infty} \wt{\CL}^{n}_t(2)=\delta_1\ee
and there
exists $r(a_n)$ such that $\intl^1_0 x \wt{\CL}^{n}_{r(a_n)}(2)dx
=\frac{1}{2}$. In fact by duality we see that $r(a)$ is strictly
increasing and continuous in $a$.

Then define \be{angr8d2}
(\CL^{n}_t(2))_{t \geq -r(a_n)}  = (\wt{\CL}^{n}_{t+r(a_n)}(2))_{t
\geq -r(a_n)}. \ee Since  for any $t>0$,
\be{Ex3}{} \lim_{n\to
\infty} \wt{\mathcal{L}}^{n}_t(2) = \delta_0, \ee it follows again
by duality that $\lim_{n\to\infty}r(a_n)=\infty$.
 The value at $t=-r(a_n)$ tends to $\delta_0$ as $n \to
\infty$.

We note that if for a subsequence $(n_k)_{k \in \N}$ the sequence
$(\CL^{n_k}_s)_{k \in \N}$
converges as $k \to \infty$ for some
$s \in \R$, then it converges for all times $t \geq s$, as we
immediately read off from the duality which gives the Feller property
(continuity in the initial state) of the McKean-Vlasov process. Therefore by considering
for example $s= -k, k \in \N$ we can obtain a further subsequence
converging for all times $t \in \R$. Now we take such a convergent
subsequence of $(\CL^n(2))_{n \in \N}$ and define $\CL^\ast(2)$ as
the $(n \to \infty)$ limit. By construction $\int x \wt \CL^n_0(2)
(dx) = \frac{1}{2}$ and therefore $\CL^\ast_t(2) \to \delta_1$ as
$t \to \infty$. It remains therefore
to argue that $\CL^\ast(2)$ is a solution and is
converging to $\delta_0$ as $t \to -\infty$.

The fact that we
obtain a solution of the equation is a consequence of the
construction together with the continuity in the initial state
which is immediate using the duality.
Namely for every $n$ we have a solution for $t \geq -\tau(a_n)$.
Then using duality we see that the limit is a solution as well,
since it satisfies the duality relation between times
$s < t; s,t \in \R$.

The convergence to $\delta_0$ as $t \to -\infty$ follows by
contradiction as follows.
Assume we converge to a  strictly
positive mean. Then applying the duality to calculate the mean at
time 0 starting at times  $t\downarrow -\infty$  we see, that this mean
converges to one and therefore we obtain
convergence to $\delta_1$ (by (\ref{11201})) at time 0. This
contradicts the assumption that the mean is $\frac{1}{2}$ at time
0. This completes the proof of existence of a solution on
$(-\infty,\infty)$ satisfying (\ref{Ex1b}).

{\bf (d, e)} Here we proceed in three steps, first we establish existence
of solutions whose mean curve has the claimed growth
behaviour at $-\infty$ and then we show uniqueness of the
mean curve with this growth behaviour and subsequently the
uniqueness of the solution for given mean curve.
\bi

{\bf Step 1} $\;$ {\em Growth of mean-curves}

The first step is to consider the growth of  solutions of
the McKean-Vlasov equation, $\CL_t(2)$, at $-\infty$. We will show
two things,
\begin{itemize}
\item first there exists a solution with mean curve $m(\cdot)$ which
behaves like $m(t)  \sim A e^{\alpha t}$ for $t \to -\infty$ with $A
\in (0,\infty)$,
\item secondly that
all solutions satisfying the growth condition posed in (d) have
mean curves satisfying $m(t)\sim Ae^{\alpha t}$ as $t\to -\infty$
for some $A \in [0,\infty)$.
\end{itemize}
In fact we shall combine these two points in our proof.

First assume that $(\CL_t(2))_{t \in \R}$ is a solution (and entrance law)
with mean-curve
$(m(t))_{t \in \R}$ which satisfies:
\be{angr8e} m(t)>0, \qquad 0< \liminf_{t\to
-\infty}e^{-\alpha t}m(t)<\infty. \ee
We want to show that $e^{-\alpha t} m(t)$ converges to some $A \in \R^+$.

We can chose a sequence $t_n$ with $t_n \to - \infty$ so that
for some $A \in (0, \infty)$ we have
\be{angr8f}
 e^{\alpha |t_n|}m(t_n)=e^{\alpha |t_n|}\int_0^1 x  \CL_{t_n}(2)(dx) \to A
 \mbox{ as } n \to \infty.
 \ee
We now abbreviate
\be{angr8f1}
\mu_{t_n} := \CL_{t_n}(2).
\ee
What we know about $\mu_{t_n}$ is its mean and from this we have to
draw all conclusions on the convergence, namely that the limit $A$
does not depend on the choice of $(t_n)_{n \in \N}$.

On the other hand we can start with such a measure
$\mu_{t_n}$ as in (\ref{angr8f1}) and try to construct a solution with a mean-curve
$(m(t))_{t \in \R}, \quad m(t)=(A +o(1)) e^{\alpha t}$ as $t \to -\infty$
as the limit point of a sequence of solutions starting at time $t_n$ in
$\mu_{t_n}$. A solution starting in $\mu_{t_n}$ at time $t_n$ is denoted:
\be{angr8f3}
(\wt \CL^n_t)_{t \geq t_n}.
\ee
Then we know, that if
$(\CL_t(2))_{t \in \R}$ is an entrance law with $\CL_{t_n}= \mu_{t_n}$, then
by uniqueness of the solution starting at $t_n$ we  must have
\be{angr8f4}
\CL_t (2) = \wt \CL^n_t (2), \mbox{ for } t \geq t_n.
\ee

Therefore our strategy now is twofold. First we verify that in fact the mean curve
$(m(t))_{t \in \R}$ of $(\CL_t(2))_{t \in \R}$ has the form:
\be{angr8d3} m(t) =
\exp (\alpha (t +\frac{1}{\alpha} \log A)) + o(e^{\alpha t}). \ee
Then we argue that in fact the $(\wt \CL^n_t)_{t \geq t_n}$ with mean
curves satisfying (\ref{angr8f}) converge to an entrance law with
mean curve satisfying (\ref{angr8d3}). We begin with the first point and
come to the second one (which is then simpler by what we have prepared
for the first point) after (\ref{rc5c2}).

In order to carry this program out we
must consider the different possibilities we have for $\mu_{t_n}\in
\mathcal{P}([0,1])$ with mean curve  $m(t_n)\sim \frac{A}{e(n)}$
where we use the abbreviation
\be{angr8f2} e(n)= e^{\alpha|t_n|}. \ee
This is somewhat complicated by the possible interplay of
probabilities and sizes in producing a given mean.

We begin by considering two special cases
that illustrate the main idea. Let
\be{rc1}
\mu_{t_n}
=(1-\frac{p_1}{e(n)})\delta_{\frac{a_1}{e(n)}}+\frac{p_1}{e(n)}\delta_{a_2},\quad
 \ee
so that $A= a_1+p_1a_2$.
We can distinguish now different cases. The first case arises in
our example if (small probability for an order one value $a_2$)
\be{angre1}
p_1=0,\;a_1>0,\ee
in which case the mean curve satisfies
\be{Ca2}
e(n)\int x \mu_{t_n} dx \la a_1,\ee
the second case arises if (large probability for very small values)
\be{angre2}
p_1=1,\; a_1=0, a_2>0,
\ee
in which case the mean-curve satisfies
\be{Ca1} e(n)\int x \mu_{t_n}
dx \la \text{const}
\ee
given in (\ref{rc5a2}).

We shall now consider these two cases each with a distinct flavour
which then allows us, since they catch the key features,
to see how we can treat the {\em general} case,
which we call case 3 and for which we need to use
an additional trick.

A key tool in all cases is the duality relation and the fact that the
dual process consists of a collection of birth and death processes
acting independently so that with generating functions we can
analyse the dual expectation.

{\bf Case 1}

Here we first treat the deterministic case $\mu_{t_n}
= \delta_{\frac{a_1}{e(n)}}$. We will use  the  dual
representation for the McKean-Vlasov limit
together with (\ref{ang2b}) (giving below
$W_n\Rightarrow W$) to compute the expected value of $x_2$ mass at
time $t$, i.e. the mean of $\CL_t(2)$. Note that then we consider the time stretch
$|t_n| +t$ from the time point $t_n$ where we have information about
the state.

Observe that we know that the dual process at time $T$ is a
collection of $K_T$ many independent birth and death processes,
which have evolved for a time $T$ and which  in the large time
limit $T \to \infty$ satisfy the {\em stable age distribution}. Therefore the
state of the empirical measure of the birth and death processes
at times $T +t$ is
given by the stable size distribution, which arises as \be{agre60}
\int \CU(\infty, ds) \nu_s \quad , \ee where $\nu_s$ is the law of
the state of a single birth and death process at time $s$ starting
with one particle at the occupied site.

We now fix some time $t > t_n$ and suppress the $t$-dependence in
the notation. We need the generating function
\be{agre61}
G_n (z) = E[z^{\zeta (t-t_n,i)}]
\ee
of the state of the birth and death process at a randomly chosen
site $i$. To calculate we introduce the generating function
\be{detca1}
 G(z)=E[z^{\zeta}],\qquad 0<z\leq 1,
 \ee
 where $\zeta$ represents the number of dual factors at a dual site under
the stable size distribution and $(\zeta(i), i \in \N)$ is an i.i.d. collection of
such variables.  We know that
\be{agre62}
G_n \la G \quad \mbox{ as } n\to\infty,
\ee
\be{agre63}
G_n^\prime \la G^\prime \quad \mbox{ on the unit circle as } n\to\infty.
\ee

Then noting that $G^\prime(1)=B=\frac{\alpha +\gamma}{c}$  and
using (\ref{stableage}):
\bea{Ex4}
\int_0^1
(1-x) \wt \CL^n_t (dx) &&=E\left[\prod^{K_{t-t_n}}_{i=1}
(1-\frac{a_1}{e(n)})^{\zeta(t-t_n,i)}\right]
= E \left[ \prod^{K_{t-t_n}}_{i=1} G_n (1- \frac{a_1}{e(n)}) \right]
\\
&& = E\left[\left( G(1-\frac{a_1}{e(n)})\right)^{W_{n} e(n)e^{\alpha t}(1+o(1)) }\right],
     \mbox{ as } n \to \infty \nonumber\\
&& \ntoo E[e^{-Wa_1 B e^{\alpha t}}]\nonumber\\
&& \sim 1- E[W]a_1 B e^{\alpha t} +O(e^{2\alpha t}), \mbox{ as } t
\to -\infty. \nonumber \eea
 From the first to the second line we
have approximated $G_n$ by $G$ and used at the same time the
convergence to the stable age distribution. To go from the second
to the third line we use that $W_{n} \la W$ as $n \to \infty$.
\bigskip

Therefore as $t \to -\infty$:
\be{Ex8}
m(t) \sim [E[W]a_1 B]  e^{\alpha t} +O(e^{2\alpha t}).
\ee

{\bf Case 2} \quad

In case 2 the number of dual sites occupied in the dual process
which have the value $a_2$ at time $t_n$,  is denoted
$Z_n(t,a_2)$, and is Bin$(W_n e(n) e^{\alpha t},\frac{1}{e(n)})$
distributed.  Hence we get
\be{rc3}
 \int_0^1 (1-x)\wt \CL^n_t(dx)=
E[\prod_{i=1,\dots,K_{t-t_n}} (1-a_2)^{\zeta(t-t_n,i})].
\ee
The r.h.s. of (\ref{rc3}) is approximated as above as in case 1 (see (\ref{Ex4}))
for $n \to \infty$
by \be{rc3a1}
 E\left[(G(1-a_2))^{Z_n(t,a_2)} \right].
\ee
Then $Z_n(t,a_2)$
converges as $n \to \infty$ in distribution to a Poisson distribution with parameter
${We^{\alpha t}}$.

For the calculation it is convenient to {\em condition on $W$} and read
expectations as expectations for given $W$ and write $\wt E$ for
this conditional expectation suppressing the condition in the
notation. We conclude that as $n \to \infty$ and then $t \to -
\infty$: \bea{rc4a1}
\wt E[\int_0^1 (1-x) \wt \CL^n_t (dx)]
  &&\la e^{[ G(1-a_2)-1]{We^{\alpha t}}}, \mbox{ as } n \to \infty\\
  && \sim 1-[ (1-G(1-a_2))]{We^{\alpha t}}+O(e^{2\alpha t}), \mbox{ as } t \to -\infty.\nonumber
\eea Hence as $n \to \infty$ and then for $t \to - \infty$ we get:
\be{rc5a2}
E[\int_0^1 x \wt \CL^n_t(dx)]\sim [(1-G(1-a_2))]E(W)]
  e^{\alpha t} + O(e^{2 \alpha t}).
\ee
In other words also in this case $m(t)=Const \cdot e^{\alpha t} + O(e^{2 \alpha t})$.

{\bf Case 3}\label{AM3}

We now turn to the general case. The key idea is to focus on the components
where we do have type-2 mass even as $t_n \to -\infty$, which is achieved
using the Palm measure. Consider the  family $\{\wh
\mu_{t_n}\}_{n \in \N}$ of Palm measures on $[0,1]$: \be{dd50} \wh
\mu_{t_n} := \frac{x\mu_{t_n}(dx)}{\int_0^1x\mu_{t_n}(dx)} \ee and
denote the normalizing constant again as in case 1,2 by
\be{dd51} e(n) =
(\int_0^1x\mu_{t_n}(dx))^{-1}. \ee
Since the family $\{\wh
\mu_{t_n}\}_{n \in \N}$ is tight (by construction the measures are supported on
$[0,1]$) we can find a subsequence and this sequence we denote again
by $(t_n)_{n \in \N}$ such that
$t_n\to -\infty$ and $\wh \mu_{t_n}$ converges. This means we can find
a measure $\wh \mu$ on $[0,1]$ (with $\wh \mu(\{0\})=0$) such
that
\be{dd52}
\frac{x\mu_{t_n}(dx)}{\int_0^1x\mu_{t_n}(dx)} \La
a\delta_0+ \wh \mu, \quad \mbox{ as } n \to \infty. \ee Moreover
there exists a scale $s_n$, with $s_n >0$ and $s_n\to 0$ such that
\be{dd53}
\frac{\int_0^{s_n}x\mu_{t_n}(dx)}{\int_0^1x\mu_{t_n}(dx)}\to a.
\ee We shall study the contributions due to the part $a \delta_0$
and due to $\wh \mu$ first separately in (i) and (ii) below and then join
things to get the total effect in (iii).

{\em (i)} We begin with the effect from $a \delta_0$.
For this purpose consider the measures $(\nu_n)_{n \in \N}$ on $[0,1]$
 defined by setting for every interval $(c,d)$, with $0\leq c<d \leq 1$:
 \be{dd54}
\nu_n((c,d))=\frac{\int_{ds_n}^{cs_n} x \mu_{t_n}(dx)}{\int_0^{s_n}x\mu_{t_n}(dx)}.
\ee
This allows us to consider the contribution due to the terms leading to
$a\delta_0$ in the expectation in the duality relation. Define
\be{dd55}
(\wt \nu_n)_{n \in \N}
\ee
as the empirical measure from
a sample of size $e^{W_n\alpha(t-t_n)}$ sampled from $\nu_n$.
Now we consider the state of the dual process at time $t -t_n$
and acting on the state at time $t_n$ but considering only the
contributions from sites with values in the interval $s_n$,
which means that we can represent the state at time by an empirical
measure as follows.  As $n \to \infty$ we have:
\be{dd56}
\int_0^1 (1-x) \wt \CL^n_t (dx) \sim E\left[
\exp\left(W_{n} e^{\alpha t} \int_0^1\log G(1-ys_n)  \wt\nu_n(dy)
\right)\right].
\ee

By (\ref{grev90a}) $G^{\prime\prime}(1)<\infty$ and therefore the
error term in (\ref{dd56}) is $O(s_n)$. Then taking the limit as $n\to\infty$ we
get using relation (\ref{dd53}) that:
\bea{dd57}
&&E\left[ \exp\left(W_{n} e^{\alpha t} \int_0^1\log G(1-ys_n)
\wt\nu_n(dy) \right)\right]\nonumber\\&&\to E\left[ \exp\left(-
G^\prime(1)We^{\alpha t}a \right)\right]\qquad\text{as
}n\to\infty\\&&\sim 1-G^\prime(1)E[W]e^{\alpha t}a\qquad\text{as
  }\;\;t\to -\infty.\nonumber
\eea
Hence in this case we get if there are only contributions from
$a \delta_0$, i.e. $\wh \mu\equiv 0$, that:
\be{dd58}
m(t) = (G^\prime (1) E[W]a) e^{\alpha t}
+ O(e^{2 \alpha t}),
\ee
i.e. $A = G^\prime (1) E[W]a$.

{\em (ii)}
Next we consider the contribution if they are only coming from $\wh \mu$.
We have to ''return'' here from the Palm $\wh \mu$ to a measure $\mu$
which need not be a probability measure since $x^{-1}$ diverges at $x=0$.
Note that for
any $\ve >0$ the set  of finite measures given by \be{Case3}
\{e(n) (1_{x>\ve} ) \mu_{t_n}(dx)\}_{n\in\N},
\ee is
tight in the space of finite measures on $\mathcal{M}([0,1])$ and
we can choose a convergent subsequence. We can then obtain a
convergent sequence, we call again $(t_n)_{n \in \N}$ such that:
\be{dd59}
e(n) (1_{x>\ve}) \mu_{t_n}(dx)
\Rightarrow \mu, \mbox{ as } n \to \infty
\ee where $\mu$ is a
measure on $(0,1]$ with
\be{dd59b}
\int_0^1 x\mu(dx)<\infty, \mbox{ i.e. in particular }
\mu ([\ve,1]) <\infty \mbox{ for all } \ve>0.
\ee
We denote
\be{ag51}
\mu_\ve = \mbox{ restriction of } \mu \mbox{ to  } [\ve,1].
\ee

By the dual representation we then have under the assumption
\be{Case3c}
\mu_{ \ve} \mbox{ is discrete with atoms at } x_1, x_2, \cdots, x_M,
\ee
that with $x_{g_n(i)}$ being the value at site $i$ at time $t_n$
the representation:
\be{rc3a}
 \int_0^1 (1-x)\wt \CL^n(dx)=
E[\prod_{i=1,\dots,K_{t-t_n}} (1-x_{g_n(i)})^{\zeta(t-t_n,i)}].
\ee
Denote by $Z_n(t,x_j)$
the number of dual sites occupied in the dual process which have at time $t_n$
the value $x_j$.
The r.h.s. of (\ref{rc3a}) we approximate again for $n \to \infty$ by
\be{rc3b}
E\left[\prod_{j}G(1-x_j)^{Z_n(t,x_j)} \right].
\ee
We know
\be{rc3c}
\CL [(Z_n(t,x_j))_{j=1,\cdots,M}] =  Mult \left(W_{n}
e(n) e^{\alpha t},(\frac{\mu(x_j)}{e(n)})_{j=1,\cdots,M}\right).
\ee

The assumption (\ref{Case3c}) is removed as follows.
For a general $\mu_{\ve}$, we consider the random measure
\be{rc3d}
Z_n(t,dx) \mbox{ (the number of dual sites  at
which  } x_2\in dx)
\ee
which converges as $n \to \infty$ in distribution to an {\em inhomogeneous Poisson
random measure} with the intensity measure
\be{rc3e}
\mu(dx){We^{\alpha t}}.
\ee

For the calculation it is
convenient to condition on $W$ and read expectations as
expectations for given $W$ and write $\wt E$ for this object.
We conclude that as $n \to \infty$ and then $t \to - \infty$
with the same reasoning as in the other cases:
\bea{rc4}
\wt E[\int_0^1 (1-x) \wt \CL^n_t (dx)]
  &&\to e^{\int^1_0 (G(1-x)-1)\mu(dx){We^{\alpha t}}}\qquad,
  \text{as   }n\to\infty\\
  && \sim 1-[\int (1-G(1-x))\mu(dx)]{We^{\alpha t}}+O(e^{2\alpha t})\qquad,
  \text{as   }t\to -\infty.\nonumber
\eea
Hence taking first $n \to \infty$ we then get for $t \to - \infty$ the expansion:
\be{rc5}
E[\int_0^1 x \wt \CL^n_t (dx)]\sim [\int^1_0(1-
G(1-x))\mu(dx)]{E[W]e^{\alpha t}} + O(e^{2\alpha t}),
\ee
which then implies that the mean curve satisfies if $a=0$, that
\be{rc5a}
m(t) = A e^{\alpha t} +O(e^{2\alpha t}),
\ee
with
\be{rc5a1}
A = E[W] \int^1_0 (1-G(1-x)) \mu (dx). \ee

{\em (iii)} Finally we have to use the results in the two different cases
where $\wh \mu \equiv 0$ respectively $a=0$, to
get the answer in full generality.  Combining the expressions for
the two contributions we obtain as $n \to \infty$ and then considering
the expansion as $t \to -\infty$ that

\be{rc5+}
\int_0^1 (1-x)\wt \CL^n_t (dx) \ntoo  E\left[ e^{\int
(G(1-x)-1)\mu(dx){We^{\alpha t}}} e^{-W G^\prime(1)e^{\alpha t}a}\right]
+ O (e^{2 \alpha t}),
\ee
and as a consequence for $t \to -\infty$ the r.h.s. is
asymptotically given by
\be{rc5b}
 \sim \left([\int^1_0 [G(1-x)-1]\mu(dx)
+aG^\prime(1)] E[W]\right) e^{\alpha t} +O(e^{2 \alpha t})
= A \cdot e^{\alpha t} + O(e^{2 \alpha t})
\ee
with
\be{rc5c}
A = \left( \int^1_0 [G(1-x) -1] \mu (dx)
+ a G^\prime (1) \right) E[W].
\ee

Finally we have to exclude the possibility that the choice of
different subsequences  in (\ref{angr8f}), (\ref{dd52}), (\ref{dd59}) (and hence  of different
values for $A$) would give a different mean curve. However this would
lead to a contradiction since this would give two values to the
same quantities.  Therefore the resulting mean curve is
independent of the choice of subsequence $(t_n)_{n \in \N}$ of
starting times.

This means that we have proved that an entrance law with at most
exponentially growing mean has the property that the mean curve
has asymptotically as $t \to -\infty$ the form
\be{rc5c2}
m(t) = A \cdot e^{\alpha t} +O(e^{2 \alpha t}),
\ee
with $A$ given by the formula (\ref{rc5c}) which completes the first of
the two points we specified around (\ref{angr8d3}).

We next need to argue that in fact such an entrance law exists by showing
that actually $\wt \CL^n$ (from (\ref{angr8f3})) converges along a subsequence.
Due to the continuity of $\CL_t$ in the initial state at time $t_0$, the
tightness is straightforward. We have to show that first of all
the limit is a solution and secondly the mean curve has the desired
asymptotics. To show that we obtain in the limit a solution we argue
as in the proof of part (c) of the proposition. It remains to verify
the asymptotics of the mean curve.

For this purpose we use the same line of argument as above for the first
moment of $x_2(t)$ to show that also  all higher
moments of $x_2(t)$ converge as $n \to \infty$ for every fixed $t$, so that we have
convergence of the processes using the Feller property. For that we observe that
we then simply have to start the dual process with $k$-initial particles and
then we carry out the same calculations, since the CMJ-theory works also for
the CMJ-process started with a different internal state of the starting site.
We omit the straightforward modification.

{\bf Step 2}  $\;${\em  Uniqueness of solution with the $t \to -\infty$ asymptotics}

We next show that if  the mean curves of two entrance
laws  differ at most by $o(e^{\alpha t})$ as $t \to -\infty$,
then they are identical.
This then allows to conclude that if we prescribe the value $A$ then
there is exactly one entrance law with this asymptotics as $t \to -\infty$.

For that purpose we consider two solutions
$(\CL^1_t)_{t \in \R}, (\CL^2_t)_{t \in \R}$ such
that the mean curves $(m^\ell(t))_{t \in \R}$ satisfy for
$\ell=1,2$ that $m^\ell(t) \sim A \exp (\alpha t)$ as $t \to -\infty$.  vbbbbbb
We consider a sequence $(t_n)_{m \in \N}$ with $t_n  \to -\infty$ as
$n \to \infty$ and we fix a time point $t \in \R$. Then we define
the two processes
we apply the duality relation at time $t$
\be{ag54}
(\wt \CL^{n,1}_s)_{s \geq t_n} \quad , \quad (\wt \CL^{n,2}_s)_{s \geq t_n},
\ee
which have at time $t_n$ the initial distributions given by
\be{ag55}
\CL^1_{t_n} \mbox{ resp. } \CL^2_{t_n}.
\ee
We apply the duality relation at time $t$.
Our goal is to show that the mean curve satisfies
\be{ag56}
m^1(t) = m^2(t), \quad \forall \; t \in \R,
\ee
by proving that for every $t \in \R$
\be{ag57}
\wt m^{n,1}(t) = \wt m^{n,2}(t) + o(1) \mbox{ as } n \to \infty.
\ee

We first verify this in Case 1. Let
$\wt{\CL}^{n,1},\;\wt{\CL}^{n,2}$ be two solutions as in case 1
but such that $m^1(t_n)-m^2(t_n) =o(\frac{1}{e(n)})$. Then by
(\ref{Ex4}),
\bea{Ex4x}
&&\left |\int_0^1 (1-x) [\wt \CL^{n,1}_t
(dx) -\wt \CL^{n,2}_t (dx)]\right|
\\&&\leq E\left[\left |\prod^{K_{t-t_n}}_{i=1}
(1-\frac{a_1}{e(n)}+o(\frac{1}{e(n)}))^{\zeta(t-t_n,i)}-\prod^{K_{t-t_n}}_{i=1}(1-\frac{a_1}{e(n)})^{\zeta(t-t_n,i)}\right|\right]\nonumber
\\&&= E \left[\left| \prod^{K_{t-t_n}}_{i=1} G_n (1-
\frac{a_1}{e(n)}+o(\frac{1}{e(n)}))- \prod^{K_{t-t_n}}_{i=1} G_n
(1- \frac{a_1}{e(n)})\right| \right]\nonumber
\\&&\leq E[G^\prime(1)W_{n}e(n)e^{\alpha t} o(\frac{1}{e(n)})] \to 0 \text{ as }n\to\infty, \nonumber \eea
where we have first used that we can split the expectation in the duality
relation into one over the initial state and then over the dual dynamic,
furthermore for $z\in [0,1]$, $G^\prime_n(z)\leq
G^\prime_n(1)\leq  G^{\prime}(1)$ (the latter by coupling).
Hence we have $m^1=m^2$ in this case.

In case 2 there iis  $o\left(\frac{1}{e(n)}\right)$ change in the mean curve resulting
from the change of $a_2$ to $a_2+ o\left(\frac{1}{e(n)}\right)$.
Moreover given
a change of the frequency $\frac{p_1}{e(n)}$ to $\frac{p_1}{e(n)}+
o\left(\frac{1}{e(n)}\right)$ we have $Z_n(t,a_2)$
is Bin$\left(W_ne(n)e^{\alpha
t},\frac{1}{e(n)}+o\left(\frac{1}{e(n)}\right)\right)$ and the
limiting Poisson distribution is not changed by the
$o\left(\frac{1}{e(n)}\right)$ perturbation.

Again, the general case follows by a combination of the arguments
of cases 1 and 2 exactly as in case 3 in Step 1.
\sm

{\bf Step 3}  $\;${\em Uniqueness given the mean curve}

To complete the proof for uniqueness of the entrance law with specified
behaviour as $t \to -\infty$, it
remains to show that for a given mean curve $(m(t))_{t \in \R}$ there is a
unique solution $(\CL^\ast_t)_{t \in\mathbb{R}}$ to the McKean-Vlasov equation
having this mean curve. Since the McKean-Vlasov
dynamics is unique for $t\geq t_0$ it suffices to prove uniqueness
of  $\CL^\ast_{t_0}$ for arbitrary $t_0$.

Now consider the tagged component in the McKean-Vlasov system
$(x_2^\infty(1,t))_{t \geq t_0}$ assuming that
the mean curve $(m(t))_{t \in \R}$ satisfies
\be{rc5d}
\limsup_{t\to -\infty}m(t) =0.
\ee
Then $\CL^\ast_t$ is the law of a tagged component
conditioned on the realization of the mean curve and this conditioned
law is a weak solution of the  wellknown SDE
\be{SDE+}
dx(t) = c(m(t)-x(t)) dt + sx(t)(1-x(t))dt
+ \sqrt{x(t)(1-x(t)} dw(t),
\ee
(cf. (\ref{angr7b0b})) which has a
unique weak solution given an {\em initial} value at time $t_0$.
However we must prove the existence and uniqueness of a solution
with time running in $\R$ and satisfying \be{Gre33}
 \lim_{t\to -\infty} x_2^\infty(1,t)=0.\ee

We  construct a ``minimal solution''  $\tilde x_2^\infty(1,t)$ for
this problem by starting with  a sequence of solutions $x_{(n)}(t)$ of
(\ref{SDE+}) (driven by  independent Brownian motions $w_n(t)$)
starting with value $x^n(t_n)=0$ at times $(t_n)_{n \in \N}$ with
$ t_n\to -\infty$. We then construct a coupled version
$\{\wh x^m:m\geq m_0\}$ by coupling paths $x^{m}$ and $x^n,\, n\leq m$ at
a time where they collide. It can be verified that the coupled
system $\{(\wh x^n(t))_{t \in \R, n \in \N}\}$  is a stochastically monotone increasing
sequence as $n\to\infty$ with limit $\wh x_2^\infty(1,t)$.
If we assume {\em existence} of a  solution of (\ref{SDE+}) with time
index $\R$ and mean curve $m(t)$ we then know it will be stochastically
bounded below by $\wh x^\infty_2$.

Next we need uniqueness in law of the solution and we shall show that
it must agree with $\wh x_2^\infty$, which then automatically must also be a solution.  We
observe that any solution to (\ref{Y777}) must have zeros at
$-|t|$ for arbitrarily large $|t|$ (cf. classical property of
Wright-Fisher diffusions \cite{S1}). If it has a zero at
$t^\ast>t_n$, then we can couple it with $\wh x_n$ and therefore we
can couple the diffusion with the minimal one therefore the
considered diffusion must agree with a version of $\wh x^\infty_2$
after the coupling time.
Combined with the uniqueness property of the
McKean-Vlasov equation this also proves uniqueness
 of $\CL^\ast_{t},\;t\geq t_0$. Since  $t_0$ is arbitrary, we obtain
 the {\em uniqueness} for all $t \in \R$.
This completes the proof of (e) of the Proposition
\ref{P.MVentrlaw}.

{\bf (f)} Finally,  we consider a {\em random} solution
$\{\CL_t:t\in\R\}$ to the McKean-Vlasov equation.  Since $\CL_t$
is a.s. a solution and since by assumption it satisfies a.s. the growth
condition with $A \in (0,\infty)$, it follows that it is given a.s. by a (random) time
shift of the standard solution.

\subsubsection{The structure of the dual process and a Crump-Mode-Jagers process}
\label{sss.structure}

The key tool for the proof of the emergence results  is the dual
representation for the system of $N$ exchangeably interacting
sites, namely,  the  dual process $(\eta_t,\CF^{+}_t)$
(note the dependence on $N$ is suppressed in the notation).

\sm

Since in the case of two types, $\mathbb{K}=\{1,2\}$ the two
frequencies of type one and two add up to 1, it suffices to
determine the law of $x_1(\frac{\log N}{\alpha}+t)$ and
$\Xi_N (\frac{\log N}{\alpha}) +t,\{1\})$ in the limit
$N \to \infty$.  This will, as
we shall see in Subsubsection \ref{sss.proofldy}, involve
computing $(\ell,k)$-moments, i.e. the product of $k$-th moments
at $\ell$-different sites for the $N$-site system and for the random
McKean-Vlasov limit. The calculation is carried out using the dual process started with $k$
particles at each of $\ell$ different sites. We denote the number
of particles in the dual process at time $t$ with such an initial condition
by \be{grev89aa2}
\Pi^{N,k,\ell}_t, \mbox{ we often write } \Pi^N_t \mbox{ if }
(k,\ell) \mbox{ are fixed}.\ee

First, to prove that the mass of the inferior type at a tagged
site goes to zero in time scale $\alpha^{-1} \log N+t$ as $N \to \infty$ and
then $t \to \infty$, it suffices to work with
$E[x_1(\cdot)]$. The calculation of this expectation involves only one initial particle
for the dual process. To carry out the calculation we assume that
$x_1(0,i)=1$ for $i=1,\cdots,N$ and
take the initial element for the dual process
\be{grev89aa}
 \eta_0 = \delta_{i_0},
\ee
i.e. one particle at a tagged site in $i_0 \in
\{1,\dots,N\}$.  Note that because of migration the law of this process depends indeed
heavily on $N$ if $t$ is large enough.

The number of factors at each site form a system of dependent
birth and death processes with immigration, i.e. a process with
state space $(\N_0)^N$, denoted \be{grev89ac} (\zeta^N(t))_{t \geq
0} \mbox{ with } \zeta^N(t) = \{ \zeta^N (t,i), i =1, \cdots,N\}.
\ee In other words

\be{grev89ad} \Pi^{N,k,\ell}_t = \suml^N_{i=1}
\zeta^N (t,i), \ee and $\zeta^N (0)$ has exactly $\ell$-non zero
components each containing  exactly $k$ factors.

Next we describe the dynamic of $\zeta^N$. Consider first a single component $i$.
If the current state of
$\zeta^N_t$ at a site $i$ is $x \in \N_0$, then
\be{grev89a} \begin{array}{l}
\mbox{the immigration rate of one new individual to this site is } c (N^{-1}
(\Pi^{N,k,\ell}_t -x)),\\
\mbox{where every individual at the other sites has the same chance to be
selected for }\\
\mbox{migration to } i,
\end{array}
\ee
\bea{grev89a1}
&& \mbox{the emigration rate for one individual at this site is }
c \frac{N-1}{N} x,\\
&& \mbox{the individual moves to a randomly chosen different site},
\nonumber \eea

\bea{grev89b}
&& \mbox{the death rate, i.e. rate for the death of one individual at
this site is }\\&& \quad (d/2)x(x-1), \mbox{ and}
\nonumber\eea
\be{grev89c}
\mbox{the birth rate for one new individual is } \quad sx. \ee The
dependence of the birth and death processes
$(\zeta(t,i))_{t \geq 0}, \quad i=1, \cdots, N$
at different sites
arises from the migration transition in (\ref{grev89a}),
namely the immigration is coupled with the emigration at another site.

Note that because of the exchangeability of $X^N$ and the initial
state $X^N(0)$ it suffices for the dual $\zeta^N$ to keep track of
\begin{itemize}
\item the number of occupied sites,
\item the number of individuals, i.e. factors, at each occupied site.
\end{itemize}

We now introduce the notation needed to  describe  the dual
population. First the process
\be{Y13}
\{ K^N_t\}_{t\geq 0} \ee
that counts the number of {\em occupied sites} at time $t$. Secondly the
process recording the {\em age and size} of each site (age is counted
from the time of
``birth'', that is, the time when the site first became occupied the last time).
Note that each site is uniquely identified by its birth time.
Equivalently, due to exchangeability it suffices to keep track of
the number of occupied sites $ K^N_t$ and the empirical age and
size distribution at time $t$, a measure on $\R^+ \times \N$ denoted
\be{Y13b} \Psi^N(t,ds,dy),\ee where
$\Psi^N(t,(a,b],y)$ denotes the number of sites in which the birth
time, that is, the time at which the initial immigrant arrives,
lies  in the time interval $(a,b]$ and the current size is $y\in
\mathbb{N}$.

We will use the abbreviation
\be{ang0b}
\Psi^N (t, ds) = \Psi^N (t, ds,\N).
\ee

Consider next for a fixed $t \in \R^+$ the collection of processes
 \be{angr15ab}
\{(\zeta^N_s(u))_{u \geq s}, \; s \in \mbox{ supp }
\Psi^N(t,\cdot,\mathbb{N})\},
\ee
which denotes the size of the
population at a site having birthtime $s$, that is, one particle
arrives at time $s$. In the absence of collisions these processes
for different birth times are independent birth and death
processes and their distributions depend only on their current
age.  However in the presence of {\em collisions} these  are {\em coupled} (by
migration) and no longer independent.

The number of occupied sites fluctuates due to migration to {\em new}
sites respectively jumps from sites with one particle to other
occupied sites. Note that the probability of a jump of the
latter type is non-negligible only if the
number of occupied sites is comparable to the total number of
sites. Hence up to the time that the number of occupied sites
reaches $O(N)$  the number of occupied sites is (asymptotically)
non-decreasing and the birth and death processes at different
sites are asymptotically independent.

We obtain an upper bound for the
growth of $( K^N_t)_{t \geq 0}$ if we suppress all collisions and
assume that always a new site is occupied. In order to construct
this process, we enlarge the geographic space from$\{1,\cdots,N\}$
to $\N$ and drop in (\ref{grev89a}) all jumps to occupied sites
replacing them by jumps to a new site, namely the free one with
the smallest label. We refer to this process as the {\em McKean-Vlasov dual}
since it arises by taking in our model of $N$ exchangeable sites
the limit as $N \to \infty$ of the dual process. This process is
denoted (if we want to stress how it arises we add the
superscript $\infty$):

\be{agr1}
(\zeta_t) = (\zeta_t(1), \zeta_t(2), \cdots)
\ee
and the process of the total number of individual by
\be{angr15ac}
(\Pi^{(k,\ell)}_t)_{t \geq 0},
\ee
and $\zeta$ has the following markovian dynamic.

If  we have at a site $i \in \N$ occupied with
$x \in \N$ individuals then we have
\be{angr15ad}
\mbox{ one death at rate } \frac{d}{2} x(x-1),
\ee
\be{angr15ae}
\mbox{one birth at rate } sx
\ee
and if $j$ is the smallest index with
$\zeta_t(j)=0$ then
\be{angr15ad1}
\mbox{at rate } c \zeta (i) \mbox{ an individual moves from } i \mbox{ to } j.
\ee

Since in the duality relation we only need the occupation numbers
of occupied sites,
we can {\em ignore migration steps of single factors} to
another unoccupied site. Therefore in the migration rate in
(\ref{angr15ad1}) we take
\be{ddy1}
c \zeta_t(i)  1_{(x \geq 2)} \mbox{ instead of }  c \zeta_t(i) 1_{(x \geq 1)}.
\ee

This way we obtain,  completely analogous to (\ref{Y13}) -
(\ref{angr15ab}) new processes which we denote \be{Y144b} (
K_t)_{t \geq 0}, (\Psi(t,ds,dj))_{t \geq 0}, \{(\zeta_s(u))_{u
\geq s}, s \in \R^+\}.\ee
Then we know that for every $t \in \R^+$ the collection
\be{ang0}
\{(\zeta_s (u))_{u \geq s}, s \in supp \Psi (t, \cdot, \N)\}
\ee
consists of {\em independent} birth and death processes.

We can show that \be{grev90d2} K^N_t \leq K_t \quad, \quad
\mbox{ stochastically for all } t \geq 0 \mbox{ and }  N \in \N.  \ee
Namely on the r.h.s. we
suppress collisions and  at each jump  occupies a {\em new} site.
Since after a collision we can have coalescence with one of the
particles already at that site, a coupling argument shows that we
have fewer particles as in the model with collisions, in fact strictly fewer
with positive probability.

Then the dual $(\eta_t, \CF^{+}_t)$ on the $N$-site model
evolves exhibiting the following features.
The dual population consists of
an increasing number, $\Pi^N_t$, of $(1_{\{1\}})$ factors (due to
birth events related to selection). Each factor $(1_{\{1\}})$ (and therefore the product)
can jump to $0$ by rare mutation  at rate $\frac{m}{N}$
(because the factor jumps to $(1_{\{2\}})$, which
becomes $ 0$ in the dual expression for $E[x^N_1(t)]$ by (\ref{dual1} with $n=1$)
since $x^N_2(0)=0$).
Therefore asymptotically as $N \to \infty$ to have non-zero probability of
such a jump to zero to occur, we need $\Pi^N_t \sim O(N)$. This means that
asymptotically the probability that a mutation event occurs
becomes positive only for sufficiently large times such that the
number of $(1_{\{1\}})$ factors in $\CF^{+}_t$ reaches $O(N)$.

In order to see the meaning of this behaviour for the original
process recall the {\em first moment  duality relation} (assuming
$x^N_1(0)=1)$ rewritten in terms of the quantities introduced
(\ref{Y13})-(\ref{angr15ab}):
\bea{angr15a}
 E[x^N_1(t)]&&=
E[\exp(-\frac{m}{N}\int_0^{t}\Pi^{N,1,1}_u du)]\nonumber\\&&
=E\left[\exp\left(-\frac{m}{N}\int_0^{t}\left(\int_0^{u}\zeta^N_v(u)
\Psi^N(u,dv)\right)du\right)\right].\eea  \sm

In order to prove Proposition \ref{P6.1b},(\ref{Y12b}), we therefore need to
show that asymptotically the probability that a mutation event
occurs by time $T_N+t$ goes to $0$ as $N \to \infty$ and $t\to -\infty$.

A complication in the analysis of the dual process
arises in that if a migrant lands at a previously
occupied site, that is, a {\em collision occurs}, then the
different birth and death processes are not independent due to
coalescence with its quadratic rate and we must take this into
account when $ K^N_t$ reaches $O(N)$ since there is then a
non-negligible probability of a collision.

We will therefore study in Subsubsections \ref{sss.dualcfr} and \ref{sss.dualcfrprop}
first the collision-free regime of the dual dynamic and in
Subsubsection \ref{sss.proofkv} its consequences for the process and then
in Subsubsections \ref{sss.collreg}-\ref{sss.limrandmarg} the
regime with collisions.

\subsubsection{The dual in the collision-free regime: the exponential growth rate}
\label{sss.dualcfr}

Due to coalescence (if $d>0$) the total number of $(1_{\{1\}})$
factors is locally stochastically bounded but due to migration the
number of occupied sites can grow and with it the total number of
factors. We need to find $\alpha$ such that the number of factors
is at time $\alpha^{-1} \log N+t$ of order $N$ for $t \geq t_0$
and $N$ large and the integrated lifetime of all factors up to
that time is $O(N)$ as $N \to \infty$, with smaller and smaller
number of factors as $t \to -\infty$ and a diverging number as
$t \to \infty$. We therefore
will fix a candidate $\alpha$ and then show lower and upper bound on the
number of factors. We will introduce this $\alpha$ as the
exponential growth rates of the number factors as $N \to \infty$ as follows.

If  $c>0$ and $N=\infty$  (i.e.
``mean-field migration''), then the number of factors (i.e.
$\Pi^\infty_t$) would grow exponentially fast as we prove below. We will now first work
with this scenario of mean-field migration in Step 1 providing a
candidate for $\alpha$ and an upper bound. In the next subsubsection
we return to the effect of
collisions, which can arise if $N<\infty$ and we have migration to a
previously occupied site to get a lower bound.

In order to verify (\ref{Y12b}), first note that by Jensen and (\ref{angr15a})
\be{d6.1}
E[x_1^N(t)]\geq \exp\left(-\frac{m}{N}E\left[ \int_0^t[\int_0^u
 \zeta_v^N(u) \Psi^N(u,dv)] du  \right]\right).
\ee
In the next three steps we analyse first $\zeta^N_v (u)$ given
$ K^N_\tau$ for $\tau$ from 0 up to $u$ to get an estimate for the
growth rate of the integral in the r.h.s. in terms of the process $K^N$ of
(\ref{d6.1}) and then in the second step we focus
on the growth of $ K^N_u$ in $u$ and in third step will then be simple, we combine
the results to get that the $\alpha$ we have chosen is an upper bound for the exponential growth
rate (in $t$) of the double-integral in the r.h.s. of (\ref{d6.1}). \bi

{\bf Step 1} $\;$ {\em Sufficient condition for (\ref{Y12b})}

We next show in order to control the r.h.s. of (\ref{d6.1}) that
\be{angr15b}
\limsup_{N\to\infty} E[\zeta^N_v(u)|\sigma\{K^N_s:s\leq u\}] \leq L
\ee
and $L$ is non-random, i.e. a constant.
Given $\alpha$ and recalling (\ref{d6.1}) in order to establish the
upper bound on the emergence time, it suffices once we have (\ref{angr15b})
to show that:
\be{Y144}
\lim_{t\to\infty}\lim_{N\to\infty} {\frac{1}{N}\int_0^{\frac{\log
N}{\alpha}-t}}E[ K^N_s]du
 =0, \ee
which we will do in the Step 2.

To prove (\ref{angr15b}) consider the collection of auxiliary processes
(here $\iota \in \R^+$)
\be{d6.1b}
(\zeta_{0,t}(\iota))_{t \geq 0},
\ee
which as it turns out bounds our process from above in a site colonized
at time 0 at time $t$ for proper choice of the parameter $\iota$.
It has the dynamic of a birth and death process
with {\em linear birth rate $sk$, quadratic death rate} $\frac{d}{2}k(k-1)$,
{\em emigration} at rate $c k1_{(k \geq 2)}$ and
{\em immigration} of a particle at rate $\iota$.  The   parameter $\iota$ we choose later.
The case $\iota=0$ corresponds to the McKlean-Vlasov dual, i.e. the dual in the collision-free regime.
The process has a unique equilibrium state
\be{Y144c}
\zeta_{0,\infty}(\iota),
\ee and is
ergodic (with exponential convergence to equilibrium).

We have (using a standard coupling argument) in stochastic order:
\be{b81}
\zeta_{0,\infty} (0) \leq \zeta_{0,\infty} (\iota).
\ee
A standard coupling argument can then be used to show that for
all $t \geq 0$ and every $s$ the random variable (recall we start in one particle)
\be{grev90a2}
\zeta_{s,t}(\iota) \mbox{ is stochastically dominated by }
\zeta_{0,\infty}(\iota).\ee

Moreover
writing down the differential equation for the second moment of $\zeta_{0,t}$ it is
easy to verify that for every $\iota \geq 0$:
\be{grev90a}
E [(\zeta_{0,\infty} (\iota) )^2] < \infty.\ee

With our McKean-Vlasov dual process we can associate a mean-field birth and death
process where the parameter $\iota$ is time-dependent, namely
the immigration $\iota_t$ is the {\em current mean} of the
process. This process arises as the limit $N \to \infty$ of the  process with migration
between $N$ sites according to the uniform distribution and starting with
a positive intensity of particles.
The mean-field birth and death process has a unique equilibrium which we
reach if we start with one particle. Namely it can be verified that the fixed point
equation

\be{grev90a33}
 \iota_\ast = E[\zeta_{0,\infty}(\iota_\ast)]\ee
has a finite solution which gives the expectation for the equilibrium of the
McKean-Vlasov dual and that in fact determines the McKean-Vlasov equilibrium.

To show this we consider the nonlinear Kolmogorov equations for
the process $(\zeta_{0,t}(\iota))_{t \geq 0}$ which are spelled out in all detail later on in
(\ref{Grev39j}) where we set in (\ref{Grev39j}) $\alpha(t)u(t)\equiv \iota$.
More precisely let
$p_{k,\ell} (t) =p_{k,\ell} (\iota,t)$
denote the transition probabilities and set
\be{koll1}
m_t (\iota) = \sum_{j=2}^\infty jp_{1,j}(\iota, t),
\ee
\be{ddy2}
m(\iota)=\lim_{t\to\infty}
\sum_{j=2}^\infty jp_{1,j}(\iota,t).
\ee
To look for equilibria we look for
fixed points $m(\iota)=\iota$. Using coupling we can show that
$m(\iota)$ is monotone increasing, continuous, ultimately
sublinear in $\iota$ and $\lim_{\iota\to 0}m(\iota)
>0$ since the smallest state is 1 if we start in 1. Therefore there exists a unique
\be{ddy2a}
\mbox{ largest solution } \iota^\ast \mbox{ of the equation } m(\iota)=\iota
\mbox{ with } m_t\leq m(\iota^\ast).
\ee

If we insert this $\iota^\ast$ as parameter the mean remains constant and
the equilibrium of the McKean-Vlasov dual process, is the unique equilibrium of the Markov process
$(\zeta_{0,t}(\iota^\ast))_{t \geq 0}$.

We see again by coupling that for $\iota \leq \iota^\ast$:
\be{b82}
\zeta_{0,\infty} (0) \leq \zeta_{0,\infty} (\iota) \leq \zeta_{0,\infty} (\iota^\ast).
\ee

We now show that this provides the  bound required in (\ref{angr15b}).  Note that  at
a tagged site in the exchangeable system of $N$ sites the number
of particles is given by

\be{ddy3}
\zeta^N_{0,t}=\zeta_{0,t}((\iota^N_u)_{u\leq t}) \mbox{ with }
\iota^N_u=\frac {\Pi^N_u}{N} \leq \frac{\Pi_u}{N},
\ee
where $\Pi_t$ is the number of particles in the McKean-Vlasov dual  process
(without collisions), which can be constructed on the same probability
space as $\Pi^N$, such that the inequality on the r.h.s. holds and such
that
\be{agdd60}
\zeta^N_{0,t} \leq \zeta_{0,t} ((\frac{\Pi_u}{N})_{u \leq t})
\leq \zeta_{0,\infty} (\iota^\ast).
\ee
 A standard mean-field limit argument
then shows that  $(\zeta^N_{0,t}, \Pi^N_t)$ converges in law as $N\to\infty$
to $(\zeta_{0,t}(\{m_{u}\}_{u\leq t}); K_t)$, where in the latter the
two components are independent.

Therefore using (\ref{b82}) and the construction of $\iota^\ast$
we see that:
\be{grev90a3}
\lim_{N\to\infty}E[\zeta^N_{0,t} |\Pi^N_s : s \leq t] =
E[\zeta_{0,t}(\{m_{u}\}_{u\leq t})]\leq E[\zeta_{0,\infty}
(\iota_\ast)] =\iota^\ast
\ee
and this proves (\ref{angr15b}) .

{\bf Step 2} $\;$ {\em A Crump-Mode-Jagers process gives a Malthusian parameter $\alpha$}.

To show (\ref{Y144})  and to obtain a candidate for $\alpha$ we can use
the collision-free mean-field dual $\Pi_t$ in particular its ingredients $K_t$ and $\zeta_t$.

The process $ K$ in (\ref{Y144b}) has some important structural
properties. Recall the process counts the number of occupied
sites, each of which has an internal state, which consists of the
number of particles at these sites. Each of the particles
exceeding the first particle
can migrate and produce a new occupied site. Whenever a new site
is colonized then an independent copy of the basic one site process starts
and evolves independently of the rest.

To analyse this object we recall the concept and some basic
results on \emph{supercritical Crump-Mode-Jagers branching
processes} and in particular its \emph{Malthusian parameter}
$\alpha$. Such a process is defined by the following properties.

The process counts the number of individuals in a branching
population whose dynamics is as follows.
Individuals can die or give birth to new individuals
based on the following ingredients:
\begin{itemize}
\item individuals have a lifetime (possibly infinite), \item for
each individual an independent realization of a point process
$\xi(t)$  starting at the birth time specifying the times at which
the individual gives birth to new individuals, \item different
individuals act independently, \item the process of birth  times is
not  concentrated on a lattice.
\end{itemize}

Let $(K_t)_{t \geq 0}$ be a process with the described structure.
The corresponding  {\it Malthusian parameter}, $\alpha
>0 $ is obtained as the unique solution of
\be{ang1}
\int_0^\infty
e^{-\alpha t}  \mu (dt) =1 \mbox{ where } \mu([0,t]) =
E[\xi([0,t])],
\ee
with
\be{ang2}
\xi(t) \mbox{ counting the number of births of
a single individual up to time } t,
\ee
(see for example \cite{J92}, \cite{N} equation (1.4)).
\sm

In our case we have
 \be{ang1ab} \mu([0,t]) = c \int_0^t
E[\zeta_0(s) 1_{(\zeta_0(s) \geq 2)}]ds. \ee

\begin{remark}\label{DMP}
Note that in our model particles at singly occupied sites do not
migrate. We can suppress the emigration step if only one particle
is left, since in that case everything would simply start afresh
at some other site.  We have done this so that in the
pre-collision regime the number of occupied sites equals the
number of migrating particles. In this case the rate of creation
of new sites due to migrations from a given site is given by
$\zeta_0(s)1_{\zeta_0(s)\geq 2}$ and the death rate when
$\zeta_0(s)=1$ is $0$.

However we can also consider the original
model in which particles at singly occupied sites do migrate. Then
sites have {\em finite lifetimes}  and the rate of production of
new sites  by a fixed site  at time $s$ is given by
$\wt\zeta_0(s)$ where now $\wt\zeta_0(s)$ has death rate $c$ when
$\wt\zeta_0(s)=1$. Then we have instead of (\ref{ang1ab}) the equation
\be{ang1abx} \mu([0,t]) = c
\int_0^t E[\wt \zeta_0(s)]ds.
\ee
However the Malthusian parameter
describing the exponential growth of the number of sites occupied
by at least one particle is the same in both cases as is shown by
explicit calculation.
\end{remark}

We therefore define $\alpha$ by (recall for $\zeta_0$
(\ref{Y144b}) and above): \be{ang1ab1} c \intl^\infty_0 e^{-\alpha
s} E[\zeta_0(s) 1_{(\zeta_0(s) \geq 2)}] ds =1. \ee

The first thing we want to know is that in our case \be{angr15x}
\alpha \in (0,s)\quad\text{if}\quad d>0 \mbox{ and } \alpha \in
(0,s] \mbox{ if } d \geq 0. \ee In the case of our dual process $
K_t \leq N_t$ where $N_t$ is a rate $s$ pure birth process and
therefore $\alpha \leq s$.

To show that the Malthusian parameter
satisfies: \be{grev90c} \alpha >0 \ee note first that when $d>0$
still no deaths (in the number of occupied sites) occur and that a
lower bound to the growth process is obtained by considering only
those migrants which come from the birth of a particle at sites
containing sofar only one particle and which then migrates before
coalescence or the next birth occurs. This process of the number of
such special sites is a classical birth and death process with
in state $k$ birth rate $s(\frac{c}{c+s+d}) k $ and death rate zero
which clearly has a positive
Malthusian parameter (if $c>0$) and hence (\ref{grev90c})
holds.

Next to see that if $d>0$ then $\alpha < s$, note that a
non-zero proportion of the new particles generated by a  pure birth
process at rate $s$ dies due to coalescence
and therefore $\alpha < s$ in this case.

If we know the Malthusian parameter, we need to know that it is
actually equal to the almost sure growth rate of the population.
It is known for a CMJ-process $(K_t)_{t \geq 0}$ that (Proposition 1.1
and Theorem 5.4 in \cite{N}) the following basic growth theorem
holds. If
\be{ang1b}
E[X \log (X \vee 1)] < \infty, \ee
where here $X =\intl^\infty_0 e^{-\alpha t} d \xi(t) $, then
\be{ang2b}
 \lim_{t\to\infty} \frac{K_t}{e^{\alpha t}} = W,  \mbox{ a.s. and in } L_1, \ee
where $W$ is a random variable which has two important properties, namely
\be{ang4}
W>0,\;a.s.$ and $E[W]<\infty. \ee
We now have to check that (\ref{ang1b}) holds in our case.

We know using (\ref{ang1ab1}) that $\intl^\infty_0 u e^{-\alpha u}
\mu(du) < \infty$, since we can bound the birth rate by the one in
which the number of particles  at each site is given by the
equilibrium state $\CL(\zeta_\infty)$, which has a finite
mean. This also implies that assumption (\ref{ang1b}) holds.
Namely, in equilibrium we can estimate
\be{grev90d}
E[X^2] \leq \intl^\infty_0 \intl^\infty_0 e^{-\alpha(s+t)}
 E[d \xi(s) d \xi(t)]  < \infty,
 \ee
by using Cauchy-Schwartz and using (\ref{grev90a}).
\sm

\begin{remark}
In the above we  considered the dual $\mathcal{F}^{+}$ starting
with one factor (particle) at one initial site.  We can also
consider the case with $k$ initial particles at one site.  In this
case the same exponential growth occurs since (\ref{ang1ab}) still
holds but the corresponding random growth factor $W^{(k)}$ in
(\ref{ang2b}) now depends on $k$.  Similarly we can start with $k$
particles at each of $\ell$ different sites and again get
exponential growth with the same Malthusian parameter $\alpha$ but
now with random variable $W^{(k,\ell)}$. This will be frequently
used in the sequel of this section.
\end{remark}

{\bf Step 3} $\;$ {\em Completion of argument}

To complete the argument we
note that (\ref{Y144})  follows from(\ref{ang2b})
(which implies that $K_t =e^{-\alpha t} W_t$ with $W_t \to W$) and
(\ref{grev90d2}) since

\be{grev90.a}
E(\int_0^{\alpha^{-1}\log N -t} K^N_s ds) \leq
Const \cdot \frac{1}{\alpha} E(W) e^{{\log N}-{ \alpha} t}.
\ee We saw earlier
that (\ref{Y144}) is sufficient for (\ref{Y12b}) and the proof is
complete that emergence does not happen before time $\alpha^{-1}
\log N$ as $N \to \infty$. \sm

\subsubsection{The dual in the collision-free regime: Further properties}
\label{sss.dualcfrprop}

We shall need for the analysis of the regime with collisions some
further properties concerning the dual in the collision-free
regime. First, we obtain further properties of the growth constant
$W$, more precisely its moments, and then the stable age
distribution which among other things allows to give an explicit
expression for the growth rate $\alpha$.
The random variable $W$ is defined in (\ref{ang2b}) in terms of
the process $(K_t)_{t \geq 0}$ introduced in (\ref{Y144b}), that
is, we consider  the dual process without collision.
\medskip

{\bf Step 1}  $\;${\em Higher moments of $W$.}

In the sequel we will perform moment calculations which use:

\noi \beL{L6.11}{(Higher Moments of $W$)}

For all $n$, \be{ga40}
E[W^n]<\infty.  \qquad \square \ee
\end{lemma}

\begin{proof}{\bf of Lemma \ref{L6.11}}  By  \cite{BD},
Theorem 1, $W$ has a finite $n$th moment if and only if
\be{ga40b}
 E \left[\left( \int_0^\infty e^{-\alpha s} d\xi(s)\right)^n\right] < \infty,
 \ee
where $\xi( t)$ denotes the number of
particles that have emigrated by time $t$ from a fixed site starting with 1
particle.  We will verify this condition for the second moment.
\bea{ga40c}
&&E[\left[\int_0^\infty e^{-\alpha s}
d\xi(s)\right]^2=2E\left[\int_0^\infty\int_s^\infty e^{-\alpha
s}e^{-\alpha (s+s^\prime))}d\xi(s^\prime)d\xi(s)\right]\\
&& =
c^2E\left[\int_0^\infty\int_s^\infty e^{-\alpha (s+s^\prime)}
\zeta_0(s)E[\zeta_0(s^\prime)|\zeta_0(s)]ds^\prime ds\right]\nonumber\\
&& =
c^2E\left[\int_0^\infty\int_s^\infty e^{-\alpha s} \zeta_0(s)
\zeta_0(s^\prime)ds^\prime ds\right]. \nonumber
\eea
Recall that here $\zeta_0$
is a birth and death process with birth rate $sk$ and death rate
$ck 1_{k \geq 2} + \frac{1}{2}dk(k-1)$. But then by explicit calculation
$\sup\limits_s E(\zeta_0(s)^2)<\infty$,
and therefore
\be{ga40d}
2E \left[\int_0^\infty\int_s^\infty e^{-\alpha s} \zeta_0(s)
\zeta_0(s^\prime)ds^\prime ds\right]<\infty.
\ee
This completes the proof  that $E(W^2)<\infty$.

Recalling that $\zeta_0(t)$ is for all $t$ stochastically
dominated by the equilibrium distribution $(p_k)_{k \in \N}$ of this birth and
death process and that the latter satisfies for some $C<\infty$,
\be{ga40e}
 p_k\leq \frac {C^k}{k!}, \qquad\text{for all } k=1,2,\cdots
 \ee it follows that
 for every $n\in\mathbb{N}$,
 \be{ga40f}
 \sup_s E[(\zeta_0(s))^n] <\infty.\ee
Then a similar argument to the second moment case argument
given above verifies that $n$th moments of $W$ are finite.

\end{proof}

\begin{remark}

Alternatively to the above we can argue as follows.
In the pre-collision regime we can also represent the particle
process as a branching Markov chain, more specifically a {\em
branching birth and death process} (linear birth rate and
quadratic death rate), that is a counting measure-valued process,
on $\mathbb{N}$, \be{ga43a} \Lambda_t\in \CM(\mathbb{N}), \ee
where $\Lambda_t(k)$ counts the number of sites occupied by
exactly $k$ partition elements. In other words the individuals of
this new process $\Lambda$ correspond to the occupied sites of the
dual particle system and the location of individuals is now the
size of the population at the sites of the dual particle system.

The dynamics of $(\Lambda_t)_{t \geq 0}$ are therefore given by
\begin{itemize}
\item The branching. Individuals branch independently with a
location-dependent branching rate, namely,  an individual at $k$ dies at
rate $ck$ and produces two new individuals with new locations, namely
one individual at $k-1$ and a second one at
 $1$. Since sites of size $0$ that would be produced if a site
of size $1$ dies would not contribute to the production of new
sites, they can be ignored.  For this reason we suppress these
deaths and set the death rate at $k$ to be $ck1_{k\ne 1}$.

\item The motion:  the motion between locations is given by a
Markov process $(\wt \zeta(t))_{t \geq 0}$ on
$\mathbb{N}\cup\{\infty\}$, namely, a birth and death process
$(\zeta(t))_{t \geq 0}$, where the birth rate is $sk$ and the
death rate $d k(k-1)/2$.
\end{itemize}

In other words the CMJ process can be
embedded in the branching birth and death process and
the number of individuals in that process is the same as the number of occupied sites
in our original CMJ-process.

We work with the probability generating function, $M_t$, of
branching birth and death process $ {\Lambda}_t$ and use this to
compute the distribution of the number of  occupied sites as
follows.

 Define the generating function:
\be{dd15}
M_t(k,z) =
E_{\delta_k}[z^{\Lambda_t(1_{\{\mathbb{N}\}})}]=\sum_{j=1}^\infty
p(t,k,j)z^j,\quad 0<z\leq 1, \ee where we write
$\Lambda_t(f)=\sum_{i=1}^\infty f(i)\Lambda_t(k)$.

Then $M_t$ satisfies the standard functional equation

 \be{dd16}\begin{array}{ll}
 M_t(k,z) =& E_k[ e^{-c\int_0^t(\wt\zeta(s))ds}]z\\
   &+E_k\left[\int_0^t e^{-c\int_0^s\wt\zeta(u)du} \cdot
   c \wt\zeta(s) \cdot M_{t-s}(\wt\zeta(s)-1,z)M_{t-s}(1,z)ds\right],
   \end{array}\ee
with $M_t(0,z)\equiv 1$. Here $\wt \zeta$
is a birth and death process with birth rate $sk$ and death rate
$\frac{1}{2}d k(k-1)$. The above equation can be solved recursively for the
coefficients of $z^j$.

This relation  can also be used to obtain moment formulas by
evaluating the appropriate derivatives at $z=1$.

Consider the semigroup $(T_t)_{t \geq 0}$ of the migration defined by
\be{dd16z}
 T_tf(k)=E_k[e^{-c\int_0^t\wt \zeta(s)ds}f(\wt\zeta(t))], \quad k \in \N \ee
and let \be{dd16y}  f_\ell(j):= j^\ell,\; \ell=1,2,3,\dots.\ee

Let for $\ell=1,2,3,\cdots$ and $k \in \N$: \be{dd16c} m_\ell(t,k)
= E_{k} ((\Lambda_t (1_{\{\mathbb{N}\}}))^\ell). \ee Then the
following moment relations hold:

\be{dd17}
  m_1(t,k)= T_t f_0(k) + c\int_0^t T_s \{f_1(\cdot)[m_1(t-s,1) + m_1(t-s,\cdot -1))]\}(k)ds.
\ee

\bea{dd18}
 m_2(t,k)=&& T_t f_0(k) + c\int_0^t (T_s\{ f_1(\cdot) \cdot[(m_2(t-s,\cdot-1) +m_2(t-s,1))\\&&\qquad +(m_1(t-s,1)\cdot
 m_1(t-s,\cdot -1))]\}(k)ds.\nonumber
 \eea

 Also note that $m_1(t,k+j)\leq m_1(t,k)+m_1(t,j)$ so that
$m_1(t,k)\leq k m_1(t,1)$. Therefore we get from (\ref{dd17})
\be{dd20}
 m(t,1)\leq T_t1 + c\int_0^t (T_sk)(1)(m_1(t-s,1)ds
+c\int_0^tT_s(k(k-1) m_1(t-s,1))(1)ds.\ee Let $M(\lambda)$,
$F(\lambda)$, $G(\lambda)$ and  denote the Laplace transforms of
$m_1$, $T_s1$, $T_sk+T_s(k(k-1))$ as functions of $s$, respectively.

Then \be{dd21} M(\lambda)(1-c G(\lambda)) \leq F(\lambda).\ee We
need $\lambda > \lambda_0$ so that this is invertible. This gives
an upper bound for the Malthusian parameter $\alpha$.
\end{remark}

{\bf Step 2}  $\;${\em Stable age and size distribution}

Return to the McKean-Vlasov version of the dual process (which is
the dual in the collision-free regime in the limit $N \to
\infty$). If we consider for every time $t$ for each individual
(i.e. occupied site in our case) currently alive its age and size,
then we can introduce a random probability measure, the normalized
empirical age and size distribution of the current population,
which we denote by \be{ang4b} \CU(t,du,j)=\frac{\Psi (t,du,j)}{
K_t}.\ee The marginal random measure $\CU(t,du,\mathbb{N})$
converges in law (we use the weak topology on measures) as $t \to
\infty$  to a {\em stable age  distribution} \be{ang4c}
\CU(\infty,du, \N) \mbox{ on } [0,\infty), \ee according to
Corollary 6.4 in \cite{N}, if condition 6.1 therein holds.
The condition 6.1 in \cite{N}  or (3.1) in \cite{JN} requires
that (with $\mu$ as in (\ref{ang1})):
\be{ang4c1} \int_0^\infty e^{-\beta t}\mu(dt)<\infty \quad
\mbox{ for some } \beta\geq 0. \ee The condition 7.1 in
\cite{N} holds for any $\beta
>0$ since the  local (at one site) population of the dual McKean-Vlasov particle system
$(\wt \eta_t)_{t \geq 0}$ given by
$(\zeta_{0,t})_{t \geq 0}$ goes to a finite mean
equilibrium. \bi

Since the
distribution of size at a site depends only on the age of the
site, it follows that $ \CU(t,du,\cdot)$ converges to a {\em stable  age
and size} distribution, i.e.
\be{ang4c2}
\CU (t,\cdot,\cdot) \La \CU(\infty,\cdot,\cdot), \mbox{ as } t\to \infty \mbox{ in law}.
\ee

We can now also calculate $\alpha$  and  the asymptotic density of
the total number of individuals $B$ as follows unsing a law of large number effect for
the collection of independent birth and death processes at the occupied sites. Namely
the frequencies of specific internal states stabilize in the stable size distribution,
while the actual numbers diverge with the order of $K_{t_N}$ as $N \to \infty$. Given site $i$ let
$\tau_{i}\geq 0$ denote the time at which a migrant (or initial
particle)  first occupies it. Noting that we can verify Condition
5.1 in \cite{N} we have that
\be{stableage}
\lim_{t\to
\infty}\frac{1}{ K_t}\sum_{i=1}^{ K_t} \zeta_{\tau_i}(t-\tau_i)
=\intl^\infty_0 E [\zeta_0(u)] \CU(\infty, du)=  B \text{  (a
constant)}, a.s.,
\end{equation} by Corollary 5.5 of
\cite{N}. The constant $B$ in (\ref{stableage}) is in our case given by the
average number of particles per occupied site and the growth rate $\alpha$
arises from this quantity neglecting single occupation. Namely define

\be{ang4c4}
\alpha =c\sum_{j=2}^\infty j
\mathcal{U}(\infty,[0,\infty),j)<\infty ,\qquad
\gamma=c \; \mathcal{U}(\infty,[0,\infty),1).\ee
Then
\be{ang4c3}
B = \frac{\alpha+\gamma}{c}.
\ee

Furthermore  the average birth rate of new sites (by arrival of a
migrant at an unoccupied site) at time $t$ (in the process
in the McKean-Vlasov dual) is equal to
\be{ang4c5}
\alpha = c \; B - \gamma.
\ee

\subsubsection{Dual process in the collision regime: macroscopic emergence\\
Proposition \ref{P6.1b}(a)} \label{sss.collreg}

We have obtained in Subsubsection \ref{sss.dualcfr} a lower bound on the emergence time by an upper bound
on the number of dual particles. Now we need an upper bound on the
emergence time via a lower bound on the number of dual particles.

In order to prove
(\ref{Y12c}) return to the dual particle system specified below
(\ref{grev89aa2}). Then after the lower bound for the emergence
time in Subsubsection \ref{sss.dualcfr}, where we ignored
collisions in the dual particle system, we derive here
a lower bound on the growth of the number
of occupied sites in this system starting with one occupied site
incorporating the effect of {\em collisions}.
We proceed in five steps, first we consider lower bounds on the number
of particles in the dual process, then refine this in a second step by
constructing a multicolour particle system which is an enrichment of the dual
particle system, then in a third step we prepare the estimation of the
difference between collision-free and dual system, and in the
fourth step we turn this into an upper bound for the emergence
time. In Step 5 we then show the convergence of the hitting times for
reaching $\lfloor \ve N\rfloor$ dual particles.

{\bf Step 1: $\;$ Number of sites occupied by  the dual process: preparation}

Here we must take into account the
effect of collisions and show that (\ref{ang1ab1}) indeed gives
the correct $\alpha$ which describes the growth of the number
of factors in the dual process.

If we have collisions, that is, migration to occupied sites,  we
have two effects. (1) We have {\em interaction} between sites. (2) If a
particle migrates from a site containing only one particle to an
occupied site, this results in a {\em decrease} by one in the number of
occupied sites. For these two reasons the process counting the number of occupied
sites, denoted by $ K^N_t$, is no longer approximated by a
nondecreasing Crump-Mode-Jagers process (which is based on independence of internal
states) when it reaches size $O(N)$ and at this moment it is also no
longer non-decreasing since now singly occupied sites can dissappear
by a jump to an occupied site.

To handle these two problems the key idea  is to carry out the analysis
of the dynamics separately in three time intervals defined as
follows. Set first
\be{grev90a4}
\tau_{\log N}=\inf \{t: K_t= \lfloor \log N\rfloor\},\quad
 \tau^N(\ve)= \inf \{t: K_t^N= \lfloor \ve N \rfloor \}.
 \ee
Then define the three time intervals as:
\bea{grev90a5}
&&[0,\tau_{\log N}), \nonumber \\
&&[\tau_{\log N},\tau_{\log N}+\frac{1}{\alpha}(\log N-\log\log N)), \\
&&[\tau_{\log N}+\frac{1}{\alpha}(\log N -\log\log N), \tau^N(\ve)]. \nonumber
\eea

What happens in the {\em first interval}? Note (since $\tau_{\log
N} <<\tau^N(\ve)$ for $N \to \infty$) that asymptotically as $N
\to \infty$, in a single migration step the probability of a collision is
$O(N^{-1} \log N)$. In particular over a time horizon of length
$o(N/\log N)$ the probability to observe a collision ever goes to zero
as $N \to \infty$. Hence
looking at the formulas we obtained for $\alpha$
and with the properties of exponentially growing populations
(occupation measure in current size)
we can replace our dual with the mean-field
dual without collisions so that with (\ref{ang2b}) we can conclude
as $N \to \infty$: \be{grev90a5b}
 \tau_{\log N}= \frac{\log\log N}{\alpha}-\frac{\log W}{\alpha}
 + o(1) \ee
and  we know furthermore that this process has
(asymptotically) reached the stable age distribution at time $\tau_{\log
N}$  according to the analysis of Subsubsection \ref{sss.dualcfrprop}, Step 2 therein.

Now consider the {\em second interval}. From the previous upper
bound result still no collisions occur comparable to the dual population
in this interval (the collision
probability is now at most $O((\log N)^{-1})$ and hence at most only finite many
collisions may occur over the whole time interval, essentially at the end and we can
assume that asymptotically as $N\to\infty$ we can replace our dual
process again by the collision free mean-field dual at the left end
point but now in the {\em whole time span} the stable age distribution is in effect.
Hence we can in the second interval continue working with the
mean-field dual without collisions and we can even assume the stable age
distribution in effect.

Returning to our dual process we conclude (recall (\ref{grev90a5})) that
at the end of the second interval we
have  as $N \to \infty$: \be{grev90a6}
  K^N_{\tau_{\log N}+\frac{1}{\alpha}(\log N-\log\log N)} \sim W \frac{N}{\log N}.
\ee

\begin{remark}
Note that this means that $\tau^N_{s_N}-\tau^N_{\log N}$, with
$s_N>> \log \log N$ but $s_N <<\log N$, is
asymptotically as $N \to \infty$ deterministic.
\end{remark}

We now consider the {\em third interval}, i.e. the regime in which
by (\ref{grev90a6}) above
\be{Grev37} K^N_t \geq W \frac{N}{\log N}. \ee
 In this
interval two new effects must be considered, namely,

\begin{itemize}
    \item[(A)] a decrease in the number of occupied sites when a lone
particle at an occupied site jumps to another currently occupied
site
    \item[(B) ] a particle moves from a site occupied by more than
one particle and jumps to another occupied site.
\end{itemize}

We note that effect (A) also increases the age distribution of
occupied sites since in a site the process is stochastically increasing
and the young sites can disappear therefore easier
(and hence the age of the total population becomes stochastically larger). Effect
(B) tends to increase the size of occupied sites of a given age
compared to the case without collision.
Therefore  the result of (A) and (B) is to tend to increase the
average size of occupied sites compared to a system without
collision (recall the older sites are stochastically larger),
but decreases the number of occupied sites.

We note that larger sites
increase the number of migration steps. Altogether
this means that we expect that the rate of growth of $ K^N_t$
denoted $\beta_N(t)$ (i.e. $\beta_N(t)=(K^N_{t_{+ \triangle t}} - K^N_t)/K^N_t)$
satisfies in the third time interval (for $N
\geq \log\log N$) if the present state of $K^N_t$ is $k$:
\be{Grev38}
\beta_N (\cdot)\geq\alpha(1-\frac{k}{N}).\ee
In order to control these effects precisely we use the technique of
coupling constructed from an enriched, i.e. multicolour particle system, which we introduce next.
\sm

{\bf Step 2: A multicolour particle system}

To examine rigorously the growth in the time span, in which the number of occupied sites
reaches $O(N)$ we construct a coupled system of {\em multicoloured particles}
such that we can compare the new system with collisions effectively with the
simpler one without collisions, and
are able to carry out estimates on the difference. This way we can obtain the lower bound
on the number of occupied sites and the total number of particles.

First an informal description. We shall have {\em white, black} and {\em red}
particles. We start with white particles only. \ In intervals 1
and 2 the process of the white particles grows as before, i.e.
without collisions. However as soon as collisions occur we
consider a {\em modified system}. The colouring of the system
allows us to keep track of events like (A) and (B) given above.
For example in the case of (A) we will mark the lost site by
placing there a black particle and the  particle that jumped
and collided with other particles on the new location is
marked by giving it a red colour. We handle the second effect (B)
 by using bounds from below on the rate of founding of new
 occupied sites.

Formally proceed as follows.
The multicolour comparison system has black, white and red particles. It has
white  and  red particles
located at the sites $\{1,2,\cdots,N\}$ and black particles located at a
site in $\mathbb{N}$ where $\{1,2,\cdots,N\}$ and $ \mathbb{N}$
are disjoint finite and countable sets respectively.
In other words the geographic space of this new system is
\be{Grev38a}
\{1,2,\cdots,N\} + (\N),
\ee
and the state space is
\be{agrev89}
(\{1,2,\cdots,N\} + (\N))^{\N^3_0}.
\ee

The initial state is given by having only white particles, which
are located at sites in $\{1,2,\cdots,N\}$ such that occupation
numbers are exchangeable on this part of the geographic space.

The dynamics of the new system is markovian more precisely it is a
pure Markov jump process and we specify the transitions and their rates.
Instead of writing the generator we describe this more intuitively in words.
This runs as follows:

\begin{itemize}
\item {\em white} particles at a site follow the same local dynamics as
the dual particle system $\eta$ as far as birth (of white particles)
and coalescence (of white particles) goes, changes occur for migration.

Let $k$ denote the number of sites having currently at least one white particle.
Each migrating  white particle moves with probability
$1-\frac{k}{N}$ to a new  site in $\{1,2,\cdots,N\}$ which prior
to this event did  not contain any white particles and with
probability $\frac{k}{N}$ changes to a black particle now located
at a new unoccupied site in $\mathbb{N}$ and at the same time also
a red particle is produced at an occupied site in
$\{1,2,\cdots,N\}$ chosen at random among the $k$ occupied sites.

Hence a migrating white particle produces:
\be{Grev38b}
\mbox{a new site occupied with a white particle with probability }
(1-\frac{k}{N}),
\ee and two  particles,
\be{Grev38b2}
\mbox{ one black and one red,  with probability } \frac{k}{N},
\ee
\be{Grev38b3}\begin{array}{l}
\mbox{the red particle is placed at a randomly chosen occupied site in $\{1,2,\cdots,N\}$,}
\\
\mbox{the black one at the smallest free site in $\N$}.
\end{array}
\ee

\item \label{Grev38b4}{\em Red} particles have the same dynamics as the
dual particle system $\eta$ on $\{1,2,\cdots,N\}$ (newborn particles are
also red) and in addition
when a red and white at the same site coalesce the outcome is
always white.

\item {\em Black} particles follow the same dynamics as the white except
that migrating black particles move on $\N$ and always go to a new,
so far unoccupied site in $\mathbb{N}$.

\end{itemize}

This means that no collisions (by migration!) occur among white particles and
furthermore by the collision convention the white particles are
not influenced by the presence of red particles, nor are they influenced
by the black particles.
Hence the key observation about the new system is that:
\begin{itemize}

\item The number of occupied sites in the {\em union} of the {\em black and white} particles follows the
dynamics of the process without collisions, i.e. the number of
sites they occupy is a version of $(K_t)_{t \geq 0}$,

\item the number of occupied sites in the {\em union} of the {\em white and red} particles follows the exact
dual dynamics, i.e. the number of sites they occupy is producing
a version of $(K^N_t)_{t \geq 0}$.

\item
The process of white and black particles per site follow the dynamics
of $\zeta$ and the one of white and red particles per site that of $\zeta^N$.

\item We have a coupling of $\zeta$ and $\zeta^N$ given by the
embedding in the multicolour system and the difference process
$\zeta_t - \zeta^N_t$ can be represented as the number of black
minus the number of red particles at the various sites. In
particular also $K, K^N$ and $K-K^N$ can be represented in terms
of the multicolour system.
\end{itemize}

This construction allows to give upper and lower bounds on the time
to reach with the dual process $\lfloor \ve N \rfloor$ occupied sites, which then
in turn gives bounds on the number of factors in the dual process.

Since we know from the last subsubsection how
the growth of the population of both black and white  particles works,
we use this as an upper bound on the number of white particles (conditioned
on $W$).

We denote the {\em number of sites occupied only by black
particles} (recall ''only'' is by the construction not a  constraint),
respectively the number of sites occupied by either
{\em white or black particles} at time $t$ by \be{bw2} \bar Z^N_t,
\mbox{ respectively } \wh K_t. \ee
Note that $\wh K$ is a version of $K$.

We shall get below an
upper bound on the number of sites occupied by black particles
and then we can use this upper bound to get
a {\em lower bound} on the number of sites occupied by white
particles and therefore a lower
bound for the number of factors of the dual (red plus white
particles).

Let \be{bw3} \wt\tau^N(\ve) =\inf \{t:\wh K_t\geq \ve\ N\}. \ee
This means in particular that
\be{w4}
\wt \tau^N(\cdot) \mbox{ is nondecreasing}.
\ee
(Note the difference between $\wt \tau^N(\ve)$ and $\tau^N(\ve)$
is that we use in the first case white and {\em black} particles and in
the second case  white and {\em red} particles. In addition $\wh K$
is a version of $K$.)

We just saw above that we have to estimate the number of black particles
to compare the dual and the dual of the McKean-Vlasov limit which differ by the difference
of black and red particles. The latter is estimated above stochastically
by the number of black particles. This holds since if we associate the black
and red upon their creation, we see that the red ones have an extra
chance to disappear once they coalesce with a white particles or
collide with red particles from another black-red creation event.

{\bf Step 3: Estimating the number of black particles.}

An upper bound on the number of black sites (recall black particles sit on
$\N$, while white and red sit on $\{1, \cdots, N\}$)
at time $\wt \tau^N(\ve)$, is
obtained as follows. Bound the number of migration steps in the model
by a Crump-Mode-Jaegers process which is in distribution given by
the process of black and white particles. This means we construct up to
terms of order $o(N)$ a stochastic
upper bound for the number of black particles by realizing {\em independently}
a Poisson stream of potential collision events according to the law of a driving CMJ-process
generating potential new sites and
with a collision probability given by the current number of white sites divided by $N$.

Therefore an upper bound
on the production rates of new black populations is obtained
by integrating the production rate of new sites generated by the process of
white and black particles at a given time (which is for large
times converging to $\alpha$ by (\ref{ang2})). Furthermore note that in
state $k$ (meaning the number of sites with at least one white
particle)  with a migration step of a white particle, the
probability of a collision and hence the probability
to become black is $k/N$. This gives at time $s$ with the CMJ-theory a
production rate of collisions which is stochastically bounded by
\be{w4b}
\alpha \frac{W_s e^{\alpha s}}{N} W_s e^{\alpha s}.
\ee
Hence we get a stochastic bound on the
intensity of migration steps resulting in collisions
of the form $W^2_s \alpha (e^{\alpha s}/N)e^{\alpha s}$.
This expression therefore bounds the rate of creation of {\em founders} of new
black families of particles.

We can replace for $s \geq \tau_{\log N}$, the quantity $W^2_s$ in the limit
$N \to \infty$ by $W^2$, due to the convergence theorem for CMJ-processes.
Therefore we can  estimate the rate at which new black families are created
asymptotically as $N \to \infty$ and work with the  expression
\be{w4b1}
W^2 \alpha e^{2 \alpha s}/N \mbox{ for } s \geq \tau_{\log N}.
\ee
The contribution to the stream of production of collisions arising before
time $\tau_{\log N}$ goes to zero in probability as $N \to \infty$ and therefore
the probability that black families are founded before time $\tau_{\log N}$ tends
to zero as $N \to \infty$, more precisely at rate $(\log N)/N$.

In principle the time we consider, i.e. from $\tau_{\log N}$
to $\wt \tau_N(\ve)$ seems random on first sight, however it is
asymptotically deterministic.
Namely
$\wt \tau_N(\ve) - \tau_{\log N}$
becomes deterministic in the sense that
\be{w4b2}
(\wt \tau_N(\ve) -\tau_{\log N}) - (\frac{1}{\alpha} \log (\frac{\ve N}{W})
- \frac{1}{\alpha} \log (\frac{\log N}{W}))
\Ntoo 0, \mbox{ a.s.}
\ee
so that conditioned on $W$ as $N \to \infty, \wt \tau_N(\ve) -\tau_{\log N}$
can be replaced by a deterministic quantity.
This allows us to generate as upper bound on the number of black particles
the following process.

(1) Realise a Poisson point process with intensity
measure
\be{intmeas1}
(W^2 \alpha e^{2 \alpha s}/N) ds.
\ee

(2) Then let from each point of this Poisson point process evolve independent
families of black particles
with the usual dynamics but independent of the Poisson point process.

Next observe that a new black founding particle generated by the Poisson point process
has a {\em descendant black population} which forms by construction
a Crump-Mode-Jagers process with growth rate $\alpha$, i.e. it grows
like $\wt W_{(\wt \tau^N (\ve)-s)} e^{\alpha(\wt \tau^N (\ve)-s)}$,
if $s$ is the birth time of the black particle and provided we
observe up to the final time $\wt \tau^N (\ve)$.
Here $\{\wt W_s\}_{s\in \mathbb{R}}$ are  independent copies of $W$.

The newly founded black populations created at the Poisson point process jump times observed at time
$\wt \tau_N(\ve)$ have sizes
\be{w4c}
\{\wt W_{(\wt \tau_N(\ve)-s),N} \exp (\alpha(\wt \tau_N(\ve)-s)), \quad s \in \{s^N_1, s^N_2,\cdots,s^N_{M(N)}\}\},
\ee
(we will suppress the subscript $N$ for $s^N_i$ in formulas below)
and are independently distributed conditioned on the process of birth-times
and $\wt \tau_N(\ve)$. For these copies we have a law of large numbers acting.
The claim is more precisely that if we take the expectation $\wt E$
over $\{\wt W_{\wt \tau_N(\ve)-s_i,N} , \quad i \in \N\}$ then as $N \to \infty$ we have the
following law of large numbers effect:
\be{w4d}
\suml^{M(N)}_{i=1}  \wt W_{\wt \tau_N(\ve)-s_i,N} \exp (\alpha(\wt \tau_N(\ve)-s_{i,N}))
\sim \intl^{\wt \tau_N(\ve)}_0 \alpha W^2
(e^{2 \alpha s}/N) \wt E[\wt W_{(\wt \tau_N(\ve)-s),N}] e^{\alpha(\wt \tau_N(\ve)-s)}ds
\ee
and the r.h.s. is asymptotically as $N \to \infty$ equal to
\be{w4e}
E[W] \intl^{\wt \tau_N(\ve)}_0 \alpha W^2
(e^{2 \alpha s}/N) e^{\alpha(\wt \tau_N(\ve)-s)}ds,
\ee
using the convergence theorem for the Crump-Mode Jagers process.

Then combining (\ref{w4c}) and (\ref{w4d}) together with (\ref{w4e}) we
obtain a {\em mean conditioned on W}
of the l.h.s. of (\ref{w4d}) which is asymptotically as $N \to \infty$ equal to:
\be{w4f}
[\intl^{\wt \tau^N(\ve)}_{\tau_{\log N}} E[W] \alpha W^2
e^{\alpha \wt \tau_N(\ve)} \cdot e^{\alpha s} ds]
\sim  E[W] W^2 \frac{1}{N} (\frac{1}{W^2} \ve^2 N^2)
= E[W] \ve^2 N.
\ee

To justify the law of large numbers given in (\ref{w4d}), we observe that the number of birth events
of black populations goes to zero for times $t_N$ with
\be{w4h}
t_N - \frac{\log N}{2 \alpha} \Nto -\infty
\ee
and for these times the order of magnitude of a descending population
is at most of order $\sqrt{N} = o(N)$. Therefore we get a diverging
number of contributions which are all $o(N)$.

More precisely, the law of large numbers is verified by showing that
\be{w4h2}
\wt Var \left[\sum \wt W_{(\wt \tau_N(\ve)-s_i),N} e^{\alpha(\wt \tau_N(\ve) - s_i)}\right]
= O(\frac{N^2 \log N}{N^3}) = O(\frac{\log N}{N}).
\ee
Here we use the explicit representation and the fact that
$Var (W_s)< \infty$ and $Var(W_s)\to Var (W)$ as $s \to \infty$
together with the asymptotics of
$\wt \tau_N(\ve)$. This concludes the proof for the law of large numbers.

Hence we can bound the number of black particles by the asymptotic growth in $N$ of
its expectation over the randomness in the process of growth, i.e. the
$\wt W_{s^N_i}$ in the black population.

This reasoning gives
the following asymptotic stochastic
upper bound on the number of black sites, namely the expression:

\be{deadsites}
W^2 \int_{\tau_{\log N}}^{ \wt\tau^N(\ve)}
\alpha\frac{e^{\alpha s}}{N} e^{\alpha s}
E [\wt W_s] e^{\alpha(\wt\tau^N(\ve)-s)}ds\;
{\stackrel{<}{ \sim}}  \; E[W] N \ve^2, \mbox{ as } N \to \infty,
 \ee
using the upper bound on $\wt \tau^N (\ve)$ implied by (\ref{ang2b}).

Therefore we conclude that the {\em proportion of black sites} among all sites
in the black and white system at time
$\wt \tau^N(\ve)$ is at most (in the limit $N \to \infty$) equal to
\be{w4g}
\frac{E[W]}{W} \ve.
\ee
In particular by letting $\ve \to 0$ this relative frequency goes to 0.

In particular this would imply that the
calculation below of $\tau^N(\ve)$ (the time for the number
of white plus red sites to reach $\ve$)  will give a stochastic
lower bound for
\be{angr17o}
\tau^N(\ve- const \cdot \ve^2),
\ee
if we condition on $W$.

\bigskip


{\bf Step 4: Upper bound on emergence time}

Recalling (\ref{Grev38}) we see that given the current state $k$
of $K^N_t$, the time of the next birth in $ K^N_t$ is a
positive random variable denoted  $\tau_k$ for which we  are
waiting for the next migration step of a white or red
particle to an unoccupied site. This waiting time is {\em bounded above} by
the waiting time of a white particle migration step to an
unoccupied site.

Observe that due to the uniform migration distribution, the property
to jump to an unoccupied site is an experiment independent of everything else and the
jump is successful with probability $1-k/N$ if $k$  is the current number of
occupied sites. Therefore we can view the formation of new occupied colonies
as a {\em pruning} of the point process of white particle migration steps.

What can we say about the point process of the attempted migration
steps by white particles?
We note first that this migration step is  {\em given all current internal states} the infimum
of $k$, independent exponential waiting times, each of which
is the minimum of exponential waiting times for the internal transitions
like birth, death and emigration. Once an internal transition other than
migration occurs, the experiment starts over with this new internal states.
The rates of all this exponential waiting times depend on the internal state
in the corresponding site. It is best to view this as trials
consisting of migration steps arising from $k$ waiting times built from
exponentials according to the internal states, waiting for one which leads to
a new colony.

We know from the CMJ-theory that the number of migration steps in the white plus
black system, with $k$ occupied at time $t$, in the time interval
$[t, t+ \triangle t]$   behaves asymptotically as  $\alpha k \triangle t$ as $k \to \infty$
a.s. and in $L_1$  as long as $k$ has not reached too large values.
We now use the observation  that the relative frequency of black colonies
is asymptotically as $N \to \infty$ negligible so that asymptotically the
white migration steps are as frequent as the ones of the white and black system
as long as $k$ is large but $k/N$ remains small.

When the current number of occupied sites is $k$ the waiting time for the next colonization  has  mean {\em bounded above} by
that of the waiting time for the birth of  white site.
To obtain information on this random time we note that the next migration of a particle is given by an exponential random clock which
occurs with a hazard function given by the  rate $c\sum_{i=1}^k\zeta_i(\cdot)$  and the probability
that a migration produces a new white site is $1-\frac{k}{N}$  (compare  (\ref{Grev38}) ).  To get an upperbound on the time $\tau_{\log N}$ we note that $\zeta_i(\cdot)\geq 1$ and note that the number of migration events needed to obtain a
new white site is geometric with mean $(1-k/N)^{-1}$.

Recall that $\tau_{\log N}$ is given by (\ref{grev90a5b}).
Then recalling the convergence to the stable age and size distribution for $t \to \infty$
for $k \to \infty$ and therefore $\frac{c}{k}\sum_{i=1}^k\zeta_i(\cdot) \sim \alpha$,
we argue below in detail that the expected time, given the internal states of the occupied sites
to obtain the next white site for $\log N\leq k\leq \ve^\prime N$, behaves for
$N \to \infty$ as follows:
\be{angr17}
\frac{1}{k(\frac{1}{k}\sum_{i=1}^k1_{\zeta_i>1}\zeta_i(\cdot))(1-\frac{k}{N})}
\sim  \frac{1}{k\alpha(1-\frac{k}{N})}
=\frac{1}{\alpha k}
+\frac{1}{\alpha N}\frac{1}{1-\frac{k}{N}} \leq \frac{1}{\alpha k}+\frac{1}{\alpha N(1-\ve^\prime)}.\ee

To make more precise the first step, note first that
\be{angr17b}
\lim_{N\to\infty}\big(\frac{1}{K^N_{\tau_{\log N}+j}}
\sum_{i=1}^{K^N_{\tau_{\log N+j}}}1_{\zeta_i>1}\zeta_i(\cdot)\big)
 = \lim_{N\to\infty} \big(\sum_{k=2}^\infty k U^N(\tau_{\log N+j},k)\big).
\ee
We know from the CMJ-theory that:
\be{angr17c}
 |\lim_{N\to\infty} \sum_{k=2}^\infty k U^N(\tau_{\log N+j},k)-\sum_{k=2}^\infty k U(\infty,k)|=0.
 \ee
Therefore
\be{angr17d}
 |\lim_{N\to\infty} \sum_{k=2}^\infty k U^N(\tau_{\log N+j},k)-\alpha|=0\quad\text{ uniformly in } j.
\ee

Secondly we now have to deal with the fact that the quantity we handled above
appears in the denominator and is random, so that taking {\em expectation} needs some care.
Note that $ \sum_{k=2}^\infty k U^N(\tau_{\log N}+j),k)\geq 1-U^N(\tau_{\log N+j},1)$ and
\be{angr17e}
\lim_{N\to\infty}(1-U^N(\tau_{\log N+j},1))>\frac{1}{2}(1- U(\infty,1)) >0.\ee
We treat below the event $(1-U^N(\tau_{\log N+j},1))\leq a$ for small $a$ as a rare event and
obtain a large deviation bound for it allowing us to proceed for $N \to \infty$ with
assuming we are in the complement.

We break the details of the argument down in two parts, namely, the probability that more than $\frac{\log N}{2}$ of the individuals have age less than or equal to $B$, that is, $\frac{1}{2} \log N$ particles produce $\frac{1}{2}\log N $ new particles in time $\leq B$. We choose $B$ so that the expected number of particles produced by one site in time $B$ is $\leq \frac{1}{3}$.  The second part is to calculate that more than a fraction $(1-a)$ of the sites of age $>B$ have only a single particle where $a>0$ is chosen so that  if a site is of age $\geq B$ then $P(\zeta \geq 2) >2a$.  Since the events that $\frac{\log N}{2}$ particles produce $\frac{\log N}{2}$ particles in time $B$ and the event that more than
a fraction $(1-a)$ sites of age $B$ have $\zeta =1$ are both ``rare'' events, we can verify  the  exponential bound:
\be{angr17f}
P[(1-U^N(\tau_{\log N+j},1))\leq a]\leq e^{-C(\log N +j)}\text{ for large }N,
\mbox{ for some } C>0.
\ee

We now take expectations {\em conditioned} to be on the good event,
\be{angr17g}
E_\ast [\cdot] = E[\cdot | (1-U^N(\tau_{\log N+j},1))>a \mbox{ for } j=1,\dots,\ve N-\log N].
\ee
We can then conclude (by bounded convergence) that:
\be{angr17h}
\lim_{N\to\infty} E_\ast[\frac{1}{ \sum_{k=2}^\infty k U^N(\tau_{\log N+j},k)}]\leq \frac{1}{\alpha}, \quad \text{uniformly in }j.\ee

This means that we now can conclude with $T_k$ denoting the waiting time for the
next migration step if we have $k$ occupied sites satisfies:
\be{angr17j}
E_\ast [T_k] \leq \frac{1}{\alpha k} + \frac{1}{\alpha N(1-\ve^\prime)},
\quad k \in [\log N, \ve^\prime N].
\ee

\begin{remark} For the case considered in Section 8 it is convenient to note the following simple generalization. Consider a system of $k$ independent exponential clocks $C_1,\dots,C_k$ with
exponential waiting times with parameters $c_1,\dots,c_k$.  Then
\bea{angr17i}
&&P(\min_i C_i \leq \frac{x}{k})=  (1- \prod_i P(C_i>\frac{x}{k})) \\&&
= (1-\prod_i e^{-c_ix/k})= 1- e^{-(\frac{1}{k}\sum_{i=1}^k
c_i)x})\qquad \text {as }k\to\infty.\nonumber
\eea
In our case $c_i= \zeta_i\in \N$ and this leads to $1- e^{-(\sum_k k\cdot U^N(\tau_{\log N+j},k))}$, and therefore mean $\frac{1}{k\cdot \sum_k k\cdot U^N(\tau_{\log N+j},k)}$ and second moment $2 \left(\frac{1}{k\cdot \sum_k k\cdot U^N(\tau_{\log N+j},k)}\right)^2$ and so the variance is $ \left(\frac{1}{k\cdot \sum_k k\cdot U^N(\tau_{\log N+j},k)}\right)^2$.

\end{remark}

In order to conclude the argument we next choose $\ve > 0$. Then for sufficiently large $N$, the mean time,
$E_\ast[\tau^N(\ve)]$, for $ K^N_t$ required to reach
$\ve {N}$ (cf. (\ref{angr17j}) ) satisfies
\bea{b-dtimes}
E_\ast[ \tau^N(\ve)] \leq && E[\tau_{\log N}]+
\frac{1}{\alpha}[\frac{1} {\log N}+\dots+ \frac{\log N}{N}]\\
&&+\frac{1}{\alpha}[\frac{\log N}{N}+\dots+\frac{1}{\ve
N}]+\frac{\ve}{\alpha(1-\ve)}.
\nonumber \eea

\noi Hence
\be{Grev39}
E_\ast[\tau^N(\ve)|W] \leq \frac{\log N}{\alpha}
+c(\ve),
\ee
where $c(\ve)$ is a constant depending on $\ve$ and $W$.

We also have from the CMJ-theory the  lower bound

\be{add50}
E[\tau^N (\ve)|W]\geq E[\wt{\tau}^N (\ve)|W] =\alpha^{-1} \log N - \left(\frac{|\log \ve^\prime| +\log W}{\alpha}\right).\ee
Therefore it follows that $\tau^N (\ve)$,  conditioned on $W$, is of the form
$\alpha^{-1} \log N +O(1)$, as claimed. We also see that the
difference between $\tau^N (\ve)$ and $\wt \tau^N (\ve)$ is
a  random variable satisfying
\be{meandiff} E_\ast[|\tau^N (\ve)-\wt{\tau}^N (\ve)|] <\infty.\ee

We can also estimate now the time till fixation as follows.
The coupling construction
now implies that if we choose $\ve$ small enough, then at
time $\alpha^{-1} \log N$ we will reach $\ve N$ dual particles
which will then in finite random time produce a mutation jump so
that we get indeed (\ref{Y12c}) and thus completing the argument.
Namely since
\be{Grev39.1}
\Pi^N_t \geq  \wh K_t - \bar Z^N_t,\quad
\bar Z^N_t \leq \mbox{ r.h.s. (\ref{deadsites})}, \ee after
time $\wt \tau^N(\ve)$ there is a mutation rate at least
\be{grev39a}
m(\ve-  E[W] \ve^2)^+ \ee uniformly in large $N$ and therefore the
probability (recall
$\tau^N(\ve) - \wt \tau^N(\ve) = O(1)$ as $N \to \infty$)
that a mutation does not occur before $\tau^N(\ve) +t$
goes to zero as $t\to\infty$ uniformly in $N$ if we choose $\ve$ small
enough. This completes the
proof of (\ref{Y12c}). \bi

{\bf Step 5: Extension: convergence of hitting times}

However we can get even more mileage from the above argument to make
it clear that the randomness in the evolution of spatially macroscopic
variables is created  only in the very early stages of the growth of the dual particle system.

We next estimate the variance of $\tau^N(\ve) -\tau^N_{\log N}$,
a random variable of which we know from the analysis of Step 4 above that it is asymptotically
close to $\wt \tau^N(\ve) - \wt \tau^N_{\log N}$ as $N \to \infty$ and $\ve \to 0$. If we now condition on all the internal
states of all occupied sites we can calculate the variance of the
waiting time $T_k$ for the next migration step explicitly since
we know it is exponential. Hence we can, denoting by
$\wt Var(\cdot)$ this conditional variance, estimate
$\wt Var(\tau^N(\ve) -\tau^N_{\log N})$.
Conditioned on $W$ and restricted to the  complement of the rare event introduced in
Step 4 (which is asymptotically negligible (\ref{angr17f})) we can verify by explicit calculation
(following the same argument as for the means in (\ref{b-dtimes}) now for the variances, in this case
we get the sum of squares of the terms we had in the previous mean calculation) that
with $\wt Var_\ast$ denoting this conditional variance {\em restricted to the complement
of the rare event}, we get the bound
\be{Grev39b} {\limsup_{N\to \infty}}\,  \wt Var_\ast[(\tau^N(\ve)
-\tau_{\log N})] \leq \limsup_{N\to\infty}
\left(\rm{const}\cdot \left(\frac{1}{\log N}-\frac{1}{\ve N}\right)\right)=0.
\ee
Now we note that the term in the limsup on the r.h.s. is independent of the internal states and we
conclude that conditional on $W$ the (restricted) variance of the difference between
$\tau^N(\ve)$ and $\tau_{\log N}$ goes to zero and hence
that is, the variance of the sum conditioned on $W$ of the
waiting times (with means given in (\ref{angr17})) in the second and third time intervals  given in
(\ref{grev90a5}) goes to $0$
as $N\to\infty$ since the rare event (and hence the restriction) become negligible as $N \to \infty$.
The analogous result also holds for $\wt{\tau}^N(\ve)-\wt\tau_{\log N}$.

This means that the growth in the second
and third time interval we have specified in (\ref{grev90a5}) is {\em deterministic}
and the randomness occurred earlier in the first time interval.

Note that  during the first time interval as $N \to \infty$ the expected number of collisions is $O(\frac{(\log N)^2}{N})$ so that
no collisions occur with probability tending to 1. Therefore  as $N \to \infty$,
\be{Grev39.2}
\wt \tau^N_{\log N} - \tau^N_{\log N} \Ntoo 0 \quad \mbox{ in probability.}
\ee
Together with the fact that as $N\to\infty$ the limiting probability of the rare event is $0$
it follows that  there exists a deterministic sequences, namely
\be{add176}
(\wt{\delta}^{\ve }_N=E[\wt \tau^N_{\ve N} - \wt\tau^N_{\log N}])_{N \in \N}$, $(\delta^{\ve
}_N=E_\ast[ \tau^N_{\ve N} - \tau^N_{\log N}])_{N \in \N}
\ee
such that for $\eta>0$,
\be{Grev39.3}
 P[| \wt \tau^N (\ve) - \wt \tau^N_{\log N}-\wt{\delta}_N^{\ve}| > \eta]\to 0,
\ee
and
\be{Grev39.3x}
 P[|\tau^N (\ve) -  \tau^N_{\log N} - {\delta}_N^{\ve} | >\eta]\to 0.
 \ee

We now collect the above estimates to prove the following.
\beP{P.C-10}{(Convergence of hitting times of level $\ve$)}

There exists a non-degenerate $\R$-valued random variable
$ \wt\tau(\ve)$ such that

\be{Grev39c} \wt\tau^N(\ve)-\frac{\log N}{\alpha} \Rightarrow
\wt\tau(\ve), \mbox{ as } N \to \infty,
\ee
 \be{Grev39c2} \tau^N(\ve)-E_\ast[\tau^N(\ve)] \text{ converges in distribution}
\ee
and there exists a constants $\{C(\ve)\}_{\ve >0}$ such that $\lim_{\ve\to 0}C(\ve)=0$ and
\bea{Grev39c3}
&& \lim_{N\to\infty}P\Big(\frac{\log N+\log \ve -\log W}{\alpha}-C(\ve) < \tau^N(\ve)\\
&&
\hspace{3cm}
 <\frac{\log N+\log \ve -\log W}{\alpha}+2 C(\ve)\Big|W\Big)=1.
\qquad \square \nonumber
\eea
\end{proposition}
\begin{proof}
We begin with some preparatory points.
Note that asymptotically as $N \to \infty$ by (\ref{ang2b}),
\be{Gre39e}
\wt \tau^N_{\log N}-\frac{\log\log N-\log W_{N}}{\alpha}=o(1)\ee
where $W_{N}\to W$ a.s., and by (\ref{Grev39.2}) it follows that
\be{Gre39e2}
\tau^N_{\log N}-\frac{\log\log N-\log W_{N}}{\alpha}=o(1).\ee

But by (\ref{Grev39.3}), (\ref{Grev39.3x}),
\be{Gre39f}
 P[|(\tilde\tau^N(\ve)-\wt\tau_{\log N}-\wt\delta_N^\ve) -(\tau^N(\ve)-\tau^N_{\log N}-\delta_N^\ve)| >2\eta]\to 0.
 \ee
 But then (\ref{meandiff}) implies that \be{meandiff2}\sup_N|\wt\delta^\ve_N - \delta^\ve_N|<\infty.\ee

 We also have that  $\wt\delta^\ve_N-(\frac{\log{(N\ve)}-\log\log N}{\alpha})\to 0$.
Therefore  as $N \to \infty$:
\be{Gre39g}
\tilde\tau^N(\ve)-\frac{\log N}{\alpha}\to \frac{\log \ve-\log
W}{\alpha}.
\ee

Having finished the preparation we already see that (\ref{Gre39g})
implies the weak convergence claimed in (\ref{Grev39c}). Now we come to the proof of
(\ref{Grev39c2}) and (\ref{Grev39c3}).

The relation (\ref{Grev39c3}) then follows from  using (\ref{Grev39f}) and inserting
(\ref{Gre39e2}), (\ref{Gre39e}) together with (\ref{meandiff2}) and the fact that
$\supl_N |\delta^\ve_N - \wt \delta^\ve_N|$ tends to zero as $\ve \to 0$.
In order to see now (\ref{Grev39c2}) we have to use that given $W$ the
$\tau^N(\ve) - \tau^N_{\log N}$ is as $N \to \infty$ asymptotically deterministic
as we saw above. Then the claim follows with (\ref{Gre39e2}).
\end{proof}

\medskip

\subsubsection{Dual process in the  collision regime: nonlinear dynamics}
\label{sss.growth}

 In order to obtain more detailed results on the emergence regime
 and on the process of fixation
  via the dual process we will need more precise information about
the growth of the number of sites occupied by the dual process,
respectively the states at these sites,
once we reach the regime where {\em collisions}  play a role.

The behaviour in this regime is qualitatively very different from the
collision-free situation  and is the result of three properties of
the evolution mechanism:
\begin{enumerate}
\item Randomness enters in early stages of the growth of the dual
process, i.e. times in $[0,O(1)]$, here the dual behaves
as a CMJ branching process described by {\em linear} (i.e.
Markov) but {\em random} dynamics, \item in the time regime
\be{Gre39g3}
[s(N),\frac{\log N}{s(N)}), $ with $s(N)\to\infty,\;
s(N)=o(\log N)\ee the law of large numbers is in effect leading to
{\em deterministic} but still {\em linear} dynamics
\item in the regime ($s(N)$ as above)
\be{R3}[\frac{\log
N}{s(N)},\frac{\log N}{\alpha}+t)\ee the dynamics is {\em deterministic}
but the evolution equation becomes {\em nonlinear} due to collisions.
\end{enumerate}
The problem is now that these regimes separate as $N \to \infty$
and we have to connect them through a careful analysis.
To integrate all three phases into a {\em deterministic nonlinear}
dynamic with {\em random} initial condition, the full dual process for the
finite system with $N$ exchangeable sites is analysed.

Indeed in this
and the next subsubsections we show that in the limit $N \to
\infty$ the initial asymptotically collision-free evolution
stage in (1) and (2) produces {\em a random initial
condition} for the nonlinear dynamics arising in (3) where the
randomness arises from stage (1). A key step in linking the
random and nonlinear aspects is the analysis of the evolution in the
time interval in (\ref{R3}) in the
limit $N\to\infty,\;t\to -\infty$.

In this subsection we consider the time interval
\be{ba10dxx}
\frac{1}{\alpha}[\log N-\log\log N, \log N+T)
\ee
in which the dynamics are deterministic.
Recall that
\be{ba10f5}
u^N(\frac{1}{\alpha}(\log N-\log\log
N))=W_N\frac{N}{\log N} =W\frac{N}{\log N} + o(1) \mbox{ as } N
\to \infty,
\ee
since $W_N \Rightarrow W$ by the equation (\ref{ang2b}).

We have to describe the dual process in the above time interval in the limit
$N \to \infty$ and to obtain a limit dynamic which allows to draw the
needed conclusions for the original process. We proceed in three steps. It turns out that
the main properties of the dual process needed for that purpose can be captured
in a triple of functionals of this process which we introduce in
Step 1. In Step 2 we formulate the corresponding limiting
(for $N \to \infty$) objects
which allows us to finally state and prove in Subsubsection \ref{sss.asycolreg}  the
key convergence relation for the dual process.
Step 3 proves that the defined evolutions have the desired properties.
\bi

{\bf Step~1} $\;${\em  Some functionals of the dual in the time span up to fixation}

First, we observe that  to determine the quantity
of interest for emergence namely
\be{ba10dx}
\frac{\Pi^N_{T_N+t}}{N} \quad , \mbox{ for } T_N=
\frac{\log N}{\alpha}  \ee as a function of $t$ we need to know
the number of sites occupied at time $T_N+t$ and relative
frequencies  of the occupation numbers of these sites.  Since the occupation
number of a site depends on the age of the site, we will keep
track of the number of occupied sites and the relative frequencies
of both age and occupation number for these sites. As $N\to\infty$  a law of large
number argument is then used to obtain the desired information.

We consider the number of sites occupied at time $t$ denoted $
K^N_t$ and the corresponding measure-valued process on
$[0,\infty)\times \N_0$ giving the unnormalized
number of sites of a certain age $u$ in $[a,b)$ and occupation size $j$:
\be{ba10d2}
\Psi^N(t,[a,b),j)) = \int_{(t-b)}^{(t-a)}1_{( K^N_{u}>
K^N_{u-})}1_{(\zeta^N_u(t)=j)}d K^N_u,
\ee
where $\zeta^N_u(t)$
denotes the occupation number at time $t$ of the site born at time
$u$, that is, a site first occupied the last time at time $u$, which is
therefore at time $t$ exactly of age $t-u$.

The {\em normalized empirical age and size distribution} among the occupied sites
is defined as: \be{ba10f1} U^N(t,[a,b),j) =\frac{1}{K^N_t} \Psi^N
(t,[a,b),j),\quad t\geq0,\; j\in \{1,2,3,\dots\}. \ee

Thus for $(\zeta^N_t)_{t \geq 0}$
we have a growing number $K^N_t$ of different interacting birth and
death processes which each have in state $k$ linear birth rate $s k$, linear death rate $c
(1-\frac{1}{N})k $, quadratic  death rate $d k(k-1)$.  In addition immigration of one
additional element into a site of size  $k$ occurs with rate
\be{Grev39f}
c\left(\frac{ K^N_t c^{-1} (\alpha_N(t) +\gamma_N(t))-k}{N}\right),
\ee
where $\alpha_N(t)$ is $c$ times the mean size of sites with more than one
particle,  $\gamma_N$ is $c$ times the frequency of singletons; these  will be
defined precisely below in (\ref{Grev39e}). Note that $U^N(t,dx,j)$
is a pure atomic measure which evolves with jumps of size
$\frac{1}{ K^N_t}$ one for each migration event to a yet unoccupied site,
or $-\frac{1}{K^N_t}$ when a  singly occupied site migrates and a collision occurs.

\begin{remark}
Note that the $\zeta^N_u(\cdot)$ interact if we are in a
regime where collisions  occur, that is at times after emergence
has occurred which happens at times of order $\frac{\log N}{\alpha}$.
The marginal distribution of age is the analogue of the empirical
age distribution in (\ref{ang4b}).  However now we take into account
 possible collisions, that is, a migration to a previously occupied
site when we consider the system of $N$ sites with symmetric
migration among these sites and in this case the distribution of
$\zeta^N_u(t)$ may depend not only on the age $t-u$  but also on $t$  due to  changes in the collision rates.
\end{remark}

Set for convenience to describe the pair given by the number of sites and size-age distribution:
\be{Grev39c2a}
 u^N(t):= K^N_t.
\ee

We have obtained with (\ref{ba10f1}), (\ref{Grev39c2a}) now a pair which is
$\N \times \CP([0,\infty) \times \N)$-valued, denoted
\be{ba10f4}
(u^N(t), U^N(t,\cdot,\cdot))_{t \geq 0}
\ee
and which completely describes  our dual particle system
up to permutations of sites.  We can recover  the functional $\Pi^N_t$ as

\be{ba10f3}
 \Pi^N_t= u^N(t)  \sum_{j=1}^\infty jU^N(t,[0,t],j).
\ee However since we work with exchangeable initial states for our original
population process, the pair $(u^N, U^N)$  provides a sufficiently complete  description of the dual process to allow us to calculate the laws of our original processes.

In the interval $\alpha^{-1} [\log N- \log \log N,
\log N+T]$ the process $(u^N(t))_{t \geq 0}$  increases by one,
respectively decreases by one, at rates

\be{Grev39d}
 \alpha_N (t)(1-\frac{u^N(t)}{N}) u^N(t), \mbox{ respectively }
\gamma_N (t) \frac{( u^N(t))^2}{N}, \ee
where $\alpha_N
(t),\gamma_N (t)$ are defined: \be{Grev39e} \alpha_N (t) = c
\intl^t_0 \suml^\infty_{j=2} jU^N (t,ds,j) ds, \quad \gamma_N (t)
= c \intl^t_0 U^N (t,ds,1). \ee
These rates of change of $U^N (t,\cdot,\cdot)$ follow directly
from the dynamics of the dual particle system $\eta$.

We note that if we start with $k$ particles at each of $\ell$
different sites we obtain similar quantities, which we denote by
\be{Grev39e2} (u^{N,k,\ell}_t, U^{N,k,\ell}_t)_{t \geq 0}. \ee

The dynamics of the the system
\be{ba10f7}
(u^N(t), U^N(t,ds,dx))
\ee
evolves in two main regimes corresponding to times
$o(\log N)$ and those of the form
$\alpha^{-1} \log N+t$, with $t=O(1)$. They are
asymptotically as $N \to \infty$ {\em linear}
and {\em random} for times of order $o(\log N)$ as we have
demonstrated in the previous subsubsection  with the randomness
captured by the random variable $W$ given in (\ref{ang2b}) respectively
(\ref{ba10f5}) above.
On the other hand the dynamics become {\em deterministic} and {\em
nonlinear} in the post-emergence regime, that is, for times of the
form $(T_N+t)_{t\in (-\infty,\infty)}$ where $T_N=\frac{\log
N}{\alpha}$.  This latter case involves the law of large numbers in the
limit $N \to \infty$ and collisions play a decisive role.
We will introduce the limiting dynamics in Step 2.
We will later in Subsubsection \ref{sss.prep2} show how the
two regimes are linked and that the first produces a random
initial condition for the second, the nonlinear regime.
\bi


\bigskip

{\bf Step~2}$\;$ {\em  The limiting dynamics for dual in collision regime}

The second step is to define for $k, \ell \in \N$  candidates
\be{Grev39g}
(\Pi^{k,\ell}_t,u^{k,\ell}(t),U^{k,\ell}(t,ds,dx))_{t \in \R} \mbox{ and } (\zeta(t))_{t \in \R}
\ee for the
limiting objects as $N \to \infty$ related to (here $T_N =\alpha^{-1} \log N$)
\be{Grev39h}\begin{array}{l}
\left(\frac{\Pi^{N,k,\ell}_{T_N+t}}{N},\frac{u^{N,k,\ell}_{T_N+t}}{N},
U^{N,k,\ell}(T_N+t,ds,dx)\right),\;-\infty<t<\infty,\\
\{\zeta^{N,(k,\ell)}(T_N+t,i), \quad i=1, \cdots, N\}, \quad -\infty <t<\infty,
\end{array}\ee
when the dual process initially has {\em $k$ particles} at each of {\em $\ell$ distinct sites} at
time $t=0$.

We now turn to the limiting objects as $N \to \infty$ which we denote
\be{Grev39m}
(\Pi^{k,\ell}_t, u(t), U(t))_{t \in \R}, \quad (\zeta^{k,\ell} (t))_{t \in \R}.
\ee

We remark that it will turn out that conditioned on
$(\Pi^{k,\ell}_{t_0},u^{k,\ell}(t_0),U^{k,\ell}(t_0,ds,dx))$, the dynamics of $(\Pi^{k,\ell}_t,u^{k,\ell}(t),U^{k,\ell}(t,ds,dx))_{t\geq t_0}$ is independent of $k,\ell$.  For this reason we often suppress the superscripts and later indicate how the initial conditions for  $k,\ell >1$ are determined.

  Since
\be{agr2}
\Pi^{k,\ell}_t = u^{k,\ell} (t) \cdot \suml^\infty_{j=1} j \cdot U^{k,\ell}
(t, \R^+,j),
\ee
and since the process $\zeta$ will depend on this pair from
the triple in (\ref{Grev39g})
it suffices to specify the pair $(u(t),U(t))_{t\in \R}$.  This is done by first
establishing that they satisfy a pair of deterministic equations and then determine the entrance law at   $-\infty$ as follows.

We first define  two ingredients (parallel to
(\ref{Grev39e})):
\be{ba10d}
\alpha(t) =\alpha(U(t)):= c
\intl^\infty_0 \sum_{j=2}^\infty jU(t,ds,j),\quad
\gamma(t)=\gamma(U(t)):=c\intl^\infty_0 U(t,ds,1). \ee
 As before
we first consider the evolution equation for times $ t \geq t_0$ and then we
try to characterize an entrance law by considering $t_0 \to - \infty$. For the
moment, for notational convenience let $t_0=0$.

Then to specify
\be{ba10d1}
(u,U) = (u(t), U(t))_{t \geq 0}, \mbox{ with } (u(t), U(t)) \in \R_+\times \CM_1(\R_+\times \N),
\ee
we introduce the (coupled) system of nonlinear evolution equations:

\be{grocoll2}
\frac{du(t)}{dt}= \alpha(t)(1-u(t)) u(t) -\gamma(t) u^2(t),
\ee

\bea{ba10g}
\frac{\partial U(t,dv,j)}{\partial t}&&\\
  =&&- \frac{\partial U(t,dv,j)}{\partial v} \nonumber
\\&&
+s(j-1)1_{j\ne1}U(t,dv,j-1)-sjU(t,dv,j)\nonumber
\\&&
  +\frac{d}{2}(j+1)jU(t,dv,j+1)-\frac{d}{2}j(j-1))1_{j\ne1}U(t,dv,j)\nonumber
\\&&
  +c(j+1)U(t,dv,j+1)-cjU(t,dv,j)1_{j\ne 1} \nonumber
\\&&
  -cu(t)U(t,dv,1)1_{j=1}\nonumber
\\&&
  + u(t)(\alpha(t)+\gamma(t))[1_{j\ne1}U(t,dv,j-1)-U(t,dv,j])   \nonumber
\\&&
  + \left ((1- u(t))\alpha(t)\right )1_{j=1}\cdot \delta_0(dv) \nonumber
\\&&
 -\Big(\alpha(t)(1-u(t)) -\gamma(t)u(t)\Big)\cdot U(t,dv,j).
\nonumber\eea
Note that if $U(0,\cdot,\cdot)$ is a probability measure on $\mathbb{R}_+\times \N$, then this property is preserved under the evolution.
(Recall here that $\alpha$ and $\gamma$ contain the constant c of migration.)

We note that the right side equation (\ref{ba10g}) is only defined in terms of
$\CM_{\mathrm{fin}} (\R^+ \times \N)$ elements, if we restrict to
measures for which  $\int_{0}^\infty\sum_j j^2 \cdot U(dv,j)<\infty$, so we have
to restrict to a subset of  $\CM_1 (\R^+ \times \N)$.
Define  $\nu$ as the measure on $\N$ given by
\be{agre65a}
\nu(j)=1+j^2.
\ee
and let $L^1 (\nu,\N)$  denote the corresponding space of integrable functions.

We can write these equations in compact form as
\be{ba10g1}
\frac{d}{dt} (u(t), U(t, \cdot,\cdot)) = \vec {G}^{\ast} (u(t), U(t, \cdot,\cdot)),
\ee
with the {\em nonlinear} operator
\be{ba10g2a}
\vec{G}^{\ast} \mbox{ mapping } \R_+ \otimes L^1_+(\R_+\times\N, \nu) \to \R_+  \otimes \F,
\ee
with $\F$ denoting the class of finite measurable functions $f$ on $\mathbb{R}_+\times \N$.

If $u(0)=u_0>0$, we also work with the equivalent $\CM_{\rm{fin}}(\R_+\times \N)$-valued process
\be{psidefn}
\wh\Psi_t:= u(t)\cdot U(t)
\ee
which solves the nonlinear systems
\be{psidefn2}
\frac{d\wh{\Psi}_t}{dt}= G^\ast(\wh {\Psi}_t)
\ee
equivalent to (\ref{ba10g1}) (note that it is easy to verify that if $u(t_0)=\wh\Psi_{t_0}(\R_+\times\N)>0$, then $u(t)>0$ for all $t\geq t_0$).

We think of this as a forward equation for a deterministic measure-valued process.  In order to formulate the corresponding Markov
semigroup generator  we consider functions $F_{g,f}$ on
$(\R_+ \otimes L^1_+(\R_+\times\N, \nu))$ of the form
\bea{dd6}
&& F_{g,f}(\wh\Psi_0)= g\left(\sum_{j=1}^{\infty}\int_0^\infty f(r,j)\wh\Psi_0(dr,j)\right)\\&&
g:\R\to \R,\;f:[0,\infty)\times\N\to\R,\; g\text{ has bounded first and second derivatives},\nonumber \\
&&f\text{ bounded and has bounded first derivatives in the first variable}. \nonumber
\eea
 The algebra of functions,  which contains all functions of the form (\ref{dd6}) is denoted
\be{add135}
D_0(\bf{ G})
\ee
and is a separating set in
$C_b(\CM_{\mathrm{fin}}([0,\infty) \times \N_0))$ and is
therefore {\em distribution-determining for laws} on
$\CM_{\mathrm{fin}}([0,\infty) \times \N_0)$.

\begin{remark}\label{R.x2}
In fact since we are concerned with $\wh \Psi$ which are
sub-probability measures we note that the integrals are bounded, we can include unbounded functions $g$
defined on $[0,1]$  and  in
$C^2 ([0,1], \R)$. In particular, we can consider $g(x)=x$ and $g(x)=x^2$.
\end{remark}

Then we introduce the semigroup
\be{ba10g2}
V_tF_{g,f}(\wh\Psi_0):= F_{g,f}(\wh\Psi(t)),
\ee
where $\wh \Psi(t)= u(t)U(t) $ and $(u(t),U(t))$ is the solution of the pair (\ref{grocoll2}), (\ref{ba10g}) with initial condition $\wh\Psi_0= u_0\cdot U_0$.
Let
\be{ba10g2b}
\wh\Psi(f):=\left(\sum_{j=1}^{\infty}\int_0^\infty f(r,j)\wh\Psi_0(dr,j)\right).
\ee

Then with $F_{g,f}$ as in (\ref{dd6}) we set
\be{ba10g21}
\begin{array}{ll}
& {\bf{G}}F_{g,f}(\wh\Psi)=g^\prime\left(\wh\Psi(f)\right) \left(\wh\Psi(G^{\wh\Psi}f(r,j))\right)\\
&= g^\prime\left(\wh\Psi(f)\right)\intl_0^\infty\intl_1^\infty \frac{\partial f(r,x)}{\partial r}\wh\Psi(dr,dx)  \\
& \qquad + g^\prime(\wh\Psi(f))\intl_0^\infty\intl _1^\infty   sx
\left[(f(r,x+1)-f(r,x)))\right]\wh\Psi(dr,dx)\\
& \qquad + g^\prime(\wh\Psi(f))\intl_0^\infty\intl_1^\infty
\frac{d}{2}x(x-1) \left[f(r,x-1)-f(r,x))\right]\wh
\Psi(dr,dx)\\
& \qquad + g^\prime(\wh\Psi(f))\cdot  c(1-\wh\Psi(1))\left(\int_0^\infty\int_1^\infty
[f(0,1)+f(r,x-1)-f(r,x)]x1_{x\ne 1}\wh \Psi(dr,dx)\right)
\\
& \qquad + g^\prime (\wh\Psi(f))\wh \Psi(1)
c\intl_0^\infty\intl_1^\infty \left[f(r,x+1)- f( r,x)
)\right]x\wh\Psi(dr,dx) \\
& \qquad + g^\prime (\wh\Psi(f))c(\int_0^\infty\int_1^\infty
x\wh\Psi(dr,dx))\left(\intl_0^\infty\intl_1^\infty \left[1_{\wt x\ne 1}f(\wt r,
\wt x -1)- f( \wt r,\wt x) )\right] \wh\Psi(d \wt r,d \wt x)\right)\\&
= g^\prime(\wh\Psi(f))\cdot \intl_0^\infty\intl_1^\infty G^{\wh\Psi}f(r,x)\wh\Psi(dr,dx),
\end{array}\ee
with the obvious extension to $D_0({\bf G})$.

We can then consider the martingale problem associated with $({\bf G}, D_0({\bf G}))$,
that is, to determine probability measures, $P$, on $C([0,\infty),\mathcal{M}_{\rm{fin}}(\R_+\times \N)\cap \mathbb{B})$ (where $\mathbb{B}$ is defined below in (\ref{Grev39ny}))
such that
\be{add1}
F_{g,f}(\wh\Psi(t))-\int_0^t ({\bf G}F_{g,f})(\wh\Psi(s))ds
\ee
is a $P$-martingale for all $F\in D_0({\bf G})$.

Since ${\bf G}$ is a {\em first order} operator any solution to this martingale problem is deterministic. To verify this consider $g(x) =x^2$ and $g(x)=x$ and
show that
${\bf G}(F_{x^2}) -2F_x {\bf G}F_x =0$. This holds by inspection of (\ref{ba10g21}).
This implies that the variance of the marginal distribution $\wh\Psi_t(f)$ for a solution to the  martingale problem is 0.

Then taking next $g(x)=x$ we have that any solution satisfies
\be{add2}
\wh\Psi_t(f)-\int_0^t\wh\Psi_s(G^{\wh\Psi(t)}f)ds= \wh\Psi_t(f)-\int_0^t G^\ast\wh\Psi_s(f)ds.
\ee
Taking $f\equiv 1$ and then indicator functions $f=1_j$ we obtain our nonlinear system.
In other words any solution to the martingale problem is a deterministic trajectory satisfying the
nonlinear dynamics (\ref{grocoll2}), (\ref{ba10g}). This form will play a role later on when we consider the convergence.

In order to discuss the existence and uniqueness of solutions to these equations we introduce a norm and reformulate the equations as a nonlinear evolution equation in a Banach space, namely,
\be{Grev39ny}
\wh{\mathbb{B}}:= \R \otimes L^1(\mathbb{R}_+\times \N,ds\times \nu )
\ee
which we furnish with the norm
\be{Grev39n}
\|\wh{\Psi}\|= \|(u,(a_i(\cdot))_{i \in \N}) \| = |u| + \|(a_i(\cdot))_{i \in \N}\|_1 ,\ee
where $u=\wh{\Psi}(\R\times\N),\; a_i (\cdot) =\wh{\Psi}(\cdot \times \{i\})$ and
\be{add3}
\|(a_i(\cdot))_{i \in \N}\|_1 := \suml^\infty_{j=1}
(1+j^2)|a_j(\cdot)|_{\rm{var}},
\ee
where $\rm{var}$ denotes the total variation norm for measures on $\mathbb{R}_+$.
On this space the r.h.s. of the equation (\ref{ba10g}) is well-defined.

We note that the Markov property is satisfied by the equivalent systems $(u(t),U(t,\mathbb{R},\cdot))$ and $\wh\Psi(t)=u(t)\cdot U(t,\mathbb{R},\cdot)$.  Since this all that is needed for the dual representation we often suppress the ages. In this case we use
\be{add51}
\mathbb{B}:= \R \otimes L^1( \N, \nu )
\ee
and the simpler norm
\be{Grev39nx}
\| (u,(a_i))_{i \in \N}) \| = |u| + \|(a_i)_{i \in \N}\|_1,\quad
\text{where}\quad\|(a_i)_{i \in \N}\|_1 := \suml^\infty_{j=1}
(1+j^2)|a_j|.
\ee

Existence of a solution to this system will follow from the convergence result
below in Subsubsection \ref{sss.asycolreg} (Proposition \ref{P.Grocoll}).
 The question of uniqueness given  the
initial point $(u(t_0), U(t_0))$ is considered  in Step 3 below.

Return now to the finite $N$-system and imagine
the collection of occupied sites in our dual particle process at a specific time $t$ represented in the form $\zeta^N(t,i),i=1,\dots,N_t$. To complete the description we must specify
what happens at a single tagged site for times $r\geq t$, i.e. to  $\zeta^N(r,i), i \in \N$ in
the limit $N \to \infty$.

Given the function $(u(t),U(t))$ this will be given by  the time-inhomogeneous birth, death and
immigration process
\be{Grev39i}
\{\zeta(t,i),t \geq 0\},
\mbox{ with transition probabilities }
\{p_{k,j}(t), (k,j) \in \N^2, \quad t \geq 0\},
\ee given by
the {\em forward} nonlinear Kolmogorov equation starting with one
particle at site $i$ at time $0$, i.e., $p_{1,1}(0)=1$ and which reads  as
follows:

\be{Grev39j} p^\prime_{1,j}(t)=(s(j-1)+\alpha(t)u(t))p_{1,j-1}(t)
   +\frac{d}{2}(j+1)jp_{1,j+1}(t)-((s+c)j+\frac{d}{2}j(j-1))p_{1, j}(t)
\ee
for $j\ne 0,1$,
\be{Grev39k}
   p^\prime_{1,0}(t)= cu(t)p_{1,1}(t)-\alpha(t)u(t)p_{1,0}(t),
\ee

\be{Grev39l} p^\prime_{1,1}(t)= {d}\cdot
p_{1,2}(t)-((s+\alpha(t)u(t)))p_{1, 1}(t),\ee where \be{Grev39l2}
 u(t)=1-p_{1,0}(t),\quad \alpha(t)u(t)=c\sum_{j=2}^\infty
jp_{1,j}(t).\ee
This process $\zeta (t,i)$ will describe the evolution of a
tagged site $i$ in the limiting system in the collision regime.

Next, in order to prepare for the study of the nonlinear system (\ref{grocoll2}),(\ref{ba10g}) we
consider a related  set of equations for the non-collision
regime obtained by setting  $u\equiv 0$,
namely,

\bea{ba10g9ax}
\\
\frac{\partial \wt U(t,dv,j)}{\partial t}=
  &&
  - \frac{\partial \wt U(t,dv,j)}{\partial v}\nonumber\\
&&
 +s(j-1)1_{j\ne1}\wt U(t,dv,j-1)-sjU(t,j)\nonumber
 \\&&
 +\frac{d}{2}(j+1)j\wt U(t,dv,j+1)-\frac{d}{2}j(j-1))1_{j\ne1}\wt U(t,dv,j)\nonumber
 \\&&
+c(j+1)\wt U(t,dv,j+1)-cj\wt U(t,dv,j)1_{j\ne 1} \nonumber
 \\&&
   +\wt \alpha(t) [\delta_0\times 1_{j=1}
     -  \wt U(t,dv,j)],
\nonumber\eea
where  $\wt \alpha(t):=c\sum_{\ell=2}^\infty\int_0^\infty  \ell\cdot \wt U(t,dv,\ell)$.

Since the dual representation expressions only involve the marginals $U(t,\mathbb{R}_+,j)$ and the jump
rates depend on $j$ but not the age of an occupied site it suffices to work with these marginal distributions
we denote by $U(t,j)$.  We then work with the corresponding versions of
(\ref{grocoll2}),(\ref{ba10g}) and (\ref{ba10g9ax}) and use the norm (\ref{Grev39nx}).

We observe next that the solution to the nonlinear system  (\ref{ba10g9ax}) can be obtained as
\be{add4}
\wt U(t,dv,j)=\frac{V(t,dv,j)}{\sum_j \int_0^\infty V(t,dv,j)}\ee
where $\{V(t,\cdot,\cdot)\}$  solve the {\em linear} system which we spell out explicitly
later on, namely it is obtained by deleting the
last  term in (\ref{ba10g9al}).  The resulting linear equations coincide with the mean equations
for the CMJ process and hence the convergence of $V(t, \cdot, \cdot)$ as $t \to \infty$
to the stable size and age distribution follows from the CMJ theory.
In addition we have that $\alpha(t)\to \alpha$ as $t\to\infty$, the Malthusian parameter of the CMJ process.

We define a linear operator given by the infinite matrix
\be{agre65aa}
Q^{a,\ast} \mbox{ on } L^1_+ (\N, \nu) \subseteq \CM_{\mathrm{fin}} (\N),
\ee
obtained by first integrating over ages setting $\alpha(t)=a$ on the right side of the
(\ref{ba10g9ax}).
Note that the equation resulting this way from
(\ref{ba10g9ax}) can be viewed as the forward equations for a Markov chain
$Z$ on $\N$ with transition matrix $Q^{a,\ast}$.  Let
\be{dualsemigroup}
(S^a)_{t\geq 0} =\mbox{ the dual semigroup on $\CM_1(\N)$ of this time-homogeneous
Markov chain on  } \N.
\ee

To verify that the semigroup $S_t^\alpha$ corresponding to $Q^{\alpha,\ast}$
is strongly continuous on $(\mathbb{B},\|\cdot\|)$, note that for $\mu\in {\mathcal{M}}_1(\N)$,
\be{add5}
\lim_{t\to 0}\|S_t\mu-\mu\|= \lim_{t\to 0}[E[Z^2(t)]-E[Z^2(0)]]=0,
\ee
follows from elementary properties of the Markov chain $(Z(t))_{t \geq 0}$.

Note that the generator of $S^\alpha$ is an unbounded operator but has domain containing
\be{add52}
D_0:= \{(a_i)_{i\in\N}\in \mathbb{ B}:\sum_{j=1}^\infty (1+j^4)|a_i|<\infty\}\ee and in fact this is a core for the semigroup (\cite{EK2}, Ch. 1, Prop. 3.3),
since the birth and death process with linear birth and quadratic death rates has moments of all orders at positive times. (Note that the corresponding  Markov process
on $\N$ has an entrance law starting at $+\infty$.
This property is wellknown for the Kingman coalescent but extends to our
birth and death process provided $d>0$ and is easily extended  to include an additional immigration term.)

The Markov chain associated with $Q^{a,\ast}$ has a  unique
equilibrium state ${U}^{a}(\infty)$.
 Let $m(a)$ denote the mean equilibrium value
 \be{add6}
 m(a)=\sum _{j=2}^\infty j\cdot U^a(\infty,j).
 \ee
The equation  $ m(a)=a$  has a unique positive solution and it is equal to $\alpha$.
Moreover $S_t^\alpha$ satisfies
\be{add7}
\liml_{t \to \infty} (S_t^{\alpha}U(0)) = \mathcal{U}(\infty), \quad \limsup_{t\to\infty} |\alpha(t)-\alpha|=0,
\ee
where  $\mathcal{U}(\infty)= U^\alpha(\infty)$ is the unique equilibrium with mean $\alpha$ and is the stable size distribution for the CMJ process.

\bigskip

{\bf Step 3} $\;${\em Uniqueness and existence results}

We next show that the equations we gave  in Step 2 uniquely determine the
objects needed  for our convergence results.

\beL{L.UCE} {(Uniqueness of the pair $(u,U)$ and of $\zeta$)}

(a) The pair of equations (\ref{grocoll2}) - (\ref{ba10g}), for $(u(t),U(t,\mathbb{R}_+,
\cdot))$  given
an initial state from $\R^+ \times \CM_1(  \N)$
satisfying
\be{agr4}
(u(t_0), U(t_0, \cdot)) \in [0,1] \otimes \CM_1(  \N)\cap L^1 (\N,\nu)
\ee
at time $t_0$, has  a unique solution $(u(t),U(t))_{t\geq t_0}$
with values in $[0,1] \otimes (\CM_1(  \N)\cap L^1_+ (\N,\nu))$.

(b) Given $(u,U)$ the process $(\zeta(t))_{t \geq t_0}$ is uniquely
determined by (\ref{Grev39j}).

(c) There exists a solution $(u,U)$ with time parameter $t\in
\mathbb{R}$ for every $A\in (0,\infty)$ with values in
$\mathbb{B}= \R \otimes L^1(\N,\nu)$, such that
\be{S081}
u(t)e^{-\alpha t}\to A \text{  as  } t\to -\infty,\ee

\be{S082} U(t)\ttooo {\CU}(\infty),\quad\text{and   } \limsup_{t\to -\infty} e^{-\alpha t}|\alpha(t)-\alpha|=B<\infty.
\ee
Here $\CU(\infty)$ is the
stable age and size distribution of the CMJ-process corresponding to
the particle process
$(K_t, \zeta_t)_{t \geq 0}$
given by the McKean-Vlasov dual process $\eta$, which was
defined in (\ref{ang4c2}).
\sm

 (d) Given any solution $(u,U)$ of
equations (\ref{grocoll2}) - (\ref{ba10g}) integrating out the age (given in
(\ref{grocoll29})-(\ref{ba10g9a})) for  $t\in \mathbb{R}$
with values in the space $\mathbb{B}:= \R \otimes L^1(\N,\nu)$ satisfying
\be{limsup2}
u(t) \geq 0, \quad \limsup_{t\to -\infty}e^{-\alpha t}u(t)<\infty,\quad
t \to \|U(t)\|\;\text{ is bounded},
\ee
has the following property of $u$
\be{limsup2b}
A=\lim_{t\to -\infty}e^{-\alpha t}u(t)\ee exists and the solution of
(\ref{grocoll2})-(\ref{ba10g}) satisfying
(\ref{limsup2b}) for given $A$ is unique.  Furthermore $U(t)$ and $\alpha(t)$  satisfy
\be{limsup2a}
U(t)\La \CU(\infty) \mbox{ as } t \to -\infty, \qquad\text{and }\limsup_{t\to -\infty} e^{-\alpha t}|\alpha(t)-\alpha|<\infty.
\ee

(e) The solution of the system including the age distribution, that is, (\ref{grocoll2}) and (\ref{ba10g}) has a unique solution.
$\qquad \square$
\end{lemma}
The following is immediate from (\ref{agr4}) and the discussion in (\ref{ba10g9ax})-(\ref{add2}).
\beC{CorUCE}
The $({\bf G}, D_0({\bf G}))$ martingale problem is well-posed.
\end{corollary}

\begin{remark} Note that looking at the form of the equation we see that
a solution indexed by $\R$ remains a solution if we make a time shift.
This corresponds to the different possible values for the growth constant
$A$ in (\ref{S081}). In particular the entrance law from 0 at time
$-\infty$ is unique up to the time shift.
\end{remark}

\noi {\bf Proof of Lemma \ref{L.UCE}}

{\bf (a)} First define
\be{agre64}
b(t)=1+\frac{\gamma(t)}{\alpha(t)} \ee and we
express $\alpha(t)$ and $\gamma(t)$ in terms of  $U(t)$ and
and consider $L^1 (\nu,\N)$ as a basic space for the analysis
of the component $U(t,\cdot)$.
Then we obtain a system of coupled differential equations for the pair $(u,U)$:

\be{grocoll29} \frac{du(t)}{dt}= \alpha(t)(1-b(t)u(t)) u(t), \ee

\bea{ba10g9a}
\frac{\partial U(t,j)}{\partial t}
  =
&&
 +s(j-1)1_{j\ne1}U(t,j-1)-sjU(t,j)\nonumber
 \\&&
 +\frac{d}{2}(j+1)jU(t,j+1)-\frac{d}{2}j(j-1))1_{j\ne1}U(t,j)\nonumber
 \\&&
+c(j+1)U(t,j+1)-cjU(t,j)1_{j\ne 1} \nonumber
 \\&&
 +
  \left [\alpha(t) - \alpha(t)u(t)-\gamma(t)u(t))\right] 1_{j=1}
  \\&&
  + u(t)(\alpha(t)+\gamma(t))[U(t,j-1)1_{j\ne1}  - U(t,j)] \nonumber
  \\&&
   -\Big(\alpha(t)(1-u(t))-\gamma(t)u(t) \Big)\cdot U(t,j).
\nonumber\eea

We now consider the related nonlinear system obtained by setting $u\equiv 0$ in (\ref{ba10g9a}), namely,
\bea{ba10g9al}
\frac{\partial \wt U(t,j)}{\partial t}
  =
&&
 +s(j-1)1_{j\ne1}\wt U(t,j-1)-sjU(t,j)\nonumber
 \\&&
 +\frac{d}{2}(j+1)j\wt U(t,j+1)-\frac{d}{2}j(j-1))1_{j\ne1}\wt U(t,j)\nonumber
 \\&&
+c(j+1)\wt U(t,j+1)-cj\wt U(t,j)1_{j\ne 1} \nonumber
 \\&&+
   \alpha(t) 1_{j=1}
     -\alpha(t) \cdot \wt U(t,j).
\nonumber\eea
where  $\alpha(t)=c\sum_{\ell=2}^\infty \ell\cdot \wt U(t,\ell)$.

\bigskip
We will consider the system (\ref{ba10g9a}) as a nonlinear perturbation of the linear strongly continuous semigroup $S^\alpha$ on the Banach space $\mathbb{B}=L^1(\nu,\N)$ given in (\ref{dualsemigroup}).

Then we can rewrite (\ref{grocoll29}) and (\ref{ba10g9a})  in  the form
($G^\ast$ stands for adjoint of $G$)
\be{Grev39m1}
\frac{d}{dt} (u(t), U(t,\cdot))= \vec{G}^\ast (u(t), U(t,\cdot)) = (0,U(t,\cdot)Q^{\alpha,\ast})+F(u(t),
U(t,\cdot)), \ee
where the infinite matrix $Q^{\alpha,\ast}$ is the  infinitesimal generator (Q-matrix) of the
 Markov semigroup of  $S^\alpha$ and denote by $U^{\alpha, \ast}(t)$ the state if
 $S^\alpha$ acts on measures and $F: \B \la \B$ is a (locally) Lipschitz
 continuous function on $\B$, namely,
 \bea{add107}
 && F(u,U)= \left (\alpha(U)u-({\gamma(U)}+{\alpha(U)})u^2,\,((\alpha-\alpha(U))+u(\alpha(U)+\gamma(U))))Z(U)\right ),\\
&&\hspace{3cm} \text{where  } Z(U)=(Z_j)_{j\in\N}, \mbox{ with }  Z_j=U(j-1), j\ne 1, Z_1=-1.\nonumber
 \eea

Note that $F$ involves a third degree polynomial in $u$ and $\alpha(t)$.  Moreover $U\to \alpha(U)$ is a linear function and differentiable on $\mathbb{B}$.  Therefore we can conclude that $F$ is $C^2$ with first and second derivatives uniformly bounded on bounded subsets of $\mathbb{B}$.  Therefore we can use the result of Marsden \cite{MA}, Theorem 4.17
( see Appendix \ref{marsden}) to conclude that there is a unique flow $(u(t),U(t))$ satisfying the equations
and in addition  the mapping $t\to (u(t),U(t))$ and the mapping $(u_{t_0},U_{t_0})\to (u(t),U(t))$
are Lipschitz in $ \mathbb{B}$.

\bigskip
\begin{remark}
{\em (Alternate proof of (a))}

We write the equation (\ref{Grev39m}) as an a system of ODE in the form
\be{Grev39m2}
B^\prime(t) = \vec{G}^\ast (B(t)) = L(B(t)) +N(B(t)),
\ee
where $L$ is a linear operator
(generator of a Markov process) and $N$ is a nonlinear operator
on the Banach space $\B$. The latter is locally Lipschitz but not globally Lipschitz.
Furthermore the operator $L$ maps $\B$ into finite well-defined functions on
$\N$, but {\em not} into $\B$. Only if we restrict to the smaller sets of
configurations $\wt \B$, i.e. configurations  $b$ satisfying
$b \in L^1 (\N, \nu)$ and $\sum j^4 |(b(j)| < \infty$ we obtain
$Lb \in \B$. This subset of configurations forms a dense subset in $L^1 (\N,\nu)$.
This means we cannot view this as a standard evolution equation in
the Banach space $\B$ and have to bring in some extra information to make
the usual arguments work. Namely we have to show that if we start in
the set $\wt \B$ we remain there; we get  an evolution in the Banach space
$\B$. The extra information needed comes from the probabilistic interpretation.

Indeed the evolution equation
\be{Grev39m4}
\wt B^\prime(t) = Q^\ast (\wt B(t))
\ee
has a nice dual Markov semigroup (backward equation) with which
we can work. Namely we know that the corresponding  Markov process
on $\N$ has an entrance law starting at $+\infty$.
This property is wellknown for the Kingman coalescent but extends to our
birth and death process provided $d>0$.
This is easily extended furthermore to an additional immigration term.

Hence the evolution (forward equation)
can be started with ''$\wt B(0)=\delta_\infty$'' and for all solution
we obtain a smaller path in stochastic order, in fact the
entrance law from infinity has the property that in particular
$\suml_j \wt B(t,j)(1+j^4) < \infty$, for $t>0$.
Therefore we can solve the evolution equation (\ref{Grev39m4})
in $\B$ uniquely. This has to be extended to the nonlinear equations
where again starting somewhere in  $\B$ after any positive time we reach
the stronger summability condition.

However if we truncate $L$ to $L^k$ where $L^k x=1_{\{j\leq k\}} Lx$,
we can view (\ref{Grev39m}) as a nonlinear evolution equation in a
subset of the Banach space  $\R \otimes L^1(\N,\nu )$ (recall (\ref{Grev39n})).
This means now we are in the set-up of {\em ODE's in
Banach spaces} and can use the classical theory of existence and
uniqueness for such equations to obtain unique solutions $B^k$.

Then we return to the original equation observing that for the solution
$B^k$ of the truncated problem we have for initial states with
$\suml_i B(0,i)(i)^4 < \infty$ as $k \to \infty$
\be{Grev39m3}
(B^k (t))^\prime \la B^\prime (t), \quad B^k (t) \la B(t), \quad t \geq 0,
\ee
using the bound in stochastic order from the entrance law as
pointed out above.

If we have existence of the solution in $R^+\otimes L^1_+(\N,\nu
)$ for all $t$ (where $L^1_+(\N,\nu)$ denotes the set of sequences
in $L^1(\N,\nu )$ having non-negative components) it suffices to
show local uniqueness.

Turn first to the existence problem.
We shall obtain in Proposition \ref{P.Grocoll}
that we have convergence of the finite $N$ systems to a solution of
the equations and we need for the existence to know that it is in our Banach space.
It can easily be verified that  the
solution obtained in the limit must lie in $\R^+\times
L^{1}_+(\N,\nu )$ and has no explosions using  apriori estimates.
(Recall that for $\zeta^N (t,i)$ the second moments (in fact {\em all} moments) are finite
over finite time intervals since this is true for the process without collision which
provides upper bounds on $u^N,U^N$ independently of $N$.)
Using the fact that second moments remain bounded in bounded
intervals the proof is standard.
Therefore a solution arising as  a limit for $N \to \infty$ lies
in the subspace of $\R^+ \times \CM_1 (\N)$
where $\| \cdot \|$ is finite.
This gives existence of a solution with values in the specified Banach space.

We turn to the uniqueness. The
form of the equation now immediately gives the local Lipschitz
property of $F$. Therefore the solution is locally unique by the
uniqueness result for nonlinear differential equations in Banach
spaces (e.g. \cite{Paz}, Chapt. 6, Theorem 1.2) applied to the truncated
problem and together with the approximation property we get uniqueness
of our equation.
\end{remark}

\bigskip

{\bf (b)} It is standard to check that $\zeta$ is for given $(u,U)$
and for (\ref{Grev39j})-(\ref{Grev39l}) a nonexplosive time-inhomogeneous
Markov process which is uniquely determined by these equations.
We have to show that the self-consistency relation (\ref{Grev39l2})
can be satisfied. This will follow from our convergence result for
$(u^N, U^N)$ below in Proposition \ref{P.Grocoll}, Part (e).

{\bf (c)} To prove  existence of an entrance law with the prescribed properties,
we construct a solution as a limit
of the sequence $(u_n(t),U_n(t))_{t \in \R}, n \in \N$ defined as follows:
\bea{add104}
u_n(t)&&=Ae^{\alpha t_n} \text{  for  } t\leq t_n\\
U_n(t)&&=\mathcal{U}(\infty) \text{  for  } t\leq t_n \nonumber\\
(u_n,U_n)&&\quad \text{satisfy (\ref{grocoll29}) and
(\ref{ba10g9a}) for } t\geq t_n  .\nonumber\eea
 Standard arguments imply the tightness of the
sequence $\{(u_n(t),U_n(t))_{t \in \R}, n \in \N\}$ and that any limit point $(u,U)$ satisfies
(\ref{grocoll29}) and (\ref{ba10g9a}).
We now have to show that the limit points satisfy the conditions
on the behaviour as $t \to -\infty$.

We observe that the solution to the equation for $u \equiv 0, U \equiv \CU(\infty)$
gives a stationary solution of the coupled pairs of equations. We know that the growth by
$A e^{\alpha t}, t \in \R$ is an upper bound for the $u$-part of the limit, since this is the
collision-free solution of the equation and therefore
$\limsup\limits_{t \to \infty} (e^{-\beta t} u(t)) \leq A < \infty.$
We therefore must verify that the limit point
matches exactly that growth and starts from $\CU(\infty)$ at time $t =-\infty$.
We note that $(u,U)$ must then agree with a stationary solution for $u \equiv 0$,
which is unique by standard Markov chain theory and hence the point at
$t = -\infty$ must be $\CU(\infty)$. We next argue that a limit point
$(u,U)$ actually satisfies $e^{-\beta t} u(t)$ converges to $A$ as $t \to -\infty$.
We shall see below that this would actually follow if
$u(t) e^{-\beta t}$ has for $t \to -\infty$ $A$ as limes superior.
If we know that the solution
with this property is unique,  then our limit point in the construction must then
have the desired properties and in particular we must actually have convergence
of $(u_n, U_n)$ as $ u \to \infty$.  This however  we shall show below in (d).

{\bf (d)} Let $(u,U)$ be a solution satisfying (\ref{limsup2}).  Then  there exists a sequence
$t_n\to -\infty$ such that as $n\to\infty$ the  limit of
$(u(t_n),U(t_n))$  exists.

Now consider the nonlinear system $(u(t),U(t))$.  We must verify that as $t\to -\infty$  we actually converge
to the case where $u\equiv 0$ and $\alpha(t)\equiv \alpha$. First note that by (\ref{ddy2a}) we can assume that $\alpha(t)$ is bounded for all $t$.  By
tightness we can find a sequence
$(t_n)_{n \in \N}, \quad t_n \to -\infty$ as $n \to \infty$
and $A \in (0, \infty)$
such that as $ n \to \infty$, $u(t_n)e^{-\alpha
t_n}\to A$.  Since by assumption $\limsup_{t\to -\infty}
e^{-\alpha t}u(t) <\infty$ it follows from the equations that
$|\alpha(t)- \alpha| \leq const\cdot e^{\alpha t}$. Then any such
solution must satisfy
\be{agre66a2}
u(t)-Ae^{\alpha t} = o(e^{\alpha t}).
\ee

We next note that given $u(t), \alpha(t),\gamma(t)$, the equation for $U(t)$ is linear
and has a unique solution by the standard Markov chain theory.
Therefore it suffices to prove that $u(t),\alpha(t),\gamma(t)$
as functions on $\R$ are unique given the asymptotics behaviour of
$u$ as $t \to -\infty$. Remember that $\alpha (t)$ and $y(t)$ are functions
$U(t)$.

We next consider the leading asymptotic term for $U$.
Let  $(S^\alpha(t))_{t \geq 0}$ the semigroup defined in (\ref{dualsemigroup})
and $S^{\alpha,\ast}$ is action on measure, $U^{\alpha,\ast}$ its path,
$U^{\alpha,\ast}(\infty)$ its equilibrium measure. Furthermore let
$L(\cdot) = \alpha (\cdot) + \gamma  (\cdot)$.
Let $t_n,t_m\to -\infty$ so that $t_n-t_m\to \infty$.  Then by (\ref{Grev39m})
\bea{agre66a3a}
 &&{U}(t_n) = \lim_{t_m\to -\infty}
   \big[ S^{\alpha,\ast}(t_n-t_m){U}(t_m)+\int_{t_m}^{t_n}S^{\alpha,\ast}(t_n-s)(\alpha-\alpha(s))Z(U(s))ds \\
&&\qquad +\int_{t_m}^{t_n} u(s)\,S^{\alpha,\ast}(t_n-s) L(s)ds \big]\nonumber \\
&& =\lim_{t_m\to -\infty}
    \big[S^{\alpha,\ast}(t_n-t_m){U}(t_m)\big]
     +\int_{-\infty}^{t_n}S^{\alpha,\ast}(t_n-s)(\alpha-\alpha(s))Z(\mathcal{U}(\infty))ds\nonumber\\
&& \hspace{4.5cm}
   +\int_{-\infty}^{t_n} A e^{\alpha s}\,S^\alpha(t_n-s)(\alpha+\gamma) Z(U^{\alpha,\ast})ds +o(e^{\alpha t_n})\nonumber
\\&&=\mathcal{U}(\infty) +\int_{-\infty}^{t_n} S^{\alpha,\ast}(t_n-s)(\alpha-\alpha(s))Z(\mathcal{U}(\infty))ds\nonumber\\
&& \hspace{3cm}
    + \int_{-\infty}^{t_n} Ae^{\alpha s}\,S^{\alpha,\ast}(t_n-s)(\alpha+\gamma) Z(U^{\alpha,\ast}(\infty))ds+o(e^{\alpha t_n}).\nonumber
\eea
Therefore we have
\be{agre66a3}
\lim_{t_n\to -\infty} U(t_n)= \mathcal{U}(\infty).
\ee
Note that then also $\alpha(t_n) \to \alpha$ and
$\gamma(t_n) \to \gamma$ as $n \to \infty$, since these are continuum functionals
of $U(t_n)$ and we have the integrability properties.

Finally, assume that there are two solutions.
We shall show below that if two such solutions $(u_1,U_1)$ and $(u_2,U_2)$ were
distinct a contradiction would result. Hence in particular in (c) we must have convergence as claimed.

Let
\be{add8}
(v(t),V(t)):=(u_1(t),U_1(t))-(u_2(t),U_2(t)),\; \wh{\alpha}(t)=\alpha_1(t)-\alpha_2(t),\;\wh \gamma(t)=\gamma_1(t)-\gamma_2(t).\ee
>From the above (i.e. (\ref{agre66a2}) and using $(\alpha - \alpha(s)) =o(e^{\alpha s})$
in (\ref{agre66a3a}) and estimating therein the third term by explicit bound) we know that
\be{afon}
|v(t_n)| =o(e^{\alpha t_n})\quad \text{ and } \|V(t_n)\| = o(e^{\alpha t_n}).\ee

Then
\bea{add9}
&&\frac{dv(t)}{dt}=\alpha v(t)+ e^{\alpha t}\wh{\alpha}(t)-(\alpha+\gamma)e^{\alpha t}v(t)-e^{2\alpha t}(\wh\alpha +\wh\gamma),\\
&&
\frac{dV(t)}{dt}= V(t)Q^{\alpha,\ast}-\alpha V(t)+A(\alpha +\gamma)e^{\alpha t}Z((V(t))) +(\alpha +\gamma)\cdot v(t)Z(U^{\alpha,\ast})\nonumber\\
&& \hspace{2cm} +e^{\alpha t}((\wh\alpha+\wh\gamma+ Z( V(t))),\nonumber\\
&& \frac{d \wh{\alpha}(t)}{dt}= -\alpha \wh{\alpha}(t)+A(\alpha +\gamma)e^{\alpha t}\alpha(Z(V(t)))+(\alpha +\gamma)\cdot v(t)\alpha(Z(U^{\alpha,\ast})) +o(e^{\alpha t}),\nonumber
\eea
where $\wh\alpha(t)=\alpha(V(t))$.

 We then obtain by inspection from (\ref{add9}) that:
 \bea{add10}
&&\frac{d(|v(t)|+\|V(t)\|)}{dt}\leq \alpha\cdot(1+\rm{const}\cdot e^{\alpha t}) (|v(t)|+\|V(t)\|).\\
\nonumber
\eea

Then use the $v$ equation and (\ref{afon}), more precisely we use the $V$ equation to get an expression for
$\wh\alpha(t)$ and $\wh\gamma(t)$ as a linear function of $v$.
to conclude $v(t)\equiv 0$ as follows.
Namely we conclude using Gronwall's inequality,
\be{add11}
|v(t_n)|+\|V(t_n)\|\leq o(e^{\alpha t_m}) e^{\alpha (t_n-t_m)}.\ee
Letting $t_m\to -\infty$ we obtain
\be{add12}  |v(t_n)|+\|V(t_n)\|=0.
\ee
This completes the proof of uniqueness of  $(u(t),U(t))$.

\textbf{(e)}
It remains to show that not only $U(t, \R^+, \cdot)$ is unique but
that this holds also for $U(t,\cdot,\cdot)$. Note that our result
already implies that $u(\cdot)$ is uniquely determined. Note also
$\alpha(t)$ and $\gamma(t)$ depend only on $U$ integrated over the
age and are continuous functions of $t$ and are therefore also uniquely
determined. Therefore $u(\cdot), \alpha (\cdot), \gamma (\cdot)$
are uniquely determined. Hence we can insert the unique objects as {\em external}
(time-inhomogeneous) input into the equation.
Then the equation (\ref{ba10g}) for given $u,\alpha,\gamma$ is the forward
equation of a time inhomogeneous Markov process with state space
$M_1([0,\infty)\times \N)$ restricted to elements in $\R \otimes L^1(\N,\nu)$,
which has a unique solution by standard arguments.

This completes the proof of Lemma \ref{L.UCE}.

\subsubsection{The dual process in the collision regime: convergence results}
\label{sss.asycolreg}

Now we are ready to state the limit theorem for the growth dynamics of the
dual population in the critical time scale where collisions are essential, with the key result
for the application to the original process given in
e) of Proposition \ref{P.Grocoll} below, see in particular (\ref{ba7}). As usual we denote by $\La$  convergence in law.

Consider the dual process which starts with $k$ particles at each
of $\ell$ distinct sites and denote the corresponding functionals of the dual process by
$(\Pi^{N,k,\ell}_s,\, u^{N,k,\ell}_s,\, U^{N,k,\ell}_s)_{s\geq
0}$. In order to consider them in the time scale in which fixation occurs, we introduce:
\be{ts1}
 \wt u^{N,k,\ell}(t)
 = u^{N,k,\ell}((\frac{\log N}{\alpha}+t)\vee 0),\; t\in (-\infty,T],\quad \wt u^{N,k,\ell}(-\frac{\log N}{\alpha})=\ell\ee
\be{ts2}
\wt U^{N,k,\ell}(t)
 = U^{N,k,\ell}((\frac{\log N}{\alpha}+t)\vee 0),\; t\in (-\infty,T],\quad \wt U^{N,k,\ell}(-\frac{\log N}{\alpha})=
 \delta_{(k,0)}.\ee
Furthermore we need the standard solution of the nonlinear system
(\ref{grocoll2}) and (\ref{ba10g}), (see Lemma \ref{L.UCE} part (c)) denoted
\be{ts3} (u^\ast(t),
U^\ast(t, \cdot,\cdot)))_{t \geq t_0},
\ee
which is the solution
satisfying (recall $\CU(\infty)$ was the stable age and size
distribution) \be{stsol1} \lim_{t\to -\infty}e^{-\alpha
t}u^\ast(t)=1,\quad \lim_{t\to -\infty} U^\ast(t)=
\mathcal{U}(\infty). \ee
Finally we have to specify three time scales.
Let  $t_0$ stand for some element of $\R$ and let
$t_0(N)$ and $s(N)$ be as follows:
\be{stsol2}
t_0(N)=s(N)-\frac{\log N}{\alpha}, \ee with
\be{ts4}
s(N)\to\infty, s(N)=o(\log N).
\ee
Recall the time shift in (\ref{ts1}) and
(\ref{ts2}).

The behaviour of the dual particle system is asymptotically as $N \to \infty$
as follows.

\beP{P.Grocoll}{(Growth of dual population in the critical
time scale in the $N\to\infty$ limit)}

(a) Assume that for some $t_0 \in \R$ as $N\to\infty$, $(\frac{1}{N} \wt u^{N,k,\ell}
(t_0), \wt U^{N,k,\ell} (t_0))$ converges in law to the pair
$(u(t_0),U(t_0))$ automatically contained in $[0,\infty) \times L_1(\N, \nu)$.

Then as $N \to \infty$
\be{agrev60}
\CL \left[\big(\frac{1}{N} \wt u^{N,k,\ell}(t),\,
 \wt U^{N,k,\ell}(t, \cdot, \cdot) \big)\big)_{t\geq t_0}\right]
 \La \CL \left[(u(t), U(t, \cdot,\cdot)))_{t \geq t_0}\right],
\ee
in law on pathspace, where the r.h.s. is
supported on the solution of the nonlinear system
(\ref{grocoll2}) and (\ref{ba10g}) corresponding to the initial state
$(u(t_0), U(t_0))$. (Note that the mechanism of the limit dynamics
does not depend on $k$ or $\ell$, but the state at $t_0$ will.
That is why we suppress $k,\ell$ on the r.h.s.)

(b) The collection  $\{(\wt u^{N,k,\ell},\wt U^{N,k,\ell}), N \in \N\}$ can be
constructed on a common probability space such that for
$t_0 (N)$ as in (\ref{stsol2})  we have:
\be{agrev61}
(e^{-\alpha s(N)}\wt u^{N,k,\ell}(t_0(N)),\wt
U^{N,k,\ell}(t_0(N)))\to (W^{k,\ell},\mathcal{U}(\infty)),\;
a.s., \mbox{ as } N\to\infty.
\ee

(c) The scaled occupation density converges (recall (\ref{ts1}) and (\ref{ts2})
and (\ref{stsol1})):
\be{grocoll1}
 \frac{\wt u^{N,k,\ell}}{N}({t})\Nto u^{k,\ell}(t),\; t\in
 (-\infty,T],
\ee
in the sense that for $\ve >0$ (with $t_0(N)$ as in (\ref{stsol2})),
\be{agrev62}
P\left(\sup_{t_0(N)\leq t\leq T} e^{-\alpha t}|\frac{\wt u^{N,k,\ell}}{N}({t})- u^{k,\ell}(t)| >\ve
\right)\to 0 \text{ as }N\to\infty.
\ee

Also the age-size distribution of the dual converges:
\be{ba10i2}
\wt U^{N,k,\ell}(t,\cdot,\cdot) \Nto
U^{k,\ell}(t,\cdot,\cdot), \; t\in (-\infty,T],
\ee
in the sense that for $\eta>0$,
with $\| \cdot \|$ denoting the variational norm,
\be{ba10i2a}
\lim_{N\to\infty}P\left(\sup_{t_0(N)\leq
t\leq T}\|\wt U^{N,k,\ell}(t, \cdot,\cdot) -
U^{k,\ell}(t,\cdot,\cdot)\|>\eta \right) =0.
\ee

(d)
The limits $(u^{k,\ell}(t), U^{k,\ell}(t, \cdot,\cdot)))_{t \geq
t_0},$  can be represented as the unique solution of the nonlinear system (\ref{grocoll2})
and (\ref{ba10g}) satisfying
\be{agrev63}
\lim_{t\to -\infty}e^{-\alpha t}u^{k,\ell}(t)=W^{k,\ell},\; \lim_{t\to -\infty} U^{k,\ell}(t)= \mathcal{U}(\infty),
\ee
with $W^{k,\ell}$ having the law of the random variable appearing as
the scaling (by $e^{-\alpha t}$)
limit for $t \to \infty$  of the CMJ-process $K^{k,\ell}_t$ started with $k$ particles at each of $\ell$ sites.

These solutions are random time shifts of the standard solution (\ref{ts3}), (\ref{stsol1}), namely:
\be{agrev64}
 u^{k,\ell}(t)= u^\ast(t+\frac{\log W^{k,\ell}}{\alpha}),\quad  U^{k,\ell}(t)= U^\ast(t+\frac{\log W^{k,\ell}}{\alpha}).\ee

(e) Let $(\Pi^{N,k,\ell}_u)_{u \geq 0}$ denote the number of dual
particles and $T_N = \alpha^{-1} log N$ and $t_0(N), s(N)$ as in (\ref{stsol2}). Then for each
$k\in\mathbb{N}, \ell \in \N$ we have:
\be{angr24}
\CL [\{N^{-1}
\Pi^{N,k,\ell}_{T_N+t},\;t_0(N)<t<\infty\}]  \Nto \CL[
\{\frac{u^{k,\ell}(t)(\alpha(t) + \gamma(t))}{c}\}_{t\in\mathbb{R}}]],
\ee
\be{agr24A}
\CL [(N^{-1} \Pi^{N,k,\ell}_u \quad , \quad 0 \leq u \leq s(N))]
\Nto \delta_{\underline{0}},
\ee
with $\underline{0}$ the function on $\R^+$, which is identically zero.

Furthermore for every $t \in \R$, there exists (a deterministic) $\nu_{k,\ell}(t)\in
\CP([0,\infty))$ such that \be{ba7} \CL [\frac{1}{N}
\intl^{T_N+t}_0 \Pi^{N,k,\ell}_u du] \Nto \nu_{k,\ell}(t).\qquad
\square \ee
\end{proposition}

\begin{remark}
Note that the dual at time $T$ describes the genealogy of a
typical sample from the original population drawn at the time $T$.
Here time runs backwards, in particular the evolution of the dual
with total time horizon $T = \alpha^{-1} \log N+t$
in the time described by varying $t$
corresponds in the genealogy of the tagged sample to the early
moments, while the emergence phenomena of the dual population
(arising collisions) reflect in the original
population later times, where law of large number effects rule
the further expansion of the fitter type on macroscopic scale
(i.e. observing the complete space).

In the case in which the rare mutant has succeeded our result  tells
us that we have  a situation where the mutations which occur  very
early in the typical ancestral path have generated $O(N)$ possible
choices for a fitter type to occur and to then prevail at the
observation time at a fixed observation site.
  The
randomness enters since these very early mutations during the
first moments of the evolution have in the limit $N \to \infty$ a Poisson structure. This
produces  a random time shift of emergence and takeover which then
follows a deterministic track. This is the global picture, if we look
at it locally we find that also a collection of tagged sites
follows  a random evolution in small time scale.
\end{remark}

{\bf Proof of Proposition \ref{P.Grocoll} }

We begin with some preparations needed for all of the proofs of the
different claims of the proposition, this is part 0, followed by
three further parts explained at the end of part 0.

{\bf Part 0 (Preparation)}

For the proof we represent the process differently, namely we
rewrite the problem in measure-valued form.
We will need the following functional of the process
$\Pi^N_t$, denoted $\wh \Pi^N_t$ respectively
$\wh \Pi^{N,k}_t, \wh \Pi^{N,(k,\ell)}_t$ if we
indicate the initial state in the notation (here $s_i = s^N_i)$:
\be{ba10i4}
\wh \Pi^{N}_t =
\suml^{ K^N_t}_{i=1} \zeta^N_{s_i} (t) 1(\zeta^N_{s_i}(t) >1), \ee
which counts at time $t$  only those dual particles which are not the only ones
at their site, i.e. they are the ones which can generate new
occupied sites. Here $\zeta_s^N(t)$ is the single site birth and
death process starting at time $s$, the time where the site was first occupied, with mean-field immigration at time $t$ given by
$c\frac{\Pi^{N,1}_t}{N} -o(1)$ as $N \to \infty$ from the other
sites.

Note that we have now made the $N$-dependence explicit.
Also note that by definition of $U^N, u^N$ and $\wh \Pi^N$ we have
\be{dd1b}
  \wh\Pi^{N,1}_t = u^N(t) \int_0^\infty \sum_{j=2}^{\infty} j U^N(t,ds,j).\ee

To continue we introduce a new random object from which we can
read off $(u^N, U^N)$ by taking suitable functions. Namely we
consider $\Psi^N(t,ds,dx)$ the
$\CM_{\mathrm{fin}}([0,\infty)\times \mathbb{N}_0)$- valued
random variable (finite measures on $[0,\infty) \times \N_0$) defined by
\be{dd2b}
 \Psi^N (t,ds,dx)=\sum_{s_i\leq t} \delta_{(t-s_i,\zeta^N_{s_i}(t))},
\ee
where $\{s_i\}$ (recall $s_i=s^N_i$)  denote the times of birth of new sites and $i$ runs from
1 to $K^N_t$ (this set of times depends on $N$).
Note that here {\em empty sites are not counted}.
We know from the CMJ-theory that
 $\Psi^N(t,\cdot,\cdot)$ is bounded by
 $W^\ast e^{\alpha t}$. Since we are interested in times
 $\alpha^{-1} \log N + r$ we introduce the scaled object:

\be{dd5}
\wh\Psi^N_t=\frac{\Psi^N_t(\cdot,\cdot)}{N}. \ee Note
that for fixed $N$ this is an element of $\CM_{\mbox{\tiny fin}}
([0,\infty) \times \N_0)$.

Next note that indeed we obtain our pair $( u^N, U^N)$ as a
functional of $\Psi^N$. Observe:
\be{dd3bb} u^N (t)
=\Psi^N_t ([0,\infty)\times \mathbb{N}), \quad U^N(t,\cdot,\cdot) =
\frac{\Psi^N_t  (\cdot,\cdot)}{u^N(t)}
= \frac{\wh \Psi^N_t  (\cdot,\cdot)}{u^N(t)/N}
\ee and the stochastic
integral below gives:
\be{dd4}
\Pi^N_t= \int_0^t\zeta^N_u(t) d
K^N_u=\int_0^t\int_0^\infty x \Psi^N_t  (ds,dx).
\ee

Therefore from (\ref{dd3bb}) we obtain the convergence of $(N^{-1} u^N,U^N)$ in the time scale
$\alpha^{-1} \log N+t$ with $t \in \R$ if we show the convergence
of $\wh \Psi^N$. This will also allow us later on  to argue that
the mass and time scaled process
\be{dd5aa}
 \bar \Pi^N (t)
 = N^{-1} \Pi^N_{\frac{\log N}{\alpha}+t}
 \ee
converges as $N \to \infty$ to a limit.

The proof is broken into three main parts. The first two parts concern the pair
$(u,U)$ the third one the quantity $N^{-1} \Pi^{N,k,\ell}_{T_N
+t}$. The parts can be outlined as follows (where
$(a), \cdots$ refers to the parts of the proposition):

\begin{enumerate}\item Assuming the marginals
of $(N^{-1} \wt u^N, \wt U^N)$ converge in distribution at $t=t_0$, prove (a)  that the
limiting evolution during $t\geq t_0$ follows the deterministic
nonlinear dynamics specified in the system of coupled equations
(\ref{grocoll2}), (\ref{ba10g}) and that pathwise convergence holds. \item Then prove (b) that the
marginal distributions converge for $t_0 = t_0 (N) =
\frac{1}{\alpha} \log N +s, s \in \R$ and show (c), (d)
on the corresponding entrance behaviour.
\item Show the convergence
of the process $\{N^{-1} \Pi^{N,k,\ell}_{T_N+t}, t \in \R\}$ and
of its time integrals to complete the proof of (e).
\end{enumerate}
\bi

{\bf Part 1 (Convergence of dynamics of $\wh \Psi^N$)}


Let  $P^N\in \mathcal{P}(D([t_0,\infty),\mathcal{M}_{\rm{fin}}(\R_+\times\N)))$ denote the law of the Markov process
\be{dd5a}
\big\{\wh\Psi^N_{T_N+t}\big\}_{t\geq t_0}$ (recall (\ref{dd5})) with $T_N=\frac{\log N}{\alpha}.\ee

In this part we assume that the marginal distributions  $P^N_{t_0}$ converge as $N\to\infty$.
Our next goal  is then to show that the processes converge in law on path space, that is,  $P^N\Rightarrow P$ as $N\to\infty$ where $P$ is the law of a deterministic dynamics.

Recall that we have defined the candidate for the limiting object by
a nonlinear equation, which we can view as forward equation for a nonlinear
Markov process on $\R_+ \times \N$. This means we want to show that if
we consider the marginal laws, that is, elements $P^N_t$ of
$\mathcal{P}(\mathcal{M}_{\rm{fin}}(\R_+ \times \N))$ they converge to an element
$\delta_{\wh\Psi_t}$ with $\wh\Psi_t \in \mathcal{M}_{\rm{fin}} (\R_+ \times \N)$.
We will now take the viewpoint of backward
equations and their generators, that is equations for expectations
of measures rather than actual weights.

In order to do
so in this present Part 1, we proceed as follows.
\begin{itemize}
\item we obtain moment bounds on $\wh\Psi^N_{T_N+t}$ uniform for $N\in\N$ and $t_0\leq t\leq t_0+T$,
\item  we formulate the martingale problem for $\wh\Psi^N$  with generator ${\bf G}_N$ given by  a linear operator on the space
$C_b(\CM_{\mathrm{fin}}(\R_+ \times \N_0))$) with common domain
$\mathcal{D}$ given by a  certain dense subset of $C_b (\CM_{\mathrm{fin}}(\R_+ \times \N_0))$,
\item we verify the tightness of the laws $P^N$,
\item we prove the weak convergence of $P^N$ to $P$ by showing
that ${\bf G}_N$ converges pointwise to  ${\bf G}$  on these functions in the domain $\mathcal{D}$ where
${\bf G}$ is the operator defined in (\ref{ba10g21}),
\item verify that limit points of $(P^N)_{N \in \N}$ are characterized as solutions to the ${\bf G}$ martingale problem,
\item conclude the convergence since  we
have proved that the ${\bf G}$ martingale problem has a unique solution given by the nonlinear system (\ref{grocoll2}), (\ref{ba10g}),
\item the convergence for fixed $t_0$ is obtained later in Part 2.
\end{itemize}

Abbreviate for bounded measureable test functions $f$ on $[0,\infty) \times \N_0$
\be{dd6a}
\wh \Psi (f) = \intl^\infty_0 \intl^\infty_\infty f(r,x)
\wh \Psi (dr,dx)
\ee
and in particular $\wh\Psi(1)=
\wh\Psi([0,\infty)\times (\mathbb{N}_0))$.

The process  $\wh\Psi^N$ is a Markov process with state space  $\mathcal{M}_{\rm{fin}}([0,\infty) \times \N)$.  We  calculate the action of the generator ${\bf G}_N$ acting on the
domain $\mathcal{D}$ of functions $F_{g,f}$ of the form  (\ref{dd6}) as follows. Let
$\wh\Psi \in \CM_{\mathrm{fin}}([0,\infty) \times \N)$. Then:

\be{dd7}\begin{array}{lll}
&&({\bf G}_NF_{g,f}) (\wh\Psi)=\\&&
c(1-\wh\Psi(1))\left(\int_0^\infty\int_0^\infty N
\left[g(\wh\Psi(f) + \frac{f(0,1)+f(r,x-1)-f(r,x)}{N})
-g(\wh\Psi(f))\right]x1_{x\ne
1}\wh \Psi(dr,dx)\right)
\\&&
+\intl_0^\infty\intl_1^\infty g^\prime\left(\wh\Psi(f)\right)
\frac{\partial f(r,x)}{\partial r}\wh\Psi(dr,dx)  \\
&&+ \intl_0^\infty\intl_1^\infty   sx N
\left[g(\wh\Psi(f)+\frac{f(r,x+1)-f(r,x)}{N})-g(\wh\Psi(f))\right]\wh\Psi(dr,dx)\\
&&+ \intl_0^\infty\intl_1^\infty \frac{d}{2}x(x-1)N
\left[g(\wh\Psi(f)+\frac{f(r,x-1)-f(r,x)}{N})-g(\wh\Psi(f))\right]\wh
\Psi(dr,dx)\\
&&+ c\wh\Psi(1)\intl_0^\infty\intl_1^\infty \intl_0^\infty\intl_1^\infty
N\left[g(\wh\Psi(f)+\frac{f(\tilde r,\tilde x+1)-f(\tilde r,\tilde
x)}{N}+\frac{1_{x\ne 1}f( r, x -1)-f( r, x)}{N}
)-g(\wh\Psi(f))\right]\\
&&\hspace{10cm} x\wh\Psi(d r,d
x)\frac{\wh\Psi(d \tilde r,d\tilde x)}{\wh\Psi(1)}.
\end{array}\ee
Note that in order for the right side to be well-defined we require $\int x^2 \wh\Psi(\R,dx)<\infty$.
This condition is automatically satisfied if $\wh\Psi\in \CM_{\mathrm{fin}} ([0,\infty) \times \N_0)\cap \mathbb{B}$.
This means that $\mathbf{G}_N F_{g,f}$ is a welldefined function on $\B$, where it is even
continuous. Note however
$C_b(\B, \R)$ is {\em  not} mapped in $C_b(\B, \R)$.
This means we now have to investigate next that under our evolution we stay with
$\wh \Psi$ being $\B$. Therefore we turn next to bounds on the empirical moments
of the number of particles at the occupied sites.

We now establish some apriori bounds which will guarantee that the empirical moments
remain stochastically bounded over the considered time interval.
By definition we have
\be{add13}
\int_0^\infty x^m\wh\Psi^N_{T_N+t}(\R,dx)=\frac{1}{N}\sum_{i=1}^{u^N(T_N+t)}\sum_{k=1}^\infty k^m1_k(\zeta^N_{s_i}
(T_N+t)).
\ee
Then we have with $s_i$ denoting the time of birth of the $i$-th colonized site
(recall $s_i = s^N_i$) and letting $S$ be the birth time of a randomly picked site:

\bea{add14}
&& E\left[\int_0^\infty  x^m\wh\Psi^N_{T_N+t}(\R,dx)\right]=\frac{1}{N}
E\left[\sum_{i=1}^{u^N(T_N+t)}\sum_{k=1}^\infty
k^m1[(\zeta^N_{s_i} (T_N+t)=k)]\right]
\\
&& =  E\left[\frac{1}{N}\sum_{i=1}^{u^N(T_N+t)}(\zeta^{N}_{s_i})^m (T_N+t)\right]
\leq E\left[(\zeta^N_S(T_N+t))^m\right]
\nonumber\eea
and by Cauchy-Schwarz applied to the $m$-th process of $\zeta^N_{s_i}$ we get

\bea{add15}
&& E\left[\left(\int_0^\infty  x^m\wh\Psi_{T_N+t}(\R,dx)\right)^2\right]
\\
&&\leq  \left(E\big[\frac{1}{N}\sum_{i=1}^{u^N(T_N+t)}(\zeta^N_{s_i}(T_N+t))^m\big]\right)^2
+ E\big[\frac{1}{N}\sum_{i=1}^{u^N(T_N+t)}(\zeta^N_{s_i}(T_N+t))^{2m}\big].\nonumber
\eea

We can obtain an upper bound for $\int_0^\infty  x^m\wh\Psi^N_{T_N+t}(\R,dx)$
by comparing with the collision free case (recall construction in Step 2  of Subsubsection \ref{sss.collreg} that is, the total number of particles is bounded by the collision free case) to obtain

\bea{b1}
&&\lim_{L\to \infty} P\left[\sup_{t_0\leq t\leq t_0+T}\int_0^\infty x \wh\Psi^N_{T_N+t}(\R,dx) >L\right]\\
&&\leq \lim_{L\to \infty} P\left[\frac{1}{N}\sum_{i=1}^{u_N(T_N+t)}\sup_{t_0\leq t\leq t_0+T}
\zeta^N_{s_i}(T_N+t)>L\right]\nonumber\\
&&\leq \lim_{L\to\infty}
 \frac{E\left[{\sup_{t_0\leq t\leq t_0+T}}(\zeta_{s}(t))^2\right]}{L^2}=0,\nonumber
\eea
using the independence of the $\{\zeta_{s}\}$ in the collision-free case and the martingale problem
together with  a standard martingale inequality for the birth and death process $\zeta_{s}(\cdot)$.

We have seen that the first empirical moment bound (\ref{b1}) holds by comparison
with the collision-free regime. We now have to establish analogous results for higher empirical
moments and incorporating collisions which involves additional
immigration from the other colonies at each occupied site.
Let  $K>0$.
Let
\be{add136}
\zeta^{N,K}_{s_i}(t)
\ee
denote the birth and death process with immigration at rate $K$
which starts at time $s_i$. This process is ergodic and has a unique equilibrium. Let $m_n(K)$
denotes the corresponding $n$-th equilibrium moment.

\begin{remark}
We used above the fact that the birth and death process with linear birth and quadratic death rates has moments of all orders at positive times. This can be verified by comparing with a subcritical branching process with immigration. Moreover the corresponding  Markov process
on $\N$ has an entrance law starting at $+\infty$.
This property is wellknown for the Kingman coalescent but extends to our
birth and death process provided $d>0$ and is easily extended  to include an additional immigration term.
\end{remark}

Then up to the time
\be{add137}
\tau_K:=\inf\{t: \int_0^\infty x\wh\Psi_{T_N+t}(\R,dx) \geq K\}
\ee
we can verify with a stochastic montonicity argument that we get a bound uniformly in $N$:
\be{add16}
E[(\zeta^{N,K}_{s_i}(\cdot))^n]\leq m_n(K)<\infty,\;\text{  for  }n,N\in\N.
\ee

Now recall that (writing $F_{x,f}$ if $g(x)=x$)
\be{add17}
M^N_t:=\wh\Psi^N_{T_N+t}(f)-\wh\Psi^N_{T_N+t_0}(f)-\int_{T_N+t_0}^{T_N+t} G_NF_{x,f}(\wh\Psi^N_{T_N+s}(\R,dx))ds
\ee
is a martingale. Using the above moment inequalities, (see (\ref{add15})),
we can show that this is an {\em $L^2$ martingale} and we can then use a martingale inequality to show that for $f$ with finite support
\bea{add18}
&& \lim_{L\to \infty} \sup_N P(\sup_{t_0\leq t\leq (t_0+T)\wedge \tau_K} \wh\Psi^N_{T_N+t}(f) >L)\\
&& \leq \lim_{L\to\infty}\sup_N\frac{4}{L^2}
\left\{E\left[\left(\wh\Psi^N_{T_N+t_0}(|f|)+\int_{T_N+t_0}^{T_N+T} |G_NF_{x,f}(\wh\Psi^N_{T_N+s}(\R,dx))|ds\right)^2\right]\right\}\nonumber\\
&&+ \lim_{L\to\infty}\sup_N\frac{4}{L^2}
 E[(M^N_{(t_0+T)\wedge \tau_K})^2]\nonumber \\
 &&=0.\nonumber
\eea
For every $n\geq 2$ we can now choose a sequence $(f_m)_{m \in \N}$ with finite
support in $C(\R,\R)$ and let $f_m(x)\uparrow x^n$  verify that for given $K$

\be{apriori1}
\lim_{L\to \infty}\sup_N P\left[\sup_{t_0\leq t\leq (t_0+T)\wedge \tau_K}\int_0^\infty x^n
\wh\Psi^N_{T_N+t}(\R,dx) >L\right]=0.
\ee
We next note that equation (\ref{b1}) implies that
\be{add114}
\tau^N_K \uparrow \infty \mbox{ as } K \to \infty
\ee
uniformly in $N$.
This together with (\ref{apriori1}) implies that the evolutions
$\{ \wh \Psi^N_t, \quad t \in [T_N +t_0, T_N +t_0 +T]\}, \quad N \in \N$
and their weak limit points are concentrated on {\em path with values in $\B$}.

We can conclude now with noting that ${\bf G}_N(F_{g,f})$ is a finite and even continuous function
  on $\B$, even though ${\bf G}_N$  does not map functions in its domain $\subseteq C_b(\B,\R)$ into
  $C_b (\B, \R)$ but only into
 the set $C (\B, \R)$ of {\em continuous functions} on $\B$.
Therefore we will have to show first that for the starting points of the evolution we
work with, namely the processes
$\{\wh \Psi^N(t), t\in [t_0,T]\}$
and the limit points of their laws are concentrated
on states satisfying this extra condition to be in $\B$.
This is now immediate from the  a priori bound (\ref{apriori1}).

We now verify first tightness in path space and then in a further step
convergence of the distribution.

{\bf Tightness}
We verify the {\em tightness} in path space
that is the tightness of the laws

\be{grevx2}
\{P^N=\CL[(\wh \Psi^N_t)_{t \geq t_0}]; N \in \N\} \mbox{ as } N \to \infty,
\ee
assuming convergence of the marginals at
time $t_0$.
Here we deal with a sequence of probability  measure-valued
(on a Polish space)
jump processes. With the help of Jakubowski's criterion, (see Theorem
3.6.4 in \cite{D}) we obtain that it suffices to consider
the tightness of laws of real valued semimartingales, namely
 of $\{\CL [(F_{g,f}(\wh \Psi^N_t))_{t \geq 0}], \mbox{ admissible } f \mbox{ and } g\}$.
In fact $g(x)=x$ and
$g(x)=x^2$ (recall Remark \ref{R.x2}) alone suffice according to that theorem
and for these two functions $g$ and all $f$ with $0 \leq f \leq 1$
we shall now verify the criterion.

Recall first  the apriori bound (\ref{apriori1}).
The tightness of the laws of $F_{g,f} \circ \wh \Psi^N$
is is obtained using the Joffe-M\'etivier criterion (\cite{JM})
in the form of Corollary 3.6.7. in \cite{D},
applied to $F_{g,f}$ for $g(x)=x$ and $g(x)=x^2$.
The  bound (\ref{apriori1}) verifies the
conditions therein. This means that now
$\{\CL[(\Psi^N_{T_N+t})_{t \geq t_0}], N \in \N\}$
has limit points in the set of laws on Skorohod space.

Let $P^N_{\wh\Psi^N_0}$ denote the law of $\{\wh\Psi^N_0\}_{t\geq t_0}$ on the Skorohod space $D([0,\infty),\mathcal{M}_{\rm{fin}}(\R_+\times\N)\cap\mathbb{B})$ with initial values $\wh\Psi^N_{t_0} \Rightarrow \wh\Psi_{t_0}\in \mathcal{M}_{\rm{fin}}(\R_+\times\N)\cap\mathbb{B}$. By tightness we can choose a convergent subsequence.  Since the largest jumps decrease to $0$ the limit measure is automatically concentrated on $C([0,\infty),\mathcal{M}_{\rm{fin}}(\R_+\times\N)\cap\mathbb{B})$ (see \cite{EK2}, Chap.3, Theorem 10.2).
Convergence
in law on the Skorohod space follows if any limit point satisfies a martingale problem with a unique solution.

{\bf Convergence}
Next we turn to the convergence in f.d.d. by establishing {\em generator convergence}.
For this purpose we consider $\mathcal{M}_{\rm{fin}}(\R_+\times\N)\cap\,\B$ as the
state space for the processes $\wh\Psi_N(\cdot)$ and the  {\em function space} $C_b (\mathcal{M}_{\rm{fin}}(\R_+\times\N)\cap\,\B, \R)$
as the basic space for our semigroup action.
Recall the algebra of functions $D_0({\bf G})$ given in (\ref{ba10g21}).
The operator ${\bf G}_N$ is a {\em linear map} from
$D_0({\bf G})$ into $C(\mathcal{M}_{\rm{fin}}(\R_+\times\N)\cap\,\B,\R)$.

We consider the convergence for two  choices for the ``initial
state'' at some time $t_0 = t_0(N)$, namely we assume that
$\wh \Psi^N_{t_0}$ converges as $N \to \infty$, but either
\be{dd7a}
\wh \Psi^N_{t_0} (1) \Ntoo 0 \quad \mbox{ or } \quad \wh
\Psi^N_{t_0}(1) \Ntoo x > 0.
\ee

In the first case we would be
back in the linear regime (no collisions) we discussed in Subsubsection \ref{sss.growth} and
which we have identified in terms of the CMJ process. Next we consider the limit corresponding to initial points at some time
$t_0$ such that $\wh \Psi^N(t_0)$ does not converge to the zero measure.

Let
$f \geq 0, \quad \|f\|_\infty \leq 1, \quad \wh \Psi(1)\leq 1$,
$g$ with $g|_{([0,1], \R)} \in C^2_b ([0,1], \R)$ and $F_{g,f}$ as in (\ref{dd6}). Recall the definition of the operator ${\bf G} $ on $D_0({\bf G})$ given in (\ref{ba10g21}). Then
there exists a constant (const)  depending only on $c,d,s$ such that:
\be{dd7b2} |{\bf G}_N(F_{g,f})(\wh \Psi)-{\bf G}(F_{g,f})(\wh \Psi)|\leq
\frac{\text{const}}{N} \|g^{\prime
\prime} 1_{[0,1]} \|_\infty \int_0^\infty\int_0^\infty x^2\wh \Psi(dr,dx)
\leq \frac{\text{const}}{N} \| \wh \Psi \|_\B.
\ee
We get that for all $F_{g,f}$ of the form in (\ref{dd6})
with the above extra restrictions that
\be{dd7c} {\bf G}_N
(F_{g,f})(\wh \Psi) \Ntoo {\bf G}(F_{g,f})(\wh \Psi),
\ee
provided that the argument $\wh \Psi$ satisfies
\be{agdd63}
\int_0^\infty\int_0^\infty x^2\wh \Psi(dr,dx) < \infty,
\ee
which is the case for $\wh \Psi \in \B$. More precisely,
considering the functions ${\bf G}_N(F_g), {\bf G}(F_g)$
on any finite ball in $\B$, we have
\be{dd7c1}
{\bf G}_N(F_{g,f}) \la {\bf G}(F_{g,f}), \quad N \to \infty
\ee
pointwise uniformly on such balls.

From (\ref{dd7c1}) and the convergence of the initial laws we want to conclude that weak limit points of
 $\CL [(\wh \Psi^N_t)_{t \geq t_0}]$
conditioned on a particular initial state at time $t_0$
are $\delta$-measures on the space of $\CM_{\rm{fin}}(\R^+ \times \N)$-valued path
which must satisfy the evolution equation

\be{add115}
\frac{d}{dt} (\wh \Psi_t) = G^{\ast} \wh \Psi_t \quad,\quad t \geq t_0.
\ee
For this we have to show
\begin{enumerate}  \item the variance process of
$<\wh \Psi^N,f>$ converges to zero, and
\item  $\wh\Psi^N$ remains in large enough balls
with high probability such that (\ref{dd7c1}) can be applied,

\item  making the choice $g(x)=x$  in (\ref{dd6}),  to verify that the path solves equation
(\ref{add115}),
\item  to prove convergence it then suffices to show that the ${\bf G}$-martingale problem is wellposed.

\end{enumerate}

For the  point (1) consider $g(x) =x^2$ and $g(x)=x$ and
show that
\be{add116}
\liml_{N\to\infty}[{\bf G}_N(F_{x^2,f}) -2F_{x,f} {\bf G}_NF_{x,f}] =0.
\ee

Note that replacing $G_N$ by $G$ the relation says that $G$ is a {\em first order}
operator and this is read of from the form of $G$ directly.
This can be verified also directly using  (\ref{ba10g21}) and (\ref{dd7b2}).
This implies that the variance of a solution to the limiting martingale problem is 0.
Since the higher moments of $\langle \wh \Psi^N,f\rangle$ are bounded this implies (1).
Point (2) is given by the apriori estimate (\ref{apriori1}).

For point (3) we take $g(x)=x$ and conclude that the limit point is given by the unique trajectory in $\mathbb{B}$ which satisfies
\be{add19}
\int_0^\infty\int_0^\infty f(r,x)\wh\Psi_t(dr,dx)-\int_0^t \int_0^\infty\int_0^\infty G^{\wh\Psi_s}f(r,x)\wh\Psi_s(dr,dx)ds=0,
\ee
that is, the nonlinear system (\ref{psidefn2}) (recall Lemma \ref{L.UCE}). In other words, any limit point is a solution to the ${\bf G}$ martingale problem.
For point (4),  we have established the uniqueness and properties of the solution to the ${\bf G}$ martingale problem in Lemma \ref{L.UCE} and Corollary \ref{CorUCE}.

To complete the proof of convergence, given $\wh\Psi_{t_0}$ and $T$, find $K_0>0$ such that the solution of the nonlinear system stays in the interior of the ball of radius $K_0$ during the interval $[t_0,t_0+T]$ and let $K\geq 2K_0$, say.  Then it follows from  the above that we have weak convergence up to time $ (t_0+T)\wedge \tau_K$
for every $K>0$ and therefore  for $\ve>0$
\be{add20}
\lim_{N\to\infty}P[\sup_{t_0\leq t\leq t_0+T}\|\wh\Psi^N_t-\wh\Psi_t\|>\ve]=0.\ee

This completes the proof of part (a) of the proposition.

\bigskip
\textbf{Part 2 Convergence of $t_0$-marginals and properties of $(u^N, U^N$)}

Our goal is to prove part (b) and (c) of the proposition under the assumptions
of convergence at a fixed time $\alpha^{-1} \log N+t_0$.

\begin{proof}{\bf of (1)}

For that purpose
we prove first the convergence of the marginals at times
$\alpha^{-1} \log N +t_0$, starting with ''$t_0=-\infty$'', (\ref{agrev61}), and then
later consider $t_0 \in \R$ to finally come to the convergence of the path.

{\bf Proof of (\ref{agrev61})}
Here we have a time
scale where as $N \to \infty$ the collisions become negligible since the observation time
diverges sublogarithmically as function of $N$. Using the multicolour particle system
(see Subsubsection \ref{sss.collreg}) we can build all variables
on one probability space. From this construction it follows that
collisions become negligible and we have the claimed convergence
from the one for the collision-free system proved in
Subsubsection \ref{sss.dualcfr}, see (\ref{ang2b}), and
Subsubsection \ref{sss.dualcfrprop} see (\ref{ang4c2})
via the CMJ-theory. This proves part (b) of the proposition.
\end{proof}
\bi

{\bf Proof of (c).}
We now have to prove the convergence result for the marginal of
$(\wh \Psi^N_t)_{t \in \R}$ at times $t = \alpha^{-1} \log N
+t_0$, where we can choose $t_0 \in \R$ as small as we want since
we know that the dynamic converges to a continuous limit dynamic. This
means we have to show (\ref{grocoll1}), (\ref{agrev62}) and (\ref{ba10i2}).
Finally we have to relate the limit of $(N^{-1} \wt u^N, \wt U^N)$ as
$N \to \infty$ to the system of nonlinear equations to then prove
(\ref{agrev63}) and (\ref{agrev64}).

{\bf Proof of (\ref{grocoll1}) and (\ref{agrev62})}
We recall that we proved in the previous part the convergence of
$N^{-1} \wt u^N$ pathwise if we consider path for
$t \geq t_0$ and if we have convergence in law at $t=t_0$.
Therefore it suffices now to show that
 $N^{-1} \wt u^N (t_0)$ converges (in law) as $N \to \infty$.

To do this  we
first note  that the time-scaled and normalized (i.e. scaled by $N^{-1}$)
number of sites occupied at time
$\alpha^{-1} \log N+t_0$, denoted $N^{-1} \wt u^N$, does not go to zero as $N\to\infty$
nor does it become unbounded. This
was proved in Proposition \ref{P6.1b}. This proves tightness
of the one-dimensional marginal distributions
at times $\alpha^{-1} \log N + t$ and
that limit points are non-degenerate.

The tightness  gives us the
existence of limiting laws for $(N^{-1} \wt u^{N}, \wt U^N)$ along subsequences for time
$\alpha^{-1} \log N+s $ for fixed $s$. Since $N^{-1} \wt u^{N}\leq We^{\alpha t}$ any limit point must satisfy
\be{add21}
\limsup_{t\to -\infty} e^{-\alpha t}u(t)<\infty\text{   and   } \|U(t)\|\text{  bounded}.\ee
Therefore by Lemma \ref{L.UCE}(d) any limit point corresponds to a solution to the nonlinear
system given by (\ref{S081}) and (\ref{S082}) and any two such solutions are (random)
{\em time shifts} of each other. Hence we have to show that each pair of limit prints
belongs to the same time shift from the standard solution of the system with growth
constant 1.

Recall the definitions of $\tau^N(\ve)$ (which  is the first time where $u^{N,k,\ell}_s$
reaches $\lfloor \ve N\rfloor$ for all $\ve \in (0,1)$), $\wt\tau^N(\ve)$ from (\ref{bw3}), (\ref{w4}).
 Now assume that there are two limit points and let
the corresponding limiting passage times be denoted $\tau^{1}(\ve)$ and $\tau^{2}(\ve)$ such that $\tau^{N_{n_i},i}(\ve)\to \tau^{i}(\ve),\; i=1,2$
in law. But by (\ref{Grev39c3}),
$P[|\tau^{N_{n_2},2}(\ve)- \tau^{N_{n_1},1}(\ve)|> 2C(\ve)]\to 0$.  But since $C(\ve)\to 0$ as $\ve \to 0$, the time shift must be arbitrarily small,
which means equal to zero and therefore the two limit points coincide. This immediately proves that $\tau^N(\ve)$ converges for any $\ve>0$.

We now elaborate on this argument.
To show convergence of the marginal distributions it suffices to prove  the
convergence in law of the {\em normalized} first passage times
$\tau_{\rm{norm}}^N(\ve):=\tau^N(\ve)-(\alpha^{-1} \log N )$.
We give two arguments which verify  that convergence of the $\tau_{\rm{norm}}^N(\ve)$  is sufficient.

(i) Given the value $u(s)$ we know
exactly the first passage time for all $\ve > u(s)$
(due to the {\em unique} and also {\em deterministic} evolution  from any initial point).
Indeed using the fact that the path converges (for convergent sequences at a
fixed time), the collection is monotone non-decreasing and cadlag,
we know that the whole collection (in $\ve$) of first passage time
converges in path space. The uniqueness of the limit
points follows since
two different limiting distributions at
times $\alpha^{-1} \log N+s $ for fixed $s$ would not be
compatible with the uniqueness in distribution of the limiting
first passage times, since the limiting dynamics $(u(t))_{t \in \R}$ is
monotone for sufficiently small (i.e. negative $t$) so that it is
uniquely determined by the first passage times in that range and
hence everywhere.

(ii)
Alternatively we can use the fact that there is according to
Lemma \ref{L.UCE} up to a time shift a unique entrance law from
0 at $t \to -\infty$ which is compatible with the growth conditions
of a limit of that quantity as $N \to \infty$.
Therefore by one first passage time this random time shift
is determined and hence it suffices to prove convergence in
law of one first passage time concluding the second argument.

Having proved the reduction of the convergence statement to
one on $\tau^N_{\mathrm{norm}}$,
we have to show next that for one $\ve \in (0,1)$
and therefore by the above for every such $\ve$:
\be{ba10h1} \CL[\tau^N(\ve) - \frac{\log N}{\alpha}] \Nto
\CL[\tau(\ve)], \ee
with
\be{ba10h2} \tau(\ve) = \inf(t \in \R
|u(t)=\ve),
\ee
where $u$ is the (unique) entrance law with $e^{\alpha |t|}
u(t) \ttooo W$ (recall Lemma \ref{L.UCE}).

We have two tasks, to
identify the limit on the r.h.s. of (\ref{ba10h1}) as the quantity defined in
(\ref{ba10h2}) this however we showed in part 1 where we established in particular
the convergence of the dynamic for the frequency of occupied sites, and secondly we have to
show the convergence in (\ref{ba10h1}). So only the latter remains.

We use our results on the corresponding times $\wt \tau^N(\ve)$, assuming no collisions
occur, namely we see that it follows from Proposition \ref{P.C-10} that
$\wt \tau^N(\ve)$ converges to a limit $\wt \tau(\ve)$
defined in terms of the {\em collision-free} process.
Hence it remains to bridge from
$\tau^N(\ve)$ to $\wt \tau^N(\ve)$.

To relate $\wt \tau^N(\ve)$ and $\tau^N(\ve)$
we will use that we can use according to the above arguments arbitrarily
small $\ve >0$ and we can work with the multicolour particle
system which we introduced in the proof of Proposition \ref{P6.1b} below (\ref{Grev38}).

Consider $\ve \downarrow 0$.
We have the fact from the multicolour particles system
(recall (\ref{deadsites}) and (\ref{angr17o})) that
for sufficiently small $\ve$ asymptotically as $N\to \infty$
in probability
\be{grevy}
\wt\tau (\ve) \leq \tau^N(\ve) \leq \tilde \tau (\wt \ve),
\ee
where $\wt \ve$ solves
\be{grevy2}
\wt \ve-\mbox{ Const } \wt \ve^2 = \ve, \quad \wt \ve>0.
\ee
Note that the Const is less than 1, so that there is a unique
 $\wt \ve>\ve$ such that
\be{grevy3}
\wt \ve - \ve = O(\wt\ve^2).
\ee
In (\ref{grevy}) the right inequality holds for $\ve$ small enough and $N$ {\em sufficiently
large}, while the l.h.s. holds even for every $N$.

Denote by $\tau^N_{\mathrm{norm}}(\ve)$ the normalised $\tau^N(\ve)$,
i.e. $\tau^N(\ve) - \alpha^{-1} \log N$.
We know that $\wt \tau^N(\ve) \la \wt \tau(\ve)$ as $N\to\infty$
if we use the coupled collection of processes we constructed
in the proof of equation (\ref{grocoll1}). Hence for $\ve_0$ sufficiently small
realizations of two different limit points
$\tau^{\infty, 1}_{\mathrm{norm}}(\ve_0), \tau^{\infty, 2}_{\mathrm{norm}}(\ve_0)$
satisfy on a suitable probability space:
\be{aga1}
\wt \tau (\ve) \leq
\tau^{\infty, 1}_{\mathrm{norm}}(\ve), \tau^{\infty, 2}_{\mathrm{norm}}(\ve)
\leq \wt \tau (\wt \ve).
\ee
Therefore $|\tau^{\infty, 1}_{\mathrm{norm}}(\ve)- \tau^{\infty, 2}_{\mathrm{norm}}(\ve)|\leq |\wt \tau (\ve)-\wt \tau (\wt\ve)|$.   We furthermore know that
$\wt \tau(\ve)$ is given by $(\log \ve - \log W)/\alpha$ so that the l.h.s.
and the r.h.s. of (\ref{aga1})  differ by a {\em deterministic} quantity, which equals
$const \cdot \wt \ve_0$ as $\wt \ve_0 \downarrow 0$.  As $N\to\infty$ we
have a limiting evolution between $\tau^N(\wt \ve_0)$ and $\tau^N(\wt \ve)$
which is deterministic. This evolution is a {\em random shift} of
some deterministic curve.
Then it is not possible to have two limit points for the law
$\tau^N(\wt \ve)$, since their difference are shifted versions of a given
curve would translate into the same
difference but now bounded above by arbitrarily small $\wt \ve_0$ and therefore equal to  0.
Hence all limit points are equal and
$\tau^N_{\mathrm{norm}}(\ve)
= (\tau^N(\ve) - \alpha^{-1} \log N)$ must converge as $N \to \infty$
as claimed and we are done. This concludes the argument for the convergence
of the marginals at time $\alpha^{-1} \log N+t_0$.
\sm

{\bf Proof of (\ref{ba10i2})}
We first need some preparation, involving the collision-free regime.
Turn to the dual particle process and recall that at
time $\alpha^{-1} \log\log N$ asymptotically as $N \to \infty$ there are $W \log N$
occupied sites and that at time $\alpha^{-1} \log\log N$ the empirical age
and size distribution is given by the CMJ stable age and size distribution
$\CU(\infty, \cdot, \cdot)$,
since we are in the non-collision regime due to $\log\log N << \alpha^{-1} \log N$.
Therefore at time $\tau_{\log N}$ the age
distribution is (asymptotically as $N \to \infty$)
given by the stable age distribution.
Moreover for each realization of $W$ there is a time shift
of the limiting dynamics which is a function of $W$.
Since the
evolution from this point $\tau_{\log N}$ on is asymptotically deterministic we
conclude that the age and size distribution $U^N$ at time
$\alpha^{-1} \log N -t_N +t$,
with $t_N \uparrow \infty$ but $t_N << \alpha^{-1} \log N$, converges to a time shift of
the one determined by equation (\ref{grocoll2}) and (\ref{ba10g}) with initial
condition given by $u=0$ and $U$ given by the stable age
distribution. We now have to bridge from times $\alpha^{-1} \log N -t_N +t$
to $\alpha^{-1} \log N +t$, where collisions play a role. This runs as follows.
\sm

For a given approximation parameter
$\ve > 0$ we choose a time horizon $\alpha^{-1} \log
N+T(\ve)$ such that $(K^N_t, \zeta^N_t)_{t \geq 0}$ is, as $N \to \infty, \ve$-approximated by the
collision-free process $(K_t, \zeta_t)_{t \geq 0}$ up to time
$\alpha^{-1} \log N +T(\ve)$ if we choose $\ve$ sufficiently small
and this implies the approximation
for the functionals $U^N$ by $\CU$ as $N \to \infty$
for $\ve$ sufficiently small.
More precisely we require that we can find a coupling between
$U^N$ and $\CU$ such that  ($\|\cdot \|$ denotes variational distance)
\be{Grev39o}
\limsup_{N\to\infty}P\left(
\sup_{t\leq \frac{\log N}{\alpha} +T(\ve)}\|U^N(t,\cdot,\cdot) -
\CU(t,\cdot,\cdot)\| > \ve \right)\to 0 \text{ as } \ve\to 0. \ee
In other words we have to show that by the collisions up to time
$\alpha^{-1} \log N +T (\ve)$ at most difference $\ve$ develops between
$U^N$ and $\CU$ for very large $N$.

To show this  note that we can bound the total effect of
collisions occurring up to time $\alpha^{-1} \log N +T(\ve)$.
Namely we return to the multicolour particle system defined below
(\ref{Grev38}), which shows that we can relate the difference of
the system without and with collision through the black and red
particles. In order to show that $\|U^N(t,\cdot,\cdot) -
\CU(t,\cdot,\cdot)\|\leq \ve$ it suffices to show that the
proportion or red and black particles is $o(\ve)$ in probability
as $N\to\infty$. But this was verified for the number of black
sites with $T(\ve)=\wt \tau(\ve)$ in (\ref{deadsites}), which
immediately by using that $\zeta$ is stochastically bounded gives
the claim for the number of particles.  Since
the number of descendants of red particles grows at a slower rate
as the black and new red founding particles occur at the same rate
as black ones then the claim follows. This completes the proof of Part (c).
\sm

\begin{proof}{\bf of (d): (\ref{agrev63}) and (\ref{agrev64})}

Again we first turn to the collision-free process. We
approximate (recall $\CU$ refers to the collision-free CMJ-process)
$\CU(t, \cdot, \cdot)$ by
$\CU (\infty, \cdot, \cdot)$.
By \cite{N} (Theorem 6.3) we know that $\CU(r,\cdot,\cdot)$
converges a.s. as $r \to +\infty$ to the stable age and size
distribution $\CU(\infty, \cdot, \cdot)$ of the Crump-Mode-Jagers
process (see (\ref{ang4c})) and hence
\be{Grev39p}
\CU(\alpha^{-1}
\log N + t,\cdot,\cdot) \Ntoo \CU(\infty,\cdot,\cdot),\;a.s.. \ee

Now take collisions into account and consider the limit of  $U^N(\frac{\log
N}{\alpha}+t,\cdot,\cdot)$ as  $N\to\infty,\; t \to -\infty$.
Combining (\ref{Grev39o}) and (\ref{Grev39p}) it follows that for
$\eta>0$ there exists $\tau(\eta)$ such that
\be{Grev39p1}\limsup_{N\to\infty}P\left(\sup_{t\leq \tau(\eta)}\|U^N(\frac{\log N}{\alpha}+t),
\cdot,\cdot) - \CU(\infty,\cdot,\cdot)\|>2\eta \right) <\eta . \ee
Therefore \be{Grev39p1+}\lim_{t\to -\infty}
U(t,\cdot,\cdot)=\mathcal{U}(\infty,\cdot,\cdot).\ee

Finally we have to verify that indeed
$(e^{-\alpha t} \frac{\wt u^N(t)}{N})_{t \in \R}$
converges to a limit $(e^{-\alpha t} u(t))_{t \in \R}$
such that $e^{-\alpha t} u(t) \to W^{k,\ell}$ as
$t \to -\infty$, where $ W^{k,\ell}$ arises as
the growth constant of the collision-free process.
This was already essentially done in the proof of (\ref{agrev61}),
if we observe that the difference between the collision-free
process and $\wt u^N(t)$ is bounded by the number of black
sites, whose intensity vanishes as $t \to -\infty$, since
\be{angr41}
\limsupl_{N \to \infty} P[\tau^N (\ve) \leq \alpha^{-1} \log N+t]
\to 0 \mbox{ as } t \to -\infty
\ee
and hence the intensity of black sites is $O(\ve^2)$
at time $\tau^N(\ve)$ for $N$ large.

This completes the proof of part (d).

\end{proof}
\bi

{\bf Part 3: Occupation density of dual particle system}

\emph{Proof of (e)}

Note first  that (\ref{agr24A}) follows from what
we proved for the collision-free regime, in particular the fact that
$\sup (e^{-\alpha t} K_t| t>0) <\infty$.

{\bf (\ref{ba7})} Next
 note that it suffices to prove the
convergence in (\ref{angr24}) to obtain (\ref{ba7}) using (\ref{agr24A}). Namely we then
have  pointwise convergence on a suitable probability space and an
upper bound is given by the collision-free process with single
sites in equilibrium, which provides an integrable integrand and
is independent of $N$ and indeed using Proposition
\ref{P.Grocoll}, part (b) we see that (\ref{ba7})
holds.

{\bf (\ref{angr24})}
Using Proposition \ref{P.Grocoll} part (c) we see that for (\ref{angr24})
we have to reduce the problem
concerning $u^N U^N$ and $uU$ to what we showed there for $u^N$
and $u$ respectively $U^N$ and $U$. Namely each of them converges weakly.
However since
\be{ba10h2c} (u(t), U(t,\cdot)) \la u(t)U(t,\cdot)
\ee
is a continuous mapping this follows from the convergence properties we proved for
$(N^{-1} \wt u^N, \wt U^N)$ as $N \to \infty$.

This completes the proof of Proposition \ref{P.Grocoll}.

\begin{remark}
Consider the normalized first passage times $\bar \tau^N(\ve)$, where
$\bar \tau^N(\ve)$ is the first time where $\Pi^{N,k,\ell}_s$
(rather than $K^{N,k,\ell}_s$) reaches $\lfloor \ve N\rfloor$.

Then one can show that for every $\ve \in (0,1)$:
\be{ba10h1a}
\CL[\bar \tau^N(\ve) - \frac{\log N}{\alpha}] \Nto
\CL[\bar \tau(\ve)],
\ee
with
\be{ba10h2b} \bar \tau(\ve) = \inf(t
\in \R |u(t) U(t,\R,\N) =\ve), \ee and $u$ is the entrance law
with $\liml_{t \to \infty} e^{\alpha |t|} u(t) = W$.
\end{remark}

\subsubsection{Dual population in the transition regime: asymptotic expansion}
\label{sss.prep2}

The purpose of this subsubsection is to use and extend the techniques
we developed above to later be able to connect the early times
(studied later in Subsection \ref{ss.bounds}) and the late times studied
in Subsubsection \ref{sss.asycolreg} above.

Return first to our original population model. On the side of the original
process the limit in the first time scale, i.e. of times $O(1)$ as $N \to \infty$, is described by the branching process
$(\gimel^m_t)_{t \geq 0}$ and in the second time range, i.e. after times
$\alpha^{-1} \log N$ by the deterministic
McKean-Vlasov equation but with random initial condition. For both
these limits we have derived (respectively for $\gimel^m_t$ we will in
Subsection \ref{ss.bounds}) exponential growth at the same rate $\alpha$
and for the first we had obtained the growth constant $\CW^\ast$
and for the latter the growth constant $^\ast \CW$. The objective
below is to show that we can consider the limit simultaneously in
both time scales in such a way that the final random value arising
from the branching process describing the number of type-2 sites
at small times provides the initial value of the McKean-Vlasov
equation describing the evolution from emergence on till fixation meaning that
in fact $^\ast\CW$ and $\CW^\ast$ have the same distribution.

This section provides the analysis of the dual needed for this
purpose.
The completion of the argument is
obtained in the following subsubsection which is devoted to the
limiting branching process and the transfer back of properties of
the dual to the original process.

Recall the definition of $u^N,\;U^N$ in  (\ref{Grev39c2a})-
(\ref{Grev39e}).  In this section we deal with the relation of the
two basic limits as $N\to\infty$ of the functional of the dual population
given by
$(u^N(t),U^N(t))$, namely, the limit in the two different time scales
\be{agrev90}
t_N= o(T_N) \mbox{ and }\;\;
t_N=T_N+t,\;\; \text{  where  }\;\; T_N:=\frac{\log N}{\alpha}
\mbox{ and } t \in \R,
\ee
as well as the {\em transition} between these two
time scales which {\em separate} in the limit $N \to \infty$.

The main tool developed in this subsection is the analysis of the
dual population in this two-time scale context. For that purpose
we make use of two techniques: (1) {\em Nonlinear evolution equations}
for the dual particle systems and their enrichments, (2) {\em couplings}
of particle systems to estimate and control the effects of
collisions in the dual process by comparing it with collision-free systems,
i.e. the dual of the McKean-Vlasov process.

Heuristically we can approximate
 $u^N/N$ by  $\wt u^N$ which are the solutions of the ODE

\be{transu3} \frac{d \wt u^N(t)}{dt}=\alpha_N(t) \wt
u^N(t)(1-\frac{b_N(t)}{N} \wt u^N(t)),\;t\geq t_0(N), \ee
\be{grevz2}
b_N(t)=1+\frac{\gamma_N(t)}{\alpha_N(t)},\quad \wt
u^N(t_0(N)) = u^N(t_0(N))= W_Ne^{\alpha t_0(N)},
\ee
which allows
to get some feeling for the behaviour based on explicit calculation
for the case where $\alpha_N$ and $\gamma_N$ are constant and
explicit formulas for the solution of the ODE (\ref{transu3}) and (\ref{grevz2})
can be given.

There are two regimes depending on which of the time scales in
(\ref{agrev90}) is used: the first (linear) regime where $N^{-1} u^N(t_N)\to 0$
as $N\to\infty$ together with $\alpha_N(t_N)\to \alpha$ and
$\gamma_N(t_N)\to \gamma$ and the second (nonlinear) regime where
$u^N(t_N)=O(N)$ and where $\alpha_N(t_N), \gamma_N(t_N)$ differ
from $\alpha, \gamma$.

We have proved in Subsubsection \ref{sss.structure} that for times
$t_N\to -\infty$ the number of sites occupied by the dual process
at times $\frac{\log N}{\alpha} +t_N$ and multiplied by $N^{-1}$
converges to zero in probability. On the other hand we saw in
Subsubsection \ref{sss.growth}  that assuming that in law
\be{gra1}
\frac{1}{N} u^N(\frac{\log N}{\alpha}+t_0)\Nto u(t_0),\ee
then in law on path space
\be{gra2}
\{\frac{1}{N} u^N(\frac{\log N}{\alpha}+t)\}_{t\geq t_0}\Nto \{u(t)\}_{t\geq t_0}
\ee
 and given the state at time $t$ the limiting dual dynamics $u(t)$  is deterministic and nonlinear,
namely, $(u(t),U(t))_{t \geq t_0}$ is for all $t_0 \in \R$ the solution to the system
(\ref{grocoll2}), (\ref{ba10g}) where in the latter quantity we
integrate out the age. Moreover we proved that this system has a
unique solution satisfying
\be{gra3}
\lim_{t\to -\infty} e^{\alpha
|t|} u(t)=W,\quad  \lim_{t\to -\infty} U(t)= \CU(\infty), \ee
where $\CU(\infty)$ is the stable size distribution of the CMJ
process induced by the collection of the occupied sites of the
McKean-Vlasov dual process.

In order to later on relate $^\ast\CW$ and $\CW^\ast$ we next
focus on $^\ast\mathcal{W}$
and try to get more information on its law. This requires,
as we shall see later, when we return to the original process
from the dual to obtain  in
(\ref{gra3}) the higher order terms as $t \to -\infty$ for the
limiting equation (as $N \to \infty$) as well as in the
approximation as $N \to \infty$ of this behaviour. This means we
want to write for $t \to -\infty$ the limiting (as $N \to \infty$) intensity of
$N^{-1} \wt u^N$, denoted $u$:
\be{gra3b}
u(t) = W e^{\alpha t} - Const \cdot
W^2e^{2 \alpha t} + o(e^{2 \alpha t}) \ee and determine the
constant, but moreover we want to consider for $N \to \infty, t \to -\infty$
the quantity
\be{gra3c} u^N (\frac{
\log N}{\alpha}+t) = W_N(t) e^{\alpha t} \cdot N + C (t,N) \ee and
to estimate the order in both $N$ and $t$ of the correction term $C(t,N)$.
The latter is equivalent
to determining up to which order the expansion in (\ref{gra3b}) is approximated by the finite
$N$-system.
First we turn to the question in (\ref{gra3c}) and then to (\ref{gra3b})
for each point of view formulating a separate proposition.

We therefore refine (\ref{gra3})  by providing speed of convergence,
 in $t$ and uniformity, in $N$, results.

\beP{strongconn}  {(Order of approximation of entrance law in $N$ and  $t$)}

(a)
Let $T_N=\frac{\log N}{\alpha}$ and $0< \delta
<1$. Then $(u^N(t), U^N(t))$ and $(u(t), U(t))$ have the following three
properties:

We can couple the $(u^N (t),U^N(t))_{t \geq 0}$ with a fixed CMJ
process $(K_t, \CU(t))_{t \geq 0}$ such that
\be{agrev91}
e^{-\alpha t} K_t \la W \mbox{ a.s. as } t \to \infty,
\ee
\be{gra4}
\lim_{t \to -\infty} \limsup_{N\to\infty}P\left(|e^{-\alpha
t}\frac{1}{N} u^N(T_N+t)- W|>e^{\delta \alpha t}\right)=0 \ee and
\be{gra4b}
\lim_{t \to -\infty} \lim_{N \to \infty} \left(
\frac{e^{-\alpha t}}{N} u^N({T_N+t}), U^N (T_N+t)\right) = (W,
\mathcal{U}(\infty)). \ee

(b)  Given $W$ the system (\ref{grocoll2}), (\ref{ba10g}) has
unique solution $(u(t),U(t))$ with $e^{-\alpha t} u(t) \to W$ as
$t \to - \infty$ which can be realized together with $\{(u^N, U^N), N \in \N\}$
on one probability space such that the process
\be{gra4c}
(e^{-\alpha t} N^{-1} u^N (T_N+t))_{t \in \R} \ee
can for every $\delta \in (0,1)$  be approximated by $(e^{-\alpha
t}u(t))_{t \in \R}$ with ``initial condition'' $W$, up to order
$e^{\alpha \delta t}$ for $t \in \R$ uniformly in $N$.$\qquad
\square$
\end{proposition}

We now have established a precise relation between $(u^N, U^N)$
and the limiting entrance law $(u,U)$.  We now  return to the
solution $(u,U)$ to (\ref{grocoll2}), (\ref{ba10g}) and
$\alpha(t),\;\gamma(t)$ as defined in
 equation (\ref{ba10d}).
 The purpose of the rest of this section is to identify the
asymptotic behaviour of this nonlinear system as $t\to -\infty$ up
to higher order terms, i.e. we want to identify in particular the
correction of order $e^{2\alpha t}$ as $t \to -\infty$ which is
due to the occurrence of collisions by migration of the dual individuals
which by subsequent coalescence of the collided individuals can
change the behaviour of the limiting system.

\beP{P.Transition}{(Transition between linear and nonlinear regime
limit dynamics $(u,U)$)}

 (a) The pair $(u,U)$ and the functionals $\alpha$ and $\gamma$
 arising as the limit in (\ref{grocoll1})-(\ref{ba10i2}) via
 (\ref{ba10d}) satisfy (here we suppress the $k,\ell$):

 \be{ga30++} (\alpha(t)+\gamma(t))u(t)\leq
(\alpha + \gamma) We^{\alpha t}, \ee

 \be{grocoll214}
 \lim_{t\to -\infty}\alpha(t)=\alpha,\quad \lim_{t\to
-\infty}U(t)=\mathcal{U}(\infty),\quad   \lim_{t\to
-\infty}\gamma(t)=\gamma,
\ee

\be{grocoll3} u(t) = W e^{\alpha t}- \kappa W^2e^{2\alpha
t}+O(e^{3\alpha t}) \quad \mbox{ as } t \to -\infty, \mbox{ for some constant } \kappa >0
\ee
and furthermore as $t\to -\infty$ the functions $\alpha(\cdot)$ and $U(\cdot)$ satisfy
that
\be{gra4d}
\alpha(t)=\alpha +\wt\alpha_0e^{\alpha
t}+O(e^{2\alpha t}),\quad U(t)=\mathcal{U}(\infty)+\wt U_0e^{\alpha
t}+O(e^{2\alpha t}),
\ee
\be{gra4e1}
\gamma(t)=\gamma + \wt \gamma_0 e^{\alpha t} + O(e^{2 \alpha t}),
\ee
where $\wt \alpha_0$ is a positive number $\wt U_0$ is an $\N$-vector and  $\wt\gamma_0= \wt U_0(1)$.

Using a
multicolor system construction, we will see below that all the randomness
is given by the (non-degenerate positive)  random variable $W$, namely
\be{gra4e}
\wt\alpha_0 \mbox{ is explicitly given by (\ref{gra39}), as } Const_1 \cdot W,
\ee

\be{gra4e2}
\wt U_0 \mbox{ is explicitly given by (\ref{U0zz}) } \mbox{ and has the form } W \cdot \vec{const}, \wt \gamma_0 = Const_2 \cdot W.
\ee

(b) The total number of particles satisfies
\be{grocoll32}
(\alpha(t)+\gamma(t))u (t) = (\alpha +\gamma) W e^{\alpha t}-
\kappa^\ast W^2e^{2\alpha t}+O(e^{3\alpha t}) \quad \mbox{ as } t \to
-\infty,
\ee
for some constant $\kappa^\ast >0$.

More precisely with
\be{gra4e3}
\kappa^\ast W^2 = (\wt \alpha_0 + \wt \gamma_0) W + \kappa W^2 = W^2 (Const_1 + Const_2 +\kappa),
\ee
we have that
\be{xg4}
\limsup_{t\to
-\infty}|(\alpha(t)+\gamma(t))u (t) - (\alpha +\gamma) W e^{\alpha
t}- \kappa^\ast W^2e^{2\alpha t}|\cdot e^{-3\alpha t}<\infty. \ee

(c) Furthermore $U(t,\cdot,\cdot)$ is uniformly continuous at $t
=-\infty$ and as $t \to -\infty$

\be{gra6}
  \| U(t)- \CU(\infty)-\wt U_0 e^{\alpha t}\|_{{\rm var}} = O(e^{2\alpha t}).
\hfill \square \ee
\end{proposition}

\begin{remark}
Note that the results above show that all the relevant randomness
in the dual process sits indeed in the random variable $W$.
\end{remark}

\bigskip

The proof proceeds in four steps.
 The proof of Propositions \ref{strongconn} and
\ref{P.Transition} are based on an enriched version of the dual
particle system $\eta$  and a reformulation of the nonlinear
system (\ref{grocoll2}) and (\ref{ba10g}) which are given in Steps
1 and 2 below which are then followed by the Steps 3 and 4 giving the proof of the two
propositions using these tools. In Step 1 we focus on the finite
$N$ system in a time regime where collisions due to migration become
essential, while the Step 2 develops the analysis of the limiting
system as $N \to \infty$ in this time regime. The main focus here
is on the behaviour as the effect of collisions becomes small and
to investigate its precise asymptotics in this regime.

\textbf{Step 1:  Modified and enriched coloured
particle system}

Our goal is to realize on one probability space the CMJ-process
and the dual particle system in such a way that we can identify in
a simple way the difference between the two dynamics with a higher degree
of accuracy than before by using three colours white, black and red. In particular  we can
identify the part of the population involved in one or in more
than one collision (more than two etc.). This joint probability space
can be explicitly constructed easily, based on collections of Poisson
processes. Since no measure theoretic subtleties occur, we do not
write out this construction in all its lengthy detail and just
spell out the evolution rules.

In order to identify the
contributions  at time $\alpha^{-1} \log N+t$ of order $e^{\alpha t}, e^{2 \alpha t}, e^{3 \alpha
t}$, we want to expand the number of occupied sites at time
$\alpha^{-1} \log N+t$, using quantities $W_N(t)$ defined as $e^{-\alpha t} K^N_t$
and then analog higher-order objects, in the form:
\be{xg5}
u^N(T_N+t)\sim W_N(t) e^{\alpha t} - W^\prime_N(t) e^{2 \alpha t}
+ W^{\prime \prime}_N(t) e^{3 \alpha t},
\ee
with $T_N= \alpha^{-1} \log N$ and $N^{-1} W_N(t), N^{-2} W^\prime_N(t),
N^{-3} W^{\prime \prime}_N(t)$
converging as $N \to \infty, t \to -\infty$, more precisely,

\be{gra6AG}
\begin{array}{l}
\lim_{t\to -\infty}\lim_{N\to\infty}[e^{-\alpha t}N^{-1} u_N(T_N+t)] = W \\
\lim_{t\to -\infty}\lim_{N\to\infty} [e^{-2\alpha t}N^{-2} [(W N
e^{\alpha t}-u^N(T_N+t))]i]=   \kappa W^2\quad , \\
\lim_{t\to -\infty}\lim_{N\to\infty} [e^{3\alpha t} N^{-3}
[u_N(T_N+t))-WN e^{\alpha t}- \kappa N^2 W^2 e^{2 \alpha t}]]  < \infty,
\end{array}\ee
where the limits are needed in probability and in $L^1$. The
latter convergence follows from the convergence in probability and the
finiteness of the moments of $W$ (see Lemma \ref{L6.11}).

We use the dual particle system to obtain information on
the asymptotics as $t\to -\infty$ of the pair $(u,U)$  and also to
identify their random initial condition.  However to carry this
out we must  enrich the {\em dual particle system} to a {\em multicolour}
particle system on an {\em enriched geographic space}.

Since we want to read off more properties than
in the emergence argument based on the mean where we had white,
black and red particles we shall need now some more colours
which allow us to record particles which have been involved
in one, two, etc. collisions and with associated colours their
counterpart in a collision-free system.

Furthermore in order to realize the CMJ and the dual on one
probability space we have to couple the evolution of certain
colours where one belongs to the CMJ-part and the other to the
dual particle system part. This will also require us to  refine
the geographic space $\{1,\cdots,N\}+\N$ (where the black
particles were collision-free and moving on $\N$) we had before
for the white, red, black particles system by adding a further
copy $\N$ (or more if we need more accuracy).

More precisely we now introduce a {\em multicolour particle
system} that is obtained by modifying and enriching the coloured
particle system defined by (\ref{Grev38a}) - (\ref{bw2}) by
introducing (1) new colours and (2) new type of sites. The number of new colours
and sites necessary depends on how precisely we want to control the behaviour as $t \to -\infty$
meaning up to which higher order terms as $t \to -\infty$ we want to go.

We first need a new colour (green) which
allows us to control exactly (i.e. not just estimating from above)
the difference between the dual particle
system $\eta$ and our old white-black collision-free particle
system. Recall that in the old system we created a red-black pair
upon a collision such that the further evolution of black remains
collision-free and the red follows the true mechanism of the dual
particle system. A real difference in the total number of particles between dual and CMJ will however occur
only if the red particle coalesces with a white particle but no such loss of an
individual occurs in the CMJ-process which now will have one particle more
due to this transition. To mark this
we will use the colour {\em green} to mark the dissappearing red particle in the dual.

To allow a fine comparison between the white-red
and the white-black system we will use coupling techniques which
allow us to estimate not only the effect of {\em collisions} (and subsequent coalescents) but also
{\em recollisions}, the latter will generate parallel to the red-black construction
two new colours {\em purple} (in the dual system)
and {\em blue} (in the collision-free comparison system).
For this purpose we have to modify the evolution rules when
red particles collide or red particles collide with white
particles. For the bookkeeping of the effects from these events we
use the new further colours.

The white particles in the
old process also have the same dynamics as the white particles
in our new formulation but in  addition to white, red, and black
particles the {\em green, purple} and {\em blue} particles
 can have now locations in
\be{xg6}
\{1,2,\cdots,N\}+\N +\N,
\ee
which we will explain as we describe the evolution rules.

In order to achieve all these goals of a higher order expansion
the modified dynamics should satisfy the following requirements:

\begin{itemize}

\item for each $N$ the union of black plus white particles is
equivalent to the particle system without collisions and can be
defined on a common probability space by a single CMJ process,

\item the white plus red plus purple particles give a version of
$(u^N(t),U^N(t,\cdot,\cdot))$,

\item the red plus purple particles give a refinement of the red particle system
in the white-black-red process used
earlier, i.e. {\em both together} have the same dynamics as the
set of red particles in  (\ref{Grev38b4}),

\item the green particles describe the loss due to coalescence after collision,
i.e. after a particle (white, red or purple) migrated to an occupied site on
$\{1, \cdots, N\}$ (by a white, red or purple particles) and then coalesced with
such particle at this site,

\item the black particles are placed  and evolve on the second copy of $\N$
and blue particles are placed and evolve on the first copy of $\N$.

\item We can match certain pairs of red and black respectively purple
and blue particles when they are created to better compare the
CMJ-part and the one with the dual dynamics.

\item The number of the sites occupied by red, green and blue
particles has the same law as the number of sites occupied by
black particles.

\end{itemize}

For each $N$ the evolution rules below  define a {\em Markov pure
jump process} which describes a growing population of individuals carrying
a {\em location} and a {\em colour} and in addition  the information which pairs of
individuals are {\em coupled}.

Let $I_n$ be a copy of $\{1, \cdots, n\}$.  Denote by
\be{agdd65}
C = \{\mbox{ white, black, red, purple, blue, green}\},
\ee
Then the state space is given by the union over
$n$ of the set of maps
\be{agr20}
I_n \la C \times (\{1,\cdots, N\} +\N +\N\}) \times (I^2_n).
\ee
If we just count the number particles with a certain colour and location
we get again a Markov pure jump process, specifying
 a multitype particle system on
 $\{1,2,\dots,N\}+\mathbb{N}+\mathbb{N}$.
The state space of this  system is given by:
\be{xg7}
(\N_0)^{(\{\{1,2,\dots,N\}+\mathbb{N}+\mathbb{N}\} \times C)}. \ee
We shall work with both processes.

In order to satisfy the properties specified above, the dynamics
of the enriched coloured particle systems has the following dynamics:

\begin{itemize}
\item Given that $k$ is the number of sites in $\{1,2,\dots,N\}$
which are occupied, and a white particle migrates to another site
in $\{1, \cdots, N\}$ the outcome is given as in (\ref{Grev38b}),
(\ref{Grev38b2}) and (\ref{Grev38b3}) (a {\em red-black pair} is
created) but we record that these particles form a red-black pair.
White particles give birth to white particles and two
white particles coalesce to produce one white particle.

\item Black particles evolve as before but live on the second copy
of $\N$.

\item Red particles give birth to red particles, two red particles
at the same site coalesce to one {\em red particle}.

\item When a red particle coalesces with a white particle (at the
same site) the outcome is a white particle and a {\em green} particle at
this site. Moreover we then couple the green particle with the
black particle that had been associated with the involved red
particle.

\item Given that the current number of sites occupied by white,
red, purple  or green particles on $\{1, \cdots, N\}$
is $\ell$, when a red particle migrates it moves to
an empty site in $\{1,\dots,N \}$ with probability
$1-\frac{\ell}{N}$ and  with probability $\frac{\ell}{N}$ it dies and produces
a  {\em blue-purple pair} where the blue particle is placed at the first
empty site in in the first copy of $\mathbb{N}$ and the purple
particle is placed at a randomly chosen occupied site on
$\{1,2,\cdots,N\}$. We also associate the blue particle with the
black particle that had been associated with the involved red
particle.

\item Blue particles have the same dynamics as black particles.

\item Purple particles have the same dynamics as the dual particle
systems $\eta$ on $\{1,\dots,N\}$, in particular they migrate to a randomly
chosen site in $\{1,\dots,N\}$.

\item Green particles give birth to green particles and two green particles at one site coalesce
 giving one green particle. Green particles on $\{1,2,\cdots,N\}$ migrate
to empty sites on $\{1,\cdots,N\}$ with probability $1-\frac{\ell}{N}$
and with probability $\frac{\ell}{N}$ to the first empty site in the first
copy of $\N$ where $\ell$ is the number of occupied sites on $\{1,\cdots,N\}$.
Green particles on the first copy of $\N$ migrate to the first empty site.

\item When a green and red particle coalesce they produce a just red particle.

\begin{remark} \label{rgremark} The reason for this rule is as follows.  The red and green pair correspond to
 two red particles one of which has coalesced with a white particle.  The green particle corresponds to a particle loss due to collision and coalescence.
First note that the green particle which coalesces with a red one must have been created
at that site. Hence both red and green must descend from the same collision of a particle
arriving at the site and hence the coupled black particles share the same site.
The red and green particle are coupled with a black or blue particle each.
The two black particles can coalesce only if they descend from the same black particle
after arriving at the current site. In this case the number of red and green, respectively
black and blue is reduced both by one and the original loss of black-white versus white-red
is cancelled.
 When the green and red coalesce this corresponds to the coalescence of two red particles that would have been at a site having no white particles.  Then the coalescence reduces the number of red particles to one.  This means that the potential loss is canceled out and in both cases we end up with a single red particle.
\end{remark}

Green and white particles do not coalesce.

\item We couple for a newly created red-black pair the birth,
coalescence with offspring up to the time of the first collision
of a red particle with a site occupied by  red or white particles.
We then continue by coupling for all future times the associated
black particle with the blue particle (from the blue-purple pair)
created at the collision time.

\item When red and purple particles coalesce they
produce a red particle.

\item When a purple particle coalesces with a white particle the
outcome is a white and a green particle at the site.

\end{itemize}
\noi Note that these rules lead to a system with the six desired
properties we had listed below (\ref{xg6}).

The state space of this  system is given by: \be{xg7a}
(\{\{1,2,\dots,N\}+\mathbb{N}+\mathbb{N}\} \times C)^{\N_0}. \ee

We shall single out specific subsystems comprised of particles of certain
colours, for example,  the {\em WRP-system} (white, red, purple)
or the {\em WRGB-system} (white, red, green, blue) which have
state-spaces:
\be{xg8}
\N_0^{(\{\{1,2,\dots,N\}+\mathbb{N}+\mathbb{N}\} \times \{W,R,P,\})},\quad
\N_0^{(\{\{1,2,\dots,N\}+\mathbb{N}+\mathbb{N}\} \times
\{W,R,G,B\})},
\ee
respectively. If we need more information, then we have the analogue of
(\ref{agr20}).

Observe  that white, red, and purple
particles are located only in $\{1,\dots,N\}$,  blue particles are
located only in the first copy of $\mathbb{N}$ and green particles
can be in either $\{1,\dots,N\}$ or first copy of $\mathbb{N}$ and
all black particles live in the second copy of $\mathbb{N}$.
It is important to note that the transitions of green particles
do not depend on whether they are located in $\{1,\cdots,N\}$
or in the first copy of $\N$.

We note that the difference between the number of red plus purple
particles and the number of black particles is due to the
coalescence of red or purple particles with white particles. Each
such loss is compensated by the creation of a green particle and
we note that green particles never migrate to an occupied site.

By the  construction above
we define on a common probability space five {\em coupled}
particle systems,
\begin{itemize}
\item
a CMJ-process which is the union of black and white particles which give
a version of $(K_t, U(t))_{t \geq 0}$,
\item
the white, red and purple particles generate
 a version of $(u^N(t), U^N(t))$,
\item a CMJ process  which is the union
of white, red, green and blue particles and which generates
a version of the pair $(K_t, U(t))_{t \geq 0}$,
\item
a subset of this above system, the white system generating the pair
$(K^{N}_{(W),t}, U^N_{(W)}(t, \cdot))$,
a CMJ-process which is the union of black and white particles which give
a version of $(K_t, U(t))_{t \geq 0}$,
\item
a {\em coupled} subsystem, on the one hand  the system of red and purple
and on the other hand a subset of the
black particles.  This coupling is induced by the
convention how purple-blue and red-black pairs evolve.
\end{itemize}

As before in the analysis of $(u,U)$ it suffices for our purposes
to study a functional of the multicolour particle system by observing
at sites only how many particles of the various colours occur and
counting the number of sites of a specific colour configuration. Our
functional of a state is given by the counting process \be{xg8b}
\Psi^N_t(i_W,i_R,i_P,i_G,i_B,i_{BL}),
 \ee denoting the number of sites containing $i_W$ white, $i_R$
 red, $i_P$ purple, $i_G$ green,  $i_B$ blue and $i_{BL}$ black
 particles at time $t$.

 The process $(\Psi^N_t)_{t \geq 0}$ evolves itself as a {\em pure jump-strong Markov
 process} with state space the counting measures on $C^\N$, the set of colour
 configurations, due to the fact that the dynamic only depends on the vector
 of the numbers of particles of the various colours at a site. Similarly we define
pure Markov jump processes
 \be{xg8c}
\Psi^N_{W,t} (i_W), \quad \Psi^N_{WRP,t}(i_W, i_R, i_P), \quad
\Psi^N_{WRGB,t} (i_W, i_R,i_G, i_B). \ee

Next we want to define the coloured versions of the $(u,U)$ system
and of certain subsystems.  We therefore define the number of sites
exhibiting certain colours:

\be{KTS}
K_{t}=\sum_{i_i, \cdots, i_{BL}} 1_{i_W+i_{BL}\geq
1}\Psi^N_t(i_W,i_R,i_P,i_G,i_B,i_{BL}), \ee

\be{KTD} K_{WRP,t}^{N}=\sum_{i_W,i_R,i_P}1_{i_W+i_R+i_P\geq
1}\Psi^N_{WRP,t}(i_W,i_R,i_P),\ee

\be{KTGB}
K_{WRGB,t}^{N}=\sum_{i_W,i_R,i_G,i_B}1_{i_W+i_R+i_G+i_B\geq
1}\Psi^N_{WRGB,t}(i_W,i_R,i_G,i_B)
\ee
and occasionally we make use of
\be{KTneu}
K_{C,t}^{N}=\sum_{i_L, \cdot L \in C}
1_{\suml_{L \in C} i_L \geq 1}
\Psi^N_{C,t}(i_C),\quad
C \subseteq \{W,R,P,G,B,BL\}, \quad i_C=(i_L)_{L \in C}.
\ee

Now we can define the coloured relatives of $(u^N,U^N)$
due to our construction all on one probability space, setting

\be{xg9}
(u(t) , u^N_{W}(t), u^N_{WRP}(t) , u^N_{WRGB}(t)) =
(K_{t}, K^N_{W,t}, K^N_{WRP,t}, K^N_{WRGB,t}),
\ee

\be{URP}\begin{array}{l}
U^N_{W,t} ( i_W) = \frac{\Psi^N_{W,t}(i_W)}{K^N_{W,t}},\\
U^N_{WRP,t}(i_W,i_R,i_P) =
\frac{\Psi^N_{WRP,t}(i_W,i_R,i_P)}{K^N_{WRP,t}},\\
U^N_{WRGB,t}( i_W, i_R, i_G, i_B)
 = \frac{\Psi^N_{WRGB,t} ( i_W, i_R, i_G, i_B)}{K^N_{WRGB,t}}.
\end{array}
\ee

\bigskip

\beL{L.enriched} {(Properties of enriched coloured particle
system)}

(a) Ignoring colour the process $(u^N_{WRP,t},U^N_{WRP} (t,\cdot))_{t \geq 0}$ has
the same law as $((u_t^N,U_t^N))_{t \geq 0}$ given by (\ref{Grev39c2a})-
(\ref{Grev39e}) .

 (b) Ignoring colour, the the total number of particles or occupied sites in the WRGB-system has the same
distribution as  the number of particles respectively occupied sites in W-BL-system and
therefore we can identify the total number of occupied sites with the total number of individuals of
a CMJ process which we denote by $(K_t)_{t \geq 0}$.

(c) The total number of particles (or sites) in the WRP-system is less than
that in the WRB-system and larger that that in the WR-system a.s.

(d) Consider the coloured particle system at times \be{xg10} T_N+t,\quad T_N
=\alpha^{-1} \log N. \ee The total number of purple or blue
particles is $O(e^{3\alpha t})$, more precisely, we have uniformly
in $N$,
\bea{bbound}
&&\sum_{i_W,i_R,i_P}i_P\cdot
\Psi ^N_{WRP,T_N+t}( i_W,i_R,i_P)\\
&& \leq \sum_{i_W,i_R,i_G,i_B}i_B\cdot \Psi^N_{WRGB,T_N+t}(
i_W,i_R,i_G,i_B) \leq const\cdot N W_\ast
e^{3\alpha t}, \nonumber
\eea
where
\be{xg11b}
W_\ast=\sup_t (e^{-\alpha t}K_{t}) <\infty,\text{ a.s. }. \ee

(e) We consider again a time horizon as in (d). Define
$W_N(t) = N^{-1} e^{-\alpha t} u^N (T_N+t)$.  The process counting
the number of red particles produced by collisions of white
particles is bounded above by an inhomogeneous Poisson process
with rate
\be{xg11}
N e^{2\alpha t} (W^2_N(t))^2
\ee
where
\be{xg11c}
\lim_{N\to\infty} W_N(t)=
W \mbox{ for each } t\in\mathbb{R}.\qquad \square \ee
\end{lemma}

\beC{C.nnred}{}

The normalized (by $N^{-1}$) number of red respectively purple particles
at time $T_N+t$ conditioned on $(W, W_\ast)$ is uniformly in $N$
\be{agdd66}
O(e^{2\alpha t}), \mbox{ respectively } O(e^{3\alpha t}), \mbox{ as } t \to -\infty.
\qad
\ee
\end{corollary}

 \begin{remark}
By introducing further colours, for example,  following the first
collision of a purple particle with a site occupied by some other
colour introducing a new pair of colours analoguous to purple and
blue, etc., we could obtain upper and lower
bounds with errors of order $O(e^{k\alpha t})$ for any
$k\in\mathbb{N}$.
\end{remark}

\begin{proof} {\bf of Lemma \ref{L.enriched}}
(a) follows by observing that combined number of white, red and
purple  particles at a site behaves exactly like the dual particle process of typical
site and when they migrate they can have collisions
according to the same rules as the dual particle system.  This
follows since when a white particle migrates it collides with the
correct probability and then becomes red and when a red particle
migrates it can collide and become purple and when a purple
particle migrates it can collide again following the same rule
(and remains purple).

(b) Since the green particles compensate for the loss of red particles
due to coalescence after collision of a white with another white
particle and the blue particles arise upon collision of red with
red or white and they evolve as the black particles, the total
number of particles is as in the W-BL system.

(c) follows since (i) the blue particles are produced in one to
one correspondence with the purple particles and (ii) the purple
particles can suffer loss due to collision and coalescence but
this does not occur with the blue particles.

(d), (e)   We have proved earlier in this section using the
CMJ-theory that the normalized rate for the number of white particles
\be{xg12}
\frac{1}{N}u^{N}_{W}(T_N +t) = O(e^{\alpha t}), \ee if we
condition on $\{W,W_\ast\}$ since for all $N$ this number is
smaller than the number of the white and black particles together, which is
then asymptotically $W_N(t) N e^{\alpha t}$, with
$W_N(t) \to W$ as $N \to \infty$
by the CMJ-theory. When white particles migrate they collide
with another white particle with probability $u^N_{W}(T_N+t)/N$.
Moreover white migrants are at times before $T_N+t$
produced at rate at most $W_N(t)$ $N\alpha e^{\alpha
t}$. Therefore  the normalized rate for the number of red particles produced
conditioned on $\{W, W_\ast\}$ is at time $T_N+t$ is
satisfying $O(N e^{2\alpha t})$.

Finally a purple or blue particle occurs when a red particle
collides with a white or red and therefore the rate of this event
is of order \be{xg13}
O(N e^{ 2\alpha t})\times O(e^{\alpha t})=O(N e^{3\alpha t}),
\ee if we condition again on $\{W,W_\ast\}$.

This completes the proof of Lemma \ref{L.enriched} which will be
a key tool in the further analysis of the dual process.

\end{proof}

\bigskip
\begin{remark}  Since the WRGB system can be identified with  the
W-BL(by assigning the union of the red and green particles at a
white site to an empty site in the first copy of $\N$)  a system
we will primarily work with the former and ignore the BL system.
We then have both upper and lower bounds for the WRP system given
by the WRB, WR systems, respectively  and this will be our primary
object for the analysis of the $t\to -\infty$ asymptotics. This
will provide upper and lower bounds with error of $O(e^{3\alpha
t})$ and therefore determine the second order asymptotics as
$t \to -\infty$.
\end{remark}

\bigskip

\textbf{Step 2:   Reformulation of the nonlinear $(u,U), (u_{(C)}, U_{(C)})$ equations}

This step has two parts, first rewriting the $(u,U)$ equation and
then a second part where we do this for the multicolour version of
this equation.
\bi

{\em Part 1} \quad
We will start by bringing the equations of the limit dynamic
$(u,U)$ in a form suitable for the purpose of the proof of
Propositions \ref{strongconn} and \ref{P.Transition}. Recall
that $(u,U)$ solves the following system of differential
equations in the Banach space $(\R \times L_1(\N,\nu), \| \cdot \|)$
(remember for the first line that $U$ is a normalized quantity):

\be{gra15} \frac{du(t)}{dt}= \alpha(t)(1-u(t)) u(t) -\gamma(t)
u^2(t), \quad t \geq t_0, \ee

\bea{gra16} \frac{\partial U(t,j)}{\partial t}
  =&&
+s(j-1)1_{j\ne1}U(t,j-1)-sjU(t,j)
\\&&
  +\frac{d}{2}(j+1)jU(t,j+1)-\frac{d}{2}j(j-1))1_{j\ne1}U(t,j)\nonumber
\\&&
  +c(j+1)U(t,j+1)-cjU(t,j)1_{j\ne 1} \nonumber
\\&&
+\alpha(t)1_{j=1}\nonumber\\&&
  + u(t)(\alpha(t)+\gamma(t))[U(t,j-1)1_{j\ne1}-U(t,j)])   \nonumber
\\&&
   - u(t)\left(\alpha(t)+\gamma(t)\right )1_{j=1} \nonumber
\\&&
 -\Big(\alpha(t)(1-u(t)) -\gamma(t)u(t)\Big)\cdot U(t,j).
\nonumber\eea

\begin{remark}
We want to interpret these evolution equations
(\ref{gra15}),(\ref{gra16}) by a particle system of the type of
the mean-field dual but now making more explicit the role of
collision and in particular the first collisions of particles.
This will allow us to analyse the behaviour as $t \to -\infty$ of
the nonlinear evolution equation above. The analysis will be based
on the fact that $u(t)$ arises as the limit of $N^{-1} u^N
(T_N+t)$ and similarly $U(t)$ as the limit of $U^N(T_N+t)$ where
$T_N =\alpha^{-1} \log N$. This will allow us in the second
part of this Step 2 to introduce enrichments of the solution of
the nonlinear equations through coloured particle systems. Then by
taking the $N\to \infty$ limit of the normalized coloured particle
systems we obtain a coloured limiting evolution and by that information about the nonlinear original system.
\end{remark}

We next rewrite the  equation (\ref{gra16}) in a form suitable for the
analysis of the solution viewed as a {\em perturbation} of the linear
(collision-free and hence $u \equiv 0$, $\alpha(t)\equiv \alpha, \gamma(t)=\gamma$)
equation in the limit as $t\to -\infty$. We shall see below in
Lemma \ref{L.fp}, that with
$u(t)\equiv 0$, the system (\ref{gra16}) has the stable age
distribution $\CU(\infty, \R,\cdot)$ of the CMJ  as equilibrium
and $\alpha =c\sum_{k=2}^\infty k \CU(\infty,\R,k)$.
Hence we should organize the r.h.s. of (\ref{gra16}) in such a way that we isolate
the linear part on the one hand and the nonlinear perturbations of this linear part
on the other hand.

We define for every parameter $a\in (0,\infty)$
(for which we shall later choose the value $\alpha$) the triple of $\N
\times \N$-matrices  \be{angre3} (Q_0^{a}$, $Q_1, L),
\ee
by the equations:
 \be{agdw19a}
Q^{a}_0 = Q^0_0 -aI+a1_{j=1}= \left(q_{j,k}\right)_{j,k\in \N}
\left(= (q_{j,k}(a))_{j, k \in \N}\right),
\ee
where
\bea{agdw19a2}
&& q_{12}=s,\quad q_{1,1}=-s \nonumber\\&&
q_{2,3}=2s,\quad q_{2,1}={d}+2c+ a,  \nonumber
\\&& q_{j,j+1}=sj,\quad q_{j,j-1}=\frac{d}{2}j(j-1)+cj,\quad
q_{j1}=a, \quad j\ne 1,2\nonumber\\&&q_{jj}=-
sj-cj-\frac{d}{2}j(j-1)-a, \quad j\neq 1, \nonumber
\eea

\be{gra12}\begin{array}{l}
Q_1= (\wt q_{jk})_{j,k\in \N}, \mbox{ where }\\
 \quad \wt q_{jj}=-1,\\
\quad\wt q_{jk}=0\quad j\ne k,
 \end{array}\ee
and finally

\be{agdw19b} L = (\ell_{jk}),\quad \ell_{jj}=0, \text{and for
}j\ne 1,\; \ell_{j-1,j}={1}.
\ee

Note that the matrix $Q^a_0$ is for the forward equation but we consider
$\bar Q^a_0$ for the backward equation. Then for each $a>0$:
\be{xg15}
\bar Q^a_0 \mbox{ generates a
semigroup } (S^a(t))_{t \geq 0} \mbox{ on } L_\infty(\N)
\ee
corresponding to a unique (pure jump) Markov process on $\N$.
This process is as follows. The matrix $Q^a_0$ defines a Markov process on
$\N$. For $a=0$ we obtain the birth and death process which corresponds
to a colony where we have birth rate $s$ for each particle, coalescence
of two particles at rate $d$ and emigration of a particle at rate $c$,
except when there is only one particle. For positive $a$ the process
is put at rate $a$ in the state with only one particle.

Furthermore we abbreviate
\be{gra12b} \wt\alpha(t)=c\sum_{k=2}^\infty
kU(t,k)-\alpha. \ee

With these four ingredients equation (\ref{gra16}) finally reads:
\beL{L.U}{(Rewritten $U$-equation)}

\bea{gra17}
\frac{\partial U(t)}{\partial t} =
&&
U(t)Q_0^\alpha \\
&&
+\wt\alpha(t)[U(t)Q_1+1_{j=1}]+u(t)(\alpha(t)+\gamma(t))[U(t)L-1_{j=1}].
\nonumber
\qquad \square
\eea
\end{lemma}

\begin{remark}
Note that our equations for $(u,U)$ are forward equations and
hence all the operators $Q^a_0, Q_1, L$ act ''from the right'' on
$U(t)$.
\end{remark}

\begin{remark}
Note that we can consider the equation (\ref{gra17}) also for arbitrary
values $a \in (0,\infty)$. Denote the solutions as
\be{agr21}
(U^a(t))_{t \in \R}.
\ee
We shall see below how we can characterize $U$ among those
by a self-consistency property.
\end{remark}

We can later use the following information about the semigroup $S^a$ to
characterize the growth rate $\alpha$ of $(u(t))_t \in \R$ in the entrance law
$(u,U)$ as $t
\to -\infty$.

\beL{L.fp} {(Representation of $\alpha$)}

Consider  the evolution equation for $U^a(t)$ in the regime in which
$(u(t))_{t\in \R}$ is identically zero and $\alpha(t)\equiv a$.
Then for given $a>0$ there is a unique equilibrium state (positive eigenvector)
\be{xg16}
\{q^\ast_j(a)\}_{j\in\mathbb{N}} \mbox{ for the operator }
Q^a_0 (\mbox{ equalling }
Q^0_0-aI+a1_{j=1}).
\ee
Then $\alpha$ is uniquely determined as the
fixed point defined by  the self-consistency equation \be{gra18}
 \alpha = c\sum_{j=2}^\infty jq^\ast_j(\alpha).\ee
Also,
\be{xg17}
\{q^\ast_j(\alpha), j=1,2,\cdots\} =\{\CU(\infty,\R,j), j=1,2,\cdots \}=(p_1,p_2,\dots),
\ee
the r.h.s. being the stable size distribution of the CMJ process.
$\qquad \square$
\end{lemma}
\begin{proof} The existence of the unique equilibrium $q^\ast_j(a)$
follows from standard Markov chain theory, the question is whether the
self-consistency equation (\ref{gra18}) has a solution.
Note for this purpose first that $q^\ast_j(0)$ is the
equilibrium for the single site birth and death process
appearing in the McKean-Vlasov dual process for one component and
observe that
$\sum_{j=2}^\infty jq^\ast_j(0)>0$. We define:
\be{xg17b}
F:a \la c \sum_{j=2}^\infty jq^\ast_j(a).
\ee
Then we are left to show that the fixed point equation
$\alpha = F(\alpha)$ has a solution.
We saw above $F(0)>0$. We claim next that the function
 is monotone decreasing in $a$, converging to 0 as $a \to \infty$ and continuous.

 To verify the monotonicity,  let $a_2>a_1>0$. Recall that the chain for $a \equiv 0$
  starting from 1 is stochastically increasing to its equilibrium. Now consider
the two Markov chains for $a_1, a_2$ as a sequence of independent excursions away
from $1$ which end when a jump to $1$ at rate $a$ occurs but otherwise follow
the $Q^0_0$ dynamic.  To
compare the average height over the excursions, we consider two
such excursion lengths given by coupled exponentials $(a_2)^{-1}
\mathcal{E},\; (a_1)^{-1}\mathcal{E}$, respectively. Then we
observe by a simple coupling argument that the heights of the $a_1$-excursion at
times $(a_2)^{-1} \mathcal{E}<t\leq \; (a_1)^{-1}\mathcal{E}$
stochastically dominate the height at time $(a_2)^{-1}
\mathcal{E}$. Therefore the average height over an excursion for
the $a_1$-chain is stochastically greater than or equal to that for
the $a_2$-chain.

The continuity follows by noting that the
$a_2$-excursions from zero converge to the $a_1$ excursion from
$0$ if  $a_2\downarrow a_1$.

Noting finally that $\sum_{j=2}^\infty
jq^\ast_j(a)$ converges to $0$ when $a\to\infty$, we obtain the
existence and uniqueness of a solution of (\ref{gra18}) which
we read as a fixed point of the equation $\alpha =F(\alpha)$.
\end{proof}

\bi

{\em Part 2.} \quad
We now turn to the second part of Step 2, where we give the
limiting dynamics of the multicolour enrichment as represented by
$(u^N_{(WRPGB)},U^N_{(WRPGB)})$. We also need systems of other subsets
of colours in our arguments and therefore we denote the limiting objects
analog to the noncoloured system by
\be{xg17b1} (u_{(C)},U_{(C)}),
\ee
with
\be{xg17b2}
C=\{W,R,P,G,B\} \mbox{ or some other colour subset of } \{W,R,P,G,BL,B\},
\ee
where $u_C$ is a positive real number and $U_{(C)}$ is a measure on
\be{xg17b3}
(\N_0)^C.
\ee

This set-up means that we
consider the processes of {\em occupied sites}, occupied with
the various colour combinations.
The dynamics, in the $N\to\infty$ limit, of the $WRP$-system
(or the WRGB-system) at time $T_N+t$ is given by an enrichment
$(u_{(WRP)}(t),U_{(WRP)}(t))$ of the nonlinear system
(\ref{gra15}), (\ref{gra16}). Namely the latter is recovered as
follows:

\be{xg18}  u(t)=u_{(WRP)}(t),\quad
U(t,k)=\sum_{\{(i_W,i_R,i_P):i_W+i_R+i_P=k\}}U_{(WRP)}(t;i_W,i_R,i_P).
\ee

In order to write down the evolution equation for
the quantities from (\ref{xg17b}), we need the changes occuring in
the underlying multicolour system and its various subsystems which
then induce the changes in the measure on these configurations.
For this purpose we need
operators associated with the various possible transitions, their
rate and their form which are associated with a particular current state.
This state  is given by a tuple of the form $(i_W,i_R,i_P, i_B,
i_G, i_{BL})$ or a tuple for a smaller set of colours which gives the
number of particles of the various colours at a site. Therefore
transitions and transition rates are specified by matrices of the form
\be{xg18b}
Q_{\underline{i},\underline{j}} \quad , \quad \mbox{ with }
\underline{i},\underline{j} \in (\N_0)^C \quad , \quad
C \subseteq \{W,R,P,B,G,BL\}. \ee

We can now distinguish two different groups of transitions
those which concern only particles of {\em one
colour}  and then there are transitions where particles of different {\em colours interact}.
We specify these transitions and the corresponding operators in
(\ref{xg19}) and then in (\ref{coloper}), (\ref{creaops}).

For the first type the one colour operators let
\bea{xg19}
Q^{0,i}&=& \text{ operator } Q^0_0\text{ applied to
the ith colour occupation numbers - see (\ref{agdw19a})}\\
&& \text{corresponding to
birth, coalescence and emigration} \nonumber \\
&& \text{(as long as there are more than
one particle) of the ith colour individuals}.\nonumber \eea

Next we introduce the appropriate {\em inter-type
coalescence} operators (for versions with formulas for the r.h.s., see
(\ref{xg201}) and sequel below) which describe the  changes in the limiting
frequency measure $U$ based on changes in the occupation numbers which correspond
to actions of the coloured particles. There are essentially two types of intertype-transitions,
(1) changes which occur at a site with rates depending on the state at this site and
(2) the creation of a new site at rates depending on the state at the founder
site leading to the migration operators.

The two groups of matrices corresponding to the coalescence, respectively
migration operator, are the following matrices
\be{agr22}
Q^A = (Q^A_{\underline{i},\underline{j}}), \quad \underline{i}, \underline{j} \in \N^C_0
\mbox{ with } A \subseteq C \times C \mbox{ or }
A \subseteq \{a \to b | a,b \in C\}
\ee
running through the list of names of the various transitions
are, beginning with the {\em coalescence} operators:
\bea{coloper}
Q^{WR}
&&= \text{coalescence of white and red at same site}\\
&& \text{yielding a white and green pair at this site},\nonumber\\
Q^{WP}
&&= \text{coalescence of white and purple at same site}\nonumber\\
&& \text{yielding a white and blue pair  at this site},\nonumber\\
Q^{PR}
&&= \text{coalescence of  red and purple at same site}\nonumber\\
&& \text{yielding a red at this site},\nonumber\\
Q^{RG}
&&= \text{coalescence of green and red at same site}\nonumber\\
&& \text{yielding a red at this site}\\
&&= \text{coalescence of purple and blue at same site}\nonumber \\
&& \text{yielding a purple at this site},\nonumber \\
Q^{WG}
&&=Q^{WB}=0.\nonumber\eea
and the {\em migration}  operators, which are effectively creation operators
for new occupied sites, which in particular are then sites occupied by
initially only one particle.  For sites occupied by {\em only} one colour we can talk of green (G), blue (B), purple (P), red (R), white (W) sites without
ambiguity. However for multicolours, a WR site denotes a site having at least one white or one red particle, etc.


\bigskip

\bea{creaops} Q^{W\to W}&&= \text{creation of white site at an empty site
in }\{1,\dots,N\},\\
Q^{W\to R}&&= \text{creation of red site at randomly chosen occupied
site in } \{1,\dots,N\},\nonumber\\
Q^{R\to R}&&= \text{creation of red at empty site
in }\{1,\dots,N\},\nonumber\\
Q^{R\to P,B}&&= \text{creation of purple-blue pair }\nonumber\\&&
\text{ with purple at randomly chosen occupied site in
}\{1,\dots,N\}\nonumber\\&&\text{ and blue at the first empty site in
 the first copy of }\mathbb{N},\nonumber
\\Q^{P\to P}&& = \text{creation of purple site at randomly chosen
site in  }\{1,\dots,N\},\nonumber\\
Q^{G\to G}&&= \text{creation of green site at first empty site in
the first copy of }\mathbb{N}, \nonumber
\\
Q^{B\to B}&&= \text{creation of blue site at first empty site in
 the first copy of }\mathbb{N}.\nonumber\eea

The verbal description on the r.h.s above corresponds to the
following expressions which we spell out in two examples:

\be{xg201} Q^{W\to W}_{(i_W,i_R,i_G,i_B),(1,0,0,0)}=c\cdot
1_{i_W+i_R\,\ne 1}\,i_W, \ee

\be{xg202} Q^{W\to
R}_{(i_W,i_R,i_G,i_B),(i^\prime_W,i^\prime_R+1,i^\prime_G,i^\prime_B)}=c\cdot
1_{i_W+i_R\,\ne 1}\,i_W U_{WRGB}(t,i^\prime_W,i^\prime_R,i^\prime_G,i^\prime_B),\quad \text{
etc.}. \ee
\bi

We now develop in the detail the system of equations for
$U_{WRGB}(t)$ since the former provides the necessary upper and
lower bounds.

The pair $(u,U)$ arises as the limit as $N\to\infty$ of the
particle system $(u^N,U^N)$ in times \be{agre52} T_N +t \text{
  where   } T_N:=\frac{ \log N}{\alpha}. \ee  The proofs are  based on first
taking the limit as $N \to \infty$ of the
$(u^N_{(WRGB)},U^N_{(WRGB)})$ system and then identifying the order of the
terms corresponding to the different colours in the limit as $t
\to -\infty$ based on the structure  of the particle system.

Now we set for the rescaled processes (misusing earlier ~-notation)
 \be{gra21b2}\begin{array}{l}
 \wt u_{WRGB}^N(t) = N^{-1} u^N_{WRGB}(T_N+t),
\\ \wt U_{WRGB}^N (t,i_W,i_R,i_G,i_B) =U^N_{WRGB}(T_N+t),
\end{array} \ee

\be{gra21b2x}\begin{array}{l}
\wt u_{WRP}^N(t) = N^{-1} u_{WRP}^N(T_N+t) , \\
\wt U_{WRP}^N (t,i_W,i_R,i_P) = U_{WRP}^N(T_N+t,i_W,i_R,i_G,i_B).
\end{array} \ee

Then let $N \to \infty$ and we  again obtain,
the existence of the limit of the rescaled system $(u^N, U^N)$,  the
existence of the coloured versions. The limiting objects are
denoted: \be{gra21c} (u_{WRGB}(t), U_{WRGB}(t,\cdot,\cdot,\cdot)),
\quad (u_{WRP}(t), U_{WRP}(t,\cdot,\cdot,\cdot)), \ee which will
satisfy a system of equations which we give below in (\ref{WRGB3}),
and (\ref{WRGB1}) respectively in
(\ref{WRP1}). We omit the details of the convergence proof
here, which are straightforward modifications of the argument
given in the proof of Proposition \ref{P.Grocoll}.

\begin{remark}
In order to give an intuitive picture, even after taking the limit
$N\to\infty$, we shall still refer to this nonlinear system in
terms of coloured particles, but of course the quantities in
question are now  continuous quantities that arise in the limit as
$N \to \infty$ of the normalized quantities corresponding the the
well defined particle system given by  the above pure jump
process.
\end{remark}

Note first the following facts
about the limiting evolution which are important for our purposes
and which follow from the corresponding finite-$N$ properties:
\begin{itemize}
\item  since the green particles compensate for the loss of any
red particles due to coalescence, the number of red plus green
particles at a given site has an evolution which is independent of
the white particles (but not the distribution of the relative
proportions of red and green),

\item the red plus green plus blue particle system has the same
dynamics as the black particle system and since red and black
founding particles are created at the same time both systems have the same
distribution,

\item the particle system comprised of white, red and purple
particles has the same distribution as the  dual particle system.
\end{itemize}


Recalling that the WBL-system is a CMJ-system and the RGB process can be identified with the set of
black particles, it follows that $(K^N_{(W)}(t)+K^N_{(RGB)}(t))_{t \geq 0}$ is
less than a CMJ process with Malthusian parameter $\alpha$ but the difference is of
order $o(e^{\alpha t})$ and we have
for $t \in \R$, that the intensity of the occupied sites in the
various coloured systems satisfies:
\bea{xg20}
&& u_{WR}(t)\leq u_{WRP}(t)\leq u_{WRB}(t) \leq u_{\mathrm{WRGB}}\\
&&\leq
u_{W}(t)+u_{RGB}(t)= \lim_{N\to\infty}
\left[ \frac{1}{N} (u^N_{W}(T_N+ t)+u^N_{RGB}(T_N+ t))\right]
\leq We^{\alpha t}.\nonumber
\eea
The intensity of individuals in the various coloured systems satisfy that
\bea{ddxg1} &&
\sum_{i_W,i_R,i_G,i_B}(i_W+i_R)u_{WRGB}(t)U_{WRGB}(t,i_W,i_R,i_G,i_B)\\&&\leq
\sum_{i_W,i_R,i_P}(i_W+i_R+i_P)u_{WRP}(t)U_{WRP}(t,i_W,i_R,i_P)\nonumber\\&&
\leq  \sum_{i_W,i_R,i_G,i_B}(i_W+i_R+i_B)
u_{WRGB}(t)U_{WRBG}(t,i_W,i_R,i_G,i_B).\nonumber\eea

As a result of these inequalities we note that the $WRGB$ system
provides {\em lower} and {\em upper} bounds for the ''number'' of occupied sites
and ''total number of particles'' for the $WRP$ system which
corresponds to the dual particle system of interest.  We note that
the difference between the upper and lower bounds for the total
number of particles (which corresponds to the number of blue
particles) is of order $O(e^{3\alpha t})$. Since each occupied
site must contain at least one particle this provides bounds for
the number of occupied sites with the same order of accuracy. For
this reason we focus on the WRGB system and then read off the
required estimates for the WRP system up to this order of a accuracy.

We now write down the equations for the WRGB-system, i.e. for the corresponding
pair $(u,U)$. We need the following abbreviations where all sums are over
$i_W, i_R, i_G$ and $i_B$:

\be{xg21} u_{(WR)}(t)=u_{(WRGB)}(t)\cdot
\left(\sum 1_{i_W+i_R\geq 1}U(t;i_W,i_R,i_G,i_B)\right), \ee
 \be{xg21b} u_{(WRB)}(t)=u_{(WRGB)}(t)\cdot
\left(\sum 1_{i_W+i_R+i_B\geq 1}U(t;i_W,i_R,i_G,i_B)\right), \ee

\be{xg22} \alpha_W(t)= c \sum  i_W\cdot 1_{i_W+i_R+i_G\geq
2}\,U_{(WRGB)}(t;i_W,i_R,i_G,i_B), \ee

\be{xg23} \alpha_R(t)= c \sum i_R\cdot 1_{i_R+i_W+i_G\geq
2}\,U_{(WRGB)}(t;i_W,i_R,i_G,i_B),\ee


\be{xg24} \alpha_G(t)= c \sum i_G\cdot 1_{i_G\geq
2}\,U_{(WRGB)}(t;i_W,i_R,i_G,i_B),\ee

\be{xg25} \alpha_B(t)=c \sum i_B\cdot 1_{i_B\geq
2}\,U_{(WRGB)}(t;i_W,i_R,i_G,i_B), \ee


\be{xg25b} \gamma_W(t)= c
1_{i_W+i_R+i_G=1}1_{i_W=1}\,U_{(WRGB)}(t;1,0,0,0),\quad
\text{similarly for }\gamma_R(t),\gamma_G(t), \gamma_B(t).
\ee
Then we can write down the equation for $(u_{(C)},U_{(C)})$ as
coloured version of the $(u,U)$-equation. Again as in the
latter case we work in the very same Banach space
$(\R \times L_1(\N, \nu), \| \cdot \|)$ without
mentioning this  explicitly below.

Consider the WRGB-system described by the
\be{add140}
(u_{(WRGB)}(t),  U_{(WRGB)}(t,\cdot)), t \in \R)
\ee
satisfying the following nonlinear equation:

\bea{WRGB3}
 \frac{\partial u_{\mathrm{(WRGB)}}(t)}{\partial t}&=&
(\alpha_W \to \alpha_W(t) + \alpha_R\alpha_R(t))u_{\mathrm{(WRGB)}}(t)(1-u_{\mathrm{(WR)}}(t)) \\
&&+ \alpha_G \to \alpha_G(t)
u_{\mathrm{(WRGB)}}(t)+\alpha_B \to \alpha_B(t) u_{\mathrm{(WRGB)}}(t)\nonumber\\
&& -(\gamma_W \to \gamma_W(t)
+\gamma_R \to \gamma_R(t)) u_{\mathrm{(WR)}}(t)u_{\mathrm{(WRGB)}}(t).\nonumber\eea

\be{WRGB1}
\frac{\partial U_{\mathrm{(WRGB)}}(t)}{\partial t}=
U_{\mathrm{(WRGB)}}(t)Q^{\mathrm{(WRGB)}}(t),
\ee
where
 \bea{WRGB2}
 Q^{\mathrm{(WRGB)}}(t)=
 &&\sum_{i=W,R,G,B}
Q^{0,i}+Q^{\mathrm{(WR)}}+Q^{RG}\\&&+ (1-u_{\mathrm{(WR)}}(t))[Q^{W\to W} +Q^{R\to
R}]+u_{\mathrm{(WR)}}(t)[Q^{W\to R}+Q^{R\to B}],\nonumber\\
&& - \alpha(t) I,\nonumber
\eea with initial condition at time $t= -\infty$ given by

\be{agdd67}
\liml_{t \to -\infty} u_{\mathrm{(WRGB)}}(t)=0,
\ee

\bea{gra29z}
\lim_{t\to -\infty} U_{\mathrm{WRGB}}(t,i_W,i_R,i_G,i_B)&&=
\mathcal{U}_\infty(i_W)\quad\text{if  } i_R=i_G=i_B=0 \\&&=
0,\quad\text{otherwise}, \nonumber
\eea
where $\mathcal{U}_\infty$
is the stable size distribution of the CMJ process.

If we write for a better understanding of (\ref{WRGB1}) this equation out pointwise, we obtain
(we suppress on the r.h.s. the subscript WRGB in $U$):

\bea{gra29}&&\\
&& \frac{\partial U_{\mathrm{(WRGB)}}}{\partial t}(t;i_W, i_R, i_G, i_B)= \nonumber\\
 &&
  +\frac{d}{2}i_W(i_R+1)U(t;i_W,i_R+1,i_G-1,i_B)-\frac{d}{2}i_W i_R U(t;i_W,i_R,i_G-1,i_B)\nonumber
\\&&
  +\frac{d}{2}(i_W+i_R)(i_G+1)U(t;i_W,i_R,i_G+1,i_B)\nonumber\\
  &&
  -\frac{d}{2}(i_W+i_R)i_G 1_{i_R\ne1}U(t;i_W,i_R+1,i_G-1,i_B)\nonumber
  \\&&
+\alpha_W(t)(1-u_{(WR)}(t))1_{(i_{W=1,i_R=0,i_G=0,i_B=0)}}\nonumber\\&&
+\alpha_R(t)(1-u_{(WR)}(t))1_{(i_{W=0,i_R=1,i_G=0,i_B=0)}}\nonumber
\\&& + (\alpha_W(t)+\gamma_W(t))u_{(WR)}(t)
(U(t;i_W,i_R+1,i_G-1,i_B)-U(t;_W,i_R+1,i_G-1,i_B)\nonumber
\\&&
+(\alpha_G(t)+\sum_{i+j\geq 1}U(t;i_W,i_R+1,i_G-1,i_B))1_{(0,0,1,0)}\nonumber
\\&&+ (\alpha_B(t)+
\sum_{j}jU(t;i_W,i_R+1,i_G-1,i_B)u_{(WR)}(t))1_{(0,0,0,1)}\nonumber
\\&&
 -\left(\alpha(t)(1- u_{WR}(t))-\gamma(t)u_{WR}(t)\right) \cdot U_{(WRGB)}(t;i_W,i_R+1,i_G-1,i_B).
\nonumber
\eea

The process $(u_{(WRP)}, U_{(WRP)})$ satisfies for $u_{(WRP)}=u$
our original equation, and the quantity
$U_{(WRP)}(t)$ satisfies a similar set of equations as above with $B$
replaced by $P$ but in this case we must add the terms
corresponding to the migration of a purple particle to a site
occupied by white, red or purple.  Also when a purple coalesces
with red at the same site (which occurs with  rate $d$) it
produces a red-green pair and when it coalesces with a  white at
the same site it produces a white-green pair.
The equation for $U_{(WRP)}$ has therefore the form
\be{WRP1}
\frac{\partial U_{(WRP)}(t)}{\partial t}= U_{(WRP)}(t)Q^{WRP}(t),
\ee
where
 \bea{WRP2}
 Q^{WRP}(t)= &&\sum_{i=W,R,P}
Q^{0,i}+Q^{WR}+Q^{WP}+Q^{PR}\\&&+ c(1-u(t))[Q^{W\to W} +W^{R\to
R}]+cu(t)[Q^{W\to R}+Q^{R\to P}]\nonumber\\&&+cQ^{P\to P} -
\alpha(t) I.\nonumber\eea

This now completes the set of limiting $(\N \to \infty)$ equations for
the coloured particle system at time $T_N+t$.

For the further analysis it is crucial to observe the following
two facts:
\begin{itemize}
\item for calculating the asymptotics of {\em first}
moments of the McKean-Vlasov system it suffices
 to isolate dual process effects
which  in the limit $N\to\infty$ are of the order $e^{\alpha t}$
as $t \to -\infty$ and  ignore error terms of order $e^{2 \alpha
t}$ and higher, \item for calculating the asymptotics of {\em second}
moments of the McKean-Vlasov system it suffices to isolate effects which are in the limit
$N\to\infty$ of the order $e^{2\alpha t}$ as $t \to -\infty$ and
ignore error terms of order $e^{3 \alpha t}$ and higher.

\end{itemize}

It is easily verified that if we {\em condition} on the growth
constant $W$ and on $W_\ast$,
we can verify that as $t\to-\infty$, we have the
following estimates on the  ''number of sites'' occupied by
certain colour combinations:
\be{ineq1}
\suml_{i_W,i_R,i_G,i_B} 1_{i_R+i_G+i_B\geq
1}\,U_{(WRGB)}(t;i_W,i_R,i_G,i_B)= O(e^{\alpha t}),
\ee

\be{ineq2}
\suml_{i_W,i_R,i_G,i_B} 1_{i_B\geq 1}\,U_{(WRGB)}(t;i_W,i_R,i_G,i_B)=
O(e^{2\alpha t})
\ee
and therefore
\be{ineq2b}
u_{RGB}(t)=O(e^{2\alpha t}) \mbox{ and } u_B(t)=O(e^{3\alpha t}).
\ee

We conclude showing that as $t \to -\infty$ this difference
between the red plus white particle system (with WRGB dynamics)
and the true dual (WRP) is of order $o(e^{2 \alpha t})$ once we
have taken the limit $N \to \infty$. This difference is the
difference between the system of {\em purple} particles and the system
of {\em blue} particles.

Again we condition on $W$ and $W_\ast$.
We first note that the intensity (after taking
$N \to \infty$) of red particles is $O(e^{2 \alpha t})$ and
therefore the rate of production for collisions of red particles
with other occupied sites by white, red or purple is $O(e^{3
\alpha t})$. Hence we have as $t \to -\infty$
\be{xg26b}
u_{P}(t) = O(e^{3 \alpha t}), \quad u (t)-u_{(WRGB)}(t) = O(e^{3 \alpha t}).
\ee

The sites occupied by  red,  green or blue particles correspond to
black sites and therefore bound the  number of occupied sites
eventually lost (compared to CMJ) due to collisions. What is the meaning of
the three components? The population of green
particles represent asymptotically as $t \to -\infty$ the total
number of particles lost (in the dual particle system compared
with the collision-free CMJ-system) due to collisions up to order
$o(e^{2 \alpha t})$.

\begin{remark}
It is important to remember that we have seen above in (\ref{xg20}),
(\ref{ddxg1}) that we can obtain upper
and lower bounds for the exact dual using {\em white plus red} for
the lower bound and {\em white plus red plus blue} for the upper
bound and that by (\ref{xg26b}) we know that the difference
of lower and upper bound is $O(e^{3 \alpha t}) = o(e^{2\alpha t})$.
Hence suffices for the purpose of first and second moment
calculations for the random McKean-Vlasov entrance law
to work with the WRGB-system on the dual side.
\end{remark}

\bigskip

{\bf Step 3:} \begin{proof} {\bf of Proposition \ref{strongconn}}

(a) We define the coupling using the multicolour system where white and black
gives the CMJ-process and WRP the dual particle system. The evolution rules
define the processes in a standard way for all $N$ on one probability space
using Poisson stream for all potential sites, colours and transitions.

Recall that the total number of sites occupied by black plus
white particles can be identified with a CMJ process $K_t$ and our
processes for each $N$ can be coupled based on a single realization
of this CMJ process using the multicolour system. Let on this common
probability space
\be{dd08091}
W:=\lim_{t\to\infty} e^{-\alpha t}K_t.
\ee Then the existence of the limit
in (\ref{dd08091}) and therefore assertion (\ref{agrev91}) follows
immediately from the Crump-Mode-Jagers theory.

In order to obtain the assertion (\ref{gra4}) we shall estimate below the {\em mean}
and {\em variance} of
\be{xg27aa}
e^{-\alpha t} N^{-1} u^N(T_N+t)\qquad (\mbox{recall } T_N=
\alpha^{-1}\log N)
\ee
by bounding the mean and variance of the
normalized number of black particles (equivalently,  red plus
green plus blue particles) and finally put this together to get
the claim of the proposition.

{\em Part 1: Bounds on the mean}.

Define here \be{xg27}
 W_{t} =  e^{-\alpha t}
K_{t},\quad W_{N,t}=W_{T_N+t}, \quad W=W_{\infty}. \ee Recall that
by CMJ theory (\cite{N}), and finiteness of $E[W^2]$,

\be{xg27a}
\bar w=  \sup_t E[W_t]<\infty,\quad \overline{w^2} = \sup_t
E[W^2_t] < \infty. \ee

We condition now on the path
$(K_t)_{t \geq 0}$ (and then in particular $W$ is given) and give a
conditional upper bound
for the number of black particles produced by time $T_N+t$. This bound is
given by considering the following upper bound for the rate of
collisions of white particles at time $s$ which is

\be{xg27b} (\alpha_N(s) + \gamma_N (s)) K^N_s \frac{K^N_s}{N},
\ee
where $\alpha_N(\cdot)$ and $\gamma_N(\cdot)$ are defined as
in (\ref{Grev39e}).
Recall that that $\alpha_N (s) \to \alpha(s)$,
$\gamma_N(s) \to \gamma(s)$ and
$(\alpha(s)+\gamma(s))K_s \leq (\alpha +\gamma)W_se^{\alpha s}$
as $N \to \infty$.
We get a stochastic upper bound of the quantity in (\ref{xg27b})
by taking the random process
\be{xg27b2}
(\alpha + \gamma) W_se^{\alpha s} \frac{K^N_s}{N}.
\ee
Now assume that we realize a Poisson point process on $[0,\infty)$ with
intensity measure given by the density given in (\ref{xg27b}). Now
we use this Poisson point process to generate the founders of the black families.
Then we let a black cloud grow descending from
one black ancestor born at time $s$ up to time $(T_N+t-s)$ and this evolution
is independent of the Poisson point process.  This
black cloud born at time $s$ grows till time $T_N+t$ as
\be{xg27c}
\wt W_{N,s}(t) e^{\alpha(t-s)}
\ee
and is given by an independent copy of the CMJ process starting with one
particle at time $s$. With this object we obtain a stochastic upper
bound on the number of black particles in the original dual process.

Since at every time $s$ of birth of a black
particle we get a cloud independent of all other clouds and also
independent  of $(K_s)_{s \geq 0}$, we get an upper bound for the expected number of black
particles {\em given a realization of } $(K_s)_{s \geq 0}$:

\be{gra6b} \int_0^{T_N+t} E[\wt W_{N, t-s}] e^{\alpha(T_N+t-s)}
(\alpha_N (s) + \gamma_N(s)) K^N_{s} \frac{K^N_s}{N}ds. \ee

Furthermore with the Crump-Mode-Jagers theory applied to $(K_s,
\CU(s))$ we get setting $D = E[\wt W] \frac{\alpha +\gamma}{\alpha}
\in (0,\infty)$ that the quantity in (\ref{gra6b}) can be bounded above
asymptotically by (recall (\ref{xg27}) for $W_{N,s}$):
\be{gra6c}
\frac{1}{N} \int_0^{T_N+t}
e^{\alpha(T_N+t-s)} \alpha D W_{s}^2e^{2\alpha s} ds
\quad \mbox{ as } N \to \infty. \ee The {\em mean} of the
expression in (\ref{gra6c}) is bounded by
\be{gra6d}
\overline{w^2} \alpha D
\frac{1}{N} \int_0^{T_N+t} e^{\alpha(T_N+t-s)}  e^{2\alpha
s} ds. \ee

This quantity in turn is equal to \be{gra6c2} N \overline{w^2}
e^{2\alpha t}(1-e^{-\alpha (T_N+t)}) = ND\overline{w^2} e^{2\alpha
t}(1-e^{-\alpha t} \frac{1}{N}). \ee
Therefore we obtain the upper
and lower bound:
 \be{gra7}
\bar w e^{\alpha t}\geq   E[ \frac{1}{N} u^N(T_N+t)]\geq
E[W_{N,t}] e^{\alpha t}- \overline{w^2} e^{2\alpha t}(1-O(\frac{1}{N})),
\ee
where the last expression is an upper bound on the expected
number of black particles.

Therefore we get the final bound for the mean of the normalized number
of sites:
\be{gra7a} 0 \geq  E [e^{-\alpha t} N^{-1} u^N (T_N+t)
-W_{N,t}] \geq c_N \cdot e^{\alpha t} (1-O(\frac{1}{N})) \ee
where $\sup_N (|c_N|)<\infty$. \bi

{\em Part 2: Analysis of the variance}

To complete the result we now have to estimate
\be{gra7a2}
Var [e^{-\alpha t} N^{-1} u^N (T_N +t)-W_{N,t}].
\ee
We will verify
as a first step that the variance of the {\em normalized} number
of black particles converges to $0$ as $N\to\infty$ and then we
will return to the dual particle system. Recall that in our
calculations we condition on $(K_t)_{t \geq 0}$.

The first step is now to condition again on a realisation of a process which
is a stochastic upper bound on the
number of white particles (compare part 1). This bound is given by a
CMJ-process. Therefore we observe first that given
this number the evolution of the black clouds once they are founded
are all independent. Therefore we obtain an {\em upper bound} on the
normalized variance of the black particles if we use a path of
the CMJ-process, which dominates the white population.

Note that therefore the {\em birth of new black clouds} of particles as $N \to
\infty$  (which arises upon collision) can by (\ref{gra6c}) be
bounded by a time inhomogeneous Poisson process with intensity
\be{xg30}
\left(\alpha \frac{W_{N,t}^2 e^{2\alpha t}}{N}\right) ds,
\ee
where we condition on a realisation of $W(t), W_{N,t} = W(T_N+t)$,
or alternatively on the path $(K_s)_{s \geq 0}$.

Define $L^r_s$ as (for given $r$ in the variable $s$) the Laplace
transform of $e^{-\alpha (r-s)} K_{r-s}$ (starting with one
particle).

The Laplace transform $L$ of the total number of black
particles at time $r$ then is given by
\be{gra8}
L^r(\lambda)=\exp\left(-\int^r_0 D \frac{\alpha
(K^N_s)^2}{N}(1-L^r_s(\lambda e^{\alpha(r-s)}))ds \right). \ee

Note that by construction $L^r_s =L_{r-s}$. Furthermore
we have that:
\be{gra8b}
\sup_{s\geq 0}Var[e^{-\alpha s}K_s] =
\supl_{s \geq 0} L_s^{\prime\prime}(0)<\infty. \ee

Now we apply this to $r=T_N+t$ and conclude that conditioned on
$W$ (this is indicated by $\sim$ on $E, Var$, etc.):
\be{gra8a}
\wt Var [e^{-\alpha t} N^{-1}\cdot(\# \mbox{ of black
particles at time } T_N + t)]
\ee
is bounded by (Const means here a function of $W$ only)
\be{gra9}
\frac{e^{-2\alpha t}}{N^2}\int_0^{T_N+t} DW^2 \alpha
e^{2\alpha(T_N+t-s)}\frac{e^{2\alpha s}}{N}L_s^{\prime\prime}(0)ds
\leq Const \frac{N^2({\log N})}{N^3}. \ee The r.h.s. is
independent of $t$. Hence {\em uniformly} in $t \in \R$, as
$N\to\infty$:
\be{gra9a}
\wt Var [e^{-\alpha t} \frac{1}{N} (\#
\mbox{ black particles at time } T_N+t)]=o(1).
\ee
Hence {\em conditioned on $W$} the number of black and white particles
normalized by $N e^{\alpha t}$ is
{\em deterministic} in the limit $N \to \infty$. Hence by (\ref{gra9a})  the
density of white particles alone and  the normalized difference between the number of
black and white particles is deterministic in the $N \to \infty$
limit.

{\em Part 3 Conclusion of argument}

Since the dual lies between the set of CMJ particles and the set
of white particles, conditioned on the variable  $W$, in the limit
$N^{-1} e^{-\alpha t} u^N (T_N+t)$ lies between {\em two
deterministic curves}, both converging to the same constant
as $t \to -\infty$, namely $W$. Using (\ref{gra7}) and (\ref{gra7a}) we
conclude  inequality (\ref{gra4}).

 The assertion (\ref{gra4b}) is proved as follows. Since the
first component was just treated,  we focus on the second, i.e.
$U^N$. Note that the total variation distance between $U^N(T_N+t)$
and $\mathcal{U}(T_N+t)$ is bounded by the normalized number of
black particles.  Then combining (\ref{gra6b}) with (\ref{gra6c2}) we have
\be{gra9a2}
\lim_{t\to -\infty} \| U^N(T_N+t)-\mathcal{U}(T_N+t)\| =0.\ee

(b) is an immediate consequence of the analysis above.

\end{proof}

\bi

{\bf Step 4:} \begin{proof} {\bf of Proposition \ref{P.Transition}}

We proceed in six parts. Four parts prove the various claims and two parts
 are needed to prove some key lemmata at the end of the argument.

{\bf Part 1 - Proof of (\ref{ga30++})}

To verify (\ref{ga30++}) we use that the $u(t),\alpha (t),
\gamma(t)$ arise as the limit $N\to \infty$ of $\frac{u^N(t)}{N},
\alpha^N(t), \gamma^N(t)$ for which we can use the representation by the
multicolour particle system. For each $N$ the inequality
\be{gra9b}
(\alpha_N(t)+\gamma_N(t))u^N(t)\leq(\gamma+\alpha)N
W_{N,t}e^{\alpha t} \ee follows since the left side counts the
number of particles in a subset of the set of particles
on the right side (the difference is induced by the green
particles in the multicolour construction). Furthermore we know
that $W_{N,t}\to W, \;\text{a.s.}$. Hence the  inequalities in
(\ref{gra9b}) are therefore inherited in the limit as
$N\to\infty$ and give
\be{gra9bb}
(\alpha(t) + \gamma(t)) u(t) \leq (\alpha+\gamma) W e^{\alpha t}.
\ee

\bi

{\bf Part 2 - Proof of (\ref{grocoll3}), (\ref{grocoll32}) }

To obtain  (\ref{grocoll3}), (\ref{grocoll32}) we use the
coloured particle system introduced  in Step 1 of this subsection.
We show first that $\kappa$ and $\kappa^\ast$ are strictly positive
and then later on we identify these numbers by a more detailed analysis.

Recall that the total number of white, red, green and blue
particles at time $T_N+t$ equals the total number,
$W_{N,t}Ne^{\alpha t}$  in the CMJ process and the collection of
white, red and purple particles is a version of the actual dual
particle system $\eta^N$. We have shown that if we consider these
particle systems at time $T_N+t$ and let $N \to \infty$ we obtain
a limit dynamic.
The analoguous statement holds for the limiting dynamic as $N \to \infty$,
now with $t$ as the time variable for the multicolour system.
The techniques are the same as used in Subsubsection \ref{sss.asycolreg}
and we do not write out here the details again.
We are now interested in the
expansion as $t \to -\infty$ in this multicolour limit dynamics.

In contrast to the proof
of Proposition \ref{strongconn} the argument here now involves the green
and blue particles.

The purpose of the green particles is to keep
track of the particles {\em lost due to collisions} of white with
red and purple particles (lost by coalescence).  The purpose of the blue particles
is to obtain with the WRB-system an upper bound for the total number of particles in
the $WRP$-system by suppressing the loss of particles that could
occur if a red or purple particle migrates to an occupied site.
We observe that therefore  that since the number of blue
particles is $O(e^{3 \alpha t})$ (recall (\ref{ineq2b})), then in order to identfiy the
$e^{2 \alpha t}$ term in the total number of dual particles as
$t \to -\infty$, it is sufficient to control the asymptotics as
$t \to -\infty$ of the green particles. Therefore we now have to
study the green population, for which we first need more information
about the red particles.
\bi

The rate at which white particles migrate and collide with
occupied sites thus creating a red particle is given and
estimated as follows (recall (\ref{gra9bb})):
\be{gra29b}
(\alpha_W(t) +\gamma_W(t)) u^2_W(t) \leq (\alpha_W(t)
+\gamma_W(t)) u^2(t) \leq W^2(\alpha+\gamma) e^{2\alpha t}.\ee
 At the particle level, once a red particle is created
it begins to develop a growing cloud of red sites which grows with
exponential rate $\leq \alpha$.   Using the latter and (\ref{gra29b}) we obtain
\be{gra29b2}
\sum_{i_W,i_R,i_P}1_{i_R\geq 1} U_{WRP}(t, i_W,i_R,i_P)=
O(e^{2\alpha t}).
\ee

Return now to the production of the green particles. We will show
first that the number of green particles and sites with green
particles is of order $e^{2 \alpha t}$ as $t \to -\infty$.

We can assume that when a red particle is
created  on an occupied site the number of white particles is
given by the stable age distribution by (\ref{gra29z}). (We will verify below
in (\ref{gra17b2}) that in fact the error in (\ref{gra29z}) is $O(e^{\alpha t})$.)
We now want to study the production of sites with green particles.
We introduce the concept of special sites for this purpose.

We briefly return  to the finite $N$-system.
We call occupied sites (by a white particle) in the coloured particle system \be{xg31}
\mbox{ {\em ``special sites'' at time }}s \ee if $s$ lies between
the (random) time when a first  red particle arrives at this site until it
contains no red, purple or green particles. This random time is a.s. finite.

We transfer this concept to the $N \to \infty$ limit at time
$T_N+t$. This means that we incorporate the additional mark in the
measure-valued description.

During the  lifetime of a special site
a special site can produce red, purple  and green
migrants. Special sites are created at rate (recall  for the first
equality sign that $\alpha,
\gamma$ are the limits of $\alpha(t), \gamma(t)$ as $t \to -\infty$,
i.e. the ones characterizing the collision-free regime and then use
(\ref{ineq2b}) to get $o(1) O(e^{2 \alpha t})$ as bound for the
difference of both sides):
\be{gra20}
(\alpha_W(t) + \gamma_W(t)) u^2_W(t) =
(\alpha+\gamma)u^2(t) + o(e^{2 \alpha t})
= W^2 (\alpha + \gamma) e^{2\alpha t}
+o(e^{2\alpha t}) \ee and therefore by integration from $-\infty$
to $t$ the number of special sites created up to time $t$ is bounded below  by
\be{gra20b} W^2\frac{\alpha+\gamma}{2\alpha}e^{2\alpha
t}+o(e^{2\alpha t}). \ee

Then green particles at a special site are produced at rate
\be{gra29c}
{d}\suml_{i_W,i_R,i_B} i_W(i_R + i_P) \,U_{WRP}(t;i_W,i_R,i_P)\geq d
\suml_{i_W,i_R,i_P}1_{i_W(i_R+i_P)\geq 1}U_{WRP}(t;i_W,i_R,i_P).
\ee
Since $i_W (i_R+i_P) \geq 1$ automatically at a special site we use (\ref{gra20}) and (\ref{gra20b}),
to get that the expected number of green
particles at time $t$  is bounded below by  \be{gra29d}
d\cdot W^2\frac{\alpha +
\gamma}{2\alpha}e^{2\alpha t}+o(e^{2\alpha t}).
\ee
Provided the limit exists (see below), this implies that (recall the green particles
describe the loss of the WRP-system compared to the W-BL-system)
\be{gra20c}
\lim_{t\to -\infty}e^{-2\alpha t}[ (\alpha
+\gamma)We^{\alpha t}-(\alpha(t)+\gamma(t))u(t)]= \kappa^\ast >0.
\ee

In order to now get  also information on $\kappa$, we need information on sites
rather than on particle numbers. In order to get that the number of
sites with only green particles is of order $e^{2 \alpha t}$ as
$t \to -\infty$, we argue as follows.  Since
there is a positive probability that  a green particles migrates
before coalescing with a red or white, this implies
(if the limit exists, see below) that
\be{gra20c2}
\lim_{t\to -\infty}e^{-2\alpha t}[ (We^{\alpha t}-u(t)]=
\kappa
>0.\ee
Hence we now know (provided that the limits taken exist) that:
\be{gra20c3}
\kappa, \kappa^\ast >0.
\ee

Next,  in order to identify the constant $\kappa, \kappa^\ast$ we need to obtain
the 2nd order asymptotics.
It remains therefore to determine the actual value of
$\kappa, \kappa^\ast$ using this information.

We first consider the production of blue particles.  Since there are
$O(e^{2\alpha t})$ special sites and a blue particle is created
only when a migrant comes from a special site and hits the set
of size $u(t)$ of occupied sites, it follows that the number of
blue particles produced is of order $O(e^{3\alpha t})$ and
hence indeed the {\em blue} particles and therefore also the {\em purple}
particles play {\em no role} determining the
$e^{2\alpha t}$-term. In particular
the production of green particles by the loss of purple
particles is of order $O(e^{3\alpha t})$ and can be omitted.

We first obtain an expression for
$\kappa^\ast$. Since the green particles represent the particles lost
in the $WRP$-system due to coalescence of red or purple with white
particles, this is obtained by  considering the growth of the
green particles in more detail. We note  that the loss of purple
particles is of smaller order than the loss of the red particles
as mentioned above.
As a result, in order to to identify the constant $\kappa$ or
$\kappa^\ast$ in the expressions for $u(t)$ (number of sites) or $(\alpha(t)+\gamma(t))
u(t)$ (number of particles)  we can work with the WRGB-system instead of the WRP-system.
Hence in carrying out the analysis using the WRGB-system we obtain
a lower bound for the number of particles lost but as shown
above the resulting error is $O(e^{3\alpha t})$.

Once a new green particle is produced by a white-red coalescence
at a special site we are interested in the number of green
particles that migrate before the end of the life time of the special
site.  As noted above we can assume that when the red arrives the
number of white particles is given by the stable size
distribution. The number of green particles at the special site
and the process  of producing green migrants can then be obtained
by the analysis of a {\em modified birth and death} process where we now
have two types, namely, red and green with birth and death rules
inherited from the WRGB dynamics.  In terms of the limit process
this involves the forward Kolmogorov equations for this modified
birth and death process which serves as a source of migrants for
the green particle system which then evolves by the CMJ dynamics.

To summarize, in order to obtain the distribution of green mass up to an
error term of order $O(e^{3\alpha t})$ we consider
\begin{itemize}
\item
the WRGB-system instead of the WRP-system, ignore blue particles
\item  the production of special sites by red-white collision,
\item the production of green particles and green migrants at a
special site,
\item a CMJ process with immigration, with immigration source given by the
green migrants from special sites.
\end{itemize}

Recall from the derivation of a lower bound on $\kappa$ above
that in the $WRGB$-system  red always migrates to a new
(unoccupied) site or otherwise a purple-blue pair is created so that the number
of {\em white plus red and blue}  particles gives an {\em upper
bound} to the original interacting dual particle system in the
limit $N \to \infty$. The number of the {\em green particles}
produced only by white-red coalescence (omitting purple-white
coalescence) gives up to an error  of order $O(e^{3\alpha t})$ a
{\em lower bound} to $W(\alpha +\gamma) e^{\alpha
t}-(\alpha(t)+\gamma(t))u(t)$.

We make the following definitions in order to turn bounds into precise
asymptotics. We say below ''expected'', to distinguish from the usual expected value, when we calculate quantities
of the form
$\sum i_G U (t,i_W,i_R,i_G,i_B,i_P)$,
where the sum is overall $i_A, A \in C$.  Let
\be{gra9b1}
g_1(t,s)
\ee
be the ''expected'' number of green particles at time $t$
at a special site created at time $s$. Note that a special site has a finite lifetime (with finite expected value) since due to coalescence it will revert to a single white particle at some finite random time after its creation.  Therefore the function $g_1(t,s)=g_1(t-s)$ is bounded.

Let
\be{gra9b2}
g_2(s,r)dr
\ee
be the rate of production of green migrants at time $r$ at a special site
created at time, i.e. the $c$ times the ''expected'' number of green particles at special sites.
This function is also bounded as above.

Let
\be{xd3}
g_3(r,t)\leq \rm{const}\cdot e^{\alpha (t-r)}
\ee
be the ''expected'' number of green {\em sites} produced
at time $t$ from a founder at time $r$
 and finally
\be{gra9b3}
g^\ast_3(r,t)\leq \rm{const}\cdot e^{\alpha (t-r)}
\ee
is the ''expected'' number of green {\em particles} produced at time $t$ from a
founder at time $r$. These four functions determine the numbers
$\kappa$ and $\kappa^\ast$ uniquely as we shall see next.

Now the creation of the first green particles at a site and then subsequently green sites
occurs from a red-white  coalescence at a site. At such an event a
green particle arises and if there is no green particle yet a new green site
is created at this moment. From these founders now a cloud of green particles
and sites develops. The evolution of these clouds is independent of the
further development of the number of white particles and white sites. Therefore
conditioned on $W$ we get that the {\em ''expected'' number of green
particles at time $t$}, denoted $\kappa^\ast_t$ is asymptotically as $t \to -\infty$
\be{gra9b4}
W^2\kappa^\ast e^{2\alpha t} +O(e^{3\alpha t}),
\ee
where $\kappa^\ast = \liml_{t \to -\infty} \kappa^\ast_t$ and
 \be{kappa*}
 \kappa^\ast_t= e^{-2\alpha t}\int_{-\infty}^t  (\alpha + \gamma)
e^{2\alpha s}\left[ g_1(t,s)+ \left(\int_s^t
g_2(s,r)g_3^\ast(r,t)dr\right)\right] ds +O(e^{\alpha t}).
\ee
Since $g_1$ and $g_2$ are bounded, the integral is finite.
Observe that $g_1 (t,s) = \wt g_1(t-s)$ and
$g_2(s,r) = \wt g_2(r-s), g^\ast_3(r,t) = \wt g^\ast_3(t-r)$
by construction of the dynamic of the green particles, which do not coalesce
with white or red particles. Therefore the first term in
(\ref{kappa*}) is independent of $t$. The second part is given by
an integral depending again only on $t-s$ and hence the complete
term again does not depend on $t$. This implies the convergence of $\kappa^\ast_t$ as $t\to-\infty$.

We now turn  to the identification of $\kappa$ that is we turn from
particle numbers to number of sites. As in the
identification of $\kappa^\ast$ we can show that there exists
$\kappa$ such that
\be{xd1}
u_{WR}(t)=We^{\alpha t}-\kappa
W^2e^{2\alpha t}+O(e^{3\alpha t}).\ee
 We obtain the constant $\kappa$ by
counting {\em green sites}
(i.e. sites which are not also occupied by only red or white)
instead of green particles. Conditioning on
$W$ we get the {\em ''expected'' number of green sites},
at time $t$ is asymptotically as $t \to -\infty$:
\be{xd2}
W^2\kappa e^{2\alpha t} +O(e^{3\alpha t}),
\ee
where $\kappa = \liml_{t \to -\infty} \kappa_t$ with
\be{kappa}
\kappa_t = e^{-2\alpha t}\int_{-\infty}^t  (\alpha + \gamma)
e^{2\alpha s}\left[  \left(\int_s^t
g_2(s,r)g_3(r,t)dr\right)\right] ds +O(e^{3\alpha t}).
\ee
The existence of the limit follows as above.

{\bf Part 3 - Proof of (\ref{grocoll214}), (\ref{gra4d}), (\ref{gra6})}

Next we return to our original nonlinear equation, which
we now relate  with the quantities of our multicolour particle
system. We will need to collect in (\ref{eo2})-(\ref{eo4}) some
facts on this system used in the proof.

We first note the relation between the original dual and the WRP
system and the possibility to replace them with the WRGB-system
for the asymptotic as $t \to -\infty$, namely as $t \to -\infty$:

 \be{eo1} u(t)=u_{\mathrm{WRP}}(t) =
u_{\mathrm{WR}}(t)+O(e^{3\alpha t}), \ee

\bea{eo2} U(t,k)
&&
=   \suml_{i_W,i_R,i_P} 1_{i_W + i_R +i_P =k}
  U_{\mathrm{WRP} }(t, i_W,i_R,i_P)\\
&&
= \frac{u_{\mathrm{WRGB}}(t)}{u(t)}{\sum_{i_W,i_R,i_G,i_B} 1_{i_W+i_R
=k}U_{\mathrm{WRGB}}(t;i_W,i_R,i_G,i_B)} + o(1).\nonumber \eea
Here the first factor arises since the normalization (by the number of
sites) for the WRP
system is $u(t)$ and it is $u_{\mathrm{WRGB}}(t)$ for the WRGB system. We
also recall that by construction

\be{xd3b}  u_{\mathrm{WR}}\leq u = u_{\mathrm{WRP}}\leq u_{\mathrm{WRB}}.\ee

Recall furthermore that the sites occupied by black particles are in
one-to-one correspondence with the sites occupied by the red,
green and blue particles and the number of white plus black sites
is $W_t e^{\alpha t}$ with $W_t \to W$ as $t \to \infty$.
Therefore the difference in the number of occupied
sites in the (W-BL)-system and the WRGB-system arises from white
particles sharing a site with the coloured particles, but those sites are represented
as two sites in the (W-BL)-system. Hence

 \bea{xd3c} \frac{u_{\mathrm{WRGB}}(t)}{We^{\alpha t}}&&=\frac
{\sum_{i_W,i_R,i_G,i_B} [1_{i_W+i_R+i_G+i_B\geq 1}] U_{\mathrm{WRGB}}(t,i_W,i_R,i_G,i_B)}
{\sum_{i_W,i_R,i_G,i_B} [{1_{i_W\geq 1}+1_{i_R+i_G+i_B\geq 1}}] U_{\mathrm{WRGB}}(t,i_W,i_R,i_G,i_B)}.\eea

Furthermore the estimates on the number of green, red and purple particles
in the part 2 of our argument for Proposition \ref{P.Transition} imply as well
\be{xd3d} \frac{u_{\mathrm{WRGB}}(t)}{u_{\mathrm{WRP}}(t)}=1+\kappa We^{\alpha t}
+O(e^{2\alpha t})\ee
and
\be{eo3}
\alpha(t)=
\frac{u_{\mathrm{WRGB}}(t)}{u_{\mathrm{WRP}}(t)}(\alpha_W(t)+\alpha_R(t)) +O(e^{2\alpha t}),
\ee

\be{eo4}
\gamma(t) =
\frac{u_{\mathrm{WRGB}}(t)}{u_{\mathrm{WRP}}(t)}\sum_{i_W,i_R,i_G,i_B}
(1_{i_R+i_W=1}U_{\mathrm{WRGB}}(t,i_W,i_R,i_G,i_B)) +O(e^{2\alpha
t}). \ee

Using the equations (\ref{eo1}-\ref{eo4}) we shall  show below that (conditioned on $W$)
we have the following approximation relations of $t \to -\infty$ for
$u,\alpha,\gamma$ and $U$, which finish the proof of (\ref{grocoll214}),
(\ref{gra4d}) and (\ref{grocoll32}):
 \be{dif1} |u(t) -
We^{\alpha t}| = O(e^{2 \alpha t}), \ee
\be{dif2}
|\alpha (t) - \alpha| = O(e^{ \alpha t}),
\ee

\be{dif3}
|\gamma(t) -  \gamma| = O(e^{ \alpha t}),
\ee
\be{dif4}
\sum_{k=1}^\infty | \sum_{i_W,i_R,i_G,i_B}
1_{i_R+i_W=k}U_{\mathrm{WRGB}}(t,i_W,i_R,i_G,i_B)-\mathcal{U}(\infty,k))| =
O(e^{ \alpha t}).
\ee
It therefore remains now to verify (\ref{dif1})-(\ref{dif4}) in the sequel.

The bound (\ref{dif1}) follows from the fact that the difference between
the W-BL-system giving $W e^{\alpha t}$ and the WRP-system giving
$u(t)$ is bounded by the RGB-system which satifies
$u_{RGB}(t) = O(e^{2 \alpha t})$.

\medskip

Turn next to the proof of (\ref{dif2}). First
note that the bound of the WRP-system (for which $\alpha(t)$ is the
rate of colonization of new sites by multiple occupied sites)
from below and above
by the WR-system respectively the WRB-particles from the WRGB-system
yields by dividing the inequality by $u_{WRP}(=u)$:

\bea{xg33} && \frac{u_{\mathrm{WRGB}}(t)}{u_{\mathrm{WRP}}(t)} (\alpha_W (t)
+\alpha_R(t))\leq \alpha(t) \leq \frac{u_{\mathrm{WRGB}}(t)}{u_{\mathrm{WRP}}(t)}
(\alpha_W (t) +\alpha_R(t) +\alpha_B(t))\eea

The difference between the first and third expressions in
(\ref{xg33}) is $O(e^{2\alpha t})$ (bounding $\alpha_B(t)$ by
$O (e^{2\alpha t})$ and using (\ref{xd3d})). Therefore we get using the representation of
$\alpha(t)$ by the WRP-system (compare (\ref{eo3})):

\bea{xg33b}
\alpha(t) &= &
 \suml_{i_W, i_R, i_G,i_B} \left[
1_{i_W>1}\,(i_W) +1_{i_W=1,i_R>0}(i_W+i_R)
+1_{i_W=0,i_R>1}(i_R)\right]\\
&&\hspace{4cm} \cdot U_{\mathrm{WRGB}}(t,i_W,i_R,i_G,i_B)\cdot
\frac{u_{\mathrm{WRGB}}(t)}{u_{\mathrm{WRP}}(t)}+O(e^{2\alpha t}).\nonumber\eea

We can obtain the expression for $\alpha$, which is suitable for the comparison
with $\alpha (t)$, as follows. This constant $\alpha$ arises
from the (W-BL)-system. Furthermore the RGB-particles are in correspondence
with the black particles, only sit always on the copy of $\N$
instead of $\N$ or $\{1, \cdots,N\}$. (Note that we are interested here on particle numbers!)
Therefore since $W-BL$ is a CMJ-system in the stable age-type distribution
\be{xg33c}\begin{array}{ll}
\alpha =& u_{WRGB}(t) (W e^{\alpha t})^{-1}\nonumber\\[2ex]
&  (\suml_{i_W, i_R,i_G,i_B}
[1_{i_W>1} i_W + 1_{i_W=0} 1_{i_R+i_G+i_B>1} (i_R+i_G+i_B) \nonumber\\[2ex]
&
\phantom{(\suml_{i_W, i_R,i_G,i_B}}
+ 1_{i_W=1, i_R+i_G+i_B>0} (i_W+i_R+i_G+i_B) \nonumber\\[2ex]
&
\phantom{(\suml_{i_W, i_R,i_G,i_B}}
+ 1_{i_W>1} (i_R+i_G+i_B)] U_{WRGB} (t,i_W, i_R, i_G,i_B)).
\end{array}
\ee

Therefore we can represent $\alpha(t) = \alpha +\wt \alpha(t)$ and by
combining (\ref{xg33b}) and (\ref{xg33c}) together with the fact that blue
particles are $O(e^{3 \alpha t})$, we get an expression for
$\wt \alpha(t)$ which together with (\ref{xd3d}) results as $t \to -\infty$ in
\be{xg33d}
\alpha(t)= \alpha + O(e^{\alpha t}),
\ee
which proves (\ref{dif2}). (See Lemma \ref{L.dw18} below for the characterization of $\wt\alpha(t)$.)

Similarly we can proceed with the remaining claims (\ref{dif3}) and
(\ref{dif4}).
This completes the proof of  (\ref{grocoll214}), (\ref{gra4d}) and (\ref{gra6})
as pointed out below (\ref{eo4}).
\sm

{\bf Part 4 - Proof of (\ref{gra4e}) and (\ref{gra4e2})}

 For this purpose we must investigate the difference between
 $\alpha(t),\;U(t),\;\gamma(t)$ and $\alpha,\;\mathcal{U}(\infty),\;\gamma$, respectively
 more accurately than in the bounds above.
Namely we need an expression with error terms of order
$O(e^{3\alpha t})$ resp. $O(e^{2\alpha t})$ in (\ref{dif1}) resp.
(\ref{dif2})-(\ref{dif4}). For that purpose we use again the multicolour
representation.
We begin by establishing the order of these differences
rather than only upper bounds. This is based on
the coloured particle system.  From  the above discussion we
expect that $\wt\alpha(t)=\alpha(t)-\alpha\sim \wt\alpha_0
e^{\alpha t}$ for some $\wt\alpha_0>0$ which we have to identify.

Recall that the intensity of green particles at time $t$
corresponds to \be{gra34b} (\alpha +\gamma) e^{\alpha t}-
(\alpha(t)+\gamma(t))u(t) \ee and the number of green particles is
given $W$ of order $O(e^{2\alpha t})$ so that conditioned on $W$ we have \bea{gra34a}
&&(\alpha +\gamma)W e^{\alpha t}-(\alpha(t)+\gamma(t))u(t) \geq
\delta e^{2\alpha t},\quad -\infty<t\leq t_0,\\&& \text{ for some
positive constant } \delta.\nonumber\eea

We want to sharpen this inequality above to a precise second order expansion.
To do this and to thereby identifying $\wt\alpha_0,\;\wt\gamma_0$ we now return to the
analytical study of the nonlinear system (\ref{gra15}),
(\ref{gra16}) and prove the following.

 \beL{L.dw18} {(Expansion of $\alpha(t)$ and $\gamma(t)$)}

(a) Let
\be{gra34c}
 \wt \alpha(t):= \alpha(t)-\alpha, \quad \wt\gamma(t)=\gamma(t)-\gamma.
 \ee
 Then as  $t \to -\infty$ we have
\be{dw18} \wt\alpha(t) =  \wt\alpha_0 e^{\alpha t} +O(e^{2\alpha
t}), \ee
\be{dw188}
\wt\alpha_0>0  \text{ is
given explicitly by (\ref{gra39}) and is of the form $\vec{Const} \cdot W$. }
\ee
Moreover,
\be{dw18b} \wt U(t) = U(t)-\mathcal{U}(\infty)= \wt U_0 e^{\alpha t}
+O(e^{2\alpha t}),\ee where
$\wt U_0$ is given explicitly in (\ref{U0zz}) and is of the form $Const \cdot W$.
Then $\wt \gamma_0= \wt U_0(1)$.

(b) Define $ u^\ast$ as the solution to the nonlinear equation
(\ref{gra15}) with initial condition at $-\infty$
\be{gra34c3}
 \lim_{t\to
-\infty}e^{-\alpha t}u^\ast(t)=1.\ee

Then as $t\to -\infty$ we have the second order expansion:
\be{gra34d}
u^\ast(t)= e^{\alpha t}-(b-\frac{\wt\alpha^\ast_0}{\alpha})e^{2\alpha t}
+O(e^{3\alpha t}),
\ee
where $\wt \alpha^\ast_0$ is obtained from $\wt \alpha_0$ by setting $W=1$
in the formula  (\ref{gra39}). $\qquad \square$

\end{lemma}

With these results we obtain immediately (\ref{grocoll32}). This
would complete the Proof of Proposition \ref{P.Transition}.

It remains to prove  Lemma \ref{L.dw18}. To prove this we  need
some preparation  we do in the next part.

{\bf Part 5 - Statement of proof of Lemma \ref{Ltransu} }

We and formulate  a statement that gives an
explicit representation of $u^\ast$ (and hence of the limit of
$u^N$ suitably shifted in terms of $\alpha(\cdot)$ which is exact
up to the third order error terms:

\beL{Ltransu} {(Identification of $u$ in terms of $\wh u$)}

Consider the  solution to the ODE

\be{ga30--}
 \frac{d u^\ast(t)}{dt}=\alpha(t) [u^\ast(t)-b(t) (u^\ast (t))^2],\;\;-\infty< t<\infty,
 \ee
with boundary condition at $-\infty$ given by
$\liml_{t\to -\infty}e^{-\alpha t}u^\ast(t)=1$
and with the abbreviations
\be{ga30--2}
b(t)= 1+\frac{\gamma(t)}{\alpha(t)}.
\ee
(Recall that the general case is
obtained by a time shift by $\frac{\log W}{\alpha}$ if 1 is
replaced by $W$ and $\alpha(\cdot), b(\cdot)$ are shifted accordingly.)

Define with $b=1 + (\gamma/\alpha)$ the function $\wh u$ on $\R$:

\be{transu7+} \wh u(t)=
\frac{e^{\alpha
t}e^{\int_{-\infty}^t(\alpha(s)-\alpha)ds}}{1+be^{\alpha
t}e^{\int_{-\infty}^t(\alpha(s)-\alpha)ds}}, \quad
-\infty<t<\infty.\ee

(a) Then \be{gra36} |u^\ast(t) - \wh u(t)| =O(e^{3\alpha t}).
\ee

The function $\wh u$ satisfies (here again $b=1+\frac{\gamma}{\alpha}$)
\be{transu7++}
\wh u(t) \sim \frac{e^{\alpha
t}e^{\int_{-\infty}^t(\alpha(s)-\alpha)ds}}{1+be^{\alpha t}}
\mbox{ as } t \to -\infty.
  \ee

(b) We have  for $T_N=\alpha^{-1} \log N$ that in distribution:
\be{agrev71} \frac{1}{N}  u^{N}(T_N+t) \Nto u^\ast(t+\frac{\log
W}{\alpha}) = \wh u(t+\frac{\log W}{\alpha}) +O(e^{3\alpha t}).
\qquad \square \ee
\end{lemma}

{\bf Part 6 - Proof of Lemma \ref{Ltransu} }

{\bf (a)} We start with the following observation.
Since by the coloured particle calculation we obtained earlier in this proof, see
(\ref{dif2}), (\ref{dif3}),
\be{cpp}
|\gamma(t)-\gamma|+|\alpha(t)-\alpha|\leq
\text{const}\cdot e^{\alpha t}\ee and $\alpha >0$, then
\be{grevz5} |b(t)-b|=|b(t)-(1+\frac{\gamma}{\alpha})| \leq
\text{const}\cdot e^{\alpha t}. \ee

When $b(\cdot)$ is not constant we cannot obtain a closed form
solution of (\ref{ga30--}). However we can obtain an approximation
that describes the asymptotics as $t\to -\infty$ accurate up to
terms of order $O(e^{2\alpha t})$ as follows.

Let  \be{agrev70} (\wh u_{t_0}(t))_{t \in \R} \ee denote the
solution of the modified equation (\ref{ga30--}) where we replace our {\em function $b(t)$ by
the constant $b$} and put the initial
condition $u_{t_0}$ at time $t_0$.

The solution of the modified equation is given by the
formula
\be{ga30---}
\wh u_{t_0}(t)
= \frac{\wh u_{t_0}(t_0)e^{\int_{t_0}^t\alpha(s)ds}}{1+b \wh
u_{t_0}(t_0)(e^{\int_{t_0}^t\alpha(s)ds}-1)},\quad t\geq t_0. \ee

Now let
\be{ga30-5}
u_{t_0}=e^{\alpha t_0},
\ee
so that
\be{ga30----} \wh
u_{t_0}(t) = \frac{e^{\alpha t}e^{\int_{t_0}^t(\alpha(s)-\alpha)ds}}{1+b  e^{\alpha
t}e^{\int_{t_0}^t(\alpha(s)-\alpha)ds}-b e^{\alpha t_0}},\quad t\geq
t_0. \ee

We then get with (\ref{transu7++}) that
\be{agrev70b}
\wh u(t)=\lim_{t_0\to -\infty} \wh u_{t_0}(t)=
\frac{e^{\alpha t}e^{\int_{-\infty}^t(\alpha(s)-\alpha)ds}}{1+be^{\alpha
t}e^{\int_{-\infty}^t(\alpha(s)-\alpha)ds}}, \quad
-\infty<t<\infty,
\ee
noting that the integrals are well defined by (\ref{cpp}).
Hence $\wh u$ satisfies a differential equation
((\ref{ga30--}) with $b(t) \equiv b$) which we now use
to estimate $\wh u-u^\ast$.

Let  $v(t):=(\wh u(t)-u^\ast(t))$.  Then $v(t)$ satisfies
\be{gra37}
  \frac{dv(t)}{dt}=\wt \alpha(t)v(t)-b(\wh
u(t)+u^\ast(t))v(t)+(b(t)-b)(u^\ast(t))^2.
\ee

Using  (\ref{cpp}), (\ref{grevz5}) and the fact that
$u^\ast(t)\leq \text {const} \, \cdot e^{\alpha t}$ we then obtain
that \be{gra37b} \frac{dv(t)}{dt} \leq \text{const}\cdot e^{\alpha
t} v(t) +O(e^{3 \alpha t}). \ee
Therefore $|v(t)|\leq
\text{const}\cdot e^{3\alpha t}$ and hence

\be{agrev74} |  u^\ast(t) - \wh u(t)| \leq \text{const}\cdot e^{3\alpha
t} \mbox{ as } t \to -\infty. \ee

We conclude with an estimate for the integral term in $\wh u$
\bea{transu11} &&
(\alpha(t)+\gamma(t))\wh u (t)
\\&& = \frac{(\alpha(t)+\gamma(t)) e^{\alpha t}e^{\int_{-\infty}^{t}
(\alpha(s)-\alpha)ds}}{1+b(e^{\alpha t}e^{\int_{-\infty}^{t}
(\alpha(s)-\alpha)ds})}.\nonumber \eea
Furthermore by (\ref{gra34a}),
\be{transu11b} (\alpha (t) + \gamma (t)) u^\ast(t) \leq
(\alpha +\gamma) e^{\alpha t}. \ee

Then solving in (\ref{transu11}) for the $e^{\int}$-term we obtain
with (\ref{gra36}) for $t \to -\infty$ the relation: \be{transu12} I(t)=
e^{\int_{-\infty}^{t}(\alpha(s)-\alpha)ds}\leq \frac{1}{1-b
e^{\alpha t}}= 1+be^{\alpha t}+o(e^{\alpha t}).\nonumber \ee

{\bf (b)} The second assertion follows from part (a). The first assertion
follows from (\ref{grocoll1}).

\noindent This completes the proof of Lemma \ref{Ltransu}.
\end{proof}

\bigskip

{\bf Part 6 - Proof of Lemma \ref{L.dw18}}

We separately show the different claimed relations in points
(1)-(3) for $\wt \alpha, \wt \gamma$ and then $u^\ast$.\\

{\em (1) Relation (\ref{dw18}) and (\ref{dw188})}

Here our strategy is to express $\wt \alpha (\cdot)$ in terms of the
rewritten $(u,U)$ using the Markov process generated by
$Q^\alpha_0$ and the operators $Q_1$ and $L$.
Using (\ref{gra17}), and (\ref{dw19})   we can rewrite the
nonlinear system defining $(u,U)$ in the following form.
Let $U(t), Q_0^a,Q_1,L$ be defined as
in (\ref{gra16}), (\ref{agdw19a}), (\ref{gra12}) and  define $\nu$
by \be{dw19a} \nu(j)=j1_{j>1}. \ee
Furthermore recall that by setting $a=\alpha$ we get the following:
\be{gra38b}
(S_t^\alpha)_{t \geq 0}  = \mbox{{semigroup with generator }}
Q^\alpha=Q^\alpha_0+\alpha1_{j=1}. \ee This semigroup has a unique
entrance law from $t \to -\infty$, since it is standard to verify
that $S^\alpha$ defines an ergodic Markov process.

Next let $(U^\alpha(t))_{t \in \R}$ solve the (forward) equation for
generator $Q^\alpha$, i.e. \be{gra38c2} \frac{\partial}{\partial
t} U^\alpha (t) = U^\alpha(t)Q^\alpha, \quad t \in \R, \ee
and define the difference process
$(\wt U(t))_{t \in \R}$  (we now suppress the superscript
$\alpha$ in $\wt U$):
\be{gra38d3} \wt U(t)=U(t)-U^\alpha(t). \ee
Then
\be{gra17b} \frac{\partial \wt U(t)}{\partial t} = \wt
U(t)Q_0^{\alpha} +\wt\alpha(t)[U(t)Q_1
+1_{j=1}]+u(t)(\alpha(t)+\gamma(t))[U(t)L-1_{j=1}]. \ee

We can now represent due to the definition of $
\alpha (\cdot), U,U^\alpha, \wt \alpha (\cdot)$
the function $\alpha(t)$ as
\be{xg36}
\alpha(t) -\alpha = <c \wt U(t), \nu>.
\ee
Therefore (recall (\ref{dw19a})) by the formula of partial integration
for semigroups (for $\nu$ see (\ref{dw19a}) with the semigroup
$S^\alpha$ of $U^\alpha$ as the reference semigroup and then
$U$ as the wanted object:

\bea{ba10g9++}
&&\wt \alpha(t)=\alpha(t)-\alpha=\langle c\wt U(t),\nu\rangle \\
&& = \int_{-\infty}^tc[ \wt\alpha (s) \langle
(U(s)Q_1+1_{j=1})S^\alpha_{t-s},\nu\rangle+
u(s)(\alpha+\gamma)\langle (U(s)
L-1_{j=1})S^\alpha_{t-s},\nu\rangle] ds.
\nonumber\eea
Next set
\be{gra38d}
  \CU(\infty)=(p_1,p_2,\dots)
  \ee
and
note that by explicit calculation (see (\ref{gra12}), (\ref{agdw19b})):
\be{gra38c}
 \CU(\infty)Q_1+1_{j=1}= (1-p_1,-p_2,-p_3,\dots),\quad \CU(\infty)
L-1_{j=1}=(-1,p_1,p_2,\dots).
\ee

Note that for $t \to -\infty$: \be{gra17b2} \lim_{t\to -\infty}
U(t)=\CU(\infty),\qquad \|{\CU(\infty)-U(t)}\|_1=O(e^{\alpha t}),
\ee where $\CU(\infty)$ is the stable size distribution of the
McKean-Vlasow dual process where the norm $\|\cdot\|_1$ is as in
(\ref{Grev39n}).
Hence we get from (\ref{ba10g9++}) that \bea{gra38da} \wt \alpha
(t) && = \int_{-\infty}^tc[ \wt\alpha (s) \langle
(\CU(\infty)Q_1+1_{j=1})S^\alpha_{t-s},\nu\rangle\nonumber\\&&+
u(s)(\alpha+\gamma)\langle (\CU(\infty)
L-1_{j=1})S^\alpha_{t-s},\nu\rangle] ds \nonumber\\
&& +O(e^{2\alpha t}). \nonumber \eea

Next note that (recall (\ref{dw19a}) and (\ref{gra38c}), (\ref{gra38d})
and the fact that the Markov process for $S^\alpha$ increases from
the initial value 1 stochastically to its equilibrium) we have the relations:
\be{gra38d2}
 -\alpha<\langle
(\CU(\infty)Q_1+1_{j=1})S^\alpha_t,\nu\rangle <0 \mbox{ and } \langle
(\CU(\infty)L-1_{j=1})S^\alpha_t,\nu\rangle >0.
\ee
Recall furthermore that
\be{gra38d2a}
\langle \CU(\infty)S^\alpha_t,\nu\rangle =
\frac{\alpha}{c} \mbox{ for all } t \geq 0.
\ee

Let now (use  (\ref{gra38d2}) for possitivity):
\be{gra38e}
A_1= c \int_0^\infty \langle
(\CU(\infty)L-1_{j=1}) e^{-\alpha s}S^\alpha_{s},\nu\rangle ds >0,
\ee \be{gra38f} A_2 =-c\int_0^\infty e^{-\alpha s}\langle
(\CU(\infty)Q_1+1_{j=1})S^\alpha_{s},\nu\rangle ds >0.
\ee
Now multiply through both sides of (\ref{ba10g9++}) by $\alpha
e^{-\alpha t}$ and set
\be{xg38}
\wh \alpha(t) = e^{-\alpha t} \wt \alpha(t),
\quad \wh S^\alpha_t = e^{-\alpha t} S^\alpha_t. \ee We
obtain the equation:
\bea{xg39}
\wh \alpha(t)& = &
c\int^t_{-\infty}[ \wh \alpha(s) \langle (\CU(\infty) Q_1 + 1_{j=1}) \wh
S^\alpha_{t-s}, \nu \rangle\\
&&\hspace{4cm}
+\wh u(s)(\alpha + \gamma) \langle (\CU(\infty) L-1_{j=1}) \wh
S^\alpha_{t-s},\nu\rangle ]ds.\nonumber
\eea
This is a {\em renewal-type  equation} of the
form (choosing suitable positive functions $f,g$, recall (\ref{gra38d2})):
\be{xg40} \wh \alpha = \wh \alpha \ast (- f) + \wh u \ast g,
\ee
where the integral over $f$ over $\R$ is equal to $A_2$, which is less than 1
due to relation (\ref{gra38d2a}) and (\ref{xg38}), furthermore
$\wh u(t) \to W$ as $t \to -\infty$ and $g$ is an integrable
positive function with integral $A_1$.

We now claim that as
$t \to -\infty$
\be{gra39b}
\wh \alpha (t) \la C \mbox{ and } |\wh \alpha (t) -C| = O(e^{\alpha t}),
\ee
where $C$ is calculated as usual in renewal theory.

Then the solution (uniqueness is verified  below) to
(\ref{xg40}) is given by \be{agrev76} \wt\alpha(t)=\wt\alpha_0
e^{\alpha t}+O(e^{2\alpha t}), \ee where $\wt\alpha_0>0$ is given
by

\be{gra39} \wt\alpha_0=\frac{A_1}{1+A_2} W.\ee

The fact that the errorterm in (\ref{agrev76}) is of the form $O(e^{2 \alpha t})$ follows
from (\ref{gra38da}).
Moreover (\ref{agrev76}) is the unique solution of order
$O(e^{\alpha t})$ to (\ref{ba10g9++}). To verify the uniqueness we
argue as follows.

The difference of two solutions $h_1,h_2$ must solve according to
(\ref{xg40}) that
\be{agr24}
(h_1-h_2)=(h_1-h_2) \ast (-f).
\ee
We can then verify that either $h_1 \equiv h_2$ or
\be{gra39a}
|(h_1-h_2) (t)| < \supl_{(-\infty,0)} |h_1(t)-h_2(t)|
\ee
and therefore indeed $h_1-h_2 \equiv 0$.

This completes the proof of (\ref{dw18}) and (\ref{dw188}).

\begin{remark} We have with (\ref{ba10d}) and (\ref{ba10g9a})
the following ODE for $\alpha(t)$ as function of $(u,U)$:
\bea{dw19}
\frac{\partial \alpha(t)}{\partial t}=
&&\frac{\partial}{\partial t}
   [c\sum_{j=2}^\infty jU(t,j)] = c \sum^\infty_{j=2} j(\frac{\partial}{\partial t} U(t,j))\\
  = && c \sum_{j=2}^\infty j   \Big\{\left [\alpha(t)(1-u(t))- \gamma(t)u(t)\right ]1_{j=1} \nonumber
  \\
  &&+ u(t){(\alpha(t)+\gamma(t))}[1_{j\ne1}U(t,j-1)  - U(t,j)]\nonumber
\\
  && +c(j+1)U(t,j+1)-cjU(t,j)1_{j\ne 1}\nonumber
 \\
  && +s(j-1)1_{j\ne1}U(t,j-1)-sjU(t,j)\nonumber
\\
  &&+\frac{d}{2}(j+1)jU(t,j+1)-\frac{d}{2}j(j-1))1_{j\ne1}U(t,j)\nonumber
  \\
  &&  -\Big(\alpha(t)(1-b(t)u(t)) \Big)\cdot U(t,j)\Big\}. \nonumber
\eea

 Collecting terms we obtain:
 \bea{dw19b}
 \frac{\partial \alpha(t)}{\partial t}
  && =  (s-\alpha(t)+\frac{d}{2})(\alpha(t)+\gamma(t))
  -\frac{cd}{2}\sum_{j=1}^\infty j^2U(t,j)
\\ && +
u(t)(\alpha(t)+\gamma(t))^2+cu(t)(\alpha(t)+\gamma(t)).\nonumber
 \eea
Hence we can determine $\alpha$ as solution of an ODE once we are given
$(u,U)$.
\end{remark}

\bigskip

{\em (2) Relations (\ref{dw18b}) and (\ref{gra4e2})}.

Recall (\ref{gra17b}) which yields

\be{gra17bz} \frac{\partial \wt U(t)}{\partial t} = \wt
U(t)Q_0^{\alpha} +\wt\alpha(t)[\mathcal{U}(\infty)Q_1
+1_{j=1}]+u(t)(\alpha+\gamma)[\mathcal{U}(\infty)L-1_{j=1}] + O( e^{2\alpha t}). \ee
Then using the expression (\ref{agrev76}) for $\wt\alpha(t)$ and $u(t)=We^{\alpha t}+O(e^{2\alpha t})$ we obtain
\be{U0zz} \wt U_0= \int_{-\infty}^0 e^{\alpha s} S^\alpha_{-s}Vds  +O(e^{2\alpha t})
\ee
where $V=W(\wt\alpha^\ast_0[\mathcal{U}(\infty)Q_1
+1_{j=1}]+(\alpha+\gamma)[\mathcal{U}(\infty)L-1_{j=1}])$.  Recalling that $S^\alpha_t$ is a  semigroup on $\mathcal{P}(\N)$, it follows that the integral is well-defined.


{\em (3) Relation (\ref{gra34d}).}

With the knowledge we have now the expression
(\ref{gra34d}) then follows by (\ref{gra36}) by substitution with
(\ref{agrev76}) in (\ref{transu7++}).

\bi

\subsubsection{Weighted occupation time for the dual process}\label{sss.limrandmarg}

In this subsubsection  we focus on the behaviour  of the quantity
$\Xi_N (\alpha^{-1} \log N +t,2)$ in the limit $N \to \infty$ as a
function of $t$. We know that this limit is $\CL_t(2)$ which
is a random
probability measure on $[0,1]$, in order to identify its
distribution it suffices to compute the moments via duality and
(\ref{ba7}). Recall also that the weighted occupation time of the
dual determines the probability that no mutation occurred
changing the value at the tagged site from type 1 to 2. Why,
what do we have to prove about the dual process?
This we now explain first.

Consider $T_N= \alpha^{-1} \log N$ and $t_0(N) \uparrow \infty, \quad t_0(N)
=o(\log N)$. Recall that

\be{xg35a} \Pi^{N,k,l}_s= \sum_{i_w,i_R,i_P}(i_W+i_R+i_P)
\Psi_s^N(i_W,i_R,i_P),
\ee
where $\Pi^{N,k,l}_s$ is given by the
dual particle system started with $k$ particles at each of $\ell$ sites.

Then we know that we only have to replace $W$ by $W^{(k,\ell)}$
and otherwise we get the same equation as we had for $k=\ell=1$. Hence we have
\be{ag32}
\int_{t_0(N)}^{(t+T_N)}[N^{-1}
\Pi^{N,k,\ell}_{s}]ds \sim \frac{1}{c}  \int_{-\infty}^t
(\alpha^{(k,\ell)} (s)+\gamma^{(k,\ell)}(s))u^{(k,\ell)}(s)ds,
\mbox{ as } N \to \infty \ee and we
need the r.h.s. in the limit as $t\to -\infty$. Integrating the
expression given by (\ref{grocoll32}) we get the following expansion:
\be{ag32+}
\int_{-\infty}^t (\alpha(k,\ell)(s)+\gamma (k,\ell)(s))u (k,\ell)(s)ds
=W^{(k,\ell)}
\frac{(\alpha+\gamma)}{\alpha} e^{\alpha
t}-\frac{\kappa^\ast}{2\alpha}(W^{(k,\ell)})^2e^{2\alpha t}+O(e^{3\alpha
t}).\ee

Recall that in terms of the {\em normalized} solution of the
equation for the dual $(u^\ast,U^\ast)$, (i.e. $e^{\alpha t}
u^\ast (t) \to 1$ as $ t \to -\infty$) we have by the above on the
one hand
\be{ag33}\begin{array}{l}
\CL \left\{[\frac{1}{N} \intl^{T_N+t}_{t_0(N)}
\Pi^{N,k,\ell}_u du]\right\}
\Nto \CL\left\{\int_{-\infty}^{\,t+\frac{\log
W^{k,\ell}}{\alpha}} \frac{1}{c}(\alpha (s)+\gamma (s))u^\ast(s)ds\right\}\\
\hspace{4cm}
=: \nu_{k,\ell}(t)\in \mathcal{P}([0,\infty))
\end{array}
\ee
and on the other hand in terms of the
original process we have (cf. (\ref{angr15a}), (\ref{dualmoment})),
\be{dd12b}
 E\left(\left[\int_0^1 x^k\mathcal{L}_t(2,dx)\right]^\ell\right)
 = \int_0^\infty e^{-my}\nu_{k,\ell}(t,dy) \mbox{ for all }
 k,\ell\in \mathbb{N}.
 \ee

 Note that the random time shift by $\frac{\log W^{k,\ell}}{\alpha}$ plays
 an essential role.  The interplay between this shift and the
 nonlinearity of $u(t)$ determines the distribution of $\CL_t (2)$.
To determine this effect we  study the properties of the r.h.s. in
the two equations above. In particular  we want to show  that
$\CL_t(2)$ is neither
 $\delta_0$ nor $\delta_1$ and we want to show in fact that it is truly random
 and has its mass on $(0,1)$.
This can be  translated into  properties of the r.h.s. of
(\ref{dd12b}) which we will study using the dual.
\bi

We next  collect some key facts needed to calculate the
probability of mutation jumps as $N \to \infty$ which appear on
the r.h.s. of (\ref{dd12b}).

\beL{L6.1}{(Properties limiting dual occupation density)}

 (a)  If  $c>0$
the limit object $\nu_{k,\ell}$ satisfies: \be{ba8} \nu_{k,\ell}
(t,({0,\infty}))=1,\ee and \be{ba88}
\lim_{n\to\infty}\nu_{n,1}(t,[K,\infty))= 1, \mbox{ for all } K.
\ee

(b) Denote by $<$ strict stochastic order of probability measures.
Then for $c>0, d>0$: \be{ba8b} \nu_{k+1,1} > \nu_{k,1}\ee \be{ba9}
\nu_{1,2}(t) < \nu_{1,1}(t) \star \nu_{1,1}(t), \ee where $\star$
denotes convolution.
\smallskip

\noi (c) Let $n_0 \in \N$ and set
\be{ba9b}
t_{n_0}(\ve,N):=
\inf\{t:\Pi^{N,n_0,1}_t= \lfloor \ve N \rfloor \} -\frac{1}{\alpha} \log N.
\ee
Then every weak limit point (in law) as $N \to \infty$, denoted
$t_{n_0}(\ve)$ satisfies: \be{angr22} t_{n_0}(\ve) \to -\infty
\mbox{ in probability as }  n_0 \to \infty. \ee Furthermore
\be{ba10neu} t_{n_0} (\ve) \longrightarrow -\infty \mbox{ in
probability as } \ve \to 0. \qquad \square \ee
\end{lemma}

\bigskip

\noi {\bf Proof of Lemma \ref{L6.1}}\\

\emph{Proof of (a)}.   \quad  The claim $\nu_{k,\ell}(t,\{0\}) =0$ follows
from (\ref{grocoll3}) since $W^{k,\ell}>0$ a.s.
and
$\nu_{k,\ell}(t,(\{\infty\})) =0 $ since for $T_N =\alpha^{-1} \log N$
and every $t \in \R$ we have:
\be{ba10n2}
\ \limsup_{N\to\infty}E\left[  \frac{1}{N}
\int_0^{T_N+t}\Pi^{N,k,\ell}_udu\right]<\infty.\ee Therefore
$\nu_{k,\ell}$ is  a probability measure concentrated on
$(0,\infty)$.

The relation (\ref{ba88}) will follow from the argument given for
the proof of (\ref{angr22}) in (c).  \bi

\emph{Proof of (b).} The first part of (b) is immediate by a
coupling argument realising the $(k+1)$ and $k$ particle system
on one probability space. The second part follows from the fact two
typical individuals picked among the descendants of the two
populations interact through coalescence with positive probability
before they jump since once they both occupy $O(N)$ sites, so that we have
a fraction of both populations overlapping. This proves (b). \bi

\emph{Proof of (c).} The main idea in the proof of (c) is that,
asymptotically as $n_0\to\infty$, starting $n_0$ particles at a
tagged site the number of particles at this site decreases to
$O(1)$ in time $1$ but during this time still produces $O(\log
n_0)$ migrants as we shall see below. Therefore if the number of
initial particles, $n_0$ increases to infinity, then an increasing
number of populations start growing like $W_i e^{\alpha t}$ and
they do so independently. Hence writing the total population as
\be{ba12.1} \suml^{m(n_0)}_{i=1}e^{ \alpha(t+\alpha^{-1} \log
W_i)}, \ee we see that if $m(n_0)$ diverges as $n_0 \to \infty$,
then the total number of particles reaches $O(N)$ at a time
$(\frac{1}{\alpha}\log N+t_{n_0})$ where $t_{n_0}$ decreases as
$n_0$ increases to $+\infty$. Hence it remains to give the formal argument for the fact
that $m(n_0)$ is of order $ \log n_0$ as $n_0 \to \infty$, which
runs as follows.

Consider Kingman's coalescent $(C_t)_{t \geq 0}$ starting with
countably many particles . Consider the number of particles in a
spatial Kingman coalescent which jump before  coalescing at the
starting site. This is the given via the rate of divergence of the
entrance law of Kingman's coalescence from 0. Hence we need the
rate of divergence  \be{ba10e} \intl^1_\delta |C_t|dt \sim
|\log\delta | \mbox{ as } \delta \to 0, \ee which follows  since
the rate of divergence is wellknown to be $\delta^{-1}$. We
connect this to an initial state with $n_0$ particles which is
obtained by restricting the entrance law of Kingman's coalescent.
We then see that  the number of migration jumps diverges as $n_0$
tends to infinity.

To get (\ref{ba10neu}) we use that $t_{n_0} (\ve,N)$ can be
bounded from above by the time where
$\Pi^{\infty,n_0,1}$ reached first $\lfloor \ve N \rfloor$.
By the analysis of $(K_t, \zeta_t)$ we know that for this it
suffices to show that $\wt \tau^{N, n_0}(\ve) \la -\infty$
as $\ve \to 0$, which follows from the exponential growth of
$K_t$ immediately.

This completes the proof of the lemma.

\subsubsection{Proof of Proposition \ref{P.measure}, Part 1: Convergence to limiting dynamics}
\label{sss.proofldy}

In the proof of Proposition \ref{P6.1b} it sufficed to work with
the expected value $E[x^N_1(\frac{1}{\alpha} \log N+t)]$ (since
$0\leq x_1 \leq 1$). However in order to prove Proposition
\ref{P.measure} it is necessary to work as well with higher moments in
order to identify not only the time of emergence but to determine
the dynamic of fixation. The strategy is to carry out the
following two steps next:
\begin{enumerate}
\item show that weak convergence of $(\Xi^{\log,\alpha}_N(t))_{t
\geq t_0}$ to the limiting deterministic McKean-Vlasov dynamics
for $t\geq t_0$ follows from the weak convergence of the marginal
distributions at $t_0$,
\item prove that the one dimensional
marginals $(\Xi^{\log, \alpha}_N(t_0))$ converge weakly to a random (not
deterministic)  probability measure on $[0,1]$.
\end{enumerate}

{\it Step~1 $\;$ Reduction to convergence of one-dimensional
marginals. }\\In order to prove that $(\Xi^{\log, \alpha}_N(t))_{t
\in \R}$ converges weakly to a random solution of the
McKean-Vlasov equation we argue first that it suffices to prove
that the one-dimensional marginals converge. Namely if one has
that, then one can use the Skorohod representation to get a.s.
convergence on some joint probability space. Then however we can use
that (using duality to obtain a Feller property)
\be{angr17d1}
\Xi^{\log,\alpha}_N(t_0)\to \mu_0, \ee implies
$\{\Xi^{\log,\alpha}_N(t)\}_{t_0\leq t\leq t_0+T}$ converges to a
solution of the McKean Vlasov equation (\ref{angr7b}) with initial
value $\mu_0$.

This follows by  noting that we claim here a standard McKean-Vlasov
limit of an exchangeable system satisfying the martingale problem
(\ref{angr2}), (\ref{angr3}) and with initial empirical measures
converging. The proof is a modification of the standard proof of
the McKean-Vlasov limit (for details see proof of Theorem 9,
\cite{DG99}). \sm

\sm

{\it Step~2 $\;$ Convergence of marginals: reformulation in terms of
the dual.}\label{marginalconvergence}\\ We now prove the required
convergence of the one-dimensional marginal distributions. To do
this we recall first that moments determine probability measures
on $[0,1]$ and that therefore the collection of ``moments of
moments'' of a random measure denoted $X$ on $[0,1]$ given by
\be{angr17d2}
E[\prod_{i=1}^m \int x^{k_i}X(dx)],\; k_i,m\in \mathbb{N},
\ee
determine the law. Compactness is automatic so
therefore it suffices to verify that for all $k,\ell\in\mathbb{N}$
the $(k,\ell)$ moments of the mass of type 1 \be{angr18}
 m^N_{k,\ell}(t)=E([\int_{[0,1]}x^k \Xi^{\log,\alpha}_N
 (t,1)(dx) ]^\ell),\quad\forall k,\ell\in \mathbb{N},
 \ee
converge as $N \to \infty$ to conclude  weak convergence of the
one-dimensional marginal distributions.

In order to analyse the empirical measure and its functionals we
return to the original system on the site space $\{1,2,\cdots,N\}$ and express
the quantity through moments of observables of the system. These
moments are obtained in terms of the dual process by considering
the dual with initial function generated by taking the $k$-product
given by $1_{\{1\}}\otimes\dots\otimes 1_{\{1\}}$ at
$\ell$  not necessarily distinct sites.

This representation can be simplified a bit. Namely note  that for
bounded exchangeable random variables $(X_1, \cdots,X_N)$ one has:
\be{gan20}
 \lim_{N\to\infty} E[(\frac{1}{N}\sum_{i=1}^{N}X^k_i)^\ell] =
E[\prod_{i_1 \ne i_2\ne\dots \ne i_\ell} X^k_{i_j}].
\ee
Therefore to compute $\lim_{N\to\infty}E[\int_{[0,1]}
x^k\Xi^{\log,\alpha}_N(t,1)(dx)]^\ell$ we use the dual process
$(\eta_t, \CF^{+}_t)_{t \geq 0}$ where the particle system $\eta$
starts in configuration $\eta_0$ with
\emph{k-particles} at \emph{each of $\ell$ distinct sites} $i=1,
\cdots, \ell$ and the function-valued part starts with
\be{and12}
\CF^{+}_0 = \bigotimes_{j=1}^\ell\bigotimes^k_{i=1} (1_{\{1\}}).
\ee
Let
\be{angr19} \Pi^{N,k,\ell}_u\ee denote the resulting number of
dual particles at time $u$.

Recall that the mutation jump to type 2 (from type 1) is given by
$1_{\{1\}}\to 0$. Recall furthermore that the selection operator
$1_{A_2} = 1_{\{2\}}$ preserves the product form (note for every
state $f$ before the first mutation $1 \to 2, f 1_{(2)} \equiv 0)$
and just generates another factor $1_{\{1\}}$ with a new variable.
Therefore the dual $\CF^{+}_t$ (before the first rare mutation
$1\to 2$ occurs) is a product of factors $1_{\{1\}}$
and this product integrated with respect to $\mu_0^\otimes$ with
$\mu_0(\{1\})=1$ is $0$ or $1$ depending on whether or not a
mutation jump has occurred. Therefore \bea{dualmoment} &&
m^N_{k,\ell}(t):= P[\text{no mutation jump $1 \to 2$ occurred by time }
T_N+t]\\&& = E[\exp (- \frac{m}{N}\intl^{T_N+t}_0 \Pi^{N,k,\ell}_u
du)].\nonumber\eea

Hence in order to prove that the following limit exists:
\be{angr20} m_{k,\ell}(t) = \liml_{N\to \infty} m^N_{k,\ell}(t),
\ee it suffices to prove that the following holds:
\be{angr20b}
\CL \left[ \frac{1}{N}\intl^{T_N+t}_0 \Pi^{N,k,\ell}_u du\right ]
\Nto \nu_{k,\ell} (t) \in \CP ([0,\infty]) \quad ,
\quad \nu_{k,\ell}(t) ((0,\infty))=1, \mbox{
for every } k,\ell \in \N
 \ee  and then we automatically have as well the formula
\be{angr20+}
m_{k,\ell}(t)=\int_0^\infty e^{-my}\nu_{k,\ell}(t)(dy).
\ee

The existence of the limit (\ref{angr20b})  is a result on the dual
process  which we proved in Proposition \ref{P.Grocoll}, part
(\ref{ba7}). This will be used in the next section to prove that
we have convergence to a random McKean-Vlasov limiting dynamic
which is truly random by showing that  the limiting variance of
the empirical mean is strictly positive. \bi

\subsubsection{Proof of Proposition \ref{P.measure}, Part 2: The random initial growth constant}
\label{sss.proofmeasure}

In the previous sections we have studied the asymptotics of
the dual process in terms of the Crump-Mode-Jagers branching
process and the pair $(u,U)$ following a certain nonlinear equation.
In this section we use
these results to establish that the limiting empirical
distribution of types is given by  the McKean-Vlasov dynamics with
random initial condition at time $-\infty$. Recall that in order to describe emergence
we assume that initially only type $1$ is present and we wish to
determine the distribution of $\bar x_2(t):=\lim_{N\to\infty}\bar
x^N_2(T_N+t)$ and its behaviour as $t \to -\infty$ subsequently.

We have so far established that
$\Xi^{\log, \alpha}_N (t,2)$ converges to $\CL_t(2)$ which is a solution to the
McKean-Vlasov dynamics and that \be{ag34} \int^1_0 x \CL_t(2)
(dx)\to 0\text{   as  } t\to -\infty.\ee Therefore  by Proposition
\ref{P.MVentrlaw} (see proof in Section \ref{sss.proofkv}),
Proposition \ref{P6.1b} and
(\ref{angr11}) the following limits exist in distribution and satisfy
\be{ag35}
\lim_{t\to -\infty}e^{\alpha|t|}\int_0^1 x\CL_t(dx)= \lim_{t\to
-\infty}\lim_{N\to \infty}e^{\alpha|t|}\bar x^N_2(\frac{\log
N}{\alpha}+t) =^\ast\CW,
\ee
so that $\CL_t$ is a random shift of $\CL^\ast$
(cf. (\ref{L.ast})), namely, $\CL^\ast_{t+\frac{\log
{^\ast\CW}}{\alpha}}$ and
\be{ag36}
\CL [^\ast \CW] =\lim_{t\to -\infty} \CL[e^{\alpha|t|}\int_0^1 x\CL_t(dx)].
\ee

We next consider the first and second moments of $^\ast\CW$.

\beP{L.dw9}{(Moments of $^\ast\CW$)}

The first and second moments of $^\ast\CW$ satisfy:
\be{dw9}
E[^\ast\mathcal{W}] = m^\ast E[W], \mbox{ with } \quad m^\ast = \frac{mb}{c},\quad
b=(1+\frac{\gamma}{\alpha}),
\ee
\be{dw10}
Var\Big[\int^1_0 x \CL_t(2) (dx)\Big] =O(e^{2\alpha t}),
\ee
and
\be{dw11}
 \text{Var}[^\ast\mathcal{W}] =m^\ast\frac{\kappa^\ast}{2\alpha} (E[W])^2.\qquad \square
\ee
\end{proposition}

\begin{proof}{\bf of Proposition \ref{L.dw9}}
Recall the equation (\ref{ag32+}) from Subsubsection
\ref{sss.limrandmarg} and note that this approximates $(u(t))_{t
\in \R}$ for $t \to -\infty$ up to terms $O(e^{3\alpha t})$:
\be{ag34b}
\int_{-\infty}^t\frac{1}{c}(\alpha(s)+\gamma(s))
u(s)ds = \frac{1}{c}\left(\frac{(\alpha+\gamma)}{\alpha}We^{\alpha
t}-\frac{\kappa^\ast}{2\alpha} W^2e^{2\alpha t}\right)
+O(e^{3\alpha t}),\quad \kappa^\ast>0.\ee
 This follows from
Proposition \ref{P.Transition} - see Lemma \ref{L.dw18} for the
proof that $\kappa^\ast>0$.

\bigskip

Therefore

 \bea{222} &&\exp\left( - \frac{m}{c}\int_{-\infty}^t
(\alpha(s)+\gamma(s))u(s)ds\right)\\&& =
\exp\left(-m^\ast\left(We^{\alpha t}-\frac{\kappa^\ast}{2\alpha}
W^2e^{2\alpha t}\right) +O(e^{3\alpha t})\right)\nonumber,
\eea where  $m^\ast=\frac{mb}{c}$.

An immediate consequence of this relation is that (recall
(\ref{angr20}) for a definition):
\be{dw15}
m_{1,1} (t)= 1- m^\ast E[W]e^{\alpha t}
+ O (e^{2 \alpha t}), \mbox{ for } t \to -\infty
\ee and therefore
\be{dw16}
\lim_{t\to -\infty} e^{\alpha|t|}
  E[\int^1_0 x \CL_t (2) (dx)]= m^\ast E[W].\ee

Now consider the second moment of the type two mass at time $t$
which we denote $ m^{(2)}(t)$. Here we consider the dual process
with one particle at each of two distinct sites.  The
corresponding growing clouds have independent random growth
constants $W_1$ and $W_2$. We use the formula \bea{dw17}
m^{(2)}(t):=E[\int x \CL_t (2) (dx)]^2
&& = E[ \int x \CL_t (1) (dx)]^2 \\
&& -2 E[ \int x \CL_t (1) (dx)] +1
=m_{1,2}(t)-2m_{1,1}(t)+1.\nonumber \eea

In order to calculate $m_{1,2}(t), m_{1,1}(t)$ we expand
(\ref{222}) to terms of order $o(e^{2\alpha t})$  for
$\wt W=W_1+W_2$ respectively $W_1$ where $W_1, W_2$ are
independent copies of $W$, to obtain

\bea{ds07}
&&\exp\left(-m^\ast\left(b \wt We^{\alpha
t}-\frac{\kappa^\ast}{2\alpha} \wt W^2e^{2\alpha t}\right)
+O(e^{3\alpha t})\right)\\&& =1-m^\ast \left(b \wt W e^{\alpha
t}-\frac{\kappa^\ast}{2\alpha} \wt W^2e^{2\alpha t}\right)
+\frac{1}{2}\left(m^\ast b  W e^{\alpha t}\right)^2 +O(e^{3\alpha
t}). \nonumber
 \eea
Therefore
 \bea{ag38}
 m_{1,2}(t)&&= 1-m^\ast bE[W_1+W_2]e^{\alpha
 t}+m^\ast\frac{\kappa^\ast}{2\alpha}E[W_1+W_2]^2e^{2\alpha t}
\\&& +\frac{1}{2}(m^\ast b)^2E[W_1+W_2]^2e^{2\alpha t}+O(e^{3\alpha
t}).\nonumber \eea
Combining this with (\ref{dw15}) and (\ref{dw17}), we see that:
 \bea{ag38+}
 m^{(2)}(t)&&= m^\ast\frac{\kappa^\ast}{2\alpha}E[W_1]E[W_2]e^{2\alpha t}
\\&& +\frac{1}{2}(m^\ast b)^2E[W_1]E[W_2]e^{2\alpha t}+O(e^{3\alpha
t}).\nonumber \eea Note that  the coefficients of the higher order
terms involve higher moments of $W$ and cancel provided that the
latter are finite. Hence to justify the cancellation it is
necessary to verify that the higher moments of $W$ are finite.
However this was established in Lemma \ref{L6.11} in Subsubsection
\ref{sss.dualcfrprop}.

Inserting (\ref{dw15}) and (\ref{ag38+}) in $Var [\bar x_2(t)] = m^2(t)-(E [\bar x_2(t)]^2$
we get:

\be{ag39}
 \text{Var}[\bar x_2(t)] = m^\ast \frac{\kappa^\ast}{2\alpha}b(E[W])^2 e^{2\alpha t}+O(e^{3\alpha t}).
\ee
If we can show that
\be{ag39b}
E [(\bar x_2 (t))^3] = O(e^{3 \alpha t}),
\ee
we get via normalizing by $e^{\alpha t}$ and letting $t \to -\infty$
indeed: \be{ag40}
 \text{Var}[^\ast\mathcal{W}]=m^\ast\frac{\kappa^\ast}{2\alpha} b(E[W])^2.\ee
Therefore we get (\ref{dw11}) once we have the bound on the
supremum over time of the third moment of
$\bar x_2(t)$, which follows from Lemma \ref{L.Himo} in Subsection \ref{ss.highermom}.

Here we can use the fact that the collision-free regime gives a stochastic upper
bound of the finite $N$ system and then that in fact according to (\ref{ag40})
we can bound the third moment of the scaled variable in $t$.

\end{proof}

\begin{remark}
The randomness expressed in $^\ast\CW$ arises as a result of the
fact that the mutation jump occurs  in this time scale at an
exponential waiting time (in $t$) together with the nonlinearity
in (\ref{grocoll2}), (\ref{ba10g}). The random variable $W$
governing the initial growth of the dual population (which
describes the limiting behaviour  of the CMJ process) arise from the initial birth events
in the dual process and has its main
influence due to the nonlinearity of the evolution in $t$. However $W$ will
enter only through its mean in the law of $^\ast\CW$ due to a law of
large number effects.

The moments of $^\ast\CW$ are determined in terms of the random
variable $W$ arising from the growth of $ K_t$ in the dual
process.  However as we have seen in the second moment calculation
above the determination of the coefficients of $e^{k\alpha t}$
depend on an analysis of an asymptotic expansion of the nonlinear
system and in particular $u(t)$ in orders up to $k$.  Although we
do not attempt to carry this out here,  we show that with respect
to Laplace transform order the law of $^\ast\CW$  lies between the
case in which the collisions are suppressed (deterministic case)
and a modified system in which the  the correction involving $\int
(\alpha(s)-\alpha)ds$ is suppressed, that is,  replacing
$\alpha(s)$ by $\alpha$. This  illustrates the role of the
collisions (nonlinearity) in producing the randomness in
$^\ast\CW$ and allows us to obtain bounds for the expected time to
reach a small level $\ve$ (cf. \cite{L1}). This we pursue further
in Subsection \ref{ss.relationsW}.

\end{remark}

\subsubsection{Completion of the Proof of Proposition
\ref{P.measure}}\label{sss.complproof}

We now collect all the pieces needed to prove all the assertions of the
proposition, we proceed stepwise.
\sm

{\it Step 1 $\;$  Completion of the Proof of Proposition \ref{P.measure}(a)}

This was proved in Subsubsection \ref{sss.proofldy} which heavily used \ref{sss.limrandmarg}.
\bi

{\it Step~2 $\;$ Completion of the proof of the Proposition
\ref{P.measure}(b)}

 In order to verify that the limiting marginal random probability measure
 $ \CL_t (j)$ has for  $j=1,2$ really a nontrivial distribution and
 is not just deterministic, it suffices to show that
\be{angr23}
 E([\int x\mathcal{L}_t(1)(dx)]^2)> (E([\int x\mathcal{L}_t(1)(dx)]))^2.
 \ee
The quantity $\int x \CL_t(1)(dx)$ arises as the limit in
distribution of $\bar x^N_1(t)$, the empirical mean mass of type
1. Note that the empirical mean mass process $(\bar
x^N_\ell(t);\ell=1,2)_{t \geq 0}$ is a random process which is not
Markov! To see this look at the second moment of mean mass, it
involves the covariance between two sites! More generally the dual
representation of the kth moment of the empirical mean is given by
starting k particles at $k$ distinct sites.

Note that no coalescence occurs until time $O(\log N)$ and then
the two clouds descending from the two different initial dual
particles \emph{interact} nontrivially. Then it is easy to
conclude that the limiting variance of $\bar x^N_2(T_N+t)$ as $ N
\to \infty$ is not zero using (\ref{ba9}). \bi

{\em Step 3 $\;$ Completion of the proof of Proposition \ref{P.measure} (c)}.

The relations (\ref{LD2}), (\ref{LD3}) follow from the dual representation
(\ref{dd12b}) of $\int^1_0 x^k \CL_t(2) dx$ from the fact that
$\nu_{k,\ell} (t, \cdot) \La \delta_0$ as $t \to -\infty$ and
$\nu_{k,\ell}(t,\cdot) \La \delta_\infty$ as $t \to \infty$.
This follows from combining (\ref{angr24}) and (\ref{ba7})
and then using (\ref{agrev63}) to get the first claim and using
$u^{k,\ell}(t) \to \infty$ as $t \to \infty$ to get the second claim.
For the proof of  (\ref{LD1}) we will use the
general fact that \be{LD6} P(\CL_t (1)(\{1\})
=1)=\lim_{k\to\infty}m_{k,1}(t). \ee The result then follows from
Lemma \ref{L6.1} equation (\ref{ba88}).

{\em Step 4 $\;$ Completion of the proof of Proposition \ref{P.measure} (e)}

The convergence of the empirical measure processes to a random
solution of the McKean-Vlasov equation follows from Proposition \ref{P.Grocoll} together
with Proposition \ref{P.measure} (b).
\sm

{\em Step~5$\;$ Proof of Proposition \ref{P.measure}(f)}

The claims (\ref{angr14a2})- (\ref{angr14d}) was proved in
Subsubsection \ref{sss.proofmeasure} in (\ref{ag35}) and
(\ref{ag36}). The  assertion (\ref{angr14a})  follows from
Proposition \ref{L.dw9}.

{\em Step 6$\;$ Proof of Proposition \ref{P.measure} (g)}

For the assertion (\ref{angr14a3}) we
refer to a calculation we do later, see (\ref{gra93}).

\bi

All these steps complete the proof of Proposition \ref{P.measure}.

\begin{remark} The same approach can be carried out on $\Omega_L$ for fixed $L$ or $\mathbb{Z}^d$.
Namely the dual can be viewed as follows.
Start with one particle, then the particle system (factors)  can
 be viewed as a stochastic Fisher-KPP equation (at least until the time of the first mutation).
 In particular in $\mathbb{Z}^1$ we expect
 linear growth according to a travelling wave solution (see \cite{CD05}) and not
 exponential growth.  In this case emergence occurs in time $O(\sqrt{N})$.
 \end{remark}

\subsection{Droplet formation: Proofs of Proposition \ref{P.DF}-\ref{CSB-longtime}}
\label{ss.bounds}

In this section we assume that $c>0,\;m>0,\;s>0$, and examine the
process in which rare mutants first appear in a finite time horizon in the population at
the microscopic level,  that is only at some rare sites they appear at a substantial level,
then later on after large times these rare sites develop into growing {\em mutant
droplets}, that is, a growing collection of sites occupied by the mutant
type but still being of a total size $o(N)$. Then finally in even much larger
times this leads  to emergence at the macroscopic level once we have
$O(N)$ sites in the droplet.

We show that in the microscopic growth regime of the droplet the
type-2 mass is described by a population growth process with
Malthusian parameter $\alpha$ and random factor $\mathcal{W}^\ast$.
In the next Section \ref{ss.relationsW} we use this structure to
investigate the relation between $^\ast\CW$ and $\CW^\ast$. Recall
that $^\ast\CW $  arose in the context of emergence and fixation
by considering time $\alpha^{-1} \log N +t$ and letting first $N
\to \infty$ and then $t \to -\infty$ whereas   $\CW^\ast$ arises
from the droplet growth at times $t$ and by letting first $N \to
\infty$ and then $t \to \infty$.

To carry out the analysis for the droplet growth we
consider the following four time regimes for our population model:
\be{ddz} [0,T_0), \text{ where
} 0<T_0<\infty,\ee
\be{ddy}
[T_0,T_N),  \text{ where }
T_N\to\infty, \mbox{ as } N \to \infty \mbox{ and }
 T_N = o(\log N),\ee
\be{ddx}
[T_N, \frac{\delta
\log N}{\alpha}+t],\text {  where  } 0<\delta<1,\; t\in \mathbb{R},
\ee
\be{ddw}
[\frac{\delta \log N}{\alpha}, \frac{\log N}{\alpha}+t],
\mbox{ where } 0 < \delta < 1, \quad t \in \R.
\ee
The first two time regimes are needed in proving the
Propositions \ref{P.DF}-\ref{CSB-longtime}, the two remaining
ones are necessary to prepare the stage for Subsection \ref{ss.relationsW}
where we shall relate $^\ast\CW$ and $\CW^\ast$.
\bi

{\bf Outline of Subsection \ref{ss.bounds}} \quad

We now give a  description of the evolution through the stages
corresponding to the time intervals given in
(\ref{ddz})-(\ref{ddw}).

The first step is to examine the sparse
set of sites at which the mutant population is of order $O(1)$ in
a fixed finite time interval and then to describe this in terms of
a Poisson approximation in the limit $N \to \infty$.
This happens in Subsubsection \ref{sss.dropletform} and proves
Proposition \ref{P.lt}.

In Subsubsection \ref{sss.proofCSB} we
analyse the consequences for the  longtime properties and
prove Proposition \ref{CSB-longtime} and in Subsubsection
\ref{sss.proofPDF} the Proposition \ref{P.DF}.

Subsequently in \ref{sss.explicitcal}
we recall some related explicit calculations and in \ref{sss.dropletgrow}
and \ref{sss.detdropgro} we continue with the time intervals (\ref{ddy})-(\ref{ddw}) to
exhibit the law and properties of $\CW^\ast$.

\subsubsection{Mutant droplet formation at finite time horizon}
\label{sss.dropletform}

Recall the definition of $(\gimel_t^{N,m})_{t\geq 0}$ in
(\ref{dd29}). Furthermore we recall the abbreviation:
\be{ddv}
\wh x^N_2(t) = \suml^N_{j=1} x^N_2(j,t)
\ee
for the total mass of type
2 in the whole population of $N$ sites. We have to prove the
convergence of $(\gimel^{N,m}_t)_{t \geq 0}$ as $N \to \infty$ to
$(\gimel^m_t)_{t \geq 0}$ and then we have to derive the
properties of this limit as $t \to \infty$, which we had stated as
three propositions which we now prove successively, but not in
order, we close with Proposition \ref{P.DF}.

\begin{proof} {\bf of Proposition \ref{P.lt}}

The strategy to prove the convergence in distribution of
$(\gimel_t^{N,m})_{t\geq 0}$ as $N\to\infty$  is to proceed in
steps as follows. We first consider only the contributions of
mutation at a typical site where we have in reality two type of rare events, (1) the
{\em immigration} of mass of type two from other sites and (2) the building
up of this mass via {\em rare mutation} at this site. Indeed if we isolate
the two effects and first suppress the immigration of type-2 mass,
this simplification leads to $N$-independent sites. We therefore first
look at one site, then at an independent collection of $N$ sites
and then finally build in the effect of immigration (destroying the independence) from other
sites. In other words we work  with the simplification in Step 1
and Step 2 and return in Step 3 to the original model. In Step 1
and Step 2 each the main point is condensed in a Lemma.

{\bf Step~1}$\;$ {\em (Palm distribution of a single site dynamic)}


Recall Lemma \ref{L.Pa2} and the definition (\ref{dd31a}) of the single site excursion measure $\mathbb{Q}$ for the process without mutation.  We must now consider the analogue for the process with mutation.
Consider the {\em single site} dynamics (with only emigration but no immigration),
which is given (misusing notation):
\bea{ssd}
d x^N_2 (1,t) & = &-c x^N_2 (1,t) dt +s\,x^N_1(1,t) x^N_2 (1,t) dt
 + \frac{m}{N}(1- x^N_2 (1,t))dt \\
 && + \sqrt{d \cdot (1- x^N_2 (1,t))(x^N_2 (1,t)} dw_2 (i,t),\nonumber\\
 x^N_2(1,0)&= &\frac{a}{N}. \nonumber
 \eea

We use size-biasing to focus on the set of sites at
which mutant mass appears in a finite time interval and
use the excursion law of the process without mutation
(with law $P$) to express the excursions of the process
with mutation (with law $\wt P$). We prove first as key tool
the following. Let
\be{q1a}
\mathbb{Q} \in \CM (C(\R^+, [0,1]))
\ee
be defined as the excursion law of the process in (\ref{onedim}) (which is(\ref{ssd}) with
$m = 0$ which then does not depend on $N$). Furthermore
let
\be{add171}
\wt {P}^\ve = \CL[(x^\ve_2(t))_{t \geq 0}] = \CL[ (x^\ve_2(1,t))_{t \geq 0}]
\mbox{ and } x^\ve_2 \mbox{ solves (\ref{ssd}) with } x^\ve_2(0) = \frac{a}{N} =0,
\frac{m}{N}=\ve.
\ee

Furthermore let
\bea{dd30im}
W_{\rm{ex}}:=&&\{w\in C([0,\infty),\mathbb{R}^+),\;w(t) =0 \mbox{ for } t < \zeta_b \mbox{ and } t > \zeta_d,\;
w(t)>0\\&&\mbox{ on the interval   } \zeta_b<t<\zeta_d \mbox{ for some }
\zeta_b<\zeta_d\in(0,\infty)\}.\nonumber
\eea
Note that then $\zeta_b$ and $\zeta_d$ are unique functions of the element
$w \in W_{\mathrm{ex}}$. We shall use these functionals
\be{add172}
\zeta_b, \zeta_d : W_{\mathrm{ex}} \la (0,\infty)
\ee
below.

Next observe that every continuous path starting at 0 and having an unbounded set of zeros we can write
uniquely as a sum of elements of $W_{\mathrm{ex}}$ with disjoint intervals of positivity and with
every point $x$ not a zero of the path we can associate a unique excursion between
\be{add138}
\zeta^x_b \mbox{ and } \zeta^x_d.
\ee

\beL{L.Pa}{(Single site excursion law and Palm distribution)}

Consider the evolution of (\ref{ssd}) and let $a=0$, $d>0$, $c \geq 0$, $t_0>0$.
Then the following properties hold:

(a) Fix $t_0>0$ and  $K\in\N$. For $k=1,\cdots,K$ set
$I_k:=x_2(\frac{k}{K}t_0)$. Let  $A_k(a_3)$ denote the event  that the continuing trajectory satisfies,
recall (\ref{dd30im}) and the sequel,
$\zeta^x_d-\frac{k}{K}>a_3$ with $x=\frac{k}{K}t_0$.
Furthermore consider the same quantities for processes $x_{2,k}$ and denote then
by $\bar I_k$, respectively $\bar A_k$.

Replace in (\ref{ssd}) the mutation term
by $\ve f_k$ and $f_k$ is the indicator of
$\left[\frac{k-1}{K} t_0, \frac{k}{K} t_0\right]$. These processes
are called $x^\ve_{2,k}$.

Then the ${\bar I_k}$ are identically distributed random variables with means $\frac{\ve}{K}t_0$ and variance asymptotically of the form $\frac{const\cdot \ve}{K}$ as $K\to\infty$.

 Then the probability of more than one such excursion is
 \be{oneex} \lim_{\ve\to 0} \wt P^\ve(\sum_{k=1}^K 1(\bar A_k) >1|\sum_{k=1}^K 1(\bar A_k) \geq 1)=0.
 \ee
This remains true for the $I_k$ and $A_k$.

(b)
The following family of measures
$\wt \Q ((a_1, a_2), a_3, \cdot)$ and a measure $\Q$
on $W_{\rm{ex}}$ exists. Choose $a_1, a_2 \in [0,\infty), a_1<a_2, a_3<0$
and $A \in \CB (C([0,\infty), \R^+))$.
Then
\bea{agrev60aa}
&& \wt \Q((a_1,a_2), a_3,A))=:
 \wt {\mathbb{Q}}(\zeta_b(w)\in (a_1,a_2), \,\zeta_d-\zeta_b > a_3,\;w\in A)\\
 &&
 =\lim_{\ve\to 0} \frac
{\wt P^\ve(\exists x \in (a_1, b_1) \zeta^x_b(w)\in (a_1,a_2), \,\zeta^x_d-\zeta_b >a_3,\;w\in A)}{\ve}.\nonumber\eea
The measures $\wt \Q$ on $W_{\rm{ex}}$ can be represented as ($\Q$ as in (\ref{q1a})),
\be{agrev60a}
\wt{\mathbb{Q}}((a_1,a_2),\,a_3,\, A)=\int_{a_1}^{a_2}\mathbb{Q}(\zeta> a_3,w( \cdot -s))\in A)ds.\ee

(c) Now consider (\ref{ssd}) which has mutation rate $\ve=\frac{m}{N}$
and put $a=0$. Then for $t_0>0$

\be{elam2} \lim_{N\to\infty} N\cdot E[x^N_2(1,t_0)]>0.\ee

(d) Let $\wh \mu^N_t$ be as defined in (\ref{Y12d}), that is, the
Palm distribution of $x^N_2(1,t)$ in (\ref{ssd}) with
$a=0$. The following limit exists:
\bea{elam2b}
&&\wh \mu^\infty_t (dx)= \liml_{N\to \infty} \wh \mu^N_t(dx)
= \frac{x\int_0^t\wt{\mathbb{Q}}(\zeta_b\in ds,\zeta_d>t,w(t)\in dx)ds} {\intl^\infty_0 {\int_0^t x\wt{\mathbb{Q}} (\zeta_b\in ds,\zeta_d>t,w(t)\in dx)ds}}\,.\eea

\noi Hence  the Palm distribution $\wh \mu^N_t$ is the law of a
random variable of order $O(1)$ which is also asymptotically
non-degenerate as $N \to \infty$.

Furthermore as a consequence of (\ref{elam2}) and (\ref{elam2b})
the first moment of the Palm distribution of $x^N_2 (1,t_0)$ has
mean satisfying:
\be{elam3}
\lim_{N\to\infty} \wh E[x^N_2(1,t_0)]=\lim_{N\to\infty}
\frac{E[(x^N_2(1,t_0))^2]}{E[(x^N_2(1,t_0))]} \in (0,\infty).
\qquad \square
\ee
\end{lemma}

\begin{proof}{\bf of Lemma \ref{L.Pa}}

We prove separately the parts a), b) and then c), d) of the Lemma.

{\em (a) and (b)}.

To prove (a) and (b) we begin by approximating the single site process $\{ x^\ve_2(t):0\leq t\leq t_0\}$ with law $\wt P^\ve$ (which satisfies (\ref{ssd}) with $\frac{m}{N}=\ve$) by a sequence of processes, where we replace the
mutation term (induced by the rare mutation) in different time intervals  by mutation terms at a grid of discrete time points
getting finer and finer.
In order to keep track of the mutations at different times
we split the type 1 into different types. Namely
let    $K\in\N$ and consider a $K+1$ type Wright-Fisher diffusion with law
\be{add139}
\wt P^{\ve,K},
\ee
which is starting with type $0$ having fitness $0$
and the other types having fitness $1>0$ (selection at rate $s$) and with mutation from type
$0$ to type $k$ at rate $\ve$ only during the interval $(\frac{(k-1)t_0}{K},\frac{kt_0}{K}]$.

We first consider a simpler process, which is built up from $K$-independent processes. Namely
consider independent processes $x^\ve_{2,k}$ representing the masses of types $k=1,\dots,K$   in which in each time interval we suppress mutation of all types except some $k$. Then  $x^\ve_{2,k}(\cdot)$ is a solution of  (\ref{ssd}) with mutation term
 $\ve 1_{(\frac{k-1}{K}t_0,\frac{k}{K}t_0]}(\cdot)$.
 Let $I_k(\ve):=x^\ve_{2,k}(1,\frac{k}{K}t_0)$. Let  $A_k(a_3)$ denote the event  that the continuing trajectory satisfies
 $\zeta_d-\frac{k}{K}>a_3$.
 Note that the $ {I_k}(\ve)$ are identically distributed random variables with means $\frac{\ve}{K}t_0$ and variance asymptotically of the form $\frac{\ve t_0^2}{2K^2}+o(\frac{\ve}{K^2})$.

Note that as $\ve \to 0$.
 \be{add53}
 \frac{1}{\ve}{P_{I_k(\ve)}(A_k (a_3))}= \frac{I_k(\ve)}{\ve} \Q(\zeta>a_3) +o(I_k(\ve)).
 \ee
 Therefore
 \bea{add54}
 && \suml^K_{k=1} 1(\frac{kt_0}{K}\in (a_1,a_2)) \frac{1}{\ve}{P_{I_k(\ve)}(A_k (a_3))} \\
 &&=\frac{1}{K}\suml^K_{k=1} 1(\frac{kt_0}{K}\in (a_1,a_2))  \frac{K I_k(\ve)}{\ve} \Q(\zeta>a_3)+o(\ve)\mbox{   as  }\ve\to 0.\nonumber\eea

 Moreover  by independence,
 $Prob(A_k\cap A_{k^\prime})= P_{I_k(\ve)} (A_k) P_{I_{k^\prime} (e)} (A_{k^\prime}) =O(\ve^2),\;k\ne k^\prime,$
 and therefore the probability of more than one such excursion in the limit $\ve \to 0$ is equal to zero.

 \bea{add108}
 \lim_{\ve\to 0}\frac{1}{K^2} \sum_{k,k^\prime}  \frac{1}{\ve}{P_{I_k(\ve)}(A_k)\cdot P_{I_{k^\prime}(\ve)}(A_{k^\prime})} =0.\eea

Then letting $K\to\infty$ we get  using the continuity of
$t \to \Q(w(t) \in \cdot)$  that the following limit exists
  \bea{add109}
  &&\wt {\mathbb{Q}}(\zeta_b\in(a_1,a_2),\zeta_d-\zeta_b>a_3, w(t_0)>x)\\
  &&=\lim_{\ve\to 0}\lim_{K\to\infty}]
  \left[\frac{1}{K} \sum_k 1(\frac{kt_0}{K}\in (a_1,a_2)) \frac{K I_k(\ve)}{\ve}\frac{1}{I_k(\ve)}{P_{I_k(\ve)}(A_k\cap \{x(t_0)>x\})} \right]\nonumber\\
  &&=\lim_{\ve \to 0}\left[\lim_{K\to\infty}\frac{1}{K}\sum_k 1(\frac{kt_0}{K}\in (a_1,a_2))  \frac{KI_k(\ve)}{\ve}
  \Q(\zeta>a_3,\, w(t-\frac{kt_0}{K})>x)+o(\ve)\right]\nonumber\\
  &&= \int_{a_1}^{a_2}\Q(\zeta>a_3,\,w(t_0-s)>x)ds.\nonumber
  \eea

 We claim that the analogous result is valid for the multitype model, that is, letting $A^\ast_k$
 denote the presence of an excursion of type $k$ of length greater that $a_3$ and having at time $t_0$
 mass bigger than $x$, starting in $(a_1, a_2)$
 \be{oneex2} \lim_{\ve\to 0} \wt P^{\ve,K} (\bigcup\limits^K_{k=1} A^\ast_k) =
 \intl^{a_2}_{a_1} \Q(\zeta > a_3, w (t_0-s) >x) ds.
 \ee

To verify that this remains true with dependence between the types as they arise in a multitype Fisher-Wright
diffusion, we return to the law $\wt P^{\ve,K}$.
The result then follows from (\ref{add109}) by  noting that

\be{add55}
\wt P^{\ve,K}(A_k\cap A_{k^\prime})\leq P_{I_k(\ve)}(A_k)\cdot P_{I_{k^\prime}(\ve)}(A_{k^\prime}).\ee
The latter holds, since

{\em (c) and (d)}

The proof of the convergence of $\wh \mu^N_t$ as $N \to \infty$
proceeds by showing that for every $k \in \N$ we have:
\be{palm5}
\lim_{N\to\infty}N\cdot E([x^N_2(i,t)]^k)=m \lim_{\delta\downarrow 0}\int_0^1 x^k\,\int_0^t \wt
\Q(\zeta_b\in ds,\zeta_d-\zeta_b>\delta, w(\cdot -s)\in dx)\ee and that the limiting variance is positive.

We first use the dual process to the process defined in equation
(\ref{ssd}) to compute the first and second moments.
Since we discuss the limit $N \to \infty$ over a finite time
horizon with finitely many initial particles, we get in the
limit no immigration term, but only emigration at the
migration rate.

We warn the reader at this point that we
use the same notation for the dual process of this single site diffusion as for the
one of our interacting system.
Let $\Pi^{(j)}_u$ denote the number of factors in
(\ref{grev62}) at time $t$ starting with $j$ factors at time $0$
where $\Pi^{(j)}_t$ is the birth and death process with birth and death rates
\be{gre71}
sk \mbox{ and } ck+\frac{d}{2}k(k-1) \mbox{  respectively.}
\ee
If we consider finite $N$ we have additional immigration
at rate $N^{-1} (\Pi^{(1),N}_t -k)^+$ if at the site considered
we have $k$ particles. This means that as $N \to \infty$ and over a
finite time horizon with
probability tending to 1 at rate $O(N^{-1})$ no immigration jumps occur.

Then for fixed $t$ we have that:
\be{grev75a}
E[x_{2}^{N}(i,t)] =1-E[x_{1}^{N}(i,t)], \ee
\be{grev75c}\begin{array}{ll}
E[x_{1}^{N}(i,t)]  &  =E[\exp(-\frac{m}{N}\int_{0}^{t}\Pi^{(1),N}_{u}du)]\\
& \sim E[\exp (-\frac{m}{N} \intl^t_0 \Pi^{(1)}_u du)], \quad \mbox{ as } N \to \infty\\
&  \sim1-\frac{m}{N}E[\int_{0}^{t}\Pi^{(1)}_{u}du], \quad \mbox{ as } N \to \infty\\
&  \sim1-\frac{mc_{1}t}{N},\quad\text{ as } N\to\infty,
\end{array}
\ee
where $c_1 = c_1(t) \in (0,\infty)$ is constant in $N$. Hence
\be{grev75d}
\lim_{N\to\infty}
N\cdot E[x_{2}^{N}(i,t)] = mc_{1}=m\int_0^t E[\Pi^{(1)}_u]ds. \ee

For fixed  $t$, as $N\to\infty$ we calculate the second moment:
\be{grev75e}\begin{array}{ll}
E[(x_{1}^{N}(i,t))^{2}]  &  =E[\exp(-\frac{m}{N}\int_{0}^{t}\Pi_{u}^{(2),N}du)]\\
& \sim E[\exp (-\frac{m}{N} \intl^t_0 \Pi^{(2)}_u du)]\\
&  \sim1-\frac{mc_{2}(t)}{N},\end{array}
\ee
for some $c_2 =c_2(t)  \in (0,\infty)$. Hence

\be{grev75f}\begin{array}{ll}
E[x_{2}^{N}(i,t)^{2}]  &  =1-2E[x_{1}^{N}(i,t)]+E[x_{1}(i,t)^{2}]\\
&
\sim 1-2E[\exp(-\frac{m}{N}\int_{0}^{t}\Pi_{u}^{(1)}du)]+E[\exp(-\frac{m}{N}
\int_{0}^{t}\Pi_{u}^{(2)}du)]\\
&  \sim2\frac{mc_{1}t}{N}-\frac{mc_{2}t}{N}, \quad \mbox{ as } N \to \infty.
\end{array}
\ee
Therefore
\be{grev75d1}
\lim_{N\to\infty} N\cdot
E[x^N_2(i,t)^2]=2mc_1(t)-mc_2(t)= m[2\int_0^t E[\Pi^{(1)}_u]du
-\int_0^t E[\Pi^{(2)}_u]du].\ee

 Therefore the l.h.s. of (\ref{elam3}) equals
\be{ag41}
\frac{(2c_1(t) - c_2(t))}{c_1(t)}. \ee  This limit is
larger than $0$ since (for every $t$)
\be{ag41b}
c_2 (t) < 2 c_1 (t)
\ee
holds because of the positive (uniformly in $N$)
probability of coalescence between the two initial particles of
the two populations before a migration step of either one. This implies
that every limit point of the sequence of laws of
$x^N_2(t)$ (solving (\ref{ssd})) has
positive variance.

Similarly we introduce
\be{ag41a}
c_k(t) = \intl^t_0 \Pi^{(k)}_u du
\ee
and we can establish writing
\be{ag41g}
E[(x^N_2 (i,t))^k] = \suml^k_{j=0} {k \choose j} (-1)^j
E[(x^N_1 (i,t))]^j]
\ee
and then using
\be{ag41h}
E[(x^N_1 (i,t))^j] = E[exp(-\frac{m}{N} \intl^t_0 \Pi^{(j)}_u du]
\ee
the convergence of all $k$-th moments as for $k=1,2$.
This then gives the weak convergence of
the Palm distribution $\wh \mu^N_t$ as $N\to \infty$.

Finally we have to verify the relation (\ref{elam2b}), that is,
the limiting Palm distribution is represented in terms of the
excursion measure.

Note that for every $k \in \N$ we know that: \bea{palm6}
&&\int_0^1 x^k \wt \Q_t(dx)\\
&&=\int_0^t \int_0^1 x^k \mathbb{Q}(x_2(t-s)\in dx) ds
 =
 \int_0^t\int_0^1
x^k(\lim_{\ve\to 0}\frac{1}{\ve}\wt P^\ve[x_2(t-s)\in dx)ds \nonumber\\
&&= \int_0^t \lim_{\ve\to 0}\frac{1}{\ve}E_\ve[(x_2(t-s))^k] ds.
\nonumber
\eea
Using the duality we calculate then (with $\wt E_\ve$ denoting the expectation
with respect to $\wt P_\ve$)
\be{palm6b}
\lim_{\ve\to 0}\frac{1}{\ve} \wt E_\ve [x_2(u)]=
\lim_{\ve\to 0}\frac{1}{\ve}E[1-(1-\ve)^{\Pi^{(1)}_u}]=E[\Pi^{(1)}_u)].\ee

\bea{palm6c}
\liml_{\ve \to 0} \frac{1}{\ve} \wt E_\ve[x_{2}(u)^{2}]
&& = \liml_{\ve \to 0} \frac{1}{\ve}\left(1-2 \wt E_\ve [x_{1}(u)]+\wt E_\ve[x_{1}(u)^{2}]\right)\\
&& = \liml_{\ve \to 0} \frac{1}{\ve}\left(1-2E[1-(1-\ve)^{\Pi^{(1)}_u}]
     +E[1-(1-\ve)^{\Pi^{(2)}_u}]\right)\nonumber\\
&&= 2E[\Pi^{(1)}_u]-E[\Pi^{(2)}_u].\nonumber \eea Substituting
this (with $k=2$) in (\ref{palm6}) and comparing with
(\ref{grev75d}) and (\ref{grev75e}) we verify (\ref{palm5}) for
$k=1,2$ and (\ref{elam3}). Similar calculations can be used to
verify the claim for all $k \in \N$.$\qquad\blacksquare$

\end{proof}

\bigskip

{\bf Step~2}$\;$ {\em (Compound Poisson limit)}

Consider a collection of processes where each component of the system
is still as in (\ref{ssd})  in which we ignore the effects of immigration
into sites and only {\em emigration} is still accounted for, namely,
 \bea{ag41a0} d x^N_2 (i,t) = -c
x^N_2 (i,t) dt
 &+& s\,x^N_1(i,t) x^N_2 (i,t) dt
 + \frac{m}{N}(1- x^N_2 (i,t))dt \\
 &+& \sqrt{d \cdot (1- x^N_2 (i,t))x^N_2 (i,t)}\, dw_2
 (i,t),\nonumber\\&& x^N_2(i,0)=\frac{a}{N}, \quad i=1,2,\cdots,N.\nonumber
 \eea
Note that this is a system of $N$ independent diffusion processes.
Furthermore we calculate for the total mass process \be{ag41a1}
d\wh x^N_2(t)\leq ((-c+s)\wh x^N_2(t)+ \frac{m}{N})dt+  dM_t, \ee
with a martingale $(M_t)_{t \geq 0}$, where \be{ag41a3} \langle
M\rangle_t \leq d \cdot \wh x^N_2(t). \ee Therefore we can bound
by a submartingale and if we start with only type 1, i.e. {\em if
$a=0$}, then by submartingale inequalities we get: \be{ag41a3a}
P[\sup_{t\leq T}|x^N_2(i,t)|>\ve]\leq \frac {\text{constant}}{N}.
\ee

We return now to the study of the independent collection of diffusions we
introduced in the beginning of this step.
Consider the atomic measure-valued process
\be{ag41a4}
(\wt \gimel^{N,m}_t)_{t \geq 0},
\ee
which is defined as the analogue of $(\gimel ^{N,m}_t)_{t \geq 0}$ but with only
excursions arising from mutation, that is ignoring immigration at
each site (emigration is still accounted for) which was defined in
(\ref{ag41a}). Similarly let
\be{ag41a5}
\{\wt x^N_2(i,t), i \in \N\}
\ee
be the corresponding collection of independent diffusions.
Then we carry out the limit $N \to \infty$.

\beL{L.step2} {(Droplet growth in absence of immigration)}

At time $t_0$,
asymptotically as $N\to\infty$ we have the following three properties
for the dynamics introduced above in (\ref{ag41a4}, \ref{ag41a5}).

(a)  There is a Poisson number of sites $j$ at which $\wt x^N_2 (j,t_0)>\ve$
where the parameter is
\be{ag41b1}
m\int_0^{t_0}\mathbb{Q}(w(t_0-s)>\ve)ds.
\ee

(b)  At time $t_0$, in the limit $N\to\infty$, there are a countable number of sites
contributing to a total mutant mass of order $O(1)$.

(c) For each $t>0$ the states of $\wt \gimel^{N,m}_{t}$
converge in the sense of the weak
atomic topology to
\be{ag41a6}
\wt \gimel^m_t \mbox{ on } [0,1],
\ee
with atomic random measure specified by the two requirements
\be{ag41b2}
\mbox{atom locations  i.i.d. uniform on } [0,1]
\ee
and atom masses given
by  a Poisson random measure on $[0,1]$ with intensity measure
\be{agrev6}
m \nu_t(dx), \mbox{ where } \nu_t (dx)
= \int_0^t \mathbb{Q}(w:w(t-s)\in dx)ds.
\ee

The Palm measure $\wh \nu_t$ of $\nu_t$ satisfies
(recall $\wh \mu^\infty_t$ from Lemma \ref{L.Pa} (d)):
\be{agrev6a}
\wh \nu_t = \wh \mu^\infty_t. \qquad \square
\ee
\end{lemma}

\begin{proof}{ \bf of Lemma \ref{L.step2}}
(a)  It follows from Lemma \ref{L.Pa} (a) that
\be{agrev6a2}
\lim_{N\to\infty}N\cdot
P[x^N_2(i,t)>\ve]={m}\wt{\mathbb{Q}}_t((\ve,1]).\ee

Therefore in the collection of $N$ sites it follows  that
asymptotically as $N\to\infty$ the number of sites $j$ with
$\wt x^N_2(j,t_0)>\ve$ (recall that in this step we suppress
immigration and we are considering only mass originating at this
site) is Poisson with parameter ${m}{\wt{\mathbb{Q}}_t}((\ve,1])$.

(b) Since $\nu_t((0,\infty))=\infty$ and
$\nu_t ((\delta,\infty)) < \infty$ for $ \delta > 0$,
there are countably many atoms in the limit as $N\to\infty$.

(c) Consider the sequence obtained by size-ordering the atom sizes
in $\wt \gimel^{N,m}_t$. Note that the limiting point process on $[\ve,1]$
is given by a Poisson number with i.i.d. sizes  and the distribution
of the atoms is given by
$\frac{\wt{\mathbb{Q}}_t(dx)}{\wt{\mathbb{Q}}_t([\ve,1])}$.  It
follows that the order statistics also converge and therefore the
joint distribution of the largest $k$ atoms converge as $N\to \infty$.

To verify that the limit is pure atomic we note that the
expected mass of the union of sites of size smaller than $\ve$
converges to $0$ as $\ve\to 0$ uniformly in $N$ by the explicit
calculations in the next section (cf. (\ref{gra45})). The
convergence of $\wt \gimel^{N,m}_t$ in the weak atomic topology then
follows by Lemma \ref{L.EK}.
\end{proof}

{\bf Step~3}$\;$ {\em (Completion of proof of convergence)}

Here we have to incorporate the migration of mass between sites
in particular the {\em immigration} to a site from all the other sites.
Note that the immigration rate at fixed site is as the rare
mutation of order $N^{-1}$ as long as the total mass is $O(1)$ as
$N \to \infty$ and therefore we get here indeed an effect
comparable to the one we found in Step 1 and 2. We can  show directly
that the first two moments of the total mass $\wh x^N_2(t)
=\gimel^{N,m}_t([0,1])$ converge to the first two moments of  $\wh
x_2(t)=\gimel^m(t)([0,1])$. These moment calculations are carried out
later on in Subsection \ref{ss.highermom}, in particular the first moment result is given by
(\ref{gra68}) and the second moment is given by
(\ref{Ma3xx2}).  This gives the tightness of the marginal
distributions.

Consider the sequence of time grids $\{\frac{\ell}{2^k},\;\ell = 0,1,2,\dots\}$
with index $k$ and  with width $2^{-k}$.  Then  the
migration of mass between time $\frac{\ell}{2^k}$ and $\frac{\ell
+1}{2^k}$ and the mutation in this interval will be replaced by immigration
from other sites and by rare mutation
at the {\em end} of the interval.
We then define using this idea an approximate evolution in discrete time by
an approximate recursion scheme of atomic measures on $[0,1]$.

Let $T_t (x)$ denote the evolution of atomic measures, where each atom follows the single site
Fisher-Wright diffusion with emigration and the initial atomic measure is $x$.
Furthermore let $Y^{N,k}_\ell$ be an atomic measure
with a countable set of new atoms such that the atoms in $Y^{N,k}_\ell$ are
produced as above in Step 1 if we observe the process at time $2^{-k}$
but now the intensity of production of a new atom is instead of just $m$ in the
previous steps now given by:
\be{gra44}
m + c\wh x^N_2(\frac{\ell}{2^k}), \ee where $\wh x^N_2$ is given in
(\ref{ddv}). Then define:
\be{agrev7} \wt
X^N_2(\frac{\ell +1}{2^k})=T_{\frac{1}{2^k}} \wt
X^N_2(\frac{\ell}{2^k}) + Y^{N,k}_\ell, \quad k \in \N,
\ee
which defines, for each $N$ and a fixed parameter $k \in \N$,
a piecewise constant atomic measure-valued process which we denote by
\be{agrev6b}
(\gimel^{N,k,m}_t)_{t \geq 0}.
\ee
We will now first investigate the weak convergence of this process, in
the parameters $k,N$ tending to infinity. First we focus on the weak convergence w.r.t.
the topology of the
weak convergence of finite measures on the state space and then we pass to the (stronger) one
w.r.t. the weak atomic topology on the state space.

Consider first $N \to \infty$.
For given $k$ and $t$, as $N\to\infty$, the random variable  $\gimel^{N,k,m}_t$,
converges  weakly even in $\mathcal{M}_a([0,1])$ to $
\gimel^{\infty,k,m}_t$ as proven in Lemma \ref{L.step2} (c).  We can
then obtain the convergence of the finite dimensional
distributions of $(\gimel^{N,k,m}_t)_{t \geq 0}$ to
those of $(\gimel^{\infty,k,m}_t)_{t \geq 0}$.
Tightness and convergence is obtained as we shall argue at the end of this proof
so that we have as $N \to \infty$ weak convergence of processes to

\be{agrev7b}
 (\gimel^{\infty,k,m}_t)_{t \geq 0}.
\ee

Next we let $k \to \infty$ and consider the weak topology on the state space.
Note first that the random variable in (\ref{agrev7b}) is by construction
{\em stochastically increasing in $k$},
 since passing from $k$ to $k+1$ we pick up additional contributions  from newly formed
 atoms. We can then show that as $k\to\infty$ the
 $\gimel^{\infty,k,m}_t$  converge
in distribution to $\gimel^m_t$. Similarly we can take the limit
$k \to \infty$ in $\gimel^{N,k,m}$ to get $\gimel^{N,m}$.

We can then use the triangle inequality for the Prohorov metric
$\rho(\gimel^{N,m}_t,\gimel^m_t)$  to get our result. Namely from the
weak convergence of $\gimel^{N,k,m}_t$ to $\gimel^{\infty,k,m}_t$ as
$N \to \infty$, the convergence as $k \to \infty$ of   $\gimel_t^{N,k,m}$ to $\gimel_t^{N,m}$
as well as the convergence as $k \to \infty$ of $\gimel_t^{\infty,k,m}$ to $\gimel_t^m$.
Using then the Markov property and a standard argument we can
verify the convergence of the finite dimensional distributions of
$\gimel^{N,m}_t$ to the finite dimensional distributions of
$\gimel^{m}_t$.

It now remains to check for
$(\gimel^{N,m})_{t \geq 0}$ as $N \to\infty$
the tightness condition of the laws and then convergence actually holds also in
$C([0,\infty),(\mathcal{M}_a([0,1]), \rho_a))$. We first prove tightness
with respect to the {\em weak topology} on the state space,
i.e. we first prove tightness in
the path space $C([0,\infty),\mathcal{M}_f([0,1]))$.

We write the state of the process as
\be{gra44b}
\gimel^{N,m}_t = \suml^N_{i=1} a^N_i (t) \delta_{x_i^N(t)}.
\ee
Then it suffices to show that the laws of
\be{gra44b2}
\{\sum_{i=1}^N f(a^N_i(t),x^N_i(t)), \quad N \in \N\}
\ee
are tight for a bounded continuous function $f$ on
 $[0,1]\times[0,1]$.  But the $x^N_i(t)$ do not change in time and
 the $a^N_i(t)$ are semimartingales with bounded characteristics and
are therefore tight by the Joffe-M\'etivier criterion.

To get to the convergence in path space based on the
{\em weak atomic topology} on the state space of the process, we next note
that since the joint distributions of the ordered
atom sizes at fixed times converge, the locations are constant in time
(as long as they are charged),  we have convergence of the
finite dimensional distributions in the weak atomic topology. The
 verification of the condition for tightness in $C([0,\infty),\mathcal{M}_a([0,1]))$
 (\ref{EKcriterion}) follows as in the proof of Theorem 3.2 in
 \cite{EK4}.
\end{proof}

\subsubsection{The long-term behaviour of limiting droplet dynamics\\
(Proof of Proposition \ref{CSB-longtime})}
\label{sss.proofCSB}

We now investigate the behaviour of the limiting $(N \to \infty)$
droplet dynamic, where again the dual process is the key tool, now
the dual process of McKean-Vlasov dynamics and certain subcritical
Fisher-Wright diffusions.

\begin{proof} {\bf of Proposition \ref{CSB-longtime}}

We prove separately the three parts of the proposition.

a)
Let $(\gimel^m_t)_{t \geq 0}$ be defined as in Propositions \ref{CSB}.  We
begin by deriving a formula for the first moment
\be{gra40}
m(t)=E[\gimel_t^m([0,1])], \ee given $\gimel^m_0(\cdot)$ with
$\gimel^m_0([0,1])<\infty$, which is accessible to an asymptotic analysis.

First we introduce a key ingredient in the formula for $m(\cdot)$
and we obtain its relation to the exponential growth rates
$\alpha^\ast$ respectively $\alpha$.

Let \be{dd33}
 f(t)=\int_{W_0} w(t)\mathbb{Q}(dw)= \lim_{\ve\to 0}\frac{1}{\ve}E_\ve[\wh x_2(t)],
\ee  where $\wh x_2(\cdot)$ denote the solution of (\ref{onedim})
and $E_\ve$ refer to the initial state $\ve$. We
now obtain an expression for (\ref{dd33}).

 To compute the first moment of
$\wh x_2(t)$ we use the dual representation for the moments of the
SDE  (\ref{onedim}). The dual process to be used then is $( \eta_t, \CF^+_t)_{t \geq
0}$ with one initial factor $1_1(\cdot)$.  The dual particle
process $(\eta_t)_{t \geq 0}$ then is effectively a birth and
death process denoted $(\wt D_0(t))_{t \geq 0}$ on $\N_0$, with a dynamic given by
the evolution rule that the process can jump up or down by 1 (birth
or death) and the rates in state $k$ are given by:
\be{dd33a}
\mbox{birth rate } sk, \mbox{ death rate
} ck +\frac{d}{2}{k(k-1)} \ee and \be{dd33b} \mbox{initial state }
\wt D_0 =1. \ee

Then the dual expression is given by \be{dd33d}
 E_\ve[\wh x_2(t)]= 1-E[(1-\ve)^{\wt D_0(t)}].
 \ee
Therefore \be{dd38}
f(t)=\int_{W_0}
w(t)\mathbb{Q}(dw)=\lim_{\ve\to 0}\frac{1}{\ve}E_\ve[\wh x_2(t)]=
\lim_{\ve\to 0}\frac{E[1-(1-\ve)^{\wt D_0(t)}]}{\ve}=E[\wt D_0(t)].
\ee

Now let $\alpha^\ast$ be chosen so that
 \be{dd34}
c\int_0^\infty e^{-\alpha^\ast r}f(r) dr =c\int_0^\infty
e^{-\alpha^\ast r}E[\wt D_0(r)] dr=1.
\ee
Recalling that by (\ref{ang1abx}), (recall the death rate of
$D_0(t)$ is zero if $k-1$)
\be{dd34a}
\int_0^\infty
e^{-\alpha^\ast r}E[\wt D_0(r)] dr= \int_0^\infty
e^{-\alpha^\ast r}E[D_0(r)1_{D_0(r)\geq 2}] dr,
\ee
it follows   that $\alpha^\ast=\alpha $ with $\alpha$ as defined in (\ref{ang1}), (\ref{ang1ab}).

Now we are ready to write down the final equation for $m(\cdot)$.
Consider
\be{dd35b} \gimel^{\ast}_0 (t)=(\sum y_i(t)\delta_{a_i})_{t \geq
0}, \ee
where $y_i(t),\;\sum y_i(0)<\infty,$ are those realizations of
the (independent) Fisher-Wright diffusions (\ref{onedim}) which represent the
atoms that were present at time 0  without immigration.
Furthermore assume that $\suml_i y_i(0)<\infty$. Then we consider
the process $(\gimel^m_t)_{t \geq 0}$ defined in Proposition \ref{CSB}
with $\gimel^\ast$ therein as given in (\ref{dd35b}).

By taking expectations  in the stochastic equation defining
$\gimel^m$, namely (\ref{ZL2m}), we get if we define
\be{gra41} g(t)=E[\suml_i y_i(t)],
\ee
the renewal equation
\bea{dd34b} m(t)=E[\gimel^m_t([0,1])]&=E[\suml_i y_i(t)]+
\int_0^t E[({m}+c\gimel^m_r([0,1]))\int_{W_0}
w(t-r)\mathbb{Q}(dw)]dr\\&  = g(t)+m \int_0^t  f(r)dr +\int_0^t
cm(r) f(t-r)dr. \nonumber
\eea

The last step is now to analyse the growth behaviour of
$(m(t))_{t \geq 0}$ using renewal theory.
Multiplying through equation (\ref{dd34b}) by $e^{-\alpha^\ast t}$ we
get an equation for $(e^{-\alpha^\ast t}m(t))_{t \geq 0}$
in terms of $(e^{-\alpha^\ast t} f(t))$ and $(e^{-\alpha^\ast t}g(t))$
as follows:

\be{dd36}
e^{-\alpha^\ast t}m(t)= e^{-\alpha^\ast t} [(g(t)+m \int_0^t
f(r)dr)] + c\int_0^t (e^{-\alpha^\ast r}m(r))(e^{-\alpha^\ast(t-r)}f(t-r))dr.
\ee
This equation in $e^{-\alpha^\ast t} f(t)$ and $e^{-\alpha^\ast t} g(t)$
has the form of a renewal equation.

We want to apply now a {\em renewal theorem} to this equation.
Define  $R$ by
\be{dd35}
 R= c\int_0^\infty re^{-\alpha^\ast r}f(r)dr.\ee
Since $f(t)$ and $g(t)$ are continuous and converge to $0$
exponentially fast, it can be verified that $a(u)= e^{-\alpha^\ast
u}\left[ g(u)+m \int_0^u f(r)dr\right]$ is directly Riemann
integrable. Therefore by the renewal theorem (\cite{KT}, Theorem
5.1) we obtain from (\ref{dd36}) that:
 \be{dd37}
 \lim_{t\to\infty} e^{-\alpha^\ast t}m(t)=
\frac{1}{R}{\int_0^\infty e^{-\alpha^\ast u}\left[ g(u)+m \int_0^u
f(r)dr\right]du}\in (0,\infty). \ee
Hence $e^{-\alpha^\ast t} E[\gimel^m(t)]$ converges as $t \to \infty$
to the r.h.s. of (\ref{dd37}) which concludes the proof of
part (a) of the proposition.

(b) Let
\be{gra43}
 \wh{x}_2(t)=\gimel^m_t([0,1]).
\ee
The convergence in distribution of $e^{-\alpha t} \wh {x}_2(t)$
as $t\to\infty$ follows since by a) the laws of
$\{e^{-\alpha t} \wh x_2(t), \quad t \geq 0\}$
form a tight family and we have to exclude only that several limit
points exist. We know that for any $\tau >0$
\be{dd39}
\lim_{t\to\infty} E([e^{-\alpha(t+\tau)}\wh{x}_2(t+\tau)-e^{-\alpha
t}\wh{x}_2(t)]^2)=0,
\ee
which is proved in Proposition \ref{P.DetGroDro-new}.
(Recall here that over a finite time horizon $[0,\tau]$ the
mutation can be ignored and hence the proposition applies).
We have to strengthen this statement (\ref{dd39}) to
\be{dg74}
\liml_{t \to \infty} (\supl_{\tau >0}
E[e^{-\alpha(t+\tau)} \wh x_2(t+\tau) - e^{-\alpha t}
\wh x_2(t))^2])=0.
\ee
This extension is provided in Corollary \ref{C.dbetatimes} of Subsubsection
\ref{sss.detdropgro} where it is proved using moment calculations.  This implies that there exists $\mathcal{W}^\ast$ such that
\be{add22}
\wh x_2(t)\longrightarrow \mathcal{W}^\ast\text{   in  }L^2.\ee
The non-degeneracy of the limit $\mathcal{W}^\ast$ follows from (c).

(c)  This claim on the variance of $\CW^\ast$ follows from Proposition \ref{P.Growthdrop2}  in
Subsubsection \ref{sss.dropletgrow}.

\end{proof}

\begin{remark}

The proof of (c) demonstrates that  the randomness arises from the
early rare mutation events.

\end{remark}

\subsubsection{ Proof of Proposition \ref{P.DF}}
\label{sss.proofPDF}

This proof builds mainly on the previous subsubsection and on calculations
we shall carry out in Subsection \ref{ss.relationsW}.

(a) This follows from Proposition \ref{P.lt}

(b) Recall that we know already that
$e^{-\alpha t} \wh x^N_2(t)$ converges in law as $N \to \infty$
to a limit $e^{-\alpha t} \wh x_2(t)$ and hence by Skorohod embedding and
our moment bounds we have convergence on some $L^2$ space, furthermore $e^{-\alpha t} \wh x_2(t)$,
converges in $L^2$ as $t \to \infty$ to a limit $\CW^\ast$ as proved in (\ref{add22})
above. Here we shall show that
\be{dd39b}
\liml_{t \to \infty}
\liml_{N \to \infty} E [(e^{-\alpha t} \wh x_2^N(t) - e^{-\alpha
t_N} \wh x_2^N (t_N))^2] =0 \ee and then we use (\ref{CSB3}) to get
the claim. The above relation follows via Corollary \ref{C.AG0} which is proven by moment calculation
which requires subtle coupling arguments which we shall develop in
Subsubsections  \ref{sss.dropletgrow} and \ref{sss.detdropgro}.

\subsubsection{Some explicit calculations}
\label{sss.explicitcal}

The entrance law from 0 for the type-2 mass can also be derived
using explicit formulas for densities of the involved diffusions
both for Fisher-Wright and branching models.
We collect this in the following two remarks.

\begin{remark}
Consider the Fisher-Wright given by the SDE \bea{ag46} &&dx(t)=
\frac{m}{N}[1-x(t)]+c[\frac{y(t)}{N}-x(t)]dt+ sx(t)(1-x(t))dt
+ \sqrt{d \cdot x(t)(1-x(t))}dw(t),\\&& x(0)=p, \nonumber\eea where
$y(t)= \wh x^N_2(t)$ is inserted as an external signal. By Kimura
\cite{Kim1} if $m=c=s=0$, then there is  no $N$-dependence and the
density of the solution at time $t$ is
\be{agrev4}\begin{array}{l}
f(p,x;t) \\
\hspace{1cm} =\sum_{i=1}^{\infty}
p(1-p)i(i+1)(2i+1)F(1-i,i+2,2,p)\cdot F(1-i,i+2,2,x)e^{i(i+1)dt/2}
\\
\hspace{1cm}  =6p(1-p)e^{-dt}+30p(1-p)(1-2p)(1-2x)e^{-3dt}+\dots,\nonumber
\end{array}
\ee
where $F(\alpha,\beta,\gamma;x)$ is the hypergeometric function.

Consider now the neutral Fisher-Wright diffusion, $s=0$ and
constant $y(t)=m$, $c>0,\; d=1$.  Crow-Kimura (see \cite{GRD},
Table 3.4, or Kimura \cite{Kim1}, (6.2)) obtained the density at
time $t$ \be{wflaw}\begin{array}{ll}
\wt f^N(y,x,t)=& x^{\mu-1}(1-x)^{\nu-\mu-1}\\
&\times\sum_{i=0}^\infty
\frac{(\nu+2i-1)\Gamma(\nu+i-1)\Gamma(\nu-\mu+i)}{i!\Gamma^2(\nu-\mu)\Gamma(\mu+i)}\times
\wt F_i(1-x)\wt F_i(1-y)\exp(-\lambda_i t),\nonumber
\end{array}\ee
where the initial value of the process is $y$ and \be{ag44} \mu
=\frac{2m}{N}, \nu=2(\frac{m}{N}+c),
\lambda_i=i[2(\frac{m}{N}+c)+(i-1)]/2,
\ee
\be{ag45}
\wt F_i(x)=F(\nu+i-1,-i;\nu-\mu;x), \ee where $F$ is the socalled
hyper-geometric function that is that solution of the equation
\be{agrev1}
 x(1-x)F^{\prime \prime}+[\gamma-(\alpha+\beta+1)x]F^{\prime \prime}-\alpha\beta F=0,
 \ee
which is finite at $x=0$.

Setting $y=0$ and taking $N\to\infty$, then for every $t>0$ and $x
\in [0,1]$, taking  $\wt f^N (x,t) = \wt f^N(0,x,t)$, we get  that
there exists a function $\wt f^\infty (x,t)$ which is for every
$t$ the density of a $\sigma$-finite measure such that:
\be{agrev1b} N\wt f^N (x,t) \la \wt f^\infty (x,t). \ee The limit
defines an entrance type law. Note that $\wt f^N$ has the same
behaviour near zero as $N\to\infty$ as the branching density below
in (\ref{G1}).

Consider now the case $s>0$. Then
selection term then leads to a new law which is absolutely
continuous w.r.t. the neutral one and the density is given by the
Girsanov factor \be{agrev5}
  L_t= \exp\left[\frac{s}{{d}}\int_0^t
\sqrt{d \cdot x(s)(1-x(s))}dw(t)-\frac{s^2}{2d^2}\int_0^t d \cdot
x(s)(1-x(s))ds \right]. \ee Therefore  the selection term does not
change the asymptotics since it involves a bounded Girsanov
density (uniform in $N$).

Alternatively to the above argument we can take the Kimura
solution (4.10) of the forward equation and take the limit
$\frac{\phi(p,x;t)}{p}$ as $p\to 0$.$\qquad \square$
\end{remark}

\begin{remark}
For comparison, consider  the  branching approximation   where there is no competition between excursions
and we obtain  for initial state 0 the Gamma density \be{G1}
f^N(x,t)=\frac{c}{N}x^{\frac{m}{N}-1}e^{-x/e^{ \wt st}},\;x\geq 0.
\ee Then $N f^N$ converges to the density $x^{-1} e^{-x/e^{\wt
st}}$ on $\R^+$ of a $\sigma$-finite measure. Moreover
\be{gra45}
N\cdot \int_0^\ve xf^N(x,t)dx \to 0 \text{ as } \ve\to 0 \text{
uniformly in } N. \qquad \square
\ee

\end{remark}

\subsubsection{First and second moments of the droplet growth constant $\CW^\ast$}
\label{sss.dropletgrow}


Continuing the study of the droplet growth we
now turn to times in $[T_0,T_N]$ as given in (\ref{ddy}) and
(\ref{ddx}). In the sequel we compute the limits as $N\to\infty$
of the first two moments of the total type-2 mass in the population
denoted $\widehat{x}_{2}^{N}(t_N)$, for
$t_N \to \infty, \quad t_N= o(\log N)$.
The present calculation allows in particular to show the existence
of $\CW^\ast$ and to determine its moments.

As before we use the dual representation to compute moments.  For
example, recall that
\be{gra45b}
 E[x^N_1(i,t_N)] = E[ \exp (-\frac{m}{N} \intl^{t_N}_0 \Pi^{N,1}_s ds)],
\ee
where $\Pi^{N,1}_t$ is of the total number of particles in the dual cloud in
$\{1,\dots,N\}$ starting with one factor $1_{\{1\}}$ at site $i$ at time $t=0$. As
developed earlier this cloud is described by a CMJ process and has the form
\be{gra45c}
\Pi^{N,1}_t = \sum_{j=1}^{W_N(t)e^{\alpha t}}\zeta^N_j(t),
\ee
where $W_N(t) \to W(t)$ as $N \to \infty$ and
$W(t)\to W$ a.s. as $t\to\infty$ and the processes at the different sites satisfy
$\CL [(\zeta^N_j(t))_{t \geq 0}] \La \CL[(\zeta_j(t))_{t \geq 0}]$ as $N \to \infty$,
where  the $(\zeta_j(t))_{t \geq 0}$ for $j \in \N$ are all evolving independently and are given by the
birth and quadratic death process at an occupied site with the same dynamic but
starting from one particle at the random time, when the site is first occupied.
 The dependence between the different occupied
sites arises through these random times even in the $N \to \infty$ limits.

If we consider second moments we have to start with two
particles and then get the analog of (\ref{gra45b}) but in
the representation of the r.h.s. of this formula given in
(\ref{gra45c}), $W_N(t)$ and $\zeta^N$ have to be replaced
accordingly and in the limit $N \to \infty$  the $W=W^{1,1}$ is replaced by $W^{1,2}$.

We have noted above
that if $t_N=o(\log N)$ then the probability that a collision
between dual particles due to a migration step
occurs up to time $t_N$ goes to 0 as $N\to\infty$. {\em The new
element that arises in this subsubsection is the fact that in the
computation of the higher moments of $\wh x^N_2(t_N)$, these
collision events must be taken into account.} The reason for this
is that we will compute moments of the mutant mass $\wh x^N_2(t)$
in terms of the moments of $\wh x^N_1(t)$.
We see below from the formula (\ref{AGr1})  that in the expansion of
$k$-th moments of the total mass powers of $N$ appear, which
forces us to consider in the expansions of these moments as well as higher
order terms in $N^{-1}$. Hence we have to analyse more carefully
the overlaps in the dual clouds corresponding to different ancestors
which  have probabilities which go to zero but not fast enough
to be discarded.

The key tool we  work with is the coloured particle systems
(WRGB) introduced in Subsubsection \ref{sss.prep2}.

\bigskip

The main result of whose proof this subsubsection is devoted to is:

\beP{P.Growthdrop2}{(Limiting growth constant of droplet)}

Let $t_N \to \infty$  and \be{AGr6b} \limsup\limits_{N\to \infty}
\frac{t_N}{ \log N}<\frac{1}{2\alpha}. \ee (a)  The family

\be{agrev29} \{e^{-\alpha t_N}\wh{x}^N_2(t_N):N\in \mathbb{N}\}
\mbox{ is tight, } \ee with non-degenerate limit points
$\mathcal{W}^\ast$ with \be{agrev30}
 0< E[\mathcal{W}^\ast] <\infty,\quad 0<\text{Var}[\mathcal{W}^\ast]<\infty.
\ee

(b) Furthermore: \be{agrev28}
  \lim_{N\to\infty}e^{-\alpha t_N}E[\wh{x}^N_2(t_N))] = m^\ast E[W_1],
\ee
\be{agre61a}
\lim_{N\to\infty}e^{-2 \alpha t_N}
E[\wh{x}^N_2(t_N))^2] = (m^\ast)^2 (EW_1)^2 + 2m^\ast (EW_1)^2
\kappa_2,
\ee
where $\alpha$ is defined as in (\ref{ang1ab1}),
$\kappa_2 >0$ (which is given precisely in (\ref{dd658++})) and
\be{agrev28b}
m^\ast =  \frac{m (\alpha + \gamma)}{c}, \mbox{
with } b=1+\frac{\gamma}{\alpha}.
\ee

(c) The third moment satisfies
\be{agrev31c}
\limsup_{N\to \infty}
E[e^{-3 \alpha t_N} (\wh x^N_2 (t_N))^3] < \infty,
\ee
so that the first and second moments of $\CW^\ast$ are given by the formulas
on the r.h.s. of (\ref{agrev28}) and (\ref{agre61a}). $\qquad \square$
\end{proposition}

The proof of the main result is given in the rest of this subsubsection
which has five parts, namely we require four tools to finally carry out the proof.

(1) First we  give a preparatory lemma \ref{generalities} that
 gives asymptotic expressions for the moments of
$\wh x^N_2 (t_N)$ in terms of the moments of $\wh x^N_1 (t_N)$.

(2) In order to calculate the moments of $\wh x^N_1 (t_N)$ we must determine the growth of dual clouds
starting with one or two initial factors. This is stated in Proposition  \ref{L.Indual}
and is then proved.

(3) Here we introduce a multicolour system to analyse the two dual clouds
and their interaction starting both from one particle (at different sites).

(4) Then we present in Lemma \ref{L.dd917} the tools to obtain
the asymptotics of the expectations of the exponentials involving the occupation times of the dual clouds.

(5) Finally we complete the proof of the above main result, the Proposition \ref{P.Growthdrop2},
using the dual expressions at the end of the section and the results from the previous Parts 1-4.

{\bf Part 1} {\em (Moment formulas)}

\beL{generalities}{($k$-th moment formulas)}

 (a) We have the $k$-th moment formula for the total type-2 mass:
 \be{AGr1}\begin{array}{ll}
 \bar m^{N,2}_k(t_N) &:= E[(\wh x^N_2(t_N))^k]=E\left[\left(
\suml^N_{i=1} x^N_2 (i,t_N)\right)^k\right]
\\& =
\sum_{\ell =1}^k\left(\begin{array}{c}  k \\  \ell \\ \end{array}
\right)(-1)^\ell N^{k-\ell} \left( \bar
m^{N,1}_\ell(t_N)\right) , \quad k=1,2,\cdots
 \end{array}\ee
where $\bar m^{N,1}_\ell$ is the $\ell$-th moment for the total type-1 mass.

(b) The term on the r.h.s. of (\ref{AGr1}) is calculated as follows.
Let $q=(q_1,\cdots,q_j)$ and let $\Pi^{N,q}$ be the dual starting
with $q_i$-particles at site $i=1, \cdots,j$.
Define

\be{AGr4}
\mu_{q_1,\dots,q_{j};q}(t_N) = E\left[\exp \left( -\frac{m}{N}
\intl^{t_N}_0 \Pi^{N,(q_1, \cdots,q_j)}_s ds\right) \right].
\ee
Then
\be{AGr2}
\bar m^{N,1}_q(t_N)  =\suml_{j=1}^{q}
\left(\begin{array}{c} N \\
  j \\
\end{array}  \right)\sum_{q_1+\dots q_j=q}\frac{q!}{q_1!\dots q_j!}
\mu_{q_1,\dots,q_j;\;q}(t_N). \qquad \square
\ee
\end{lemma}

\begin{proof}{\bf of Lemma \ref{generalities}}

(a) First recall that
\be{agrev8}
\widehat{x}_{1}^{N}(t)
=\sum_{i=1}^{N}x_{1}^{N}(i,t),\quad \widehat{x}_{2}^{N}(t)
=\sum_{i=1}^{N}x_{2}^{N}(i,t)=N-\widehat{x}_{1}^{N}(t) \ee and
therefore

 \be{AGr1x}\begin{array}{ll}
 \bar m^{N,2}_k(t_N) &:= E[(\wh x^N_2(t_N))^k]=E\left[\left(
\suml^N_{i=1} x^N_2 (i,t_N)\right)^k\right]
\\& =
\sum_{\ell =1}^k\left(\begin{array}{c}  k \\  \ell \\ \end{array}
\right)(-1)^\ell N^{k-\ell} E\left[\left(
\sum_{i=1}^Nx^N_1(i,t_N)\right)^\ell\right] \quad , \quad
k=1,2,\cdots
 \end{array}\ee
Recalling  the definition of $\bar m^{N,1}_\ell(t)$ the claim follows.

(b) Moreover we can express the moments of the mass of type 1
appearing on the r.h.s. of (\ref{AGr1x}) by expanding according to
the {\em multiplicity} of the occupation of sites in the mixed
moment expression, namely:
\bea{Gr1bx}
&&\bar m^{N,1}_k(t_N) = E\left[ \left(
\sum_{i=1}^{N}x_{1}^{N}(i,t)\right) ^{k}\right]\\&&
=\sum_{i_{1}\neq i_{1}\neq\dots\neq i_{k}}E\left[
x_{1}^{N}(i_{1},t)\dots x_{1}^{N}(i_{k},t)\right] +\dots+
\sum_{i=1}^{N} E[(x_{1}^{N}(i,t))^{k}]. \nonumber \eea

The moments on the r.h.s. above can be computed  as follows.
Define \be{AGr2bx} \mu_{q_1,\dots,q_j;\;q}(t_N)=E\left[
\prodl_{\ell=1}^{j} (x^N_1 (i_\ell, t_N))^{q_\ell}\right],\quad
i_\ell \text{ distinct},\quad \suml^j_{\ell=1}  q_\ell =q.\ee Then
 \be{AGr2x} \bar
m^{N,1}_q(t_N)  =\suml_{j=1}^{q}\left(\begin{array}{c}
  N \\
  j \\
\end{array}  \right)\sum_{q_1+\dots q_j=q}\frac{q!}{q_1!\dots q_j!}
\mu_{q_1,\dots,q_j;\;q}(t_N). \ee Next use the duality relation to
express $\mu_{q_1, \cdots, q_j ;q}$.
Define
\be{AGr3x}
\Pi^{(q_1, \cdots,q_j),N}_{t_N}, \ee
as the total number of dual particles at time $t_N$  starting the
particle component $\eta_t$ of the dual process $(\eta_t,
\CF^{++}_t)$ with
\be{agrev9} \mbox{ $q_i$ particles at distinct
sites }i=1,\cdots,j. \ee

Then  write \be{AGr4b}
\mu_{q_1,\dots,q_{j};q}(t_N) = E\left[\exp \left( -\frac{m}{N}
\intl^{t_N}_0 \Pi^{N,(q_1, \cdots,q_j)}_s ds\right) \right]. \ee
q.e.d.

\end{proof}


\bi

{\bf Part 2} {\em (Growth of dual clouds)}

Looking at the r.h.s. of formula (\ref{AGr1})-(\ref{AGr2}) we see that
for calculating higher moments we have to study
the asymptotic (as $N \to \infty$) behaviour of the dual system starting with {\em more} than
one particle. Of particular importance is to estimate the {\em overlap} in the dual particle
population of different families due to different ancestors at the starting time
of the evolution since this will allow to control the terms of order $N^{-k}$
which are needed to compensate for the $N \choose k$ term on the r.h.s.
of (\ref{AGr1}). This is addressed in the form needed for the calculation of first and second
moments in the following definition and proposition.

\beD{D.dcloud}{(Dual clouds for moment calculations)}

(a) To calculate first and second moments at one time we calculate the
dual particle number asymptotically as follows if we start with one or
two particles. We choose $W_a$ as:
\begin{itemize}
\item $W_a=W_1$ where $W_1=W^{1,1}=W$ defined by (\ref{ang2b}) starting with one
$1_{\{1\}}$ factor at  time $0$
\item  $W_a= W^{1,2} = W_1+W_2$ where $W_1,W_2$ are independent copies of the
random variable $W$  arising from one factor at each of two different sites,
\item $W_a=W^{2,1}$ which is the
corresponding random variable when the CMJ process is started with
the two factors $1_{\{1\}} \otimes 1_{\{1\}}$ at time $0$ at the same site.
\item
Similarly we write $W_a(\cdot)$ for the function $(W_a(t))_{t \geq 0}$
arising from $e^{-\alpha t} K^a_t$.
\end{itemize}

Let $\xi$ be a time-inhomogeneous Poisson process on $[0,\infty)$ with
intensity measure given by
\be{agr25}
\frac{(\alpha + \gamma)}{c} \frac{(W_a(s))^2e^{2\alpha s}}{N} ds.
\ee
Note $\xi$ depends on the  chosen  $a$ but only through $W_a$.

Let
\be{agr26}
 \wt g(u,s)\ee
(which evolves independently of the parameter $N$) be the size of that green
cloud of sites of age $u$ produced by one collision (creating a red particle)
at time $s$. Then $\wt g(t,s)$ has a law depending only on $t-s$, namely
\be{add23}
\wt g(t-s)=e^{\alpha(t-s)}W^g(t-s),\ee
where $W^g(t-s)\to W^g$ as $t-s\to \infty$ and $W^g$ is defined by this relation.

Furthermore introduce for a realization of $\xi$ a process $\mathfrak{g}$ depending on
the parameters $N$ and whose law depends on the choice of $a$ but only through
$W_a(\cdot)$ and which describes the size of the complete green population at time $t$:
\be{agrev16x}
\mathfrak{g}(W_a(\cdot),N,t)=\int_0^{t}\wt{g}(t,s)\xi(ds).
\ee

Analogously to $\wt g, \mathfrak{g}$ we define the pair
\be{agr27}
\wt g_{\mathrm{tot}}, \mathfrak{g}_{\mathrm{tot}}
\ee
as the number of particles in the green cloud, instead of the number of sites.

(b) To compute the joint second moment at two different times we consider
\be{agrev16x4}
W_{a}= W_{j_1}+W_{j_2,s},
\ee
which is the corresponding random variable when
the CMJ process is started with the factor $1_{\{1\}}$ located at $j_1$
at time 0 and an additional factor $1_{\{1\}}$ is added at location
$j_2$ at time $s>0$.

Let $\xi^{(1,2)}$ be an
inhomogeneous Poisson process on $[0,\infty)$ with intensity measure
\be{agrev16x6}
\left[(\alpha+\gamma)\frac{W_1(u)W_2(u)e^{2\alpha u}}{N}\right]du.
\ee

Then with
\be{agrev16x5}
t=t_N,
\ee
we consider time $t+s$. Furthermore
set (recall (\ref{agrev16x5}))
\be{dd658+}
\mathfrak{g}_{(1,2)}(W_1 (\cdot) ,W_2 (\cdot), N,t,t+s)= \int_s^{t}\wt g(t_N,u)\xi^{(1,2)}(du).
\ee

We denote by $K^{N,(j_1,0),(j_2,s)}_t$
the number of occupied sites in the dual particle system starting with a particle at
$j_1$ at time 0 and one at time $j_2$ at time $s$. $\qad$
\end{definition}
Now we can state the main result on the growth of the dual clouds.

\bigskip

\beP{L.Indual} {(Dual clouds)}

Let
\be{agr23}
t=t_N, \mbox{  where } t_N\to\infty,\quad t_N=o(\frac{\log
N}{2\alpha}),\ee
more precisely,
\be{add24}  2\alpha t_N -\log N \to -\infty \mbox{ as } N\to\infty.\ee

Then consider the growth of the
dual cloud in various initial states indicated by the subscript a.

(a) We have

\be{agr28}
K^{N,(a)}_t = W_a(t) e^{\alpha t} -\mathfrak{g} (W_a(\cdot), N,t) + \CE_a(N,t),
\ee

\be{agrev16} \Pi_{t}^{N,(a)}
= \frac{1}{c}(W_a)(t))(\alpha +\gamma) e^{\alpha t}- \mathfrak
{g}_{\mathrm{tot}} (W_a(\cdot),N,t) + \mathcal{E}_{a, \mathrm{tot}}(N,t),
\ee
where we have for some $\kappa_2 \in (0,\infty)$ \be{agrev16xx}
\lim_{N\to\infty}Ne^{-2\alpha t_N}E[ \mathfrak {g}_{\mathrm{tot}}(W_a(\cdot),N,t_N)]= \kappa_2 E[(W_a)^2]
 \ee and the error term
$\CE_a$ or $\mathcal{E}_{a, \mathrm{tot}}$ satisfies that

\be{agrev16xxx} \lim_{N\to\infty}E[Ne^{-2\alpha t_N}|\mathcal{E}_a(N,t_N)|]= 0.
\ee

(b) We have (as $N \to \infty$, recall (\ref{agrev16x5})):
\be{gra62}
 K^{N,(j_1,0),(j_2,s)}_{t+s} \sim (W_{1}e^{\alpha
(t+s)}+W_2 e^{\alpha t})- \mathfrak {g}(W_1 (\cdot), W_2 (\cdot),N,t,t+s)
+O(\frac{e^{3\alpha t}}{N^2}),
\ee
and
\be{agr34}
\Pi^{N,((j_1,0),(j_2,s))}_{t+s}
= \frac{\alpha +\gamma}{c} (W_1 e^{\alpha(t+s)} + W_2 e^{\alpha(t+s)})
 - \mathfrak{g}_{\mathrm{tot}} (W_1 (\cdot), W_2  (\cdot), N,t,t+s) + O(\frac{e^{3 \alpha t}}{N^2}),
 \ee

\noi
where $ \mathfrak {g}_{\mathrm{tot}}(W_1 (\cdot),W_2 (\cdot),N,t,t+s)$
satisfies the two conditions
\be{gra62b}\begin{array}{ll}
\mathfrak {g}_{\mathrm{tot}}(W^{1,2}(\cdot),N,t,t+s)=&
\mathfrak{g}_{\mathrm{tot}}(W_1(\cdot),N,t+s)
+\mathfrak{g}_{\mathrm{tot}}(W_2(\cdot),N,t)\\
&+\mathfrak{g}_{(1,2),\mathrm{tot}}(W_1(\cdot),W_2(\cdot),N,s,t),
\end{array}
\ee
and
\be{dd658++}
\lim_{N\to\infty}N\cdot {e^{-\alpha(2t_N+s)}}E[\mathfrak{g}_{(1,2),\mathrm{tot}}(W_1(\cdot),W_2(\cdot),N,t_N,t_N+s)]=2\kappa_2
E[W_1]E[W_2].
\ee
Similarly such relations hold for $\mathfrak{g}$ instead of
$\mathfrak{g}_{\mathrm{tot}}$ with $\kappa_2$ replaced by the appropriate
different positive constant.

(c) Now consider the case in which $N \to \infty$ but $t$ is independent of $N$ (that is,
we no longer take  $t_N\to\infty$). Then  we do not replace $W_i(\cdot)$ by $W_i$ and $\alpha(t),\gamma(t)$ by $\alpha,\gamma$.
Then in (\ref{gra62}), (\ref{agr34})
$W_1 \cdot (\alpha+\gamma)$ is replaced by $W_1(t+s)(\alpha(t+s)+\gamma(t+s))$ and $W_2 \cdot (\alpha+\gamma)$
 is replaced by $W_2(t)(\alpha(t)+\gamma(t))$.
$\qad$
\end{proposition}

\begin{remark}
Note that in the expansion $(W_1+W_2)^2=W_1^2+2W_1W_2+W_2^2$ we
can associate the term $2W_1W_2$ with collisions of a particle from
one cloud with sites occupied by the other cloud.  In order to
asymptotically evaluate these collision terms we could keep
track of the white, red, green, blue particles from the two sites
by labelling them, for example, by
\be{agrev31a}
(W_-,R_-,G_-,B_-), \;(W_+,R_+,G_+,B_+),
\ee
respectively, and then considering the
collisions of $W_+$ with $W_-$, etc.
\end{remark}

\begin{remark}
Note that the existence of $\kappa_2$ follows exactly in the same way as the argument leading up to (\ref{kappa}) (proof of existence of $\kappa$) in Subsubsection \ref{sss.prep2}. The only difference here is that it involves the growth of the number of green sites arising from either self collisions or between collisions of two different clouds.
\end{remark}

\bi

\begin{proof}{\bf of Proposition \ref{L.Indual}}
{\bf (a)} In order to justify (\ref{agrev16}) we now reconsider the
development of the dual particle system on $\{1,\dots,N\}$
starting from one factor $1_{\{1\}}$ at location 1. We have used above
the fact that for the range $0\leq t\leq t_N$ where
$\lim_{N\to\infty} \alpha \frac{t_N}{\log N} <1$, the particle
system is described asymptotically as $N \to \infty$ in first order
by a CMJ process with Malthusian parameter
$\alpha$. Correction terms are analysed using the multicolour
particles system.

Moreover if we follow two families starting with one
factor $1_{\{1\}}$ at locations $1,2$ respectively, then the
probability that they collide, that is occupy the same site in
$\{1,\dots,N\}$, in this time regime goes to $0$ as $N\to\infty$.
This is essentially equivalent to the statement that the
probability that a green particle is produced in the WRGB system
in this time regime goes to zero as $N \to \infty$. However for
finite $N$ these probabilities are not zero and we now determine
them by expansion in terms of powers of $N^{-k}$.

Recall that the system of white and black particles
is a CMJ process as described
above in which the number of occupied sites has the form
\be{agrev19d1}
W(t)e^{\alpha t} \ee where  $\lim_{t\to\infty}W(t)=W$ and this
black and white system serves as an {\em upper bound} for the system of occupied sites in the
dual process. Hence if at
time $t$ a migrant is produced by this collection following the
dual dynamics, it moves to a randomly chosen point in
$\{1,\dots,N\}$.  Then the  probability that  it hits an occupied
site is therefore at most:
\be{agrev19d2}
 \frac{W(t)e^{\alpha t}}{N}. \ee Therefore a upper
bound for the process of collisions is given by an inhomogeneous
Poisson process $\xi(s)$ with rate \be{agrev19d3}
\frac{W^2(t)e^{2\alpha t}}{N}
\ee
and number of sites that have been hit this way in
the time interval $[0,t]$ is at most
\be{agrev19d4}
\int_0^t \frac{W^2(s)e^{2\alpha s}}{N}ds \leq \frac{1}{N}
(\sup_t W(t))^2 \cdot \frac{1}{2\alpha} e^{2\alpha t}
=
\frac{1}{N} (W^\ast)^2 \frac{1}{2 \alpha} e^{2\alpha t}.
\ee
Recall that $W^\ast = \supl_{t \geq 0} W(t) < \infty$ a.s. For bounding second moments we can also
show that $\sup_{t\geq 0}E[(W(t))^2]<\infty$.
In particular the number of sites at which collisions take place by
time $t_N$ is $o(N)$, i.e. has spatial intensity zero as $N \to \infty$.

Moreover the expected number of sites to be hit more than once
in the time interval $[0,t]$ goes to 0 as $N\to\infty$ at the order $N^{-2}$.
This order holds also for red-red collision
by migration in the multicolour particle representation.

Next we recall that when a white particle hits an occupied site it
produces a red particle at this occupied site. Subsequently if the
red particle is removed by coalescence with a white particle, then
a green particle is produced. Due to the coalescence process eventually no red particles will remain at this site.  Also due to migration eventually this site will revert to a pure white site.
During the ``lifetime'' of such a multicolour site  there is a positive probability that one or more green particles will migrate and then produce a growing cloud.  (Recall that the lifetime  has finite exponential moment.)

 Once such a new site containing white and red particles is formed, a green particle
is created by coalescence unless either the white or the red
disappears before coalescence. The number of green particles to migrate from
such a special site is random with probability law that depends  on the number $\zeta(s)$ of the occupying
white  particles at the time $s$ of arrival  of the red
particle (asymptotically as $N \to \infty$ no further red particle
will arrive at this site before the coalescence or migration step
of the red particle or disappearance of all whites at the site).

 In other words  there is a positive
probability that the green particle will produce green migrants
before it coalesces with a red particle (such that only one red particle
remains and the green particle disappears (recall Remark \ref{rgremark})).
 We denote the probability $p^N(s)$ that at least one green migrant is produced during the lifetime
 of a special site due to a collision at time $s$ satisfies


\be{gra57b}
p^N(s) \la p(s) \mbox{ as } N \to \infty, \quad p(s) \la p>0 \mbox{ as } s \to \infty.
\ee
After the migration step the resulting green family, if non-empty,
grows according to the CMJ process. Combining the two facts
we see that the number of
green sites this family occupies at time $t$, denoted
$\wt g(t,s)$, satisfies:
\be{agrev19d5}
\wt g(t,s)=W^g(t,s) e^{\alpha (t-s)},
\ee
and $\lim_{t\to\infty}W^g(t,s) =W^g_s$ exists and converges to $W^g$ as $s \to \infty$.
Note that  the event $\{W^{g}=0\}$ occurs with the probability that a collision does not create a green migrant (see Remark \ref{nogr} below) and otherwise
\bea{add25}
&&W^g(t,s):=\mbox{ is the random growth constant for all the sites occupied by  }\\
&&\qquad\qquad\quad\mbox{ green migrants from a collision at a site at time $s$ in the interval }[s,t].\nonumber\eea  Also $W^g (t, r)$ is independent of $W(\cdot)$ and
$W^g (t, r), W^g (t,s)$ are independent for $r \neq s$. The analogous quantity $W^g_{\rm{tot}}(t,r)$ is defined to be the number of green particles produced in $[r,t]$ by a collision at time $r$.

Define the analog of $W^g(t,s)$ in the case where the white sites are in the stable age distribution as
\be{add26}
W^{g,\ast} (t,s).
\ee

If we assume that the white population has reached the stable size distribution (for the number of white particles at an occupied site) at time $s$, then the probability law of the corresponding process
$W^{g,\ast}(t,s)$ depends only on $t-s$.

\begin{remark}\label{nogr}

We note that at times $t \geq s_N$ with $s_N\to\infty$ the distribution of $\zeta$ in the CMJ process
population approaches the stable age distribution as $N \to \infty$
uniformly in $t \geq s_N$.
\end{remark}

Then combining (\ref{agrev19d4}) and (\ref{agrev19d5}) we get that
\be{agrev19d6}
E[\mathfrak{g}(W (\cdot),N,t)] \leq \int_0^t E[(W(s))^2]\frac{e^{2\alpha
s}}{N}E[W^g(t,s)]e^{\alpha(t-s)}ds\leq
\text{const}\cdot\frac{e^{2\alpha t}-e^{\alpha t}}{N}.
\ee
Provided the first inequality could be {\em modified} such that it
becomes up to terms of order
$o(N^{-1})$ an {\em equality},  the following limit exists
\be{agrev19d7}\lim_{N\to\infty}  Ne^{-2\alpha t_N} E[\mathfrak{g}(W(\cdot),N,t_N)]=\kappa_2 E[W^2],\quad \kappa_2=\int_0^\infty e^{-\alpha u} E[W^{g,\ast}(u)]du.
\ee

We note that the probability of  collision in $[0,t_N]$ tends to $0$ as $N\to\infty$ therefore we determine the probability of such a collision, the distribution of the collision time given that a collision occurs and the resulting number of green sites produced by this one collision. Given that a collision occurs the site at which it occurs is obtained by choosing a randomly occupied site.

In order to verify that such a modification exists and that
the assumption on the error in the first inequality in (\ref{agrev19d6})
we have to (1) estimate how
much we overestimate the collision rate by using the CMJ process
(of black and white particles) instead
of the real dual (i.e. $(W,R,P)$ particles) and (2) in the application of $W^g$ we need to control the error introduced by assuming that  collision events involve  sites with the stable age distribution for the white particles and
replacing  $\int_0^{t_N}e^{-\alpha (t_N-s)}E[W^g((t_N,s)]ds$ as $N \to \infty$
by the second part of the equation (\ref{agrev19d7}).

(1) We again note that the over estimate of the number of green particles corresponds to events  having two or more collisions and this occurs with rate  $\frac{e^{3\alpha t}}{N^2}$ and this correction as well as the exclusion of the
purple particles produces a higher order error term.

(2)  Consider the early collisions before reaching the stable age distribution.  Noting that
\be{earlyd} \int_0^{s_N} e^{2\alpha s}e^{\alpha (t_N-s)}ds = Const \cdot e^{\alpha(t_N+s_N)}\ee
therefore these contribute only a term increasing as $e^{\alpha t_N+s_N}$  whereas the latter ones
grow at a rate $e^{2\alpha t_N}$.  If  $s_N\to\infty$, $s_N=o(t_N)$, then the early collisions form a negligible contribution.  But then as $N\to\infty$, the white population approaches the stable age distribution by time $s_N$.

Note that for the late ones we cannot assume that $W^g(t,s)$ takes on its limiting value and for this reason we obtain the integral in the definition of $\kappa_2$.
This proves Part (a) of Proposition \ref{L.Indual}.

{\bf (b)}  We now consider  (\ref{agr34}) when we begin with one
$1_{\{1\}}$ factor at site $j_1$ at time $0$ and one factor
$1_{\{1\}}$ at $j_2\ne j_1$ at time $s>0$ and the total number of occupied sites at time $t$
is what we called $\Pi^{N,(1,2)}_t$.  This can be handled exactly as in (a).  However we want to keep track of the subset of sites occupied by migrants coming from collisions between the two different families. For this reason
we constructed the multicolour system.

Finally we have to show that the bound (\ref{agrev19d7zz}) gives actually
asymptotically the correct answer. Following the argument in part (a) we note that the upper bound on
the $G_{1,2}$  particles gives us a lower bound on the white
particles and this in turn gives us by (\ref{gra59})  with a
bootstrapping a lower bound on the number of $G_{1,2}$  particles.  It then follows that
the overestimate of the number of $G_{1,2}$ particles is no larger than $\rm{const}\cdot \frac{e^{3\alpha t_N}}{N^2}$.

The analogous arguments also yields  the claim (\ref{dd658++}) for $\mathfrak{g}_{(1,2),\rm{tot}}$.

{\bf (c)} Other than including the time dependence for $W_i(\cdot)$ the proof follows the same lines as above. In particular the first two terms
in (\ref{gra62}) now correspond to the CMJ process obtained when $N=\infty$ and the limiting expected value of the correction term now has the form

\be{dd658++td}
\lim_{N\to\infty}N\cdot {e^{-\alpha(2t+s)}}E[\mathfrak{g}_{(1,2),\mathrm{tot}}(W_1(\cdot),W_2(\cdot),N,t,t+s)]=2\kappa_2(t)
E[W_1(t)]E[W_2(t)].
\ee


This completes the proof of Proposition \ref{L.Indual}.
$\hfill \blacksquare$

\end{proof}

{\bf Part 3} {\em (Two dual clouds interacting)}

We want to compare
the dual populations (note the different brackets and the $\sim$ indicating the
populations instead of just the population sizes)
\be{sys1}
\wt \Pi^{N, \{j_1\}}, \wt \Pi^{N,\{j_2\}}, \wt \Pi^{N, \{j_1,j_2\}}
\ee
corresponding to the descendants of the
particles at $j_1,j_2$, respectively evolving independently and
the resulting population starting with both initial particles
and evolving jointly.
In order to carry out this comparison we construct a {\em coupling},
i.e. construct all three populations
on one probability space in such a way that the respective laws are preserved.
At the same time we want to be able to compare these populations with
the respective collision-free versions evolving as CMJ-processes.

{\em Construction of enriched multicolour system.}

To construct this coupling  we use again the WRGB-system and
we construct the particle system starting with 1 particle
at both sites $j_1$ and $j_2$, recall here the constructions from subsubsection \ref{}.
We now
modify this coloured WRGB-particle system by marking the offspring of
the initial particles at $ j_1,j_2$ by using an enriched colour system,  i.e. using the colours
\be{dg10}
\{W_1,R_1,G_1,B_1,P_1, W_2,R_2,G_2,B_2,P_2\}.
\ee
together with
\be{dg11}
\{W_1^\ast, W_2^\ast, R^\ast_1, R^\ast_2, G^\ast_1, G^\ast_2, B^\ast_1, B^\ast_2,
P^\ast_1, P^\ast_2,
R_{1,2}, G_{1,2}, B_{1,2}, P_{1,2}\}.
\ee
The second set of colours is reserved for particles involved in a collision between the two different families or their descendants.

The new features are
\begin{itemize}
\item  When a $W_1$ particle hits a $W_2$ particle the result is a $W_1^\ast$ and $R_{1,2}$ pair as well as the
original $W_2$ particle at this site.  A $W^\ast_1$ particle does not coalesce with the $W_2$ particle.
\item If the $R_{1,2}$ particle coalesces with the $W_2$ particle, the result is the pair consisting of a  $W_2$ particle and a $G_{1,2}$ particle.
\item The above holds with 1 and 2 interchanged.
\item If $W^\ast_1$ hits a $W_1$ particle it changes to a $R^\ast_1$ particle.
\end{itemize}

With this expanded colour set, we now have to define a suitable multicolour dynamic.
We define this system more precisely after outlining the goals of
this construction. This is done so that the $(W_i,W^\ast_i,R_i,G_i),\;i=1,2$ describe the families
in the absence of interactions between them.  The particles $(R_{1,2},G_{1,2})$ correspond to the
particle families resulting from the interaction between the two families and  $(W_1,W_2,R_1,R_2,G_1,G_2,R_{1,2},G_{1,2})$ describe the interacting system.
As before we want to couple the interacting system with a CMJ process.

Recall that the covariance of the type-1 mass in two locations
if we write it in the form $E[XY]-E[X] E[Y]$ consists in
the dual representation of a cloud of dual particles evolving
jointly from two initial particles at different locations and two independently
evolving clouds starting from one particle each.
Hence the idea is that the covariance of $x^N_1(i,t+s)$ and
$x^N_1(i,t)$ is determined by the particles in
\be{dg12}
\{\wt \Pi^{N,\{1\} }\cup \wt \Pi^{N,\{2\} }\} \backslash \wt \Pi^{N,\{1,2\}},
\ee
where $\wt \Pi^{N,\{1\}}, \wt \Pi^{N,\{2\}}$ denotes the two independently
evolving dual particle systems and
$\wt \Pi^{\{1,2\}}$ the jointly evolving ones.
The set $\wt \Pi^{N,\{1\} }\cup
\wt \Pi^{N,\{2\} }\backslash \wt \Pi^{N,\{1,2\}}$ corresponds to the
$G_{1,2}$ green particles.

We list now the precise evolution rules of the new coloured particle system.
In the evolution of the populations arising from the two initial particles
as before particles {\em migrate, coalesce} and
{\em give birth} to  new particles at the rates as in the dual.
One modification is required, namely, in this case we allow
singletons of one family to migrate to sites occupied by the other
family and also singletons of one family in the presence of one or
more particles of the other type can migrate. Here are the precise
rules.
\begin{itemize}
\item[(0)] The rules for the particles and
colours appearing in (\ref{dg10}) are exactly as before as long as
the two populations (i.e. particles with a colour with one index
either 1 or 2) do not interact or occupy the same site. If particles
of the two families occupy the same site singletons of one
family can now migrate.

Furthermore we have the four additional rules concerning the colours
newly introduced in (\ref{dg11}):

\item[(1)] When a white particle from family 1 and a
white particle from family 2 coalesce the result is one $W^\ast_{1}$,
one $W^\ast_2$ particle  and a $R_{1,2}$ particle.
The $W^\ast_\ell$ particles evolve according to the rules of the
dual particle system, but ignoring further collisions with the
respective other population and upon collision with the own population
a $R^\ast_\ell$-particle is created replacing the $W^\ast_\ell$-particle.
Similarly proceed with the other $\ast$-colours.
\item[(2)] $R_{1,2}$ particles that coalesce with $W_1$ (respectively $W_2$) produces a
$W_1$ (respectively $W_2$) and $G_{1,2}$ pair.\\
Two $R_{1,2}$ particles that hit an
occupied site produce blue-purple  pairs $(B_{1,2}, P_{1,2})$ where as before the $B_{1,2}$
particle is placed on the first empty site in the first copy of
$\N$ continuing according to the collision-free dynamic
and the $P_{1,2}$ particle behaves as a typical dual particle for
$\wt \Pi^{N,\{1,2\}}$.

 \item[(3)] Coalescence of $G_{1,2}$ and $R_{1,2}$
particles produce $R_{1,2}$ particles (cf. Remark \ref{rgremark}).

\item[(4)] If a $R_{1,2}$ or
$G_{1,2}$ particle more generally every colour in (\ref{dg11})
gives birth the new particle carries the same type.
\end{itemize}
These set of rules uniquely define  a pure jump  Markov process on the
state space describing the individuals with their colour and location
and with the coupling information (recall (\ref{agr20}) for how
to formalize this).

The key properties of the new coloured dynamics are the following.
On the same probability space we have the various systems we consider.
Namely
\be{agr29}
\wt \Pi^{N,\{i\}} = \mbox{ union of the } W_i,  W_i^\ast, R_i,P_i,R^\ast_i, P^\ast_i, \mbox{ for } \;i=1,2
\ee
The new system also includes our dual starting with two particles, namely,
(by construction)
\be{dg15}
\wt \Pi^{N,\{1,2\}} = \mbox{ the union of the }
W_1,W_2,R_1,R_2,P_1,P_2,R_{1,2},G_{1,2} \mbox{ particles}.
\ee
Furthermore we can again identify a CMJ-process,
\be{dg16}
\mbox{union of }
(W_1,W_2,R_1,R_2,G_1,G_2,G_{1,2}, B_1,B_2,R_{1,2},B_{1,2})\mbox{-particles}
\wh{=} \mbox{ CMJ process.}
\ee

\begin{remark}
To compute higher moments  a {\em complete} bookkeeping of the
difference would require more and more different colours whenever
such a collision with subsequent coalescence occurs. However
although these collisions do occur, we note that they do not
contribute to the calculation up to terms of order $N^{-2}$.
Therefore we again
 ignore these higher order terms instead of introducing further colours.
\end{remark}

In the sequel we shall retain  the
$W^\ast_1$ and $W^\ast_2$ particles but not  the $R^\ast_\ell, B^\ast_\ell, P^\ast_\ell$
particles since the order in $N$ of these sets of particles are of too small  to play a role for our purposes.

\beL{L.negcol}{(Negligible colours)}

Let $\CN^\ell_t$ denote the number of $R^\ast_\ell, B^\ast_\ell$ or
$P^\ast_\ell$ particles in the system $\ell$. Then
\be{agr33}
E[\CN^\ell_t] = O(\frac{e^{3\alpha t}}{N^2}).\qquad \square
 \ee
\end{lemma}

\begin{proof}{\bf of Lemma \ref{L.negcol}}
This follows since these correspond to families resulting from two or more collisions and these have order
$ O(\frac{e^{3\alpha t}}{N^2})$. q.e.d.
\end{proof}

If we count only the occupation numbers by the various colours the
state space is given by
\be{dg14}\begin{array}{l}
\mathcal{M}_c[(\{W_1,W_1^\ast,R_1,G_1,B_1\} \cup\\[2ex]
\phantom{\mathcal{M}_c} \{W_2,W_2^\ast,R_2,G_2,B_2\}
\cup  \{R_{1,2},G_{1,2},B_{1,2},P_{1,2}\}) \times (\{\{1,2,\dots,N\}+\mathbb{N}+\mathbb{N}\})]
\end{array}\ee
(where $\mathcal{M}_c$ denotes the set of counting measures).

This concludes the construction of a coupling via a new enriched
enriched multicolour particle system.
\bi

{\em Consequences of the coupling}

Using the arguments of part (2) above we can conclude  that:
\be{agrev1612}\begin{array}{l}
  \Pi_{t}^{N,\{1\}} + \Pi_{t}^{N,\{2\}}\sim
\frac{1}{c}((W^\prime +W^{\prime \prime}))(\alpha +\gamma) e^{\alpha t}
- [\mathfrak {g}_{\mathrm{tot}}((W^\prime (\cdot)),N,t)
+ \mathfrak{g}_{\mathrm{tot}} ((W^{\prime \prime} (\cdot)),N,t)]\\[2ex]
\hspace{9cm}
+ \mathcal{E}_{1, \mathrm{tot}}(N,t) + \CE_{2, \mathrm{tot}}(N,t),
\end{array}\ee

\be{agrev1613}
\Pi_{t}^{N,\{1,2\}} \sim
\frac{1}{c}((W^\prime +W^{\prime \prime}))(\alpha +\gamma) e^{\alpha t}- \mathfrak
{g}_{\mathrm{tot}}((W^\prime (\cdot)+ W^{\prime \prime} (\cdot)),N,t)
+ \mathcal{E}_{1,2;\mathrm{tot}}(N,t),
\ee
where
\be{dg17}
\CE_{\ell, \mathrm{tot}} (N,t) = o(N^{-1} e^{2 \alpha t_N}) \mbox{ as } N \to \infty,
\quad \ell=1,2 \mbox{ or } (1,2),
\ee
\be{agrev1614}
\mathfrak{g}_{\mathrm{tot}}((W^\prime (\cdot) +W^{\prime \prime} (\cdot)),N,t)
=\int_0^{t_N} \wt g_{\mathrm{tot}}(t_N,s)\xi^{(1,2)}(ds),
\ee
with $\xi^{(1,2)}$  a time inhomogeneous Poisson process on $[0,\infty)$
with intensity measure
\be{dg18}
\frac{((W^\prime (\cdot)(s)+ W^{\prime \prime}(\cdot) (s)))^2e^{2\alpha s}}{N} ds
\ee
and $\wt g_{\mathrm{tot}}(t,s)$ is the size of the
green cloud at time $t$ produced at time $s$ by a {\em newly} created red particle at a white
site.

We know furthermore that

\be{dg19}
E[ \mathfrak {g}_{\mathrm{tot}}((W^\prime+W^{\prime \prime})(\cdot),N,t)]
= \kappa_2 E[(W^\prime+W^{\prime \prime})^2]
\frac{(e^{2\alpha t}-e^{\alpha t})}{N}.
\ee

Then combining all these properties (\ref{dd658+}) follows by comparing the expression for
$\mathfrak{g}_{\mathrm{tot}} ((W^\prime+W^{\prime\prime})(\cdot),N,t)$ with
$\mathfrak{g}_{\mathrm{tot}}   (W^\prime (\cdot),N,t)+\mathfrak{g}_{\mathrm{tot}}(W^{\prime\prime} (\cdot),N,t)$.

Therefore the process corresponding to the quantity
\be{dg20}
 \Pi_{t}^{N,1} +  \Pi_{t}^{N,2}-  \Pi_{t}^{N,\{1,2\}}, \quad t \geq 0
\ee
has an upper bound which is generated by
the Poisson process on $[0,\infty)$ with intensity measure
\be{dg21}
\left(\frac{(2W^\prime W^{\prime \prime})e^{2\alpha s}}{N} + O(e^{3 \alpha s} N^{-2})\right)ds.
\ee
This corresponds to the  difference term (i.e. lost particles due to collision and coalescence between the two families) given by the $G_{1,2}$-particles.  It is
bounded above by allowing  the rate of production of collisions leading to the production of $G_{1,2}$ particles by the  CMJ-process which excludes collisions.  We can then obtain a lower bound by reducing the white families by the respective upper bounds for the sets of green particles as well as the blue particles. This replaces the white family $W(t)e^{\alpha t}$ by $W(t)e^{\alpha t} - O(\frac{e^{2\alpha t}}{N})$.  Using this we can show that the result is an error term which is of order $O(\frac{e^{3\alpha t}}{N^2})$.

We now want a more detailed description in which we keep track explicitly of the collisions of the two families and the resulting $G_{1,2}$ families.
We  therefore consider
\be{add56}
K^{N,\{1,2\}}_t,\Psi^N_t(i_{W_1},i_{W_2},i_{R_{1}},
i_{R_{2}},i_{R_{1,2}},i_{G_{1}}, i_{G_{2}},i_{G_{1,2}},i_{B_1},i_{B_2},i_{B_{1,2}},i_{B_{1,2}}),\ee
which denotes the number of occupied sites at time $t$,
respectively the number of sites having $i_{W_1}$  white particles
from the first population, $i_{W_2}$ white particles from the
second population, etc. Similarly we can define systems for subsets
of colours as we did before. We also consider the  pair $(u^N, U^N)$
describing the new coloured particle system and a new corresponding limiting system
$(u,U)$ for $N \to \infty$.   For the
$(W_1,W_2,R_1,R_2,G_1,G_2,R_{1,2},G_{1,2},B_1,B_2)$-system we modify the
equation for $U_{WRGB}$, namely, (\ref{gra29}) to include the
contribution to the  dynamics of collisions between $W_1$ and
$W_2$ particles. Similar we can proceed with other colour combinations.

Then we obtain the following representation formula.
Here we use the convention to write
\be{gra50a}
U^N_t (i_{W_j}, i_{W_j^\ast}, i_{R_j}, i_{P_j}),
U^N_t (i_{W_1},i_{W_2},i_{R_{1,2}}),\cdots
\ee
for $U^N_t(\cdot)$ all other colours not appearing in the argument summed out.
Then:

\be{gra50}\begin{array}{l}
\Pi^{N,\{j\}}_t\\
=u^N(t)\suml_{i_{W_j},i_{W_j^\ast},i_{R_j}, \cdots, i_{P^\ast_j}}
(i_{W_j}+i_{W_j^\ast}+i_{R_j}+i_{P_j} + i_{R^\ast_j} + i_{P^\ast_j})
U^N_t(i_{W_j},i_{W_j^\ast},i_{R_j},i_{P_j},i_{R^\ast_j}, i_{P^\ast_j}),\\
\hspace{8cm}  j=1,2,
\end{array}
\ee

\be{gra51}\begin{array}{l}
 \Pi^{N,\{1,2\}}_t = u^N(t) \suml_{i_{W_1}, \cdots,i_{P_{1,2}}}
({i_{W_1}+i_{W_2}+i_{R_1}+i_{R_2}+i_{R_{1,2}}+i_{P_1} +i_{P_2} + i_{P_{1,2}}})\\
\hspace{7cm}
U^N_t(i_{W_1},i_{W_2},i_{R_1},i_{R_2},i_{R_{1,2}},i_{P_1}, i_{P_2}, i_{P_{1,2}}),
\end{array}\ee
where the sum is over
$i_{W_1}, i_{W_2}, i_{R_1}, i_{R_2}, i_{R_{1,2}}, i_{P_1}, i_{P_2}, i_{P_{1,2}}$.

We now focus on the are sites at which a collision occurs between a $W_1$ particle and a $W_2$ site (or vice versa) which can then produce migrating
$R_{1,2}$ or $G_{1,2}$-particles. Therefore
we define as {\em ``special'' sites} at time $t$,  those sites
at which
\bea{agr30}&&
\mbox{one or more }R_{1,2} \mbox{  particles are present at time } t
\mbox{ and  {\em also}  }\\&& \qquad\qquad\qquad \mbox{a }  W_1\mbox{ or }W_2\mbox{-particle is present at time }t.
\nonumber \eea
Note that special sites have a finite lifetime and can produce
migrating $R_{1,2}$, $G_{1,2}$  particles during their lifetime.
Denote the number of special sites by:
\be{kappa23}
K^{N,\ast}_t=\left(\sum_{i_{W_1},i_{W_2},i_{R_{1,2}}}
(1_{i_{R_{1,2}}\geq 1}\cdot(1_{\,i_{W_1}\ne 0}\vee 1_{i_{W_2}\ne
0})U^N_{t}(i_{W_1},i_{W_2},i_{R_{1,2}})\right)\cdot u^N(t).
\ee

In the  calculation  of  second moments we can ignore errors
of size $o(N^{-2})$, we can work with the system of
$W_1,W_2,R_1,R_2,R_{1,2}$  and ignore the purple particles
i.e. the colours $P_1,P_2,P_{1,2}$. Therefore it suffices to determine the number of sites occupied by $W_1,W_2,R_1,R_2,R_{1,2}$ particles which estimates $K^{N,(j_1,0),(j_2,s)}_{t+s}$ up to errors of order $O(\frac{e^{3\alpha t}}{N^2})$.  This results in the expansion
given in (\ref{gra62}) where in (\ref{gra62b}) the three $\mathfrak{g}_{\rm{tot}}$ terms correspond to the sites occupied by only $G_1$, $G_2$ and $G_{1,2}$  particles.  It remain to estimate these terms.  The $G_1$, $G_2$ terms are similar to those considered earlier so we focus on the $G_{1,2}$ sites and particles.
We will see below that key to the determining the asymptotics of the covariance at two sites  as $N \to \infty$ is to analyse the number
of $G_{1,2}$-particles asymptotically as $N \to \infty$.
We have now all the tools to start the estimation.
\bigskip

{\em Estimates on the number of $G_{1,2}$-particles}

In order to study the asymptotics of the {\em number of $G_{1,2}$
particles} at time $t_N$,  in sublogarithmic time scales
which we denote this number by:
\be{gra52b}
 \Pi^{N,G_{1,2}}_{t_N} \mbox{ with }t_N\mbox{  as  in }(\ref{agr23}).
 \ee
We next recall that {\em independent} systems of dual particle
systems with starting particles in site $j_1$ at time $0$ and site $j_2$,
at time $s>0$ with $j_1 \neq j_2$, produce clouds as $N\to\infty$ with independent stable
size distributions and the number of sites of the two dual
populations is denoted by
\be{gra52b2}
K^{N,\{1\}}_t, K^{N,\{2\}}_t.
\ee
We then have by the CMJ-theory and the $K_t$ for $K^N_t$ approximation
$N \to \infty$,
\be{gra52b1} K^{N,\{1\}}_{t_N} \sim
W^\prime e^{\alpha t_N}, \quad K^{N,\{2\}}_{t_N}\sim W^{\prime \prime} e^{\alpha t_N}. \ee

To identify the correction terms needed for second moment calculations we need to consider
 the dual particle system starting with two particles
one at site 1, one at site 2  and estimate the effects of the {\em collisions between the two clouds}. Using our
embedding in the coloured particle system this can involve
a collision in which a  $W_1$ particle can migrate to a site already occupied by a $W_2$ particle (or vice versa). Once a collision
occurs, this produces a $W_1^\ast,R_{1,2}$ pair. The $R_{1,2}$ particle can move to another site before it coalesces with one of the $W_2$ particles of the other
colour. If it coalesces the  $W_2$ particle remains an a new  $G_{1,2}$-particle is added at this site.  The $G_{1,2}$  particles
can then reproduce and coalesce with other particles. This way they can produce more
$G_{1,2}$ particles or migrate to produce a new $G_{1,2}$ site.
The system of $G_{1,2}$ particles after the creation of a $G_{1,2}$-particle
has the property that after foundation a $G_{1,2}$-family does not
anymore depend on the other colours,  therefore the $G_{1,2}$-system
evolves as a copy of the basic
one type {\em CMJ particle system} with {\em immigration}
given by a randomly fluctuating source which is independent of the
CMJ-process.  The creation of new $R_{1,2}$ particles at time $s$ is determined by
the number of $(W_1,W_2)$-pairs at time $s$.

Suppose we are given the evolution of the $W_1, W_2$-particles. Then
the events of collisions between $W_1$ and $W_2$ particles is
given by a Poisson process and when a collision  occurs at a site there is
a positive probability of a coalescence and therefore the
production of a green particle. This then produces a growing
population of $G_{1,2}$ descendants according to the CMJ-theory
leading to the analogue expression to (\ref{agrev19d6}) which we had for $W_1, W_2$-particles.

If we replace the system of white particles of the two types by the {\em independent}
site-1 and site-2 system, we obtain an {\em upper} bound for the
mean number of $G_{1,2}$ particles given the growth constants
$W_1$ and $W_2$ of the two clouds, so that at time
$t=t_N$ provided that
\be{kappa22b}
t_N=o(\frac{\log N}{\alpha}),
\ee
by  just calculating the mean.

This results in the following asymptotic population size for the
$G_{1,2}$ particles which is then again given by the behaviour of the
conditional mean (given the $W_1,W_2$ populations):
\be{gra59}
\frac{c}{N}\int_0^{t_N} E[W]
e^{\alpha(t_N-u)} (\Pi^{N,\{1\}}_u K^{N,\{2\}}_u + \Pi^{N,\{2\}}_u
K^{N,\{1\}}_u) p^N (u) du +o({N^{-1}}).
\ee

The real system requires a correction term due to the interaction of the two clouds
leading to a reduction of the two independent white systems by the
$G_{1,2}$ but which is of lower order in $N$ than the $W_1,W_2$-particle
numbers and their effect, i.e. $o(N)$. Hence from
(\ref{gra59}) we get the asymptotic upper bound for the mean
number of $G_{1,2}$ particles as $N \to \infty$ if we condition on
the growth constants $W^\prime$ and $W^{\prime \prime}$ of the two clouds.
As bound we get  an expression which is asymptotically
$(N \to \infty)$ equal to: \be{gra60} \sim \frac{\kappa^N_2}{N}
W^\prime W^{\prime \prime} e^{2\alpha t_N}, \ee
and which differs from the real system by
an  error term of order $o(\frac{e^{2\alpha t_N}}{N})$.

In fact we can obtain
\be{agrev19d7zz}\lim_{N\to\infty}  Ne^{-2\alpha t_N} E[\mathfrak{g}_{(1,2)}(W(\cdot),N,t_N)]=\kappa_2 E[W^\prime]\cdot E[W^{\prime\prime}],\quad \kappa_2=\int_0^\infty e^{-\alpha u} E[W^{g,\ast}(u)]du.
\ee
exactly the same way as for (\ref{agrev19d7}).

\begin{remark}
In the longer time regime $\frac{\log N}{\alpha}+t$ the growing cloud of
green particles $G_{1,2}$ satisfies  a law of large number effect
 similar as in (\ref{gra9a}), that is, conditioned on $W_1,W_2$,
 the variance of the normalized number of $G_{1,2}$-particles goes to zero
 as $N \to \infty$.

\end{remark}

\bigskip

{\bf Part 4} {\em (Random exponentials)}

We next present  the facts needed about the expansion of random exponentials
in the following lemma.

\beL{L.dd917} {(Asymptotics of random exponentials)}

(a) Let $X$ be a nonnegative random variable satisfying
$E[X^k]<\infty$ for all $k\in\N$. Let
$f(\ve,X)$ be such that
\be{dd9171}
f(\ve,X)\geq 0 ,\quad f(\ve,X)\sim
\ve[ X-\ve \kappa X^2 +O(\ve ^2) X^3],
\ee
that is, \be{dd917A}
\sup_{0<\ve\leq 1}\frac{|f(\ve,x)-\ve X+\ve
^2X^2|}{\ve^3}<\infty.
\ee
Consider
\be{dd917B}
E[ 1- e^{-f(\ve,X)}].
\ee

Then the following two relations hold:
\be{dd120}
\lim_{\ve\to 0}\frac{ E[1-e^{-f(\ve,X)}]}{\ve}= E[X],
\ee
\be{dd121}
\lim_{\ve\to 0}\frac{ E[1-e^{-f(\ve,X)} -\ve
X]}{\ve^2} = (1+2 \kappa) E[X^2]. \ee

(b)  Consider now random variables $X, Y$ with
$E[Y^k+X^k]<\infty$ for all $k\in\N$.
 Let $Y_\ve$ be a superposition of a Poisson number of independent copies of a random
 variable $Y$ with intensity $\ve X^2$.
Assume  now that $f(\ve,X,Y)\geq 0$ and  that we have:
 \be{dd120B}
f(\ve,X,Y)\sim \ve[X-Y_\ve].
 \ee

Then (\ref{dd120}) and (\ref{dd121}) are satisfied with
\be{gra100}
f(\ve,X) \mbox{ replaced by } f(\ve,X,Y) \mbox{  and } 1+2 \kappa
\mbox{ by } 1 + EY. \qquad \square
\ee
\end{lemma}

\begin{proof}
(a)  Since the limits (\ref{dd120}) and (\ref{dd121}) are not
affected by the $O(\ve^2)$ term in (\ref{dd9171}) we can assume
that it is a constant  without loss of generality. Then by the
mean value theorem,
\be{dd120C}
\frac{1-e^{-\ve X+\ve^2 X^2-\ve^3X^3
}}{\ve} = D(X,\ve^\ast)\text{
 where  }D(X,\ve)=\frac{d}{d\ve}[1- e^{-\ve X+\ve^2\kappa X^2-\ve^3X^3 }] \ee
and  where $0<\ve^\ast\leq \ve$ and $\lim_{\ve\to 0} D(X,\ve)=X$.
Therefore we have for $c>0$ that
 \be{dd120D}
E \left[\frac{1-e^{-\ve X+\ve^2\kappa X^2-\ve^3X^3 }}{\ve}
\right]= E[D(X,\ve^\ast)1_{(X\leq \frac{c}{\ve})}]+
\frac{\text{const}}{\ve}\cdot P[X> \frac{c}{\ve}].
\ee
Then
\be{dd120E}
D(X,\ve)=\frac{d}{d\ve}[1- e^{-\ve X+\ve^2\kappa X^2-\ve^3X^3
}] = [(X-2\ve^2\kappa X^2+3\ve^2X^3)e^{-\ve X+\ve^2\kappa
X^2-\ve^3X^3 }].
\ee
Note that
\be{dd120F}
E[\sup_{0<\ve \leq 1,0\leq
 X \leq\frac{c}{\ve}}|D(X,\ve)|]\leq \text{const}E[X+X^2+X^3]<\infty.
 \ee
 Therefore by dominated convergence for every $c>0$:
\be{dd120G}
\lim_{\ve\to 0} E[1_{(X\leq \frac{c}{\ve})}\cdot
D(X,\ve)]=E[X].
\ee
Since  $E[X^k]<\infty$ for $k\in\N$ we have
\be{dd120H}
P[X>\frac{c}{\ve}]\leq E[X^k]\cdot (\frac{\ve}{c})^k.
\ee
Combining (\ref{dd120G}), (\ref{dd120H}) we conclude from
(\ref{dd120D}) the claim (\ref{dd120}).

Similarly for the second order expansion write ,
\be{dd120I}
\frac{1}{\ve}[\frac{1-e^{-\ve X+\ve^2\kappa
X^2-\ve^3X^3 }}{\ve}-X ]=
\frac{1}{\ve}[D(X,\ve^\ast)-D(X,0)]=D^2(X,\ve^{\ast \ast}),
\ee
where
$D^2(X,\ve)=\frac{d^2}{d\ve^2}[1- e^{-\ve X+\ve^2 \kappa
 X^2-\ve^3X^3
 }]$ and
 \be{dd120J}
 \lim_{\ve\to 0} [D^2(X,\ve)] =X^2+2\kappa X^2
 =  (1+2 \kappa) X^2.
 \ee
As above we then conclude the claim (\ref{dd121}).

(b) We have
\be{dd120K}
1-e^{-f(\ve,X,Y)}= 1-e^{-\ve X+ \ve Y}1_{Y\ne
0}-e^{-\ve X}1_{Y= 0}
\ee
and \be{dd120L}
\sup_{0<\ve\leq
1,\;Y<\frac{1}{\ve}}|\frac{d}{d\ve}(1-e^{-f(\ve,X,Y)})|\leq
\text{const}\cdot(X+Y).\ee

Therefore the relation is now obtained as follows:
\bea{dd120M}
&&\lim_{\ve\to 0}\frac{E[1-e^{-f(\ve,X,Y)}]- \ve
E[X]}{\ve}\\
&&=\lim_{\ve\to 0}\left[ \frac{Er[((1-e^{-\ve X}])+\ve
E[(1-e^{-\ve X^2})]E[Y 1_{0<Y\leq\frac{1}{\ve}}]  -\ve
E[X]}{\ve^2}\right]\nonumber\\
&& = \lim_{\ve\to 0}\left[
\frac{\ve^2E[X^2] +\ve^2 E[X^2] E[Y] +P[Y> \frac{1}{\ve}]
}{\ve^2}\right]= E[X^2](1+E[Y]). \nonumber \eea

$\hfill\blacksquare$

\end{proof}
\bigskip

{\bf Part 5:}
{\em (Proof  of Proposition \ref{P.Growthdrop2})}

To prove (a) (\ref{agrev29}) and (\ref{agrev30}) it suffices to show
that \be{agrev31}
 0< \lim_{N\to\infty} E[e^{-\alpha t_N}\wh{x}^N_2(t_N)]
 <\infty,\quad 0<\lim_{N\to\infty}\text{Var}[e^{-\alpha t_N}\wh{x}^N_2(t_N)]<\infty,
 \ee
\be{agrev31b} \limsup_{N\to \infty} E[(e^{-\alpha t_N} \wh x^N_2
(t_N))^3] < \infty, \ee so that the family
$\{\mathcal{L}({\wh{x}^N_2(t_N))}:N\in \mathbb{N}\}$ is tight and
their limit points are non-degenerate and random since they have finite
non-zero mean and variance. If we can
identify the limits in (\ref{agrev31}), then we can obtain
(\ref{agrev28}) and (\ref{agre61a}) via (\ref{agrev31b}).

The proof proceeds in steps, using duality representation of
general moments in order to obtain  first,  second and then third
moments of $\wh x^N_2(t_N)$. It follows from Lemma
\ref{generalities} that we must consider {\em correction terms} of
order $\frac{1}{N^k}$ in the calculation of kth moments of
$\wh{x}^N_1(t)$ in order to later derive the asymptotics of the
moments of $\wh x^N_2$ using the representation of formula (\ref{AGr1}).

Recall that the basic dual relation
\bea{dd659}
&& E[(x^N_2(i,t))^a]= E\left[1-\exp \left(-\frac{m}{N}
\intl^{t_N}_0 \Pi^{N,a}_s ds\right)\right]\eea
where $\Pi^{N,a}_s$ denotes the total number of dual particles if we start with $a\in\N$  particles at site $i\in\{1,\dots,N\}$.
We have obtained a representation of this set of particles in terms of a decomposition into white, red and
purple particles embedded in a multi-colour particle system of white, red, purple, blue and green particles
(WRGPB-system).  We have also coupled this system to a related CMJ process which in particular provides an upper bound,
which was the WRGB-system.
We now look at this system in more detail in the time interval $[0,t_N]$ where $t_N$ satisfies the hypothesis (\ref{AGr6b}).
We first recall that the number of occupied sites in the CMJ has the form
\be{add57}
\Pi^{CMJ,a}_s= W_a(s)e^{\alpha s},\quad  \lim_{s\to\infty}W_a^s=W_a,\;\rm{a.s.}\ee
Moreover the set of white and red particles in given by removing a set of green and blue particles from the CMJ process and we can bound the number of purple particles by the set of blue particles.
Therefore we can represent the dual population by removing the green and blue particles with an error given by the difference between the blue and purple particles.  Note that the blue and purple particles result from the event of {\em 2 or more collisions} and therefore the error term can be obtained by determining this probability.  Moreover we have seen that the collision process is given by a Poisson process with intensity bounded by $\alpha(s)W_a(s)e^{\alpha s}\frac{W_a(s)e^{\alpha s}}{N}$.  Therefore the probability of a collision in $[0,t_N]$ is $O(\frac{e^{2\alpha t_N}}{N})$ and this goes to $0$ as $N\to\infty$ as $N\to\infty$ under hypothesis (\ref{AGr6b}
). Let $\#$ denote the number of collisions in $[0,t_N]$.  Moreover the probability of the event two or more collisions is
\be{add58}
P(\#\geq 2) \leq Z_1\left(\frac{e^{4\alpha t_N}}{N^2}\right),
\ee
where $Z_1$ is a random variable having finite moments.
Conditioned on the presence of a collision at time $s\in [0,t_N]$ the resulting set of green particles is given by $W^g_{\rm{tot}}(t_N,s)e^{\alpha(t_N-s)}$ where $W^g_{\rm{tot}}(t_N,s)$ is a nonnegative random variable with all moments and bounded by an independent copy of the CMJ process.

Then using Taylor's formula with remainder we obtain
\bea{ddtay1} &&E\left[ 1-\exp \left(-\frac{m}{N}
\intl^{t_N}_0 \Pi^{N,a}(s) ds\right)-\left(\frac{m}{N}\left[
\intl^{t_N}_0 \Pi^{N,a}(s) ds\right]\right)\right]\\&&
=-\frac{m^2}{2N^2}E\left[
\left(
\intl^{t_N}_0 \Pi^{N,a}(s) ds\right)\right]^2+E[\mathcal{E}_2(N)]\nonumber\\
&&
\mathcal{E}_2(N)\leq \frac{1}{6} \left(-\frac{m}{N}\left[
\intl^{t_N}_0 \Pi^{CMJ,a}(s) ds\right]\right)^3\leq \frac{m^3}{6N^3}Z_2 e^{3\alpha t_N},
\eea
where $Z_2$ has finite moments. Calculate first the expectation of the occupation integral.

\bea{ddtay2}
&& E\left[\frac{m}{N}\left[ \intl^{t_N}_0 \Pi^{N,a}(s) ds\right]\right]\\
&&=\frac{m}{N}E\left[\intl^{t_N}_0 \frac{(\alpha(s)+\gamma(s))}{c}W_a(s)e^{\alpha s} ds-\frac{1}{N}\intl^{t_N}_0 \frac{(\alpha(s)+\gamma(s))}{c}(W_a(s))^2W_{g,\rm{tot}}(t_N,s)ds\right]\nonumber\\
&& \qquad +E[\mathcal{E}_1(N)]\nonumber\\
&& \mathcal{E}_1(N)\leq Z_1\frac{e^{4\alpha t_N}}{N^3}.\nonumber
\eea
Here the second term comes from integrating over the event $\{\#=1\}$ and corresponds to a green family produced by a single
collision. To estimate the first term on the r.h.s. of (\ref{ddtay1}) we calculate:
\bea{ddtay3}  && E\left[\frac{m}{N}\left(
\intl^{t_N}_0 \Pi^{N,a}(s) ds\right)\right]^2\\
&&=\frac{m^2}{N^2}E\left[\intl^{t_N}_0 \frac{(\alpha(s)+\gamma(s))}{c}W_a(s)e^{\alpha s} ds\right]^2 +E(\mathcal{E}_3(N))\nonumber\\
&&
\mathcal{E}_3(N) \leq Z_3\frac{e^{3\alpha t_N}}{N^3}.\nonumber
\eea

We are now ready to obtain study the first and second moments of $\wh x_2(t_N)$. Consider a sequence $s_N\to\infty$ $0\leq t_N-s_N\to \infty$.  Then using the basic properties of the CMJ process (\ref{ang2b}) and (\ref{ang4c2})  we have $W_a(s)\to W_a, \alpha(s)\to \alpha$  and $\gamma(s)\to \gamma$ (a.s. and in $L^1$)  as $s\to \infty$.  Using this we obtain

\bea{ddinitre}
&&\lim_{N\to\infty} Ne^{-\alpha t_N}E\left[\intl^{t_N}_0 \Pi^{N,a}_s ds\right]\\
&&
=\lim_{N\to\infty} Ne^{-\alpha t_N}E\left[\intl^{t_N}_0\frac{\alpha(s)+\gamma(s)}{c} W_a(s)e^{\alpha s} ds\right]\nonumber\\&&
= \lim_{N\to\infty}Ne^{-\alpha t_N}\left( E\left[\intl^{s_N}_0 \frac{\alpha(s)+\gamma(s)}{c} W_a(s)e^{\alpha s} ds
\right]+ E\left[\intl^{t_N}_{s_N}\frac{\alpha(s)+\gamma(s)}{c} W_a(s)e^{\alpha s} ds\right]\right)\nonumber\\&&
= \frac{(\alpha +\gamma)}{c}W_a.\nonumber\eea

Therefore
\bea{dd924}
&&\lim_{N\to\infty}Ne^{-\alpha t_N}\left(1-E[ \exp (-\frac{m}{N}
\intl^{t_N}_0 \Pi^{N,a}_s ds)]\right)\\
&& =\lim_{N\to\infty}Ne^{-\alpha t_N}E\left[m\int_0^{t_N}\Pi^{N,a}_sds\right]\nonumber \\
&& = m^\ast E[W_a]\eea
where  $m^\ast=\frac{\alpha+\gamma}{\alpha c}$.

\bigskip
In the same way we obtain
\bea{dd925}&&\lim_{N\to\infty}N^2e^{-2\alpha t_N}\left(1-E\left[ \exp \left(-\frac{m}{N}
\intl^{t_N}_0 \Pi^{N,a}_s ds\right)-\frac{m}{N}\int_0^{t_N}\Pi^{N,a}_sds\right]\right)\\&& =
{\kappa_2} E[(W_a)^2]  + \frac{(m^\ast)^2}{2} (E[W_a])^2. \qquad
\nonumber \eea

 {\em Step 1$\;$ First moment of $\widehat{x}_2(t_N)$ }

Note that by exchangeability \be{AGr9}
 E[\widehat{x}_{1}^{N}(t)] =NE [x^N_1 (t)].
\ee Then
\be{AGr9a}
 E[\widehat{x}_{2}^{N}(t)]=N-NE [x^N_1 (t)].
\ee

Using  the above calculations we have
\be{AGr9c}
 E[\widehat{x}_{2}^{N}(t_N)]=NE [x^N_2 (t_N)]= NE\left[1-\exp \left(-\frac{m}{N}
\intl^{t_N}_0 \Pi^{N,a}_s ds\right)\right].
\ee
Therefore using (\ref{dd924}),  we obtain

\be{AGr9b}
\lim_{N\to\infty}E[e^{-\alpha t_N}x^N_2(t_N)] = m^\ast E[W_a].
\ee

{\em Step 2$\;$ Second moment of $\wh{x}_2(t_N)$}

We will focus next on the calculation of the second moment.
We follow the method used above but leave out some details.

To compute the second moment of $\wh{x}^N_2(t)$ note that
\be{agrev13}\begin{array}{lcl}
E[(\widehat{x}_{2}^{N}(t)]^{2})
&  =&
N^{2}-2NE[\widehat{x}_{1}^{N} (t)]+E[(\widehat{x}_{1}^{N}(t))^{2}]\\
&  =& N^{2}-2N^2 E[x^N_1 (i,t)] + E[(\wh x^N_1(t))^2]
\end{array}
\ee
and with $i \neq j$
\be{agrev14} E[(\widehat{x}_{1}^{N}(t))^{2}]
=NE[(x_{1}^{N}(i,t))^{2}]+N(N-1)E[x_{1} ^{N}(i,t)x_{1}^{N}(j,t)].
\ee

Since the formula contains  $O(N^2)$ factors,   in order to
compute the r.h.s. up to order 1,  as $N \to \infty$, it is
necessary to include corrections of order $\frac{1}{N^2}$ in
both $E[x^N_1(i,t)]$ and
$E[x_{1}^N(i,t)x_{1}^{N}(j,t)]$, etc. The first moment expansion we
have done in Step 1. This means that we have to
consider corrections due to possible {\em intersections} of
different clouds (clouds meaning the descendants of one of the
initial individuals.)

 Consider two clouds starting with one  factor at different sites  in the
dual population which occupy  the following number of sites respectively
\be{agrev14b}
W_{1,N}e^{\alpha t_N}
\mbox{ and } W_{2,N}e^{\alpha t_N}, \mbox{ as } N \to \infty,
\ee
with $t_N$ as in (\ref{AGr6b}) and where
$W_{i,N}=W_i(t_N)\to  W_i$ as $N\to\infty$ for $i=1,2$.
Then the rate at
which the first cloud produces particles that migrate onto one of
the sites occupied by the second cloud and symmetrically is as $N \to \infty$:
\be{agrev15}
\frac{2(\alpha+\gamma)^2}{c}\frac{W_{1,N} W_{2,N}e^{2\alpha
t_N}}{N}+o(\frac{e^{2\alpha t_N}}{N})
=
\frac{2(\alpha+\gamma)^2}{c} \frac{W_1 W_2 e^{2\alpha t_N}}{N}+o(\frac{e^{2\alpha t_N}}{N}).
\ee
If the two clouds have
an intersection, coalescence can reduce the number of particles
compared to $(W_1+W_2) (\alpha+\gamma)e^{\alpha t_N}$ and this
produces a {\em correction term} corresponding to the $G_{1,2}$ particles.  We have analysed this correction
term in all detail  (in the proof of Proposition \ref{L.Indual})
using a coloured particle system but here using that
construction and its consequences, recall (\ref{agrev16}-\ref{agrev16xxx}),
we get a result on the behaviour of the number of
dual particles and can perform moment calculations.
We use these results to calculate covariances for $i\ne j$, namely, we
use the dual and expand the time integral using (\ref{agr34}) with $s=0$
to represent $\Pi^{N,(1,2)}_u$ for $u \in [0,t_N]$. We
estimate $E[\intl^{t_N}_0 \mathfrak{g}_{\mathrm{tot}}(W^{1,2}(\cdot), N,u) du]$ by
$\frac{\kappa_2}{2c\alpha}E[(W^{1,2})^2]\left(\frac{e^{2\alpha t_N}-2e^{\alpha t_N}}{N^2}\right)$
and (\ref{agrev16xxx})
to bound the error term.

If $t_N$ satisfies the above growth assumption then as $N \to \infty$ we can expand the
exponential and again verify  that we can take expectation in this
expansion, as in (\ref{AGr9b}). (Recall Lemma \ref{L.dd917} for the justification of these expansions.)
Then using again (\ref{dd659}) we calculate as follows.

\be{cov}\begin{array}{l}
E[x_{1}^{N}(i,t_N)x_{1}^{N}(j,t_N)]
\\[2ex]
=E\left[e^{-\frac{m^\ast}{N}\left[(W_1+W_2)(\alpha+\gamma)
(e^{\alpha t_N}-1)-\int_0^{t_N}\mathfrak{g}_{\mathrm{tot}}
((W_1+W_2)(\cdot),N,u)-e^{\alpha t_N}du
+o(N^{-1})\right]}\right]  \nonumber \\[2ex]
= 1- \frac{m^\ast}{N} E[(W_1+W_2)] (\alpha+\gamma) (e^{\alpha t_N}-1)
+ \frac{m \kappa_2}{2 c\alpha} E [(W^{1,2})^2]
  \frac{e^{2 \alpha t_N}}{N^2} \nonumber \\[2ex]
\hspace{3cm}  + \frac{(m^\ast)^2}{2N^2} (E[W_1 +W_2])^2 e^{2 \alpha t_N}
  + o \left(\frac{e^{2 \alpha t_N}}{N^2}\right).\nonumber
\end{array}
\ee

We have as $N \to \infty$ for $i \neq j$:
\bea{AGr8}
&&E[x_{1}^{N}(i,t_N)x_{1}^{N}(j,t_N)]\\
&&     \sim 1 - \frac{m^\ast}{N}E[W_1+W_2] e^{\alpha t_N}\nonumber\\
&& +
\frac{m\kappa_2}{2c\alpha}E[(W^{1,2})^2]\frac{e^{2\alpha
t_N}}{N^2}+\frac{(m^\ast)^2}{2N^2}(E[W_1+W_2])^2e^{2\alpha t_N} +
o(\frac{e^{2 \alpha t_N}}{N^2}),\nonumber
\eea
where $W^{2,1}$ is the growth constant if we start with {\em two} particles at {\em
   one} site.
Similarly we get

\bea{secmom}
&&E[(x_{1}^{N}(i,t_N))^2] \\
&&  \sim 1 - \frac{m^\ast}{N}E[W^{2,1}] e^{\alpha t_N}\nonumber\\
&& +
\frac{m\kappa_2}{2c\alpha}E[ (W^{2,1})^2] \frac{e^{2\alpha
t_N}}{N^2}+\frac{(m^\ast)^2}{2N^2}(E[W^{2,1}])^2e^{2\alpha t_N} +
o(\frac{e^{2 \alpha t_N}}{N^2}). \nonumber
\eea

Furthermore note that
\be{agrev41}
\text{Var}\left( \wh x_{2}^{N}(i,t_N)\right)
  = \sum_{i=1}^N \text{Var}(x^N_1(i,t_N))
+\sum^N_{i\ne j}\text{Cov}(x^N_1(i,t_N),x^N_1(j,t_N)). \ee

Then by (\ref{agrev14}) combined with (\ref{dd659}) and
(\ref{cov}) the total mass of type 2 satisfies for $N \to \infty$
by expanding the exponential (below $i \neq j$):

\be{agrev18}
\begin{array}{l}
E[(\widehat{x}_{2}^{N}(t_N))^{2}]= N^2-2NE[\wh x_1^N(t_N)]
+E[(\wh x^N_1(t_N))^2]
\\[2ex]
\hspace{1cm}
=N^2-2N^2E[x^N_1(i,t)]+NE[(x^N_1(i,t_N))^2]+N(N-1)E[x^N_1(i,t_N)x^N_1(j,t_N)]
\\[2ex]
\hspace{1cm}
= -N^{2}+2Nm^\ast E[W_{1}]e^{\alpha t_N}
+ N(1-\frac{m^\ast}{N}E[W^{(2,1)}]e^{\alpha t_N})-m^2E[(W_1)^2]e^{2\alpha t_N}
\nonumber\\[2ex]
\hspace{1cm}
 +N(N-1)(1-\frac{m^\ast}{N}2E[W_{1}] e^{\alpha t_N}
+\frac{(E[(W_1+W_2)^2]\kappa_2 m^\ast e^{2\alpha t_N})}
{  N^{2}})\nonumber \\[2ex]
\hspace{1cm}
+\frac{1}{2}N(N-1)\frac{m^{\ast 2}}{N^{2}}E[(W_{1}+W_{2})^{2}]e^{2\alpha t_N}. \nonumber
\end{array}\ee

Therefore we get
\bea{agrev191}
(\mbox{ r.h.s. (\ref{agrev18}) })
 \\ = && 2E [W_1 W_2] \kappa_2 m^\ast e^{2\alpha t_N}-m^\ast E[W^{(2,1)}] e^{\alpha t_N}\nonumber
  \\ && +m^{\ast 2}(E[W_{1} ])^{2}e^{2\alpha t_N} + o(1), \mbox{ as } N \to \infty.
    \nonumber \eea

This implies
\be{agrev19b}
\lim_{N\to\infty} (e^{-2\alpha t_N} E[(\wh x^N_2(t_N))^2])
  = (2\kappa_2 m^\ast + m^{\ast 2}) (E[W_1])^2,
\ee

\be{agrev19c} \lim_{N\to\infty} (e^{-2\alpha t_N} Var[(\wh
x^N_2(t_N))])
  = 2\kappa_2 m^\ast  (E[W_1])^2. \ee

We now consider the case with two times $t_N$ and $t_N+t$.
Following the same steps as above we can show that
\be{agrev19d} \lim_{N\to\infty} (e^{-2\alpha t_N} E[(\wh
x^N_2(t_N))\wh
x^N_2(t_N+t)])
  = (2\kappa_2 m^\ast + m^{\ast 2}) (E[W_1)]^2e^{\alpha t}. \ee

{\em Step 3 $\;$ (Third moments)}

We can proceed again as before, expand the $x^N_2(t)$ in terms of
$x^N_1(t)$ and apply the dual representation again. Here however
we only need an upper bound and therefore we can bound the total number of
dual particles from above
by the ones of the collision-free process. Then the claim follows again easily
by using (\ref{ga40}).

This completes the proof of Proposition  \ref{P.Growthdrop2}.
\bi

\subsubsection{Asymptotically deterministic droplet growth}
\label{sss.detdropgro}

We now return to the analysis of the beginning of the evolution,
where the droplet growth starts.
Over  time periods which are large but  remain much smaller than
$\frac{\log N}{\alpha}$ as $N \to \infty$ the total mass of type 2
is expected to grow  deterministically as an exponential with rate
$\alpha$ and constant (random) factor. The random factor developed over a short
time period at the very beginning of the evolution. This leads to
the idea to focus
here in a first step on the case where after some time $t_0$ {\em we turn off the
mutation}. The point is that the mutation which enters later than
some large time $t_0$ is not relevant for the growth of the
droplet anymore if we choose $t_0$  sufficiently large which
we show in the second step.

In this section  we therefore compute the moments of $\wh
x_2^N(t)$ conditioned on the configuration at a fixed time $t_0$,
that is, given a specific realization of the random configuration
${x}_{2}^{N}(i,t_0),\;i=1,\dots,N$, then identify the asymptotics
as  first  $N \to \infty$ and then $t \to \infty$. In the
Corollary \ref{C.dbetatimes} and Corollary \ref{C.hmdg} and
then later also in the next subsubsection in
Proposition \ref{P.CW} we will then consider the joint limit of $N
\to \infty, t \to \infty$, for the quantity $e^{-\alpha t} \wh
x^N_2(t)$, more precisely
$e^{-\alpha t_N} \wh x^N_2(t_N)$ for $N \to
\infty$ with $t_N \uparrow \infty$ and
$\limsup\limits_{N\to \infty} (t_N/\log N)<\alpha$.
Those quantities converge in probability to a random limit.

In order to mimic the situation, we consider a process, where the
rare mutation is turned of after time $t_0$ for times $t \geq 0$,
we consider a process starting in finite mass and $m=0$
for all times.

\beP{P.DetGroDro-new} {(Deterministic regime of droplet growth)}

Assume that at some time $t_0\geq 0$  the  distributions of $\{\wt
x^N_2(t_0,i):i=1,\dots,N\}$ for the system (\ref{angr2}),
(\ref{angr3}) are converging to the one of  $\{\wt x_2(t_0,i), i \in
\N\}$ with $\sum_i \wt x_2(t_0,i)<\infty$ a.s.

Moreover we assume that for every
$N$ and for the limit configuration we have that $\forall\, \ve>0$ there is
a finite random set (we suppress the dependence on $N$ in the
notation):
\be{agre61b}
\mathcal{I}(\ve) =\CI (\ve,N) \mbox{ of } k=k(\ve,N) \mbox{ sites,}
\ee
such that for every $N \geq 2$:
\be{agre62b}
\sum_{i\in
(\mathcal{I}(\ve))^c}\wt x^N_2(t_0,i)<\ve. \ee

Assume
\be{agre60b}
m=0 \quad\text{for } t\geq t_0.
\ee

We set  (recall (\ref{ba10g9a}) for $\CU(t,\cdot)$ and
(\ref{limsup2a}), (\ref{ang4c}) for $\CU(\infty, \cdot)$):
\be{agrev43ll} g(t,x)=(1-\sum_{\mathfrak{t}}^\infty (1-x)^{\ell}
\CU(t,\ell)), \quad g(\infty,x) =(1- \suml^\infty_{\ell=1}
(1-x)^\ell \CU(\infty,\ell)). \ee

Then the following four convergence properties hold:

(a)
\bea{agrev431m}
\lim_{t\to\infty} e^{-\alpha t}
E[\gimel^{m,t_0}_t([0,1])]& =&\lim_{t\to\infty} e^{-\alpha t}
\lim_{N\to\infty}E[\wh x^N_2(t)]\\
&=&E[W]\sum_{i=1}^\infty g(\infty, \wt x_2 (t_0,i)),\nonumber \eea

(b) \be{agrev43l2}
\liml_{t \to \infty} \liml_{N \to \infty}
(E[e^{-\alpha (t+s)}
  \wh x^N_2 (t+s)]-E[e^{-\alpha t}  \wh x^N_2(t)]) =0, \quad \forall s \in \R,
\ee

(c)   \be{gra63}\begin{array}{l} \liml_{t\to \infty} \liml_{N \to
\infty}
 e^{-\alpha (s+2t)}\rm{Cov}(\wh x^N_2(t+s),\wh x^N_2(t))  \\
 \hspace{2cm}
 = (E[W])^2 \left[ 2\kappa_2\suml^k_{i=1} g(\infty, \wt x(t_0,i))
   \right], \quad \forall s \in \R,
\end{array}\ee
and

(d)  \be{agrev43h}
\liml_{t\to \infty} \liml_{N \to \infty} E
\left[ \left([e^{-\alpha (t+s)} \wh x^N_2(t+s)] - e^{-\alpha t} E[\wh
x^N_2(t)]\right)^2\right]  = 0 , \quad \forall \; s \geq 0.
\ee

(e) Now assume that $m>0$ for all $t \geq t_0$ (recall (\ref{agre60b}).  Then the conclusions of (a)-(d) remain valid.
$\qquad \square  $
\end{proposition}

We next verify that we can compute the simultaneous limits
$N\to\infty$, $t\to\infty$ for times $t_N(\beta,t) =
(\beta/\alpha) \log N+t$, $0<\beta<1$.

\beC{C.dbetatimes} {(Droplets of size $N^\beta$ grow deterministically in $\beta$) }

Consider again the case where $m=0$ and the initial state is as in
Proposition \ref{P.DetGroDro-new}. Choose $\beta$ with $0 < \beta
< 1$ and  let \be{agrev43g} t_N(\beta,t) =\frac{ \beta(\log
N)}{\alpha}+t. \ee

We have \be{agrev43i}
\lim_{t\to\infty}\lim_{N\to\infty}E\Bigg[\Big(e^{-\alpha t}\wh
x^N_2(t)-e^{-\alpha t_N(\beta,0)} \wh
x^N_2(t_N(\beta,0))\Big)^2\Bigg]=0. \qquad \square \ee
\end{corollary}

This result shows that the randomness in the droplet growth is
generated in the very beginning at times of order $O(1)$ as $N \to
\infty$.
In particular we can then
show that, conditioned on the total mass of type 2 at time
$t_N(\beta_1,0)$, the total mass at a later time $t_N(\beta_2,0),
1>\beta_2 > \beta_1 \geq 0$ is {\em deterministic} in the limit
$N\to\infty$.

\begin{remark}
The limiting total droplet mass dynamic as $N \to \infty$ becomes
deterministic in the following sense.
Assume that $\{{x}_{2}^{N}(i,t_0),\;i=1,\dots,N\}$ is a measure
with total mass $N^{\beta}$ and supported on a set of size
$aN^{\beta}$ for some $a \in (0,\infty)$ and $\beta \in (0,1)$.
Then we can show that for $t>t_0$ we have:
\be{gra75}
 \frac{Var[\wh x^N_2(t)]}{(E[\wh x^N_2(t)])^2}\to 0\ee
as $N\to\infty$.

\end{remark}

\begin{remark}
Given a fixed
initial measure $\{\wt x(t_0,i),i\in\N\}$ as in the statement of Proposition \ref{P.DetGroDro-new} we can also prove that as $t \to
\infty$ \be{gra76}
 \lim_{N\to\infty}E[\sum_{i} (x^N_2(i,t))^2] = O(\lim_{N\to\infty}E[\sum_i (x^N_2(i,t))]).
 \ee
This is analogous to  the Palm picture (\ref{elam3})  and implies that the
mass of type 2 is clumped on a set of sites on which the mass is almost one.
\end{remark}

\begin{proof} {\bf of Proposition \ref{P.DetGroDro-new}. }
The proof proceeds by using  the dual representation of moments and
the multicolour particle system introduced
in Subsubsection \ref{sss.dropletgrow} to analyse the expressions.

(a)  We determine the contribution to the first moment coming only from
the initial mass in $\mathcal{I}(\ve)$.  The result then follows
by letting $\ve\to 0$.

Without loss of generality we use (recall $k=k(\ve,N),\CI(\ve)=\CI(\ve,N)$)

\be{agre62c} \CI(\ve) =
\{1,\cdots,k\}. \ee

 Note that since we are now not working with a spatially homogeneous initial
measure, we must modify the dual. In particular we eliminate the rule that single particles at a site do not move. In the calculation
of $E[\wh x^N_2(t)]$ we  compute
\be{add200} E[\wh x^N_2(t)]=\suml^N_{j=1} E[x^N_2(j,t)].\ee
To compute $E[x^N_1(j,t)]$ we start the dual with one particle at site $j\in\{1,\dots,N\}$.  However the number of particles at this site will return to $1$ after at most a finite random time and therefore after some finite random time the single particle will jump to a randomly chosen site and site $j$ will become empty.  Therefore for large times $t$ the value of $E[x^N_1(j,t)]$ is independent of $j$ (even if it is in $\CI(\ve)=\{1,\dots,k\}$).

We calculate the expectation of type 1 in a site $j$ using the dual starting with one particle at site $j$.

 \be{1-time}
 E[x^N_1(t,j)] = E \left[\prod_{i=1}^{K^N_t}
(1-\wt x^N_2(t_0,i))^{\zeta^N (t,i)}\right]. \ee

We first consider the initial configuration where we have
\be{dg60}
\wt x^N_2 (t_0,i) = 0 \mbox{ for every } i \notin \CI(\ve).
\ee
Recall that here $U^N(t,\ell)$ is the frequency of occupied sites that have
$\ell$ particles.

Hence noting that there will be a $Bin(K^N_t,\frac{k}{N})$ number
of hits of the $k$-set $\CI(\ve)$ by the dual cloud and the frequency of sites in the cloud
with $\ell$-sites is given by $U^N(t, \ell)$ and the contribution of such a site by
$(1-\wt x^N_2 (t_0,i))^\ell$, we get for sites $j \notin \{1, \cdots, k\}$ that:
\bea{gra64} E[x^N_1(t,j)] &&
=E\Bigg[(1-\frac{k}{N})^{K^N_t}+(1-\frac{k}{N})^{K^N_t-1}\frac{K^N_t}{N}
\sum_{i=1}^k \sum_{\ell=1}^\infty U^N(t,\ell)(1-\wt
x_2^N(t_0,i))^{\ell} \\
&&
+(1-\frac{k}{N})^{K^N_t -2}\frac{K^N_t(K^{N}_t-1)}{2N^2}\nonumber\\
&& \times \sum_{i_1,\,i_2=1}^k \sum_{\ell_1,\ell_2=1}^\infty
U^N(t,\ell_1)U^N(t,\ell_2)(1-\wt x_2^N(t_0,i_1))^{\ell_1}(1-\wt
 x_2^N(t_0,i_2))^{\ell_2} \Bigg ]\nonumber\\
 && +o(N^{-2}). \nonumber
\eea
Note that as $N\to\infty$,  $U^N(t,\ell)\to
U(t,\ell)$, with $U(t, \cdot)$ the size distribution of the CMJ process at time $t$.

We rewrite the formula (\ref{gra64})  as follows. Let

\be{gra65}
 g^N(t,x)=1-\sum_{\ell =1}^\infty(1- x)^{\ell}U^N(t,\ell)
 \ee
 and

 \be{gra65b}
 g^{N,2}(t,x)=1-\sum_{\ell =1}^\infty(1- x)^{\ell}(U^N(t)^\ast U^N(t))(\ell).
 \ee
Note that $g^N , g^{N,2}$ do depend as well on the initial state of the dual process
via $U^N = U^{N,(1,1)}$, this dependence  we suppress here in the notation.

Then:
\bea{1-time-2}
&&\\
&&E[x^N_1(t,j)] \nonumber\\
  &&=E\Bigg[ 1-\frac{K^N_t}{N}\sum_{i=1}^k g^N(t,\wt x^N_2(t_0,i)) \nonumber \\
  &&\hspace{1cm} +\frac{K^N_t(K^N_t-1)}{2N^2}\left\{
{\sum \sum}_{i_1\ne i_2}^kg^N(t,\wt x^N_2(t_0,i_1))g^N(t,\wt x^N_{2}(t_0,i_2))\right. \nonumber \\
  && \left. \hspace{5cm} +\sum_{i=1}^kg^{N,2}(t,\wt x^N_2(t_0,i))\right\}
  \Bigg]
  +o(N^{-2}).\nonumber
\eea
Note that because of the convergence of the dual particle system to the collision-free
McKean-Vlasov dual we have:
 \be{gra65c}
 g^N(t,x) \to 1-\sum_{\ell =1}^\infty(1-
 x)^{\ell}U(t,\ell)\text{   as   }N\to\infty,
 \ee
\be{gra65d}
 g^{N,2}(t,x) \la 1-\sum_{\ell =1}^\infty(1- x)^{\ell}U^{\ast 2} (t, \ell)
 \mbox{ as } N \to \infty.
 \ee

Therefore we can conclude from equation (\ref{1-time-2}) that for $i \in \{1,2,\cdots,k\}$:
\be{gra66}
 E[x^N_2(t,i)]=E\left[
\frac{K^N_t}{N}\sum_{i=1}^k g^N(t,\wt x^N_2(t_0,i))
\right]+o(N^{-1}). \ee

Consequently we calculate:
\be{gra67}
 E[\wh x^N_2(t)]=E\left[
{K^N_t}\sum_{i=1}^k g^N(t,\wt x^N_2(t_0,i)) \right]+o(1).
\ee
Then by $K^N_t \to K_t$ as $N\to\infty$, (\ref{gra65c}) and
dominated convergence we get:
\be{gra68}
\lim_{N\to\infty}E[\wh x^N_2(t)]=E\left[
{K_t}\sum_{i=1}^k g(t,\wt  x_2(t_0,i)) \right]
\ee
where we define $g(t,x)$ using (\ref{gra65c}) by
\be{gra69}
g(t,x)=1-\sum_{\ell =1}^\infty(1- x)^{\ell}\CU(t,\ell),
\ee
together with the fact that
$K_t,\CU(t,\cdot)$ can be identified as the number of occupied sites and the size
distribution for the McKean-Vlasov limit dual process obtained
above in Subsubsection \ref{sss.dualcfr} (Step 2).

Finally putting everything together we have:
\be{gra70}
\lim_{t\to\infty} (e^{-\alpha t} \lim_{N\to\infty}E[\wh x^N_2(t)])
=E[W]\sum_{i=1}^k g(\infty, \wt x(t_0,i)).\ee

Finally we have to remove the restriction on the initial state being 0
outside $\CI(\ve,N)$.
We note that the system where $\wt x^N (0)$ is different
from zero only at $\CI(\ve)$ and the true system differ at time
$t$ at most by $Const \cdot \ve e^{\alpha t}$ and hence we get
letting $\ve \to 0$ (in $\CI(\ve)$) finally the claim of part (a).

\noindent (b) The claim on the first moment is immediate from (a).

\noindent (c) We now consider the second moment calculation again
using the dual particle system.
As before in the proof of part (a) we consider first the system arising by setting the type
2 mass equal to zero outside the set $\mathcal{I}(\ve)$ and then
later we shall let $\ve \to 0$ to obtain our claim for the general
case with the very same argument as in part (a).

First note that

\bea{twotime} &&E[\wh x^N_2(t+s)\wh x^N_2(t)] = E[\wh x^N_2(t+s)]
E[\wh x^N_2(t)]
   + cov (\wh x^N_2(t+s), \wh x^N_2(t))\\
&& = E[\wh x^N_2(t+s)] E[\wh x^N_2(t)]+ \sum_{i=1}^N
Cov(x^N_1(t+s,i),
(x^N_1(t,i)))\nonumber\\
&& \hspace{5cm}
+\sum^N_{j\ne \ell} Cov(x^N_1(t+s,j)),x^N_1(t,\ell)).\nonumber\eea
\end{proof}

Since the last term of (\ref{twotime}) has $N(N-1)$ summands it is
necessary to consider terms of order $O(\frac{1}{N^2}) $ in the
computation of $Cov(x^N_1(t+s,j)),x^N_1(t,\ell))$.  Therefore the
calculation is organized to identify the terms up to those of
order $o(N^2)$ as $N\to\infty$, and then to calculate the
contributing terms of $O(N^{-2})$.

To handle the calculation of $Cov(x^N_1(t+s,j),x^N_1(t,\ell))$ both
for $j \neq \ell$ and for $j=\ell=i$ we use the dual process. We
introduce as ''initial condition'' for the dual particle system:
\be{dualx} \mbox{one initial dual particle at $j_1$ at time $0$
and a second dual particle at $j_2$ at time $s$} \ee and then
consider the dual cloud resulting from these at time $t+s$.

Since we shall fix times $0$ and $s$ where the $j_1$, respectively
$j_2$, cloud start evolving we abbreviate the number of sites occupied
by the cloud at time $u$ and similarly the process of the frequency distribution of
sizes of sites by
\be{cloud1}
K^{N,j_1,j_2}_u = K^{N, (j_1,0), (j_2,s)}_u \mbox{ and }
U^{N,j_1,j_2} = U^{N,(j_1,0), (j_2,s)}.
\ee

Precisely this runs as follows. There will be a $Bin(K_t^{N,j_1,j_2},\frac{k}{N})$
number of hits of the $k$-set $\CI(\ve)$, hence as $N\to\infty$ we
obtain:

\bea{2time2} &&\\
&& E[x_{1}^{N}(t+s,j_1)\cdot x_{1}^{N}(t,j_2)]\nonumber\\
&& =E[ \prod_{\ell=1}^{K^{N,j_1,j_2}_{t+s}}
     (1-\wt x^N_2 (t_0, \ell))^{\zeta_\ell^{(t+s)}}]\nonumber\\
&& =E \Bigg[ {(1-\frac{k}{N})^{K_{t+s}^{N,j_1,j_2}}}
  +\frac{K_{t+s}^{N,j_1,j_2}}{N}
    (1-\frac{k}{N})^{K_{t+s}^{N,j_1,j_2}-1}
    \sum_{i=1}^k
    \sum_{\ell=1}^\infty{U^{N,j_1,j_2}_{t+s}(\ell)}(1-\wt x_2^N(t_0,i))^{\ell}  \nonumber\\
&& \hspace{1cm}
+(1-\frac{k}{N})^{K_{t+s}^{N,j_1,j_2}-2}\frac{K^{N,j_1,j_2}_{t+s}(K^{N,j_1,j_2}_t-1)}{2N^2}\nonumber\\
&&
  \hspace{1cm} \sum_{i_1,\,i_2=1}^k \sum_{\ell_1,\ell_2=1}^\infty
   U^{N,j_1,j_2} (t+s,\ell_1)U^{N,j_1,j_2}(t+s,\ell_2)\nonumber\\&&
   (1-\wt x_2^N(t_0,i_1))^{\ell_1}(1-\wt x_2^N(t_0,i_2))^{\ell_2} \Bigg ]
 +o(N^{-2}).\nonumber\eea

 Therefore if we now use the function
 $g^N(t,x)$ defined in (\ref{gra65})
 but now for the initial state specified by the space-time points
$(j_1,0), (j_2,s)$ (note $g^N=g^{N,(j_1,0)}$, respectively $= g^{N,(j_2,0)})$,
 we can rewrite the above expression.
 However there is a small problem with the time index, since our two clouds
 now have a time delay $s$ therefore depending on which of the two clouds has the
 hit of $\CI(\ve)$ we have to use $t$ or $t+s$ as the argument
 in our function $g^N(\cdot, x)$. Here we use therefore the {\em notation} $t(\cdot)$
which indicates that actually the time parameter here is either
$t$ or $t+s$ and depends on which of the two clouds is involved.
However this abuse of notation is harmless since in the arguments
below both the times $t$ and $t+s$ will go to $\infty$.

 We get, with the above conventions,  from the formula (\ref{2time2}) that:
 \bea{gra71}
 &&E[x^N_1(t+s,j_1)x^N_1(t,j_2)] \\
 &&
 = 1- E\left[\frac{K_{t+s}^{N,j_1,j_2}}{N}\sum_{i=1}^k g^N(t (\cdot),
   \wt x^N_2(t_0,i))\right] \nonumber \\
&&
 + E\left[\frac{K^{N,j_1,j_2}_{t+s}(K^{N,j_1,j_2}_{t+s} -1)}{2N^2}
\sum_{i_1\ne i_2=1}^k g^N(t(\cdot), \wt x^N_2(t_0,i_1))g^N
(t(\cdot),\wt x^N_2(t_0,i_2))\right]+o(N^{-2}). \nonumber \eea

Now we continue and calculate the covariances of $x^N_1(t+s,j_1)$
and $x^N(t,j_2)$ based on (\ref{gra71}) and the first moment calculations.
For this purpose we can consider the dual systems starting in
$(j_1,0)$ respectively $(j_2,s)$ and evolving {\em independently}
and the one starting with a factor at each time-space point
$(j_1,0)$ and $(j_2,s)$ but evolving together and with {\em interaction}. As we saw we
can couple all these evolutions via a multicolour particles dynamic.
This now allows to calculate with the first moment formulas
obtained in (\ref{1-time}) to (\ref{gra69}) and using the notation
\be{gra71b}
s_m=0 \mbox{ if } m=1 \mbox{ and } s_m=s \mbox{ if } m=2,
\ee
 and we obtain the covariance formula
 \bea{gra72}
 &&\text{Cov}(x^N_1(t+s,j_1)x^N_1(t,j_2))
 \\
 && = 1- E\left[\frac{K^{N,j_1,j_2}_{t+s}}{N}\sum_{i=1}^k g^N(t(\cdot),\wt
x^N_2(0,i))\right] \nonumber \\
&& \hspace{0.5cm} +
E\left[\frac{K^{N,j_1,j_2}_{t+s}(K^{N,j_1,j_2}_{t+s}-1)}{2N^2}
\sum_{i_1\ne 1_2=1}^k g^N(t(\cdot), \wt x^N_2(t_0,i_1))g^N(t
(\cdot),\wt
x^N_2(0,i_2))\right] \nonumber \\
&& \hspace{0.5cm} -\prod_{m=1}^2 E\Bigg[
1-\frac{K^{N,j_m}(t+s_m)}{N}\sum_{i_m=1}^k g^N(t+s_m,\wt
x^N_2(t_0,i_m) ) \nonumber \\
&& \hspace{0.5cm}
 +\frac{K^{N,j_m}_{t+s_m}(K^{N,j_m}_{t+s_m}-1)}{2N^2}
 \sum_{i_{m,1}=1}^k\sum_{i_{m,2}=1}^k
 g^N (t+s_m,x_{i_{m,1}})g(t+s_m ,x_{i_{m,2}})\Bigg] \nonumber \\
 &&\hspace{0.5cm}  +o(N^{-2}).\nonumber
\eea
We now have to distinguish two cases, namely $j_1 \neq j_2$ and
$j_1=j_2$ and we argue separately in the two cases.

From (\ref{gra72}),  using (\ref{gra66}) and Proposition \ref{L.Indual}, we get that
if $j_1\ne j_2$, then
\bea{2timecov2} &&\text{Cov}(
x_2^N(t+s,j_1), x_2^N(t,j_2))\\&&=
\frac{1}{N^2}\left\{E[N({K^{N,j_1}_{t+s}+K^{N,j_2}_t}-
{K^{N,j_1,j_2}_{t+s,t}})\sum_{i=1}^k g^N(t(\cdot),\wt
x^N_2(t_0,i)) ]\right\}
+o(N^{-2}). \nonumber\eea

The expression for {\em $j_1=j_2$} is a slight modification of this,
since either we have as first step a coalescence and are left with
one particle or we get first a migration event and get a $j_1 \neq
j_2$ case, or finally we get a birth event first. Since we have
the time delay $s$ of course coalescence can act only with this
delay.

Therefore combining both cases up to terms of order $o(\frac{1}{N^2})$ as $N \to
\infty$, \bea{2timecov3} &&
\text{Cov}( \wh x_2^N(t+s), \wh x_2^N(t))\\
&& = \left\{E[N({K^{N,j_1}_{t+s}+K^{N,j_2}_t}-
{K^{N,j_1,j_2}_{t+s,t}})\sum_{i=1}^k g^N(t,\wt x^N_2(t_0,i))
]\right\}
+o(1).\nonumber \eea

Recall that    $g^N(t,\wt x^N_2(t_0,i_1))\to g(t,\wt x_2(t_0, i_1))$
as $N\to\infty$.  Therefore the main point is to identify the
behaviour of
\be{dg22}
N(K^{N, j_1}_{t+s} + K^{N, j_2}_{t+s} - K^{N, j_1,
j_2}_{t+s}), \mbox{ as } N \to \infty,
\ee
which is $N$ times the difference between two
independent dual populations each starting in $j_1$ and $j_2$  and
one combined interacting system. In order to evaluate this we must determine the effect of the collisions of the dual populations.

In order to analyse the dual representation of the r.h.s. of equation
(\ref{twotime}) rewritten in the form (\ref{2timecov3})
we will again work with the enriched coloured
particle system which allows us to identify the contribution, on
the one hand of the two evolving independent clouds, and on the other
hand the effects of the interaction between the two clouds by
collision and subsequent coalescence.

The construction of the
coloured particle system (observing $W_1,W_2,R_1,R_2,R_{1,2},G_1,G_2,G_{1,2}$
particles which are relevant for the accuracy needed here) was given above in the sequel of (\ref{sys1}). Recall that  the
$G_{1,2}$ particles are produced in the coloured particle system
from collision of $W_1$ and $W_2$ particles, that is, $W_1$ and
$W_2$ particles at the same site which then suffer coalescence.
Such an event reduces the number of occupied sites in the
interacting clouds by 1 and then creates a $G_{1,2}$ family. Then
the number of occupied sites in the $W_1,W_2$ interacting system
$K^{N,j_1,j_2}_t$ is given by (\ref{gra62}).

First we observe that if we evaluate the duality relation (recall the
form of the initial state we consider here) we see that the contributing terms of order
$O(\frac{1}{N^2})$ arise from either double or single hits of the
dual particle cloud of a point in $\CI(\ve)$.

We note that in general there are two cases which must be handled
differently, namely, the case in which the difference in (\ref{dg22}) is based on
a Poisson collision event with mean of order $O(\frac{1}{N})$ as
in Lemma \ref{L.dd917} (b) and the other in which there is a law of
large numbers effect as in(\ref{gra9a}) and Lemma \ref{L.dd917}(a).
However in both cases the expected value used in the calculations
below have the same form. In the present case for fixed $t$ we are
in the Poisson case with vanishing rate so that at most one
Poisson event contributes to the limit.


Using Proposition \ref{L.Indual} (\ref{gra62}) we se that the key
quantity in (\ref{dg22}) satisfies as $N \to \infty$ that:
\be{dg61}
E[N({K^{N,j_1}_{t+s}+K^{N,j_2}_t}- {K^{N,j_1,j_2}_{t+s,t}}) ]
\sim N\cdot E[\frak{g}_{(1,2)}(W^\prime (\cdot), W^{\prime \prime} (\cdot),N,s,t)]
\ee where
$\frak{g}_{(1,2)}(W^\prime (\cdot), W^{\prime \prime} (\cdot),N,s,t)$ is given by (\ref{dd658+}) and has mean
given in (\ref{dd658++}).
Hence we conclude with (\ref{dd658++}) that the limit of the r.h.s. of (\ref{dg61}) is given by
\be{addx1}
2 \kappa_2 (E[W])^2 e^{\alpha(2t+s)}.
\ee

\begin{remark}
Here we have noted that using the coupling introduced above the number of {\em single
hits} of the $k$ sites in $\CI(\ve)$ coming from either one of the two
interacting  clouds, can be
obtained by deducting from the  single hits of the "non-interacting
$W_1$ and $W_2$ particles" the single hits by $G_{1,2}$ particles
in the interacting system. Moreover the event of a {\em double hit}, that two occupied sites in the dual cloud coincide
with one of the $k=k(\ve)$ sites, can occur from either two
$W_1$ sites, two $W_2$ sites or one $W_1$ and one $W_2$ site.
(Since the number of $G_{1,2}$  sites or including sites at which
both $W_1$ and $W_2$ particles coexist are of order
$O(\frac{1}{N})$, the contribution of double hits involving
$G_{1,2}$ sites is of higher order and can be ignored.)\\

\end{remark}

Taking the limit as $N\to\infty$   in (\ref{2timecov3}), using (\ref{addx1}), that leads to
\be{2timecov3a2}\begin{array}{l} Cov (\wh x^N_2 (t+s), \wh x^N_2
(t)) \Ntoo \\ 2\kappa_2 e^{\alpha(2t+s)}
E[(W^\prime (t+s)W^{\prime \prime} (t) \suml^k_{i=1} g (t, \wt x_2(t_0,i)].\\
\end{array}
\ee

Recalling (\ref{gra69}), (\ref{earlyd})  and the fact that $\mathcal{U}(t,k)\to
\mathcal{U}(\infty,k)$ as $t\to\infty$, we get  next as $t \to \infty$:
\be{2timecov3b}\begin{array}{ll}
&\mbox{r.h.s. of (\ref{2timecov3a2}) }
\sim  2 e^{\alpha(2t+s)}E\left[ \kappa_2 W^\prime W^{\prime \prime}
\sum_{i=1}^k g(t,\wt x_2(t_0,i))\right]\mbox{   and}\\
 & \left| e^{-(2t+s)}Cov (\wh x^N_2 (t+s), \wh x^N_2
(t)) - 2(E[W])^2 \left[ \kappa_2 \suml^k_{i=1} g(\infty, \wt x_2(t_0,i))
\right]\right|\to 0,
\end{array}
\ee
uniformly in s.
To verify the uniformity we work with the expression for the expected number of $G_{1,2}$ sites
(in the $N \to \infty$ limit dynamics), recall here (\ref{add25}), to get:
\be{add59}
E\left[ \int_s^{t+s}W^\prime(u)W^{\prime\prime}(u-s)e^{\alpha(u-s)}e^{\alpha u}W^g(t+s,u)e^{\alpha(t+s-u)}du\right].
\ee
Then using that for the $N \to \infty$ limiting dynamics always
$W(u)\to W$ as $u\to\infty$ and showing that the contribution from the early collisions (in $[0,s+s_0(t)]$, $s_0(t)=o(t)$) form a negligible contribution as in the
argument following (\ref{earlyd}) we obtain
\bea{add60}
&& e^{-\alpha(2t+s)}E\left[ \int_s^{t+s}W_1(u)W_2(u-s)e^{\alpha(u-s)}e^{\alpha u}W^g(t+s,u)e^{\alpha(t+s-u)}du\right]\\&&\qquad\qquad -(E[W])^2 \left[ \kappa_2 \suml^k_{i=1} g(\infty, \wt x_2(t_0,i))
\right]\leq \rm{const}\cdot e^{-\alpha t}\mbox{   as  }t\to\infty.\nonumber
\eea
This proves (\ref{gra63}) and completes the proof of (c).

(d) follows immediately from (c) and (a).

(e) To verify this note let $\wh{x}_2^{N,m,[t_0,\infty)}(t)$ denote the contribution to the
population resulting from rare mutations that occur after $t_0$.
Then
\be{add61}
E[e^{-\alpha t}\lim_{N\to\infty}\wh{x}_2^{N,m,[t_0,\infty)}(t)] \leq\frac{m}{\alpha}e^{-\alpha t_0}.\ee
Since the above results are valid for any $t_0$, we can let $t_0\to\infty$, to conclude that they remain true if we do not turn off the rare mutation at time $t_0$.

\begin{remark}
In the above calculations we have included    double hits - that is two clouds including
the same site  and that conditioned on a hit the number of
particles at that site is given by the size distribution and
therefore converges to the stable size distribution as
$t\to\infty$.
\end{remark}

\begin{proof} {\bf of Corollary \ref{C.dbetatimes}}

Again we use the dual representation of the moments and
analyse the dual based on the enriched coloured particle system.
We  consider separately two cases $0<\beta<1/2$ and $1/2 \leq \beta<1$,
which correspond to the distinction whether two clouds of
dual particles descending from the two initial particles do not meet as $N \to \infty$ in the
$(\beta/\alpha) \log N$ time scale with a nontrivial
probability or whether they do next.

The proof in the case $0<\beta<1/2$ (where the clouds do not meet
asymptotically) follows the same lines as
for Proposition \ref{P.DetGroDro-new} only that now in
(\ref{AGr8}) the limits $N \to \infty$ and $s\to \infty$ are
carried out at once. Hence the additional step is to verify that
the error terms remain negligible when $s$ is replaced by $s_N$
with \be{x2} s_N = t_N(\beta)-t. \ee
Now we have to analyse the first and second moments in this time scale.

To prove the claims we first note that for $\beta<1$ the probability
that the particle system hits a given site still goes to $0$ as
$N\to\infty$. Therefore with (\ref{gra67}) we still can apply the
CMJ-theorem on exponential growth which takes care of the {\em first
moment} term asymptotic.

Moreover,  the remainder term in the
calculation of the {\em second order term} $s$, i.e.  of \be{x2b}
E[e^{-\alpha t}x^N_1(t,j_1)e^{-\alpha
t_N(\beta)}x^N_1(t_N(\beta),j_2)] \ee has the order \be{gra77}
O(\frac{e^{\alpha t_N(\beta)}}{N^3})=o(N^{-2}).\ee This gives the
claim.

The additional consideration in case $\frac{1}{2}<\beta <1$ is
that there is a positive probability that two clusters will
collide in this time scale.  However the question is whether the collisions
with their effects can change the occupation integral appearing in
the dual representation in the exponential. We see that
since the number of those sites
at which collisions occur among the total of the
$N$ sites  is $o(N)$ we do not need to
consider higher level occupation.

\end{proof}

\subsection{Relation between $^\ast\CW$ and $\CW^\ast$}
\label{ss.relationsW}

Let us briefly review the two limiting regimes we have considered.   In the first regime
we have considered
  the limiting droplet process $\gimel^m_t$ which is a $\mathcal{M}_a([0,1])$-valued Markov process.  We have also proved that there exists $\alpha^\ast >0$ such that
\be{add62}
\{e^{-\alpha^\ast t}\gimel^m_t([0,1])\}_{t\geq 0} \mbox{ is tight}\ee
and determined the limiting first moment and variance of $e^{-\alpha^\ast t}\gimel^m_t([0,1])$
as $t\to\infty$.  This has been obtained by considering the process $\gimel^{N,m}_t$ for times in
$[0,t_N]$ with $t_N\to\infty,\;t_N-\frac{\log N}{\alpha^\ast}\to -\infty$ and letting $N\to\infty$.
Moreover by  (\ref{agrev43h}) it follows that $\gimel^{m}_t([0,1])$ is Cauchy as $t\to\infty$
and therefore has a limit $\CW^\ast$ in $L^2$ as $t\to\infty$.

On the other hand there exists $\alpha >0$ such that the  empirical process $\int_0^1 x \Xi^{\log,\alpha}_N(t,dx)$ in time window $\{\frac{\log N}{\alpha}+t\}_{t\in \R}$ converges to $\CL_t, $ and the limit $^\ast\CW = \lim_{t\to -\infty} e^{-\alpha t}\int_0^1 x\mathcal{L}_t(dx)$ exists and is random.

Next recall that we have established  that the two exponential growth rates
$\alpha$ and $\alpha^\ast$ in the two time regimes $t_N \uparrow
\infty, t_N=o(\log N)$ as $N \to \infty$ respectively
$t_N = \frac{1}{\alpha} \log N +t$ as $N \to \infty, t\to -\infty$
are the same, so that in that respect limits interchange.

It is then natural to  conjecture that also the distributions
of $^\ast\mathcal{W}$ and ${\mathcal{W}}^\ast$ are the same even
though there is the following obstacle.  Define
\be{agrev43a2} t_N(\beta,t) = \frac{\beta}{\alpha} \log
N+t, \mbox{ with } 0<\beta<1. \ee Consider the law $\CL[\bar
x^N_2(t_N(\beta,t))]$ and note that for fixed $t$ there is a
discontinuity at $\beta =1$ since then collisions of the dual
occur  at order $O(1)$ and we know that as function of
$\beta$ we jump from type-2 mass zero for $\beta<1$ to a positive
value at $\beta =1$. Hence the matter to relate $^\ast \CW$
and $\CW^\ast$ is very subtle.

In order to investigate the relation between
the laws of $^\ast\CW$ and $\CW^\ast$ we next consider the relations
between the dual process at times corresponding to  the droplet
formation and the  dual  in the time scale corresponding to the
macroscopic emergence. In this direction we verify first that the
first two moments of $\CW^\ast$ and $^\ast \CW$ agree and then use
the techniques of proof to conclude with $L_2$-arguments that the
two random variables have asymptotically as $N \to \infty,
t \to + \infty $ respectively $N \to \infty, t \to -\infty$ $L_2$-distance zero.

Recall the notation
$t_N (\beta,t) = \beta \alpha^{-1} \log N +t \quad , \quad
t_N =\alpha^{-1} \log N , \mbox{ which we shall use in the sequel}.
$

\beP{P.CW}{(Relation between $^\ast\CW$ and $\CW^\ast$)}

(a) For $0 < \beta <1$, \be{gra79}
  \lim_{t\to
-\infty}\lim_{N\to\infty} E\left[\left(e^{-\alpha t_N(\beta,0)}\wh
x^N_2(t_N(\beta,0))-e^{-\alpha(\frac{\log N}{\alpha}+t)}\wh
x^N_2(\frac{\log N}{\alpha}+t))\right)^2\right] =0.\ee

(b) We have the relation:
\be{gra80} \CL{[^\ast\CW]}=\CL[\CW^\ast]. \qad \ee
\end{proposition}

\begin{remark}\label{R.result}
We can summarize the results in (\ref{gra79}) as follows.
Consider for $0<\beta<1$, with  $\bar x^N_2(t):= \frac{\wh
x^N_2(t)}{N}$, the quantities: \be{gra93}
 \{\wh x^N_2(t):0\leq t\leq \frac{\beta\log N}{\alpha}\},\quad
\{\bar x^N_2(\frac{\log N}{\alpha} + t):\frac{(\beta -1)\log
N}{\alpha}\leq t\leq 0\}. \ee
Then we have shown so far that these two processes converge in law to
\be{gra94}
 \{\wh x_2(t):0\leq t<\infty\},\quad
\{\bar x_2(t):-\infty< t\leq 0\}, \mbox{ respectively}. \ee

We can, by the classical embedding theorem of weakly convergent distributions
on a Polish space, construct both sequences and
their limits on one probability space $(\Omega, \CA,P)$ preserving the law
of each random variable in question and such that we have a.s.
convergence in the supnorm on compact time intervals. On this probability
space we consider $L^2(\Omega, \CA,P)$.

Then we have:
\be{gra95}
 \CW^\ast =L^2 - \lim_{t\to\infty}e^{-\alpha t} \wh x_2(t)
 = L^2 - \lim_{t\to -\infty}e^{-\alpha t}\bar x_2(t) = ^\ast \CW.
 \ee
\end{remark}

 \beC{C.hmdg} {(Equality of moments for growth constant)}
\be{agrev43f} E[(^\ast\CW)^k]=E[(\CW^\ast)^k],\quad k=1,2,\dots. \qquad
\square\ee
\end{corollary}
\begin{proof} {\bf of Corollary \ref{C.hmdg}}   This follows from Proposition \ref{P.CW}.
\end{proof}

\begin{remark}
We can calculate the moments of $^\ast\CW$ and
$\CW^\ast$ by explicit calculation. However since we have not proved
that all moments are finite we cannot use this to prove the
equality in law (\ref{gra80})).

\end{remark}

\begin{proof} {\bf of Proposition \ref{P.CW}}

We start with the rough idea. We have calculated $Var [\CW^\ast]$
in (\ref{agrev19c}) which involved the time scales smaller than
$\alpha^{-1} \log N$. An expression for  $Var [^\ast\CW]$ was
derived in (\ref{ag40}). We see that $Var [\CW^\ast]=Var
[^\ast\CW]$. We shall now prove (\ref{gra79}) which will give
(\ref{gra80}). The point is that the moment calculations now have to
include time scales up to $\alpha^{-1} \log N$ and must therefore
different from before treat collisions of the dual cloud for a given point occuring
with positive probability.

We proceed in five steps.
We again start  by considering the case without mutation but start instead with
positive but finite initial  type-2 mass.
We calculate in Step 1 first moments, in Step 2
second moments at the given time, in Step 3 we calculate the mixed
moment between the two times specified in the assertion under our assumption.
In Step 4 of our arguments we include  mutation. In Step 5 we
conclude the proof.

{\bf Step 1} We begin with the {\em first moment} calculation.
Consider the calculation of the expected value using the dual
particle system starting with one particle at $j$ and recall
$\CI(\ve)$ from (\ref{agre62c}).  We will assume
that at time 0 we have an initial distribution as in Proposition
\ref{P.DetGroDro-new} and that furthermore $m=0$. As above there
will then be a $\text{Bin}(K^N_{t_N(1,t)},\frac{k}{N})$ distributed number of hits of
the $k$-set $\CI(\ve)$ by the cloud of dual particles. To apply the Poisson
approximation recall that by (\ref{agr28})-(\ref{agrev16xxx}) we have:

\be{agr28xxx}
K^{N,(a)}_t = W_a(t) e^{\alpha t} -\mathfrak{g} (W_a(\cdot), N,t) + \CE_a(N,t),
\ee
where the error term satisfies
\be{add64}
E[\CE_a(N,t)]\leq \rm{const}\frac{e^{2\alpha t}}{N} \cdot\sup_u E[W_g(u)].\ee

 Therefore the
number of hits of the $k$-set converges as $N \to \infty$ weakly to a {\em Poisson
distribution with parameter} $We^{\alpha t}-\kappa_2 W^2e^{2\alpha
t}+O(e^{3\alpha t})$. Therefore from (\ref{2time2}) we get looking
at 0,1 or 2 hits of $\CI(\ve)$ and denoting by
\be{gra80c}
\{\xi_\ell,\quad  \ell =1,\cdots, K^N_{t_N(1,t)}\}
\ee
the states of the sites occupied by the dual at time $t_N(1,t)$
 and by $K^N_t,U^N_t$  the number
of occupied sites and the size distribution for the finite $N$
dual (nonlinear) system studied above:
\bea{gra81} &&
E[x^N_1(t_N(1,t),j)]=E\left[\prod_{\ell=1}^{K^N_{t_N (1,t)}}
(1-x^N_2(t_0,\xi_\ell))^{\zeta^N (t_N(1,t),\xi_\ell)}\right] \\
&&
=E\Bigg[(1-\frac{k}{N})^{K^N_{t_{N}(1,t)}}\nonumber\\
&&
+(1-\frac{k}{N})^{K^N_{t_N(1,t)}-1}\quad  \frac{K^N_{t_N(1,t)}}{N}
\sum_{i=1}^k \sum_{\ell=1}^\infty U^N(t_N(1,t),\ell)(1-\wt
x_2^N(t_0,i))^{\ell}\nonumber \\
&&
+(1-\frac{k}{N})^{K_{t_N(1,t)}-2}\quad
\frac{K^{N}_{t_N(1,t)}(K^{N}_{t_N(1,t)}-1)}{2N^2}\nonumber\\
&& \sum_{i_1,\,i_2=1}^k \sum_{\ell_1,\ell_2=1}^\infty
U^N(t_N(1,t),\ell_1)U^N(t_N(1,t),\ell_2)(1-\wt
x_2^N(t_0,i_1))^{\ell_1}(1-\wt
 x_2^N(t_0,i_2))^{\ell_2} \Bigg ]\nonumber\\
 && +o(N^{-2}). \nonumber
\eea

We rewrite this in the form
\bea{gra82}
&&E[x^N_1(t_N(1,t),j)]\\
&& =E\Bigg[ 1-(W_N(t_N) e^{\alpha
t}+O(W^2_N(t_N) e^{2\alpha t}))\sum_{i=1}^k g^N(t_N(1,t),\wt
x^N_2(t_0,i))\nonumber\\
&&\qquad\qquad
 +\frac{W_N(t_N) e^{2\alpha t}}{2}\sum_{i_1=1}^k\sum_{i_2=1}^kg^N
 (t_N(1,t),x_{i_1})g^N(t_N(1,t),x_{i_2})\Bigg] +o(N^{-2}),\nonumber
 \eea
with
\be{gra86}
 g^N(t_N(1,t),x)=1-\sum_{\ell =1}^\infty(1- x)^{\ell}U^N(t_N(1,t),\ell).
\ee
Therefore as $N \to \infty$, we have
(recall that $W_N (t_N) \to W$ a.s. as $N \to \infty$ and $W_N (t_N)$
are uniformly square integrable):
 \be{gra83}
 E[x^N_2(t_N(1,t),j)]=E\left[
(We^{\alpha t} +O(W^2 e^{2 \alpha t})) \sum_{i=1}^k
g^N(t_N(1,t),\wt x^N_2(t_0,i)) \right]+o(N^{-2}). \ee

Next we need
to know more about $g^N(t_N(1,t), \wt x_2^N (t_0,i))$ as $N \to
\infty$ and then $t \to -\infty$. Recall that  $
g^N(t_N(1,t),x)=1-\sum_{\ell =1}^\infty(1- x)^{\ell}U^N(t,\ell)$ (
see
 \ref{gra65}), $\lim_{N\to \infty} U^N(t_N(1,t),\cdot)= U(t,\cdot)$ (see \ref{agre60})   and $\lim_{t\to -\infty}
U(t,\cdot,\cdot)=\mathcal{U}(\infty,\cdot,\cdot)$ (see
\ref{Grev39p1+}). Therefore
 \be{gra83b} \liml_{t\to -\infty}
\liml_{N\to\infty} g^N (t_N(1,t), \wt x^N_2(t_0,i)) = g(\infty, \wt
x_2(t_0,i))=1-\sum_{\ell =1}^\infty(1-\wt x_2(t_0,i))^{\ell}\CU(\infty,\ell).
\ee

This allows to calculate next as $N \to \infty, t \to - \infty$:
\be{gra84} E[\wh x^N_2(t_N(1,t))]=E\left[ WNe^{\alpha
t}\sum_{i=1}^k g^N(t_N(1,t),\wt x^N_2(t_0,i)) \right]+o(e^{\alpha
t}N^{-1}). \ee Then \be{gra85}
\liml_{t \to -\infty}
\lim_{N\to\infty}E[e^{-\alpha t_N (1,t)} \wh
x^N_2(t_N(1,t))]=E\left[ W \sum_{i=1}^k g( \infty,\wt x_2(t_0,i))
\right]. \ee

Finally putting everything together results in
\be{gra87}
\lim_{t\to -\infty}  \lim_{N\to\infty}
E [ e^{-\alpha t_N(1,t)} (\wh x^N_2(t_N(1,t))]
=E[W]\sum_{i=1}^k g(\infty,\wt x(t_0,i)). \ee

{\bf Step 2} We now consider the {\em second moments}. Consider the time points
$t_N(\beta_1,0)$ and $t_N(\beta_2,t), \beta_1 < \beta_2$. The case
$\beta<2$ follows as in Proposition \ref{P.DetGroDro-new}.  The case
$\beta=1$ requires more careful analysis, we carry out now. Let
$0<\beta\leq1$, then
\bea{2time2mf}
&& E[x_{1}^{N}(t_N(\beta,0),j_1)\cdot x_{1}^{N}(t_N(1,t),j_2)]=\\
&& =E[\prod_{j=1}^{K^{N,j_1,j_2}_{t_N(1,t)}}
\prod_{\ell=1}^{\zeta(t_N(1,t))}
    (1-\wt x^N_{2}(t_0,j))^\ell]\nonumber\\
&& =E \Bigg[
     (1-\frac{k}{N})^{K_{t_N(1,t)}^{N,j_1,j_2}}
     +\frac{K^{N,j_1,j_2}_{t_N(1,t)}}{N}(1-\frac{k}{N})^{K^{N,j_1,j_2}_{t_N(1,t)}-1}
     \sum_{i=1}^k \sum_{\ell=1}^\infty{U^N_t(\ell)}(1-\wt x_2^N(t_0,i))^{\ell}  \nonumber\\
 && +(1-\frac{k}{N})^{K^{N,j_1,j_2}_{t_N(1,t)}-2}
      \frac{K^{N,j_1,j_2}_{t_N(1,t)}(K^{N,j_1,j_2}_{t_N(1,t)}-1)}{2N^2}\nonumber\\
 &&\sum_{i_1,\,i_2=1}^k \sum_{\ell_1,\ell_2=1}^\infty
     U^N(t,\ell_1)U^N(t,\ell_2)(1-\wt x_2^N(t_0,i_1))^{\ell_1}(1-\wt
     x_2^N(t_0,i_2))^{\ell_2} \Bigg ]\nonumber\\&&+o(N^2).
 \nonumber\eea

Proceeding as in the proof of Proposition \ref{P.DetGroDro-new} in
the case where $j_1\ne j_2$ we get: \bea{gra89} &&\\&&\text{Cov}(
x_2^N(t_N(\beta,0),j_1), x_2^N(t_N(1,t),j_2))\nonumber\\&&=
\frac{1}{N^2}\left\{E\left[N\left(K^{N,j_1}_{t_N(1,t)}+K^{N,j_2}_{t_N(\beta,0)}-
{K^{N,j_1,j_2}_{t_N(1,t),t_N(\beta,0)}}\right)\sum_{i=1}^k
g^N(t_N(\cdot),\wt x^N_2(t_0,i))\right]\right\} \nonumber\\&&
+o(N^{-2}).
\nonumber\eea

{\bf Step 3}
The next point is to prove the following claim on second moments for $N \to \infty$
and then as $t \to -\infty$. For $0< \beta\leq
\beta_2\leq 1$,
\be{gra78}
  \lim_{t\to -\infty}\lim_{N\to\infty} E\left[\left(e^{-\alpha t_N(\beta,0)}\wh
x^N_2(t_N(\beta,0))\cdot e^{-\alpha t_N(\beta_2,t)}\wh
x^N_2(t_N(\beta_2,t)))\right)\right] =\text{const.}, \ee
with a constant independent of $\beta, \beta_2$.

In order to verify (\ref{gra78}) we must now include  new
considerations that arise since $e^{\alpha t_N(1,t)}=Ne^{\alpha
t}$ leading to a positive particle density and hence collisions.
Recall our representation in terms of the {\em coloured particle
system} and a careful analysis of the number of $W_1$, $W_2$ and
finally $G_{1,2}$ particles. We begin with the $W_1$ and $W_2$
particles.

The key point is that in the analysis of the $\kappa_2$ term one
needs to include slower growth of the $W_1$ and $W_2$ particles
due to the nonlinear term.  In particular the
populations of white particles are asymptotically as $N \to \infty$  of the form $N
F(t)$ and the function $F(t)$ can be expanded as $t \to -\infty$
in terms of $k$-th powers $(k=1,2,\cdots)$ of $e^{\alpha t}$.
Hence the number of particles in the two white clouds,
the one starting in $j_1$ and the one starting in $j_2$ are each
asymptotically as $N\to \infty$ behaving as
\be{gra78b}
We^{\alpha t_N(\gamma,t)}- \frac{1}{N}\kappa_2 W^2
e^{2\alpha t_N(\gamma,t)}\sim WN e^{\alpha t}-\kappa_2 W^2Ne^{2 \alpha t}
= N e^{\alpha t} W(1-\kappa_2 W e^{\alpha t}),
\ee
where we take two independent realisations of $W$ for each cloud and with
$\gamma= 1-\beta$. The system with start of a particle at $j_1$ at time
0 and $j_2$ at time $t_N(\beta,0)$ evolves similar but now an additional
correction is needed because of the collisions between the two clouds.
This effect is represented by the $G_{1,2}$ particles.

Hence  the important term for the covariance calculation is the
one involving the $G_{1,2}$  particles. The number of $G_{1,2}$
particles is taking again $N\to \infty$ and then afterwards $t \to -\infty$
of order
\be{gra78c}
O(W^\prime W^{\prime \prime} \kappa_2 e^{2\alpha t}).
\ee
The  self-intersections of the $W_1$ and the
$G_{1,2}$ particles produce only a higher order correction to the number and
size distribution of the $G_{1,2}$  particles. As a result in the
limit as $t\to -\infty$ the correction terms are all of higher
order.

 From together with (\ref{gra89}, \ref{gra78b}) and (\ref{gra78c})
 we can verify that the covariance term is
asymptotically the same as (\ref{gra63}), that is,

\bea{x12}
&& \lim_{t\to -\infty}\lim_{N\to\infty}\text{cov}
(e^{-\alpha (t_N(1,t))} \wh x^N_2 (t_N(1,t)), e^{-\alpha
t_N(\beta,0)} \wh x^N_2 (t_N(\beta,0)))\\&&
=2(E[W])^2\left\{\kappa_2 \sum_{i=1}^k g(\infty,\wt
x_2(t_0,i))
\right\}
  .\nonumber
  \eea
 and we get the claim in (\ref{gra78}).

{\bf Step 4} We are now ready to conclude the argument for the two
   assertions of the proposition. First we
 have to show the initial masses  we used  in the calculation above can
 in fact arise as configuration at time $t_0$. This follows from (\ref{agrev2}), where
 we show convergence of the  $\gimel^{m,N}$ to $\gimel^m$ (recall (\ref{ag41e})).
 Second we have to prove that the case we treated can be
replaced by the masses which arise from mutation. For any fixed $t_0$ these
masses satisfy (\ref{agre61b}), (\ref{agre62b}) and the above
argument gives therefore the result for mutation turned of after time $t_0$.
We now have to argue that this is not relevant for large $t_0$.

The idea is that the contribution to
$\CW^\ast$ due to mutation after $t_0$ is negligible for large $t_0$,
more precisely, the expected contribution is $O(e^{-\alpha
t_0})\to 0$ as $t_0\to\infty$. The result then follows by letting
$t_0\to \infty$ in Corollary \ref{C.AG0} below.

{\bf Step 5} Next we have to go through the list of claims.

(a) This follows from (\ref{x12}) above, (\ref{agrev43i}) together with (\ref{gra63})
and (\ref{gra78}) for $\beta = \beta_2=1$ since we see that both the variances of
the two random variables and their covariance all converge to the same value.

(b) Follows from (a) by using the construction of Remark \ref{R.result}
$(\CW^\ast, ^\ast\CW)$ both on one probability space) and
noting that  with $t_N(0,t)=t$
\be{gra90}
 e^{-\alpha t_N(1,t)}\wh x^N_2(t_N(1,t))= e^{-t}\bar x^N_2(t_N(1,t)),\ee
\be{gra91}
 \CW^\ast=\lim_{t\to\infty}\lim_{N\to\infty}e^{-\alpha t_N(0,t)}\wh
x^N_2(t_N(0,t)) \ee and \be{gra92}
 ^\ast\CW =\lim_{t\to -\infty}\lim_{N\to\infty}e^{-t}\bar
x^N_2(t_N(1,t)).
\ee
This completes the proof of Proposition \ref{P.CW} if we could
show that in (\ref{gra92}) we can replace $t_N(0,t)$
and $t \to \infty$ by $t_N$ with $t_N \to \infty$ but
$t_N=o(\log N)$.
Hence the proof of Proposition \ref{P.CW} follows from the next
corollary.
\end{proof}

\beC{C.AG0} {}
Consider $t_N \uparrow \infty$ with $t_N =o(\log N)$. Then:
\be{agre64b} \liml_{t \to \infty} \liml_{N \to \infty}
E[(e^{-\alpha t} \wh x^N_2(t) - e^{-\alpha t_N} \wh x^N_2(t_N))^2]
=0 \qquad \square \ee
\end{corollary}

{\bf Proof of Corollary \ref{C.AG0}}

If we consider the case in which we exclude mutation but insert
initially some finite mass (even for $N \to \infty$) then the
result follows  due to (\ref{agrev43i}). What we want is however
the statement where initially we have no type-2 mass but where we
have mutation from type 1 to type 2 at rate $mN^{-1}$. The idea
now is to use the fact that  only early mutation counts, that is,
mutation after some late time has a negligible effect.  We can use
this idea together with coupling to prove the result  with
mutation.

We consider solutions of our basic equation for mutation rate $m$
and mutation rate 0. We denote the solutions by \be{agre65}
(X^{N,m}_t)_{t \geq 0} \mbox{ and } (X^{N,0}_t)_{t \geq 0}. \ee We
can define these two processes on one probability space so that
the following relation holds: \be{agre66} \wh x^{N,0}_2 (t) \leq
\wh x^{N,m}_2(t),
\ee
with
\be{agre67}
E[\wh x^{N,m}_2(t) - \wh x^{N,0}_2 (t)] \leq Const \cdot m \intl^t_0 e^{\alpha s} ds.
\ee
We now apply this to compare the process with mutation turned off
after  time $t_0$ with the original process with mutation running
on $[0,\infty)$.  We know $\wh x^{N,0}_2 (t_0)= \wh x^{N,m}_2(t_0) =O(e^{\alpha t_0}))$
and with (\ref{agre66}) and (\ref{agre67}) we have
\be{agre68}
0 \leq E[e^{-\alpha (t_0+t)}\wh x^{N,m}_2(t_0+t) -
e^{-\alpha (t+t_0)}\wh x^{N,0}_2 (t_0+t)] \leq Const \cdot m
e^{-\alpha t_0}. \ee Letting $t_0\to\infty $ we get that we can
approximate the systems with mutations by ones where mutation is
turned off after $t_0$ uniformly in $N$. Then the claim follows by
using the result for $O(1)$ initial mass of type 2 and no
mutation.

\subsection{Third moments} \label{ss.highermom}

We expect that the random variables $\CW^\ast (^\ast\CW)$ have moments
of all orders which determine the distribution.  We do not verify this here but in
 this section we verify that  the  third  moment
of the rare mutant type is finite given if we start with total mass of
type $1$ only using the same method as that used for the first and second moments which for
convenience we review here.

\beL{L.Himo} {(First, second and third  moments of $\CW^\ast$)}

Consider  times $t_N \to \infty$ with $t_N=o(\log N)$.
The first, second and third moment of the rescaled total type two mass satisfies:

\be{Ma3xx1}        \lim_{N\to\infty} e^{-\alpha t_N}
                   E[(\wh{x}^N_2(t_N))] = m^*E[W],
		   \ee

\be{Ma3xx2} \lim_{N\to\infty} e^{-2\alpha t_N}
E[(\wh{x}^N_2(t_N))^2] =(m^{*})^2(E[W])^2+2m^*\kappa_2(E[W])^2
\ee
and for some $0<\kappa_3<\infty$, $m^*=\frac{1}{c}(1+\frac{\gamma}{\alpha})m$ (as in (\ref{dw9})), $\kappa_2$ given by
(\ref{agrev19d7zz})
\be{Ma3}
\lim_{N\to\infty} e^{-3\alpha t_N}
E[(\wh{x}^N_2(t_N))^3]= (m^{\ast})^3(E[W])^3+6(m^*)^2\kappa_2 (E [W])^3
 +6m^*\kappa_3 (E [W])^3. \hfill \qquad \square
 \ee
\end{lemma}

\begin{proof} {\bf of Lemma \ref{L.Himo}}

The proof will be based on duality and as a prerequisite we need (in order to determine the moment
asymptotics) to evaluate the
expression which represents the moments in terms of the dual particle system. We therefore
carry out the proof in {\em two steps}, first an expansion of dual expression and then the
moment calculation.

{\bf Step 1} {\em (Expansion of $\Pi^{N}(t)$ up to terms of order $O(\frac{1}{N^3})$)}

In this subsection we extend the coloured particle system in order to derive and expansion of
$\Pi^{N,a}_{t}$ including terms of order $\frac{1}{N}$ and $\frac{1}{N^2}$ which are needed in the calculation of
the second and third moments of $\wh x^N_1(t)$.
\beL{kappa3}{(Expansion of $\Pi^{N,a}_{t_N}$)}

Consider the dual particle system starting with one, two or three particles.
To obtain a unified expression  we state the result using  $W_a$.

The corresponding CMJ
random growth constants (recall (\ref{ang2b}) $W_a=W_1$ with one initial particle, $W_a=W_1+W_2$ if two particles start at different sites, $W_a=W^{(2)}$ if two particles start at the same site,
$W_a=W_1+W_2+W_3$ is three particles start at different sites and $W_a=W_1^{(2)}+W_2$ if two particles start at one site and one particle at
a different site.

Then there exists $\kappa_3>0$ such that
\bea{agr35}
\Pi^{N,a}_{t_N} \sim &&\frac{1}{c}(\alpha+\gamma)\left[ (W_a) e^{\alpha t_N}
- \frac{2 \kappa_2}{N} (W_a^2)  \frac{e^{2 \alpha t_N}}{N}\right.
\\&&\left.- \frac{3 \kappa_{3}}{N^2}(W_a^3) e^{3 \alpha t_N} +O(\frac{e^{4\alpha t_N}}{N^3})\right].\nonumber
\eea
Then we get for the time-integral and the exponential of it:
\bea{agr36}
H_a^N(t):=\frac{m}{N} \intl^{t_N}_0 \Pi^{N}_s ds
\sim \frac{m^\ast}{N} (W_a) e^{\alpha t_N}
&-& \frac{m^*\kappa_2}{N^2}(W_a^2) e^{2 \alpha t_N}\\
&-& \frac{m^*\kappa_{3}}{N^3}   W^3_a e^{3 \alpha t_N}
+ O(\frac{e^{4 \alpha t_N}}{N^4}),\nonumber
\eea
\be{agr362} e^{-H_a^N(t)}=1- H_a^N(t)+\frac{1}{2}(H_a^N(t))^2-\frac{1}{6}( H_a^N(t))^3 +O((H_a^N(t))^4).
\hfill \square
\ee
\end{lemma}
\begin{proof} {\bf of Lemma \ref{kappa3}}
We use an extension of the multicolour system defined so far  in Subsubsection \ref{sss.prep2} (Step 1)
therein and in the above Subsection \ref{ss.relationsW} in order to identify
$\Pi^{N,a}_{t_N}$ up to terms of order $e^{3\alpha t_N}$ with an error term of order $O(\frac{e^{4\alpha t_N}}{N^3})$.
  To obtain this
we must introduce new colours.

Recall in the expansion to order $e^{2\alpha t}$ we worked with the system of white, black, red, green, purple and blue coloured particle system
and that the dual system is given by the white, red and purple
particles. The green particles are produced by coalescence of a red and white particle and
a green site is produced by the migration of a green particle. These represent particles which
are lost due to one collision (followed by coalescence). The white particles grow as
$W_a \frac{\alpha +\gamma}{c} e^{\alpha t}$ and the green particles are produced at rate
$W^2_a e^{2\alpha s}$ at time $s$. Since the green families grow with exponential rate $\alpha$ the total number of green particles grows
like $W^2_a \kappa_2e^{2\alpha t}$.
The number of purple-blue  particles is of order
$O(e^{3\alpha t})$ and the the lower bound is given by the white
and red particles obtained by subtracting the number of green
particles. The upper bound is obtained by adding the blue particles.

In order to identify the {\em next order term} we must include the purple particles (and exclude the blue particles) and
follow what happens when the purple particles coalesce or migrate. This requires to introduce two further colours,
namely {\em pink} and {\em yellow}. When a purple particle migrates and hits a white, red or purple site we now
produce a pink-yellow pair.
There are $O(e^{4\alpha t})$ pink (yellow) particles where now the pink particles provide the dual and the yellow particles now behave as the blue particles did before in obtaining the error term.

We note that the purple particles can coalesce with a
white (or red) particle and when this happens a white (or red) and
green  (to distinguish these we label them $G_2$) pair of particles is produced at the site where coalescence
occurs. Green $(G_2)$ particles are  produced at rate $\text{const}\cdot W^3_a e^{3\alpha t}$.

The two particle systems we have to compare, the dual one and the collision-free one, are now given
as follows.
\begin{itemize}
\item
The dual is given by the {\em white, red, purple} and {\em pink particles} and
\item
the collision free by {\em white, red, green} ($G\cup G_2$), {\em purple} and {\em yellow}.
\end{itemize}

The modified rules are now:

\begin{itemize}
\item  If a purple migrates and hits a site
occupied by white or red  it dies and creates a pink-yellow
pair which have subsequently coupled times of birth, migration,
\item the yellow particle is placed in the first unoccupied
site in the first copy of $\mathbb{N}$,
\item if a purple particles coalesces with a white or red particle it produces a
green $G_2$ particle
\item the pink particles behaves like a standard dual particle.
\end{itemize}

This will allow us to  approximate the dual
with an error of order as: $N \to \infty$:
\be{pink}
O(N^{-3} e^{4\alpha t_N}).
\ee

Now we have to represent the quantity $\Pi^{N,a}_t$ in terms of the numbers of particles of
the various colours. Here we follow the procedure explained in Subsubsection \ref{sss.dropletgrow},
Part 2 therein. Then we obtain the claimed expansion of (\ref{agr35}).
Recall that $\kappa_2 W^2e^{2\alpha t}$ denotes the term corresponding to the growth of the green $G$ particles in subsubsection \ref{sss.dropletgrow}.  In the same way
the $\kappa_3 W^3 e^{3\alpha t}$ in (\ref{agr35}) term corresponds to the set of $G_2$ particles, that is, the particles lost
due to particles involved in exactly two collisions.

\end{proof}

The assertions (\ref{agr36}) and (\ref{agr362}) follow from (\ref{agr35}) by explicitly calculating
the integral and then inserting the expansion in the Taylor expansion of the exponential around
0. q.e.d.

\bigskip

\noi
{\bf Step 2} {\em (Moment calculation)} \quad\\

First recall the dual identity
\be{ax1} E[(x^N_1(i,t))^k]= E\left[e^{-H^N_k(t)}\right],
\ee
where $H^N_k$ is given by (\ref{agr36}) and where $W_a$ is the CMJ random growth constant (recall (\ref{ang2b})) obtained by starting the dual
with $k$ particles at site $i$. We obtain a similar expression for $E[x^N_1(i,t)x^N_1(j,t)],\,$ where now we start one particle at $i$ and one particle at $j$, etc.

\medskip

 We now compute the first three moments following the same method.
 Note that
 \bea{agrev20xx}
 &&E[(\widehat{x}_{2}^{N}(t))]= E[N-\wh{x}_1(t)]= N- NE[x_1(t)]\\&&= N- NE[(1- m^*\frac{W_1}{N}e^{\alpha t} +O(\frac{1}{N^2}))]. \nonumber
 \eea
 Then replacing $t$ by $t_N$ we get

\be{Ma3xx1p}        \lim_{N\to\infty} e^{-\alpha t_N}
                   E[(\wh{x}^N_2(t_N))] = \lim_{N\to\infty} e^{-\alpha t_N}N[\frac{m^*E[W]}{N}e^{\alpha t_N}+O(\frac{1}{N^2})]\ee
 and (\ref{}) follows provided that $t_N=o(\log N)$.

We now consider the second and third moments:

\bea{ax2}
       &&E[(\widehat{x}_{1}^{N}(t))^{2}]=\\
         && NE[(x_{1}^{N}(i,t))^{2}]+N(N-1)E[(x_{1} ^{N}(i,t))x_{1}^{N}(j,t)], \quad i \neq j \nonumber
         \eea
       and for $i,j,k \in \{1,\cdots,N\}$  distinct
       \bea{agrev20}
       E[(\widehat{x}_{1}^{N}(t))^{3}]=
         && NE[(x_{1}^{N}(i,t))^{3}]+3N(N-1)E[(x_{1} ^{N}(i,t))^{2}x_{1}^{N}(j,t)] \\
       && +N(N-1)(N-2)E[x_1^{N}(i,t)x_{1}^{N}(j,t)x_{1}^{N}(k,t)].
       \nonumber \eea

We first illustrate the method of calculation by calculating in detail the second moment.  We substitute (\ref{agr36} and (\ref{agr362}) in (\ref{3mkk}) keeping track of
only those terms that do not go to 0 as $N\to\infty$. We get:

       \bea{2mkk}
       &&E[(\widehat{x}_{2}^{N}(t))^{2}]= E[N-\widehat{x}_{1}^{N}(t)]^2\\
       &&= N^2-2NE[\widehat{x}_{1}^{N}(t)] +E[\widehat{x}_{1}^{N}(t))^{2}]\nonumber
       \\ && =N^2-2N^2 E[x_1^{N}(t)]+NE[(x_1^{N}(t))^2]+N(N-1)E[x_1^{N}(i,t)x_1^{N}(j,t)].\nonumber
       \eea

We now insert for the appearing moments of $x^N_1$ the formula (\ref{ax1}) and then use the
expansion of the exponential from (\ref{agr362}).
       We first compute the terms coming from the $1-H^N(t)$ in the expansion of the exponential.
       This gives
\bea{ax3}
&& N^2 -2N^2 E\left[1-\frac{m^*}{N}We^{\alpha t}+\frac{m^*\kappa_2}{N^2}W^2e^{2\alpha t}+O\left(\frac{e^{3\alpha t}}{N^3}\right)\right]\\&&
       +NE\left[1-\frac{m^*}{N}W^{(2)}e^{\alpha t}+\frac{m\kappa_2}{N^2}W^2e^{2\alpha t}+O\left(\frac{e^{3\alpha t}}{N^3}\right)\right]\nonumber\\&&
       +N(N-1)E \big[1-\frac{m^*}{N}(W_1+W_2)e^{\alpha t}+\frac{m^*\kappa_2}{N^2}(W_1+W_2)^2e^{2\alpha t}\big]\nonumber\\
       &&= N^2(1-2+1)+NE\left[(2m^*We^{\alpha t}+1-1 -m^*(W_1+W_2)e^{\alpha t}) +O\left(\frac{e^{3\alpha t}}{N^3}\right)\right]\nonumber\\&&
       E \left[-2m^*\kappa_2W^2e^{2\alpha t}-m^*W^{(2)}e^{\alpha t}+m^*\kappa_2(W_1+W_2)^2e^{2\alpha t}+O\left(\frac{e^{3\alpha t}}{N^3}\right)\right]\nonumber\\
       &&= 2m^*\kappa_2 (E[W] )^2e^{2\alpha t}+m^*E[W_1+W_2]e^{\alpha t} -m^*E[W^{(2)}]e^{\alpha t} +O\left(\frac{e^{3\alpha t}}{N}\right). \nonumber\eea

       We next note that the $\frac{1}{2}(H^N(t))^2$ gives
\be{ax4}
(m^*)^2(E[W])^2e^{2\alpha t}.\ee
       Here $W_1,W_2$ are independent random variable coming from the CMJ limit starting with particles at two disjoint
       sites and $W^{(2)}$ is the random variable coming from the CMJ limit starting with two particles at the same site.

       Combining all these terms we obtain

\bea{ax5}
&& E[(\widehat{x}_{2}^{N}(t))^{2}] \\
&&= 2(E[W])^2m\kappa_2 e^{2\alpha t}+ 2m^*E[W]e^{\alpha t}-m^*E[W^{(2)}]e^{\alpha t} +(m^*)^2(E[W])^2e^{2\alpha t}\nonumber\\
&&+O\left(\frac{e^{3\alpha t}}{N}\right). \nonumber
\eea
Letting $t=t_N$, multiplying by both sides $e^{-2\alpha t_N}$ and taking the limit $N\to\infty$,
we obtain (\ref{Ma3xx2}) provided that $t_N=o(\log N)$.

\bigskip

       We now follow the same method for the third moment but omit writing out the numerous
       intermediate terms which cancel and get:
       \bea{3mkk}
        && E[(\widehat{x}_{2}^{N}(t))^{3}]\\&&=N^{3}-3N^{2}E[\widehat{x}_{1}
       ^{N}(t)]+3NE[(\widehat{x}_{1}^{N}(t))^{2}]-E[(\widehat{x}_{1}^{N}(t))^{3}]\nonumber\\ &&
        = N^3 -3N^3E[(x_1^N(t))]+3N^2E[(x^N_1(t))^2]+3N^2(N-1)E[x_1^{N}(i,t)x_1^{N}(j,t)]\nonumber\\&&
       -NE[(x^N_1)^3]-3N(N-1)E[(x^N_1(i))^2x_1(j)]-(N^3-3N^2+2N)E[x^N_1(i)x^N_1(j)x^N_1(k)].\nonumber
       \eea

The proof proceeds by substituting (\ref{agr36}) and (\ref{agr362}) via (\ref{ax1}) in (\ref{3mkk}).
One can verify by direct algebraic calculations that the resulting
coefficients of $N^3,N^2,N^1$ are zero.  The $O(1)$ terms involves powers of $e^{\alpha t}$, $e^{2 \alpha t}$ and $e^{3\alpha t}$.  To prove the
Lemma it suffices to identify the coefficient of $e^{3\alpha t}$.  The coefficient involving $(m^*)^3$
comes from the term $\frac{N^3}{6}\cdot(H_N(t))^3$ in the expansion of the exponential.
The $\kappa_2$
term comes from the term $\frac{N^3}{2}\cdot(H_N(t)^2$ and the $\kappa_3$ term comes from the term
${N^3}H_N(t)$  in the expansion of the exponential.

To illustrate this structure we consider the contribution coming from the term
\be{ax6}
  \frac{N^3}{2}(H^N_a)^2=-(m^*)^2\kappa_2e^{3\alpha t}(W_a)^3.\ee
We must apply this to each term in the following expression
\be{ax7} N^3[-3E(x_1^N(t))+3E[x_1^N(i,t)x_1^N(i,t)]-E[x_1^N(i,t)x_1^N(j,t)x_1^N(k,t)]\ee
      The corresponding terms with a $(W_a)^3$ result in
      \be{ax8} -3E[W_1^3]3E[(W_1+W_2)^3] -E[(W_1+W_2+W_3)^3]=-6E[W_1]E[W_2]E[W_3].\ee

       Note that the coefficients of
       $e^{\alpha t},e^{2\alpha t}$ depend on higher moments of $W$. This leads to

       \bea{ax9}
       && E[(\widehat{x}_{2}^{N}(t))^{3}]=m^*\kappa_3(E[W])^3e^{3\alpha t} +6(m^*)^2\kappa_2(E[W])^3e^{3\alpha t}+ 6(m^{*})^3 (E[W])^3e^{3\alpha t}\nonumber\\&& \qquad +\text{const}
        \cdot e^{2\alpha t} +\text{const} \cdot e^{\alpha t}
       +O\left( \frac{e^{\alpha 4t}}{N}\right)\nonumber
       \eea

\noi Letting $t=t_N$, multiplying by $e^{-3\alpha t_N}$ and taking the limit $N\to\infty$ we obtain (\ref{Ma3}) provided that $t_N=o(\log N)$.

\end{proof}

\subsection{Propagation of chaos: Proof of Proposition
\ref{P.chaos}} \label{ss.chaos} Here we prove Proposition
\ref{P.chaos}. We proceed in three steps.

{\bf Step~1} Fix $t_0> -\infty$ and  tagged sites $1,\dots,L$.  We
can again prove for $T_N = \frac{1}{\alpha} log N$ that
$\{x^N_2(1,T_N+t_0),\dots,x^N_2(L,T_N+t_0)\} $ converges in
distribution as $N\to\infty$ using the dual representation to show
that all joint moments converge. Here of course it  is equivalent
to do the same for $x_1(i,T_N+t_0)$ which is more convenient for
calculation.
In particular we  argue as below (\ref{and12}) using the dual
process, except that now we start  $k_{i}$-particles at site
$i,i=1,\dots,L$ to calculate the mixed higher moments of $x^N_1
(1,T_N+t_0), \cdots, x^N_1 (L, T_N+t_0)$, i.e.
\be{angr23c3}
E[\prod_{j=1}^L(x^N_1(j,T_N+t_0))^{k_j}], \quad
k_1, k_2, \cdots, k_L \in \N,
\ee
in terms of
\be{agr38}
E[exp(-\frac{m}{N} \intl^{T_N+t}_0 \Pi^{N,(k_1, \cdots, k_L)}_s ds)]
\ee
using formula (\ref{AGr1}) and (\ref{AGr4}) and denoting with the
superscript $(k_1,\cdots,k_L)$ the initial position.

Then we obtain convergence as $N \to \infty$ from Proposition
\ref{P.Grocoll} part (e). (Note the fact that the $k_i$ here are
possibly different does not change the argument). Therefore we
have  proved the convergence of the marginal distribution at time
$t_0$ of the tagged size-$L$ sample.

{\bf Step~2} Next we verify that the limiting dynamics for $t \geq
t_0$ is as specified in (\ref{Y777}). To do this we first review
some results proved earlier. Using the results from Step 1 and Skorohod
representation of weakly converging laws on a Polish space
we can assume that $\Xi^{\rm{\log},\alpha}_N(t_0,2)\Ntoo
\mathcal{L}_{t_0}(2),\;a.s.$. Note that conditioned on
$\mathcal{L}_{t_0}(2)$, the path $\{\mathcal{L}_{t}(2)\}_{t\geq t_0}$ is
deterministic.

Consider for $i\in \{1,\dots,L\}$ the system:

\bea{don6.1} d x^N_2 (i,t) = c(\bar x^N_2 (t) - x^N_2 (i,t)) dt
 &-& s\,x^N_2(i,t) (1-x^N_2 (i,t)) dt\\
 &-& \frac{m}{N} x^N_2 (i,t)dt\nonumber \\
 &+& \sqrt{d \cdot x^N_2(i,t)(1-x^N_2(i,t))} dw_2(i,t). \nonumber
\eea

Conditioned on $\mathcal{L}_{t_0}(2)$, we know already for $N \to \infty$  (on our
Skorohod probability space):
\be{ga43e}
\bar x^N_2 (t)\to \int x\mathcal{L}_t(2,dx),\; a.s.,\quad
\text{for}\; t\geq t_0.\ee Since we know from Step 1 that
\be{ga44} \CL [\{x^N_2 (i,t_0), \; i=1,2,\cdots,L\}] \Nto \CL[ \{
x^\infty_2 (i,t_0),\;i=1,\dots,L\} ], \ee then a standard coupling
argument yields
\be{ga45}
\CL [ \{x^N_2 (i,t)\}_{i=1,\dots,L;\;t\geq t_0}]
  \Nto \CL [\{ x^\infty_2 (i,t)_{i=1,\dots,L;\;t\geq t_0}\}], \ee
where $x^\infty_2(\cdot,\cdot)$ satisfies for $i=1,\dots  L, \quad t \geq t_0$.
\bea{don6.2} d x^\infty_2 (i,t) = && c(\int x \mathcal{L}_t (t,dx)
- x^\infty_2 (i,t)) dt
 - s\,x^\infty_2(i,t) (1-x^\infty_2 (i,t)) dt\\
 &&+ \sqrt{d \cdot x^\infty_2(i,t)(1-x^\infty_2(i,t))} dw_2(i,t). \nonumber
 \eea

{\bf Step~3} Now we need to show that we have for the equation
(\ref{don6.2}) a unique solution for $t \in \R$.

We can prove that, given $(\int x\mathcal{L}_t(2,dx))_{t \in \R}$,
which converges to 0 as $ t \to -\infty$,
(\ref{don6.2}) has a unique solution on $(-\infty,t_0]$ as
follows.

We can construct a minimal solution by considering a sequence
(in the parameter $n$) of solutions corresponding to starting the
process at time $t_n$. Then consider the
tagged site SDE driven by this mean curve starting at $0$ at time
$t_{n}\rightarrow-\infty$. The sequence of solutions forms  a stochastically monotone
increasing sequence as $n\rightarrow\infty$\ (use coupling). Then
observe that any solution $x(t)$ that satisfies $x(t)\rightarrow0$
as $t\rightarrow-\infty$ must have zeros (cf. classical result on
Wright-Fisher diffusions). \ Then by coupling at a zero we get
that it must agree with the minimal solution. This gives
uniqueness for a given mean path $(\CL_t(2))_{t \in \R}$.
Note that this also implies the
uniqueness of $(\CL_t)_{t \in \R}$ since we have already established
the uniqueness of the limiting mean curve.

This completes the proof.

\subsection{Extensions: non-critical migration, selection and mutation rates} \label{ss.non-critical}

The hierarchical mean field analysis allows us to investigate the
emergence times and behaviour for a wide range of scenarios
involving different parameter ranges for mutation rates, migration
rates and relative fitness of the different levels. Our main focus
in this work is the {\em critical case} in which mutation, migration and fitness all play a
comparable role and for this reason we have chosen the
parameterization
$\beta_1=0,\;\beta_2=0,\; \beta_3 =1$. \ However the basic
tools and analysis we develop here can be adapted to other scenarios. \
Although we will not carry out the analysis in detail we now
briefly indicate the main features of the different cases.

For example, consider the following ''general'' parametrization of the
migration, selection and mutation rates:
\be{agrev53}
c_{N,k} =\frac{c}{N^{(k-1)(1-\beta_{1})}},\;\;0 \leq \beta_{1}<1, \ee
\be{agrev54} s_{N}   =\frac{s}{N^{\beta_{2}}},\;\;\beta_{2} \geq 0,
\ee
\be{agrev55}
m_{N}   =\frac{m}{N^{\beta_{3}}}, \quad \beta_3 \geq 0.
\ee

The different ranges of the migration parameters have the following interpretation. The
value $\beta_{1}=0$ corresponds to euclidian space with dimension
2. \ The  values $0<\beta_{1}<1$ correspond to dimensions $d>2$. We do
not consider this case in this paper but for a detailed discussion
of the relation of random walks on the hierarchical group to
random walks on $\Z^d$ and for further references see  \cite{DGW01}.

We have discussed the case $\beta_1 = \beta_2 =0$ and $\beta_3=1$ above
and we will now briefly discuss the case (1) where  $\beta_1=0$ but $\beta_2>0$ and
(2) the cases $\beta_1=\beta_2=0$ but either $\beta_3<1$ or $\beta_3>1$,
which are the other most interesting cases.

{\em (1)}
In order to understand the behaviour
$\beta_1=0$, $\beta_{2}>0$
consider the behaviour of the dual process.
Then the dual process has the property that
most occupied sites have only
one factor (until the number of factors is $O(N)$) and the number
of factors grows like $e^{\frac{st}{N^{\beta_2}}}$. \ In
speeded-up time scale $N^{\beta_2}t$ the factors move quickly with
possibility of coalescence each time a pair collides.
If $\beta_3\leq 1$, then heuristically, emergence will occur when
\be{agrev56} e^{\frac{st}{N^{\beta_{2}}}} =O(N^{\beta _{3}}) \ee
that is, \be{agrev57}
 t=O(N^{\beta_2}\frac{\beta_3}{s} \log N).\ee
This means that emergence occurs much later than
in the critical regime considered in our work.

{\em (2)} We now restrict our attention to the case $\beta_1=\beta_2=0$ and
discuss the effect of different mutation rates.
The case of {\em
faster mutation}, $0<\beta_{3}<1$, has the effect that emergence
can occur before collisions become non-negligible. In particular,
if the mutation rate is $\frac{m}{N^{\beta_3}}$ with $\beta_3<1$
we get emergence at time $\beta_3\alpha\log N$ and the {\em
emergence is deterministic}, i.e. $\text{Var}(\bar
x^{N}_{2}(\alpha\beta_3\log N))\to0$. In this case if we look at
the dual process there are
$O(N^{1-\beta_3})$ occupied sites at time $O(1)$ and therefore we
have this number of droplets of size $N^{\beta_3}$ at time
$\alpha\beta_3\log N$. Therefore we have a law of large numbers
and deterministic emergence dynamics.

On the other hand, {\em slower mutation}, $\beta_{3}>1$,  means
that the ``wave of advance'' of the dual system must move to
distance 2 or higher depending on the value of $\beta_{3}>1$  in order to produce
the number of dual particles needed to guarantee emergence.

\section{Appendix. Nonlinear semigroup perturbations}\label{marsden}

\setcounter{secnum}{\value{section}} \setcounter{equation}{0}

We use the following result of Marsden \cite{MA}, (4.17).

\beT{T.pert}{(Perturbation)}

Let $\mathbb{B}$ be a Banach space and let $A_S$ be the infinitesimal generator of a
{\em strongly continuous semigroup},
with  $\|S_t\|\leq Me^{Ct}$ for some $C$.  Let $A_T:\mathbb{B}\to \mathbb{B}$ be a vector field on $\mathbb{B}$ such that $A_T$ is of class $C^2$ with its first and second derivatives uniformly bounded on bounded subsets and let $\{T_t\}$ be the flow of $A_T$.

Then $A_S+A_T$ has a unique flow which is Lipschitz for each $t,\; 0\leq t\leq T$, and
\be{add35}  V_tx=\lim_{n\to\infty} (S_{t/n}\cdot T_{t/n})^nx\ee
uniformly in $t$ for each $x$ on bounded sets of $t$.  If  $x\in \mathcal{D}(A_S+A_T)$, then
\be{add36}
\frac{d}{dt} V_tx=(A_S+A_T)U_t\ee
on $[0,\tau)$ where $\tau$ is the exit time from $\mathbb{B}$.
$\qad$
\end{theorem}

\newpage

\begin{center}

\end{center}

\end{document}